\documentclass[12pt,twoside,english]{article}
\usepackage[T1]{fontenc}
\usepackage{units}
\usepackage{amsmath}
\usepackage{amssymb}
\usepackage{graphicx}

\makeatletter

\providecommand{\tabularnewline}{\\}

\newcommand{\lyxaddress}[1]{
\par {\raggedright #1
\vspace{1.4em}
\noindent\par}
}

\usepackage{a4wide}
\usepackage{graphicx}
\usepackage{fancyhdr}
\usepackage{amsmath}
\usepackage{amssymb}
\usepackage{subfigure}
\usepackage{setspace}
\usepackage{verbatim} 
\numberwithin{equation}{section}

\newcommand{\vek}[1]{\mathchoice{\displaystyle\boldsymbol#1}
{\textstyle\boldsymbol#1}{\scriptstyle\boldsymbol#1}
{\scriptscriptstyle\boldsymbol#1}}
\newcommand{\mat}[1]{\mathchoice{\displaystyle\mathbf#1}
{\textstyle\mathbf#1}{\scriptstyle\mathbf#1}
{\scriptscriptstyle\mathbf#1}}

\makeatother

\usepackage{babel}
\begin{document}

\title{Higher-order meshing of implicit geometries---part II:\\Approximations
on manifolds}

\author{T.P. Fries, D. Sch{\"o}llhammer}

\maketitle

\lyxaddress{\begin{center}
Institute of Structural Analysis\\
Graz University of Technology\\
Lessingstr. 25/II, 8010 Graz, Austria\\
\texttt{www.ifb.tugraz.at}\\
\texttt{fries@tugraz.at}
\end{center}}
\begin{abstract}
A new concept for the higher-order accurate approximation of partial
differential equations on manifolds is proposed where a surface mesh
composed by higher-order elements is automatically generated based
on level-set data. Thereby, it enables a completely automatic workflow
from the geometric description to the numerical analysis without any
user-intervention. A master level-set function defines the shape of
the manifold through its zero-isosurface which is then restricted
to a finite domain by additional level-set functions. It is ensured
that the surface elements are sufficiently continuous and shape regular
which is achieved by manipulating the background mesh. The numerical
results show that optimal convergence rates are obtained with a moderate
increase in the condition number compared to handcrafted surface meshes.

\end{abstract}
\newpage{}\tableofcontents{}\newpage{}

\section{Introduction\label{sec:Introduction}}

Many challenging applications in engineering and natural sciences
are characterized by physical phenomena taking place on curved surfaces
in the three-dimensional space. There are numerous examples for \emph{transport
and flow} phenomena on biomembranes or bubble surfaces \cite{Gross_2011a,Xu_2012a}.
Examples in \emph{structures} are membranes and shells \cite{Chapelle_2011a,Blaauwendraad_2014a}.
Phenomena on surfaces may also be coupled to processes in the surrounding
volume such as in surfactant transport, hydraulic fracturing, reinforced
structures etc. As an additional challenge, the surfaces may be moving
\cite{Dziuk_2007b,Xu_2003a,Elliott_2012a}, i.e.~the domain of interest
changes. The modeling of such phenomena naturally leads to boundary
value problems where partial differential equations are formulated
on manifolds. For the solution of such models, customized numerical
methods are needed.

The first application of the finite element method for the solution
of the Laplace-Beltrami operator on manifolds is reported in 1988
by Dziuk \cite{Dziuk_1988a}. Since then, the topic has attracted
a tremendous research interest leading to a variety of numerical methods
for PDEs on surfaces existing today, see \cite{Dziuk_2013a} for an
overview. The most straightforward approach is to generate surface
meshes on the manifold and extend the finite element method in a natural
way using tangential differential calculus. That is, standard gradients
of the planar two-dimensional case are replaced by surface gradients
on the manifold. It is interesting to note that this approach has
been chosen from the beginning in the simulation of transport phenomena,
e.g.~\cite{Bertalmio_2001a,Adalsteinsson_2003a,Du_2011a,Dziuk_2013a}.
However, for the modeling of membranes and shells, a less intuitive
path using local coordinate systems and Christoffel symbols is standard
since a long time \cite{Chapelle_2011a,Blaauwendraad_2014a}. It is
rather recent that these models have been recasted in the frame of
global tangential operators \cite{Delfour_1995a,Delfour_1996a,Hansbo_2014a,Hansbo_2015a}.

Another approach is to only employ an implicit description of the
manifold and solve the model equations on \emph{all }iso-surfaces
at once \cite{Dziuk_2008a,Burger_2009a}. Then, the problem is naturally
set up in the three-dimensional space embedding the manifolds, i.e.~volumetric
elements and shape functions are employed. However, typically only
the solution on \emph{one} iso-surface, say the zero-isosurface, is
of interest. One may then restrict the surrounding domain to a narrow
band around the manifold \cite{Bertalmio_2001a,Deckelnick_2010a,Elliott_2009a}.
There are interesting similarities to phase field and diffuse interface
approaches \cite{Raetz_2006a}. A recent approach is to collaps the
narrow band to the manifold itself. Then, shape functions of the volumetric
background elements are used, however, the integration takes place
on the trace of the manifold only \cite{Olshanskii_2009b,Olshanskii_2009a,Hansbo_2014a}.
The resulting approaches are labeled TraceFEM \cite{Deckelnick_2014a,Grande_1916a,Olshanskii_2009b,Olshanskii_2016a,Reusken_2014a}
or CutFEM \cite{Hansbo_2004a,Hansbo_2014a}. Higher-order approximations
of PDEs on manifolds have been reported in different contexts before:
For explicit handcrafted surface meshes in \cite{Demlow_2009a} and
in the context of the TraceFEM in \cite{Reusken_2014a}. Adaptivity
is considered e.g.~in \cite{Demlow_2007a,Demlow_2012a}.

Herein, we propose a higher-order accurate approach for the approximation
of PDEs on manifolds. The manifold is described implicitly based on
the level-set method. A surface mesh composed by mixed higher-order
quadrilateral and triangular finite elements is automatically generated
from a background mesh and given level-set functions. A master level-set
function defines the shape of the manifold. However, as the implied
zero-isosurface may be infinite, it is restricted by additional (slave)
level-set functions. That is, several level-set functions imply the
bounded manifold being the domain of interest in the BVP. As a model
problem, we consider the Laplace-Beltrami operator and an instationary
advection-diffusion problem. Based on this, the extension of the approach
to more advanced transport problems on surfaces and in the simulation
of membranes and shells will be reported in forthcoming publications.

The automatic detection of higher-order surface elements has been
reported by the authors in \cite{Fries_2015a,Fries_2016b} in the
context of integration and interpolation. There, only a set of surface
elements is needed featuring double nodes and not necessarily fulfilling
$C_{0}$-continuity. In order to be suited for the approximation of
PDEs on surfaces as discussed herein, (1) continuity requirements
have to be fulfilled, (2) the elements must be sufficiently shape
regular, and (3) connectivity information in the usual FEM sense has
to be provided, enabling the concept of nodal degrees of freedom.
These issues are addressed herein with emphasis on higher-order accurate
approximations. Also, the concept of using several level-set functions
for the definition of the bounded manifold is new and an extension
of \cite{Fries_2016b}.

The paper is organized as follows: In Section \ref{X_Preliminaries}
we outline the geometric description of the bounded manifold based
on several level-set functions defined on a background mesh composed
by higher-order elements. The automatic generation of suitable higher-order
surface meshes is described in Section \ref{X_MeshGeneration}: The
reconstruction of surface elements approximating the zero-isosurface
of the master level-set function, the restriction by means of additional
(slave) level-set functions, the extraction of a continuous surface
mesh from the element set, and the manipulation of the background
mesh to achieve shape regular elements. Section \ref{X_FEMonManifolds}
shortly recalls the standard finite element approach for approximations
on meshed surfaces. Numerical results are presented in Section \ref{X_NumericalResults}
for curved lines in two dimensions and curved surfaces in three dimensions.
The Laplace-Beltrami operator is considered as well as instationary
advection-diffusion on manifolds. Finally, a summary and outlook is
given in Section \ref{X_Conclusions}.

\section{Preliminaries\label{X_Preliminaries}}

\begin{figure}
\centering

\subfigure[unconnected]{\includegraphics[height=4cm]{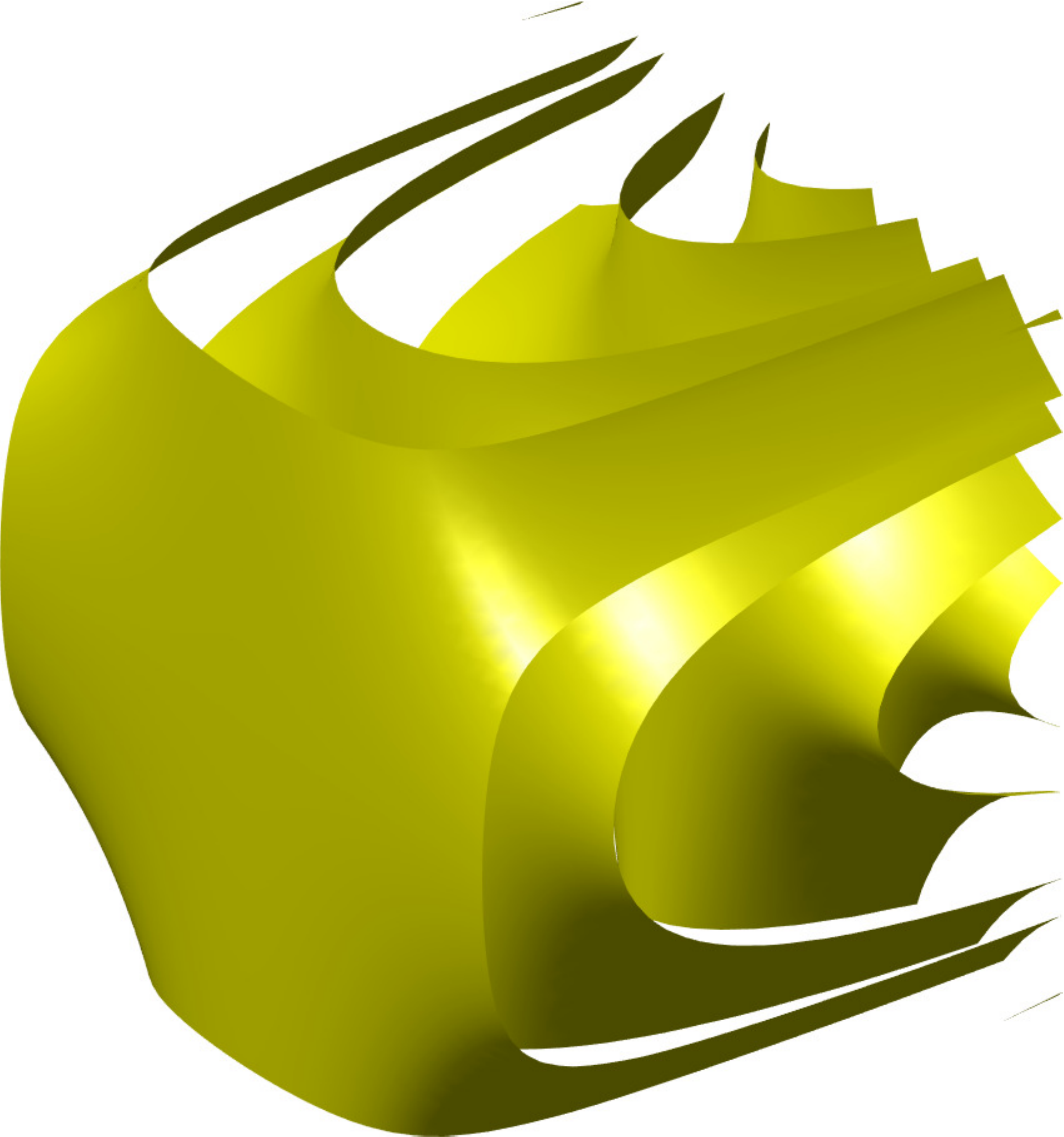}}\qquad\qquad\subfigure[orientable]{\includegraphics[height=4cm]{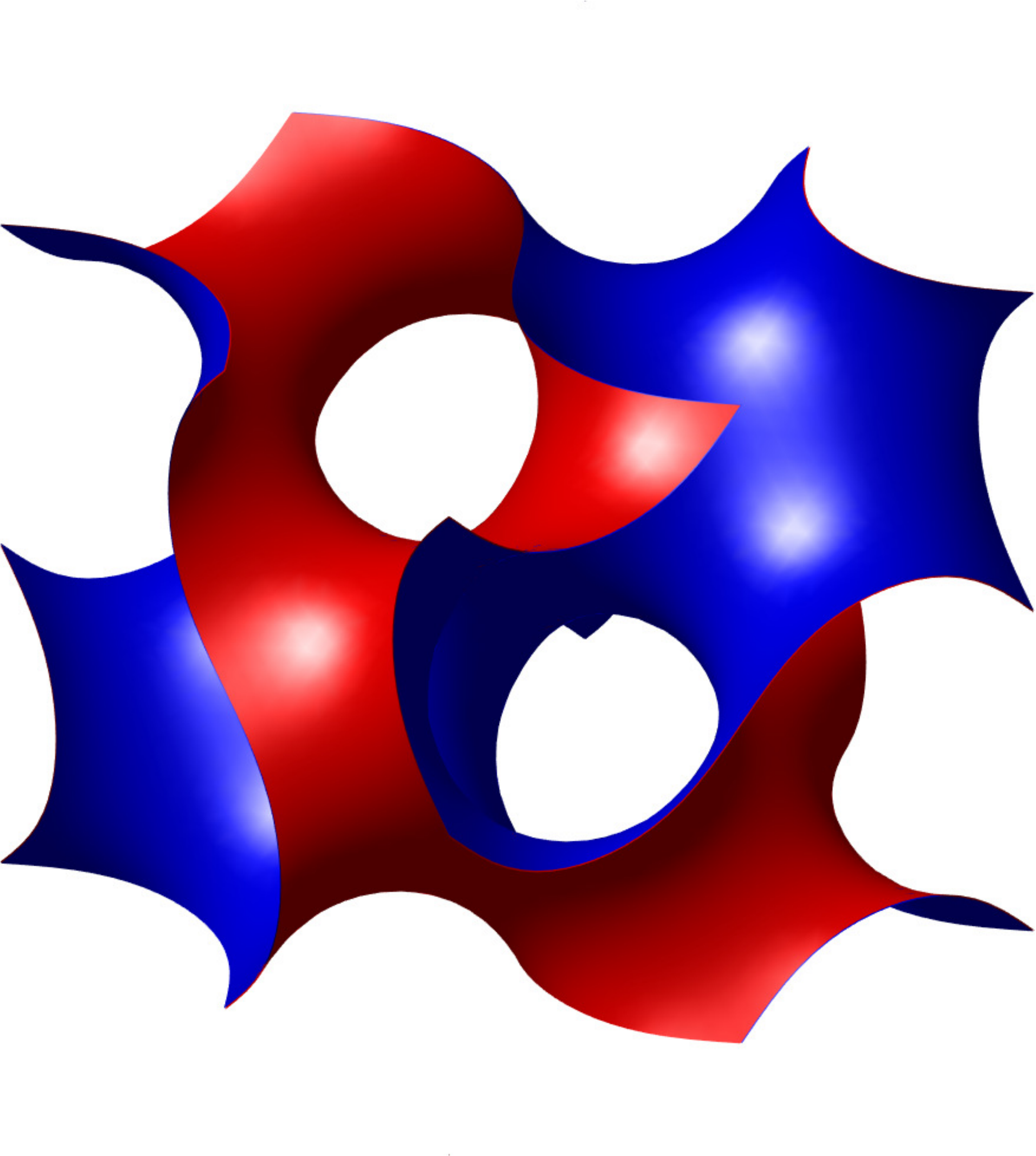}}\qquad\qquad\subfigure[closed]{\includegraphics[height=4cm]{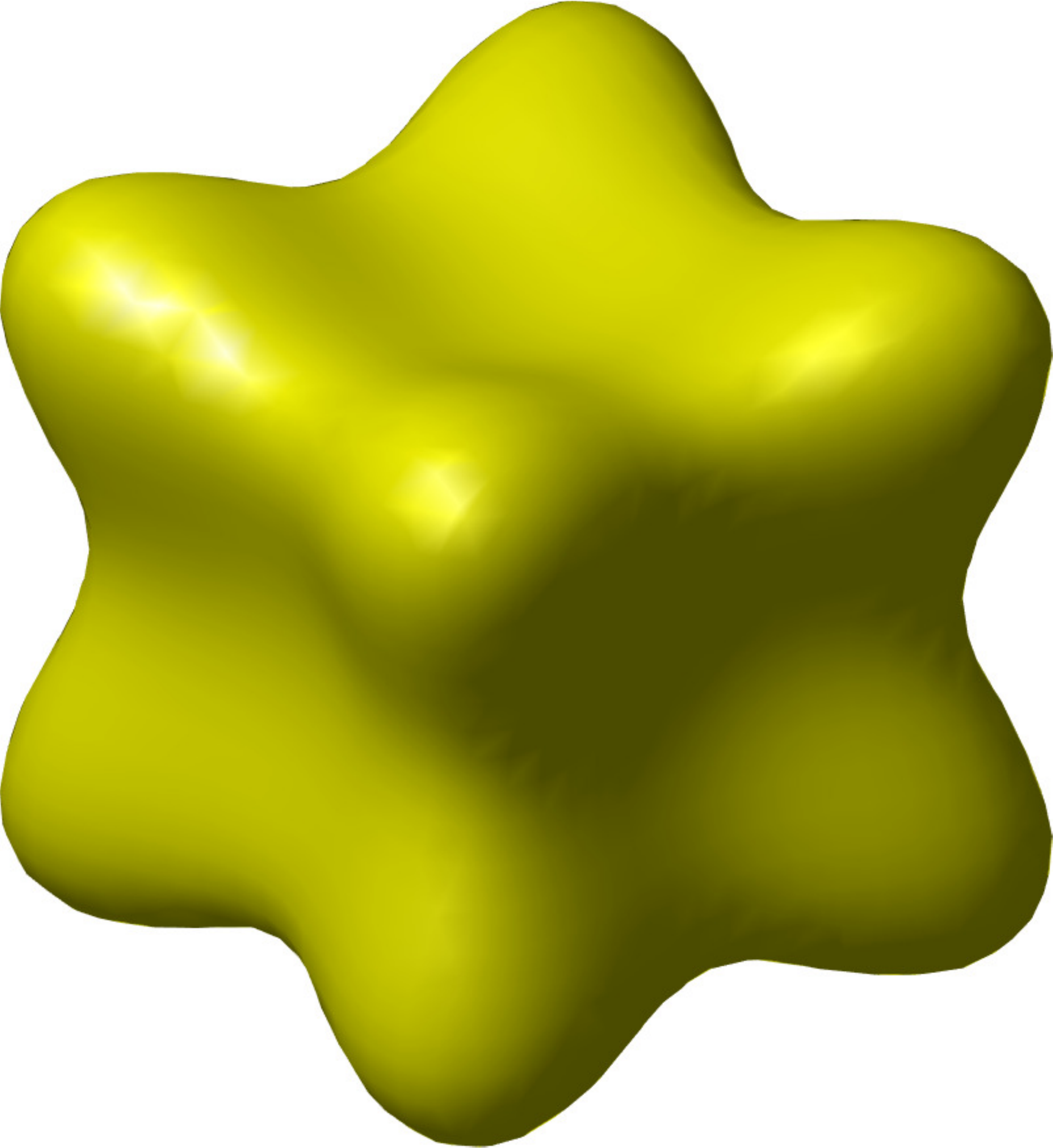}}

\caption{\label{fig:ZeroLevelSets3d}Examples for iso-surfaces: (a) shows unconnected
isosurfaces in $[-1,1]^{3},$ (b) shows that isosurfaces are orientable,
(c) shows a closed isosurface.}
\end{figure}

The task is to solve a boundary value problem (BVP) on a surface $\Gamma$
in three dimensions. Let the surface be possibly curved, sufficiently
smooth, orientable, connected (so there is only \emph{one} surface),
and feature a finite fixed area. The surfaces discussed herein are
defined following the concept of the level-set method. That is, there
is a continuous scalar function $\phi\left(\vek x\right)$ with $\vek x\in\mathbb{R}^{3}$
and its zero-level set
\begin{equation}
\Gamma_{\phi}=\left\{ \vek x\in\mathbb{R}^{3}:\,\phi\left(\vek x\right)=0\right\} \label{eq:UnboundedIsosurface}
\end{equation}
defines a (zero-)isosurface. We note that a general isosurface of
$\phi\left(\vek x\right)$ may be unconnected, see Fig.~\ref{fig:ZeroLevelSets3d}(a).
Then, it is necessary to select the \emph{one} zero-level set of interest
for the solution of a BVP, so this is not a problem. Furthermore,
isosurfaces of level-set functions are orientable by default, i.e.~there
is a consistent normal vector pointing from the negative to the positive
side, see Fig.~\ref{fig:ZeroLevelSets3d}(b).

There is no concept of a boundary when only \emph{one} level-set function
$\phi\left(\vek x\right)$ is considered with $\vek x\in\mathbb{R}^{3}$.
As a consequence, isosurfaces with \emph{finite} area must be closed
(compact) as e.g.~shown in Fig.~\ref{fig:ZeroLevelSets3d}(c), hence,
there is no boundary. Otherwise, in the more general case, they are
open and feature an \emph{infinite} area without a boundary. It is
thus seen that the manifold of interest, i.e.~the surface $\Gamma$
where a BVP is to be solved, is typically only a subset of $\Gamma_{\phi}$,
hence, $\Gamma\subseteq\Gamma_{\phi}$.

\begin{figure}
\centering

\subfigure[bounded manifold]{\includegraphics[height=4cm]{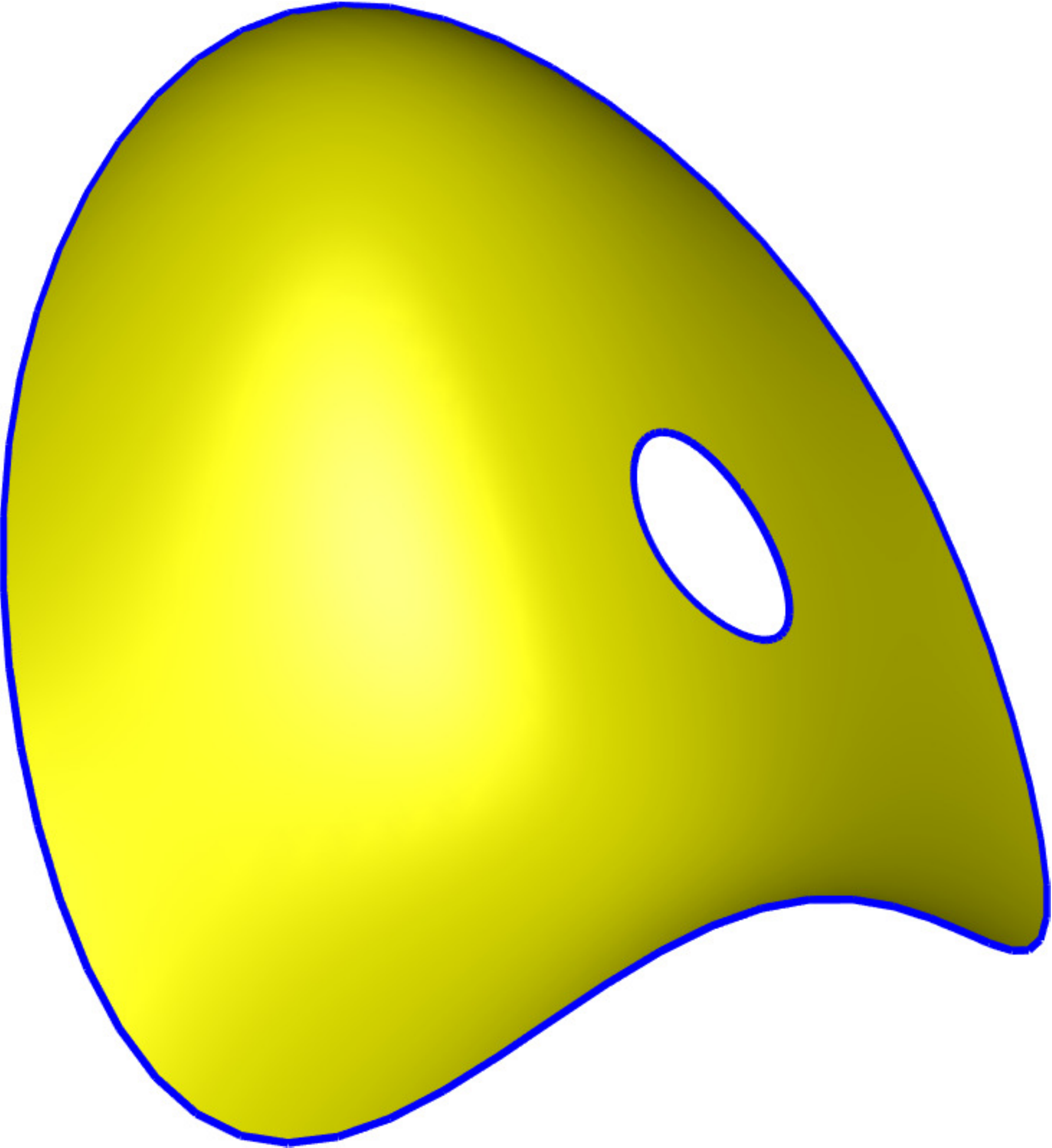}}\qquad\subfigure[subregion $\Omega$]{\includegraphics[height=4cm]{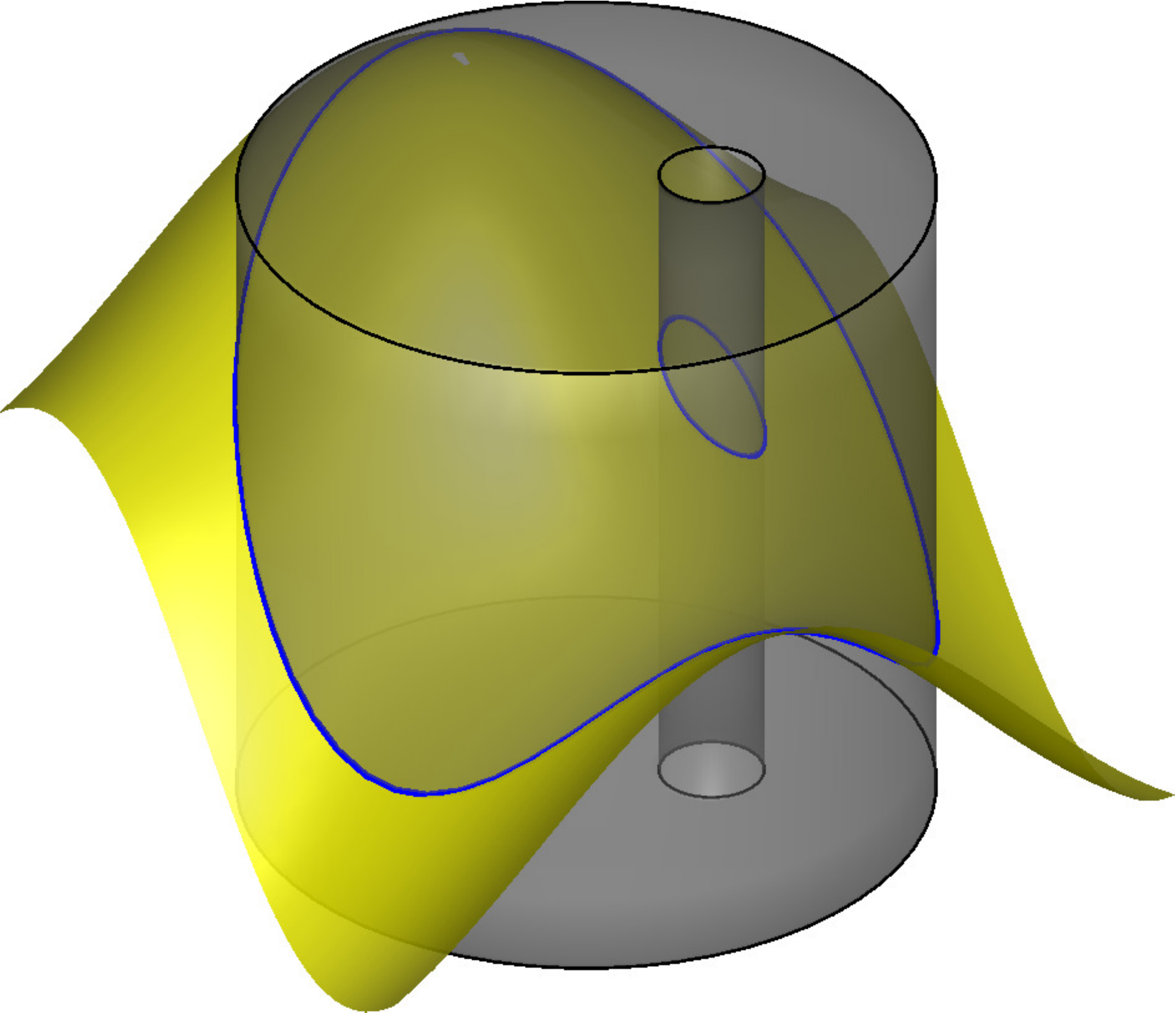}}\qquad\subfigure[several level-sets]{\includegraphics[height=4cm]{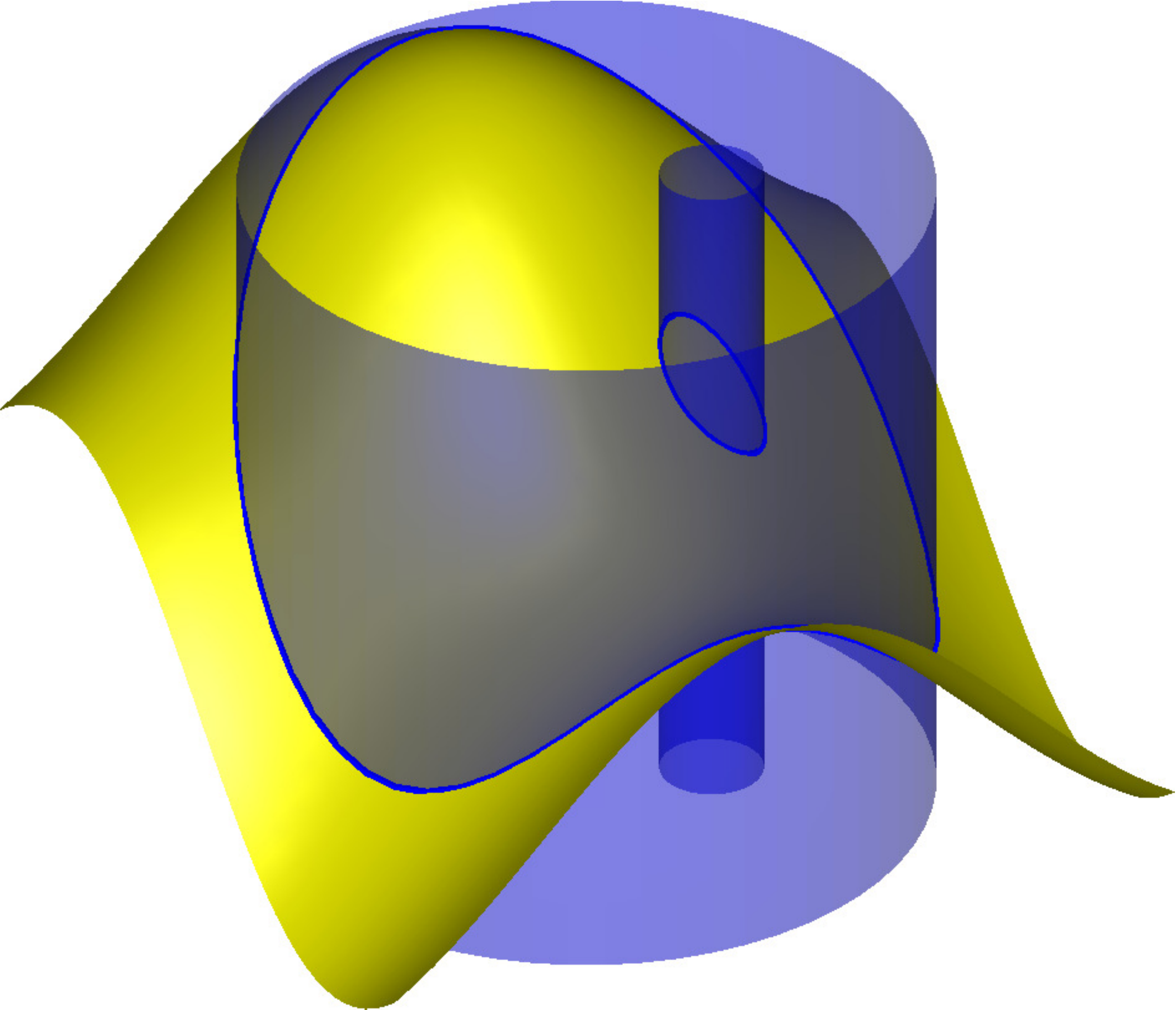}}

\caption{\label{fig:BoundedIsosurfaces}The bounded isosurface in (a) may be
described by restricting the evaluation of $\phi\left(\vek x\right)$
to the subregion $\Omega$ shown in (b) or by using additional level-set
functions whose isosurfaces are seen in (c).}
\end{figure}

There are two alternatives to define \emph{bounded} isosurfaces with
finite area, see Fig.~\ref{fig:BoundedIsosurfaces}(a), which fully
represent the manifold of interest: One is to evaluate the level-set
function $\phi$ only in a subregion $\Omega\in\mathbb{R}^{3}$, hence
\begin{equation}
\Gamma=\left\{ \vek x\in\Omega:\,\phi\left(\vek x\right)=0\right\} .\label{eq:BoundedIsosurfaceOmega}
\end{equation}
See Fig.~\ref{fig:BoundedIsosurfaces}(b) for an example, where $\phi\left(\vek x\right)$
is only evaluated in the gray subregion $\Omega$ instead of $\mathbb{R}^{3}$.
The other is to employ additional level-set functions $\psi^{i}\left(\vek x\right)$
to restrict $\Gamma_{\phi}$, hence
\begin{equation}
\Gamma=\left\{ \vek x\in\mathbb{R}^{3}:\,\phi\left(\vek x\right)=0\;\mathrm{and}\;\psi^{i}\left(\vek x\right)\leq0,\, i=1,2...\right\} .\label{eq:BoundedIsosurfaceMult}
\end{equation}
See Fig.~\ref{fig:BoundedIsosurfaces}(c) for an example, where two
additonal isosurfaces of $\psi^{1}$ and $\psi^{2}$ are shown in
blue. The same bounded isosurface from Fig.~\ref{fig:BoundedIsosurfaces}(a)
is thereby defined.

\begin{figure}
\centering

\subfigure[ex.~1, isosurfaces]{\includegraphics[width=3.5cm]{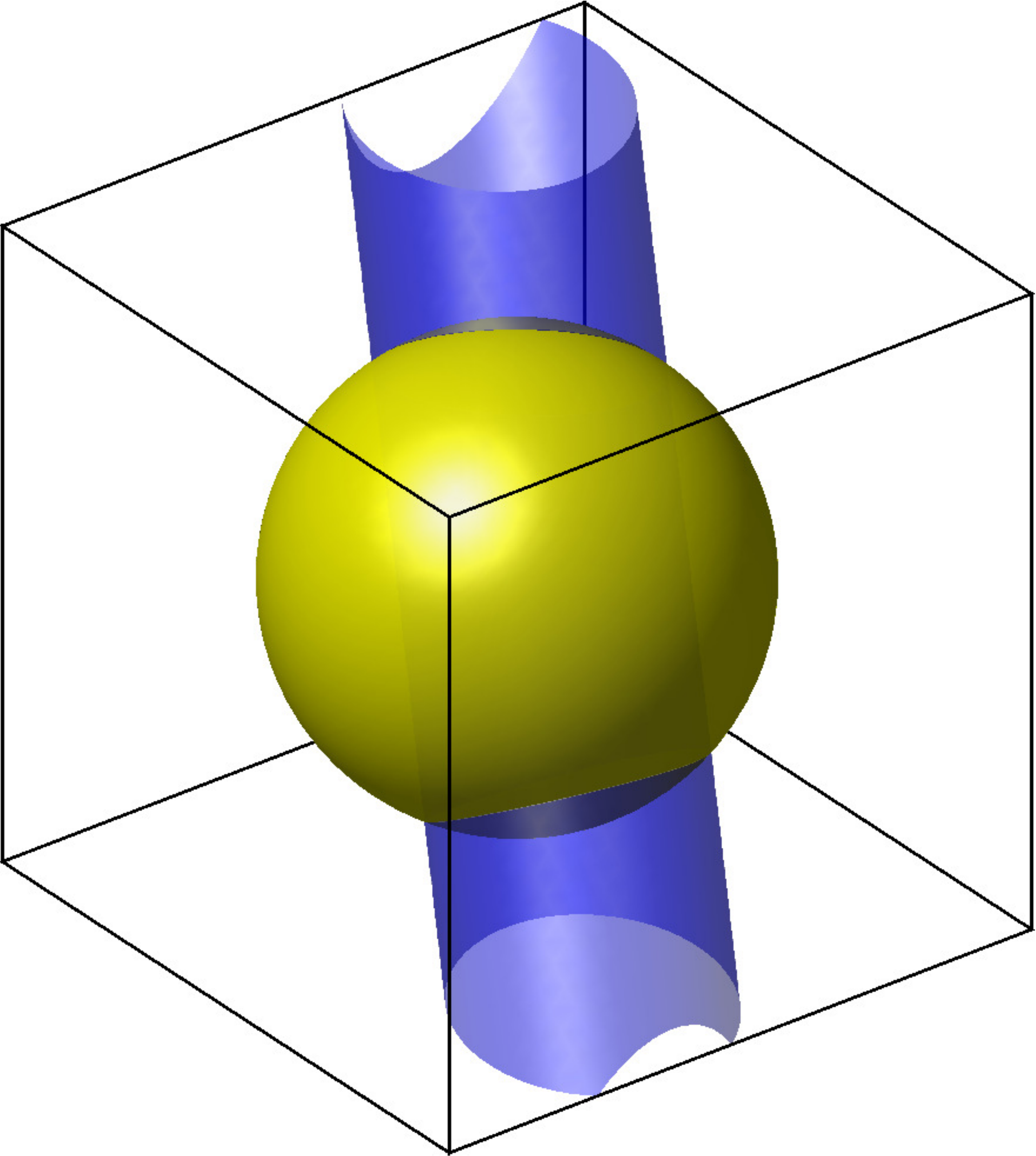}}\qquad\subfigure[ex.~2, isosurfaces]{\includegraphics[width=3.5cm]{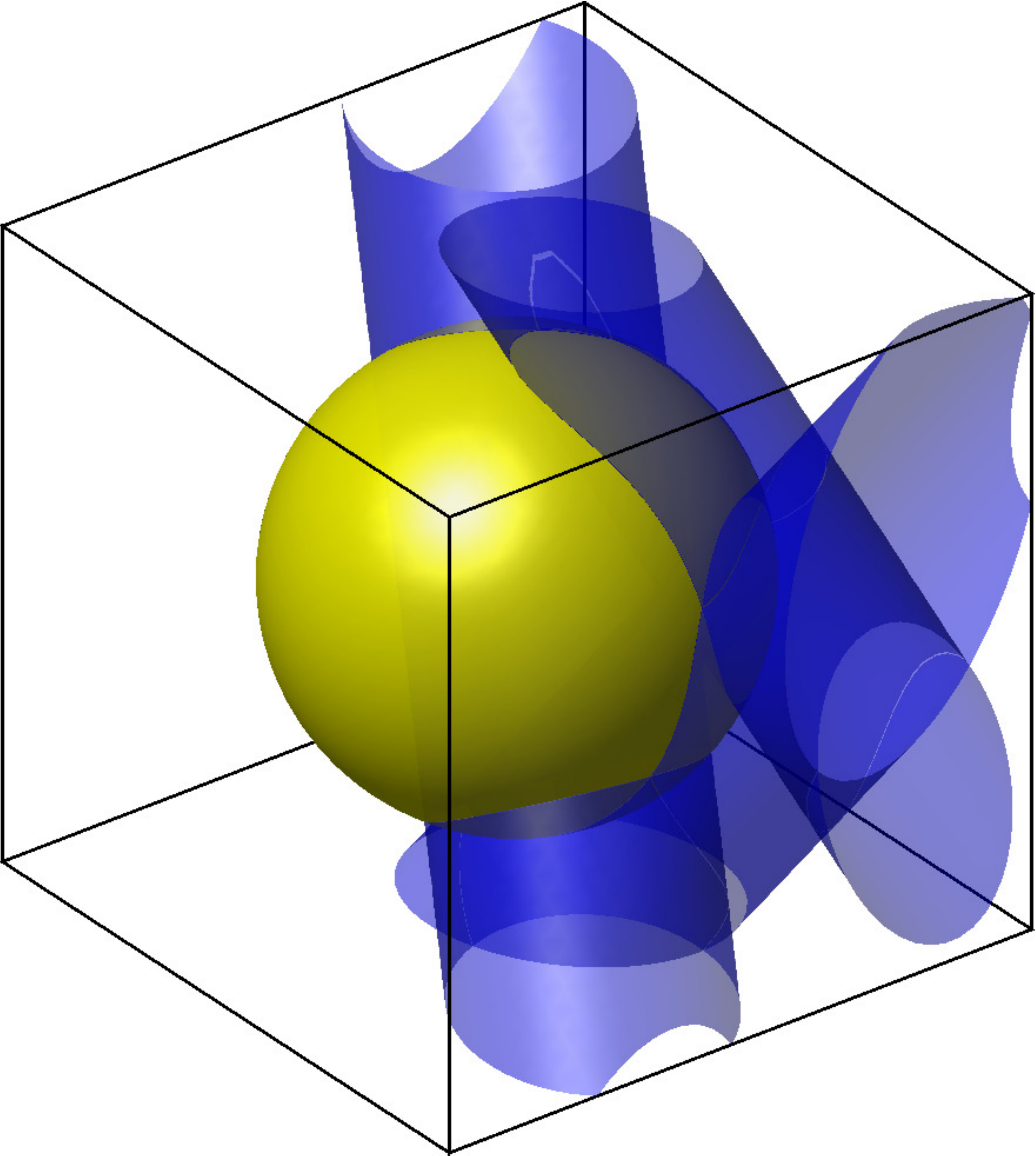}}\qquad\subfigure[ex.~3, isosurfaces]{\includegraphics[height=3.5cm]{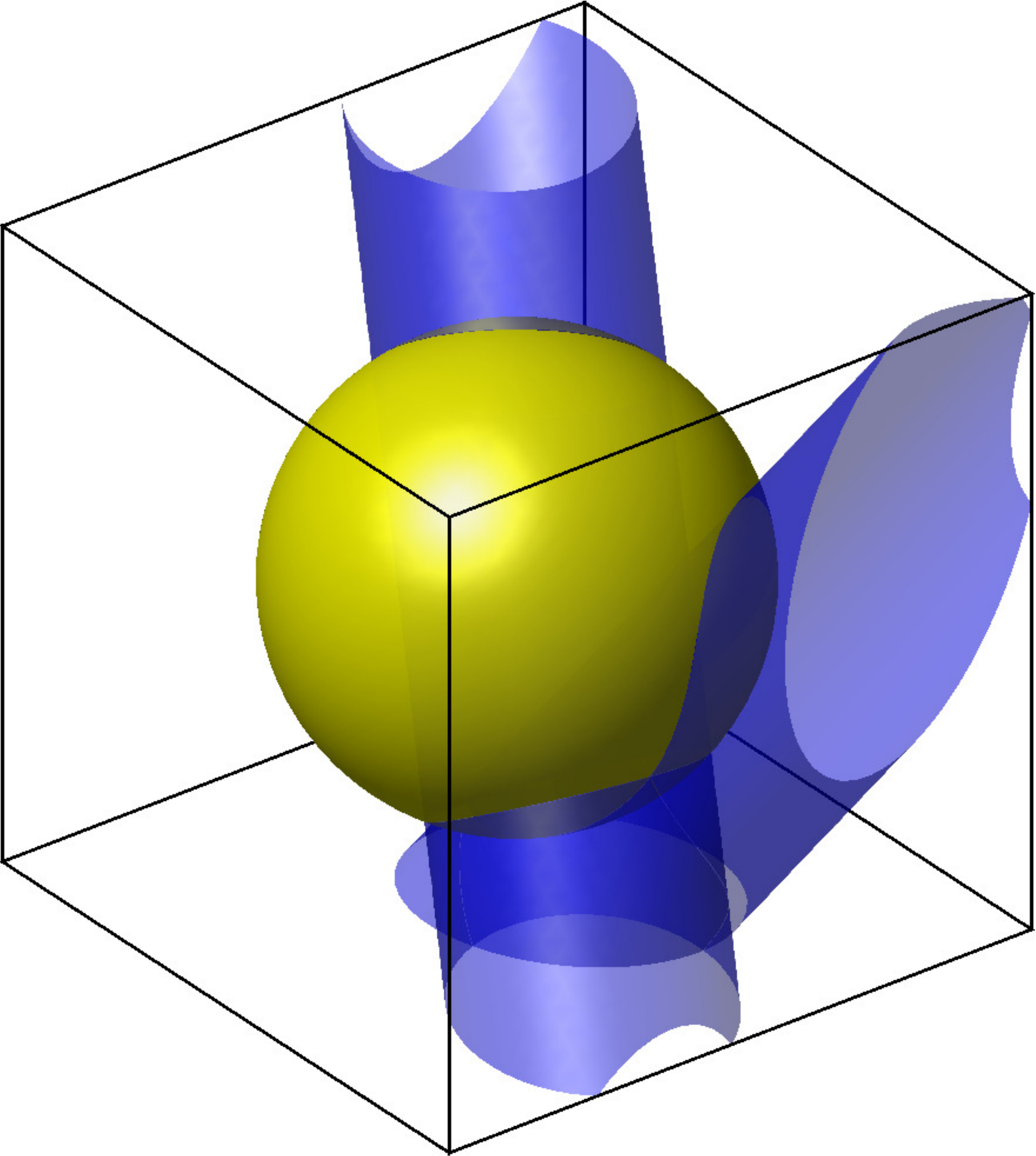}}\qquad\subfigure[ex.~4, isosurfaces]{\includegraphics[width=3.5cm]{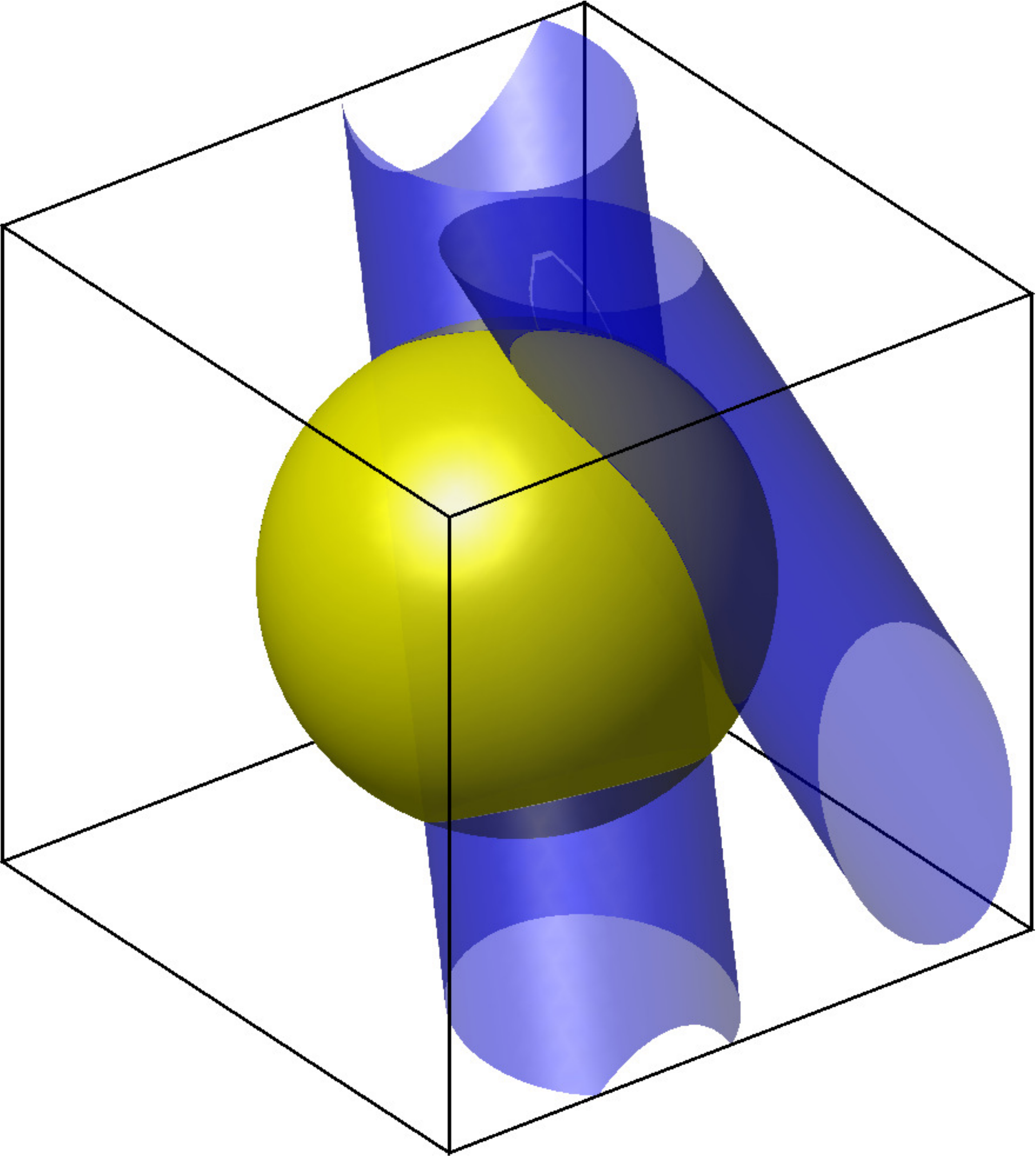}}

\subfigure[ex.~1, manifold]{\includegraphics[width=3cm]{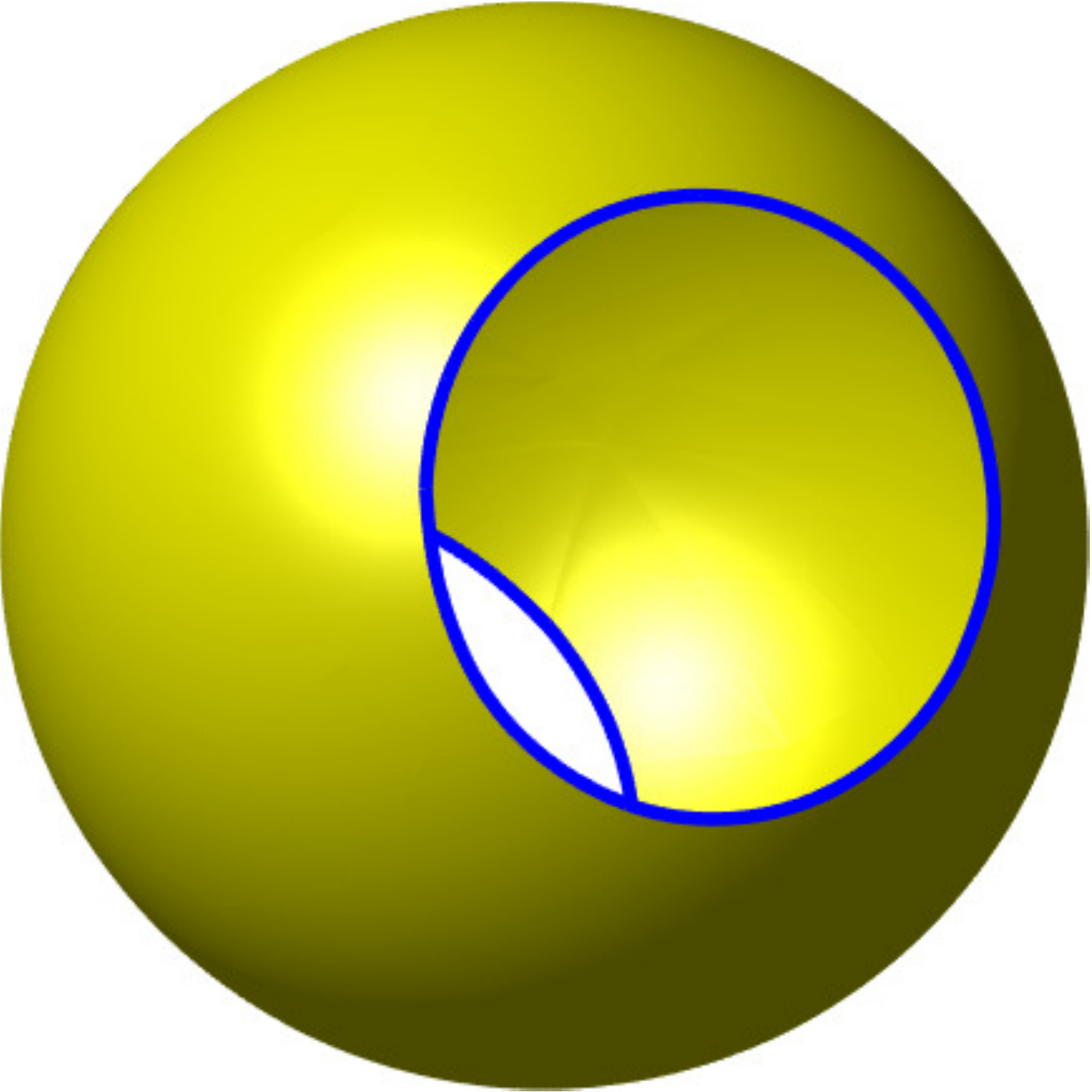}}\qquad\subfigure[ex.~2, manifold]{\includegraphics[width=3cm]{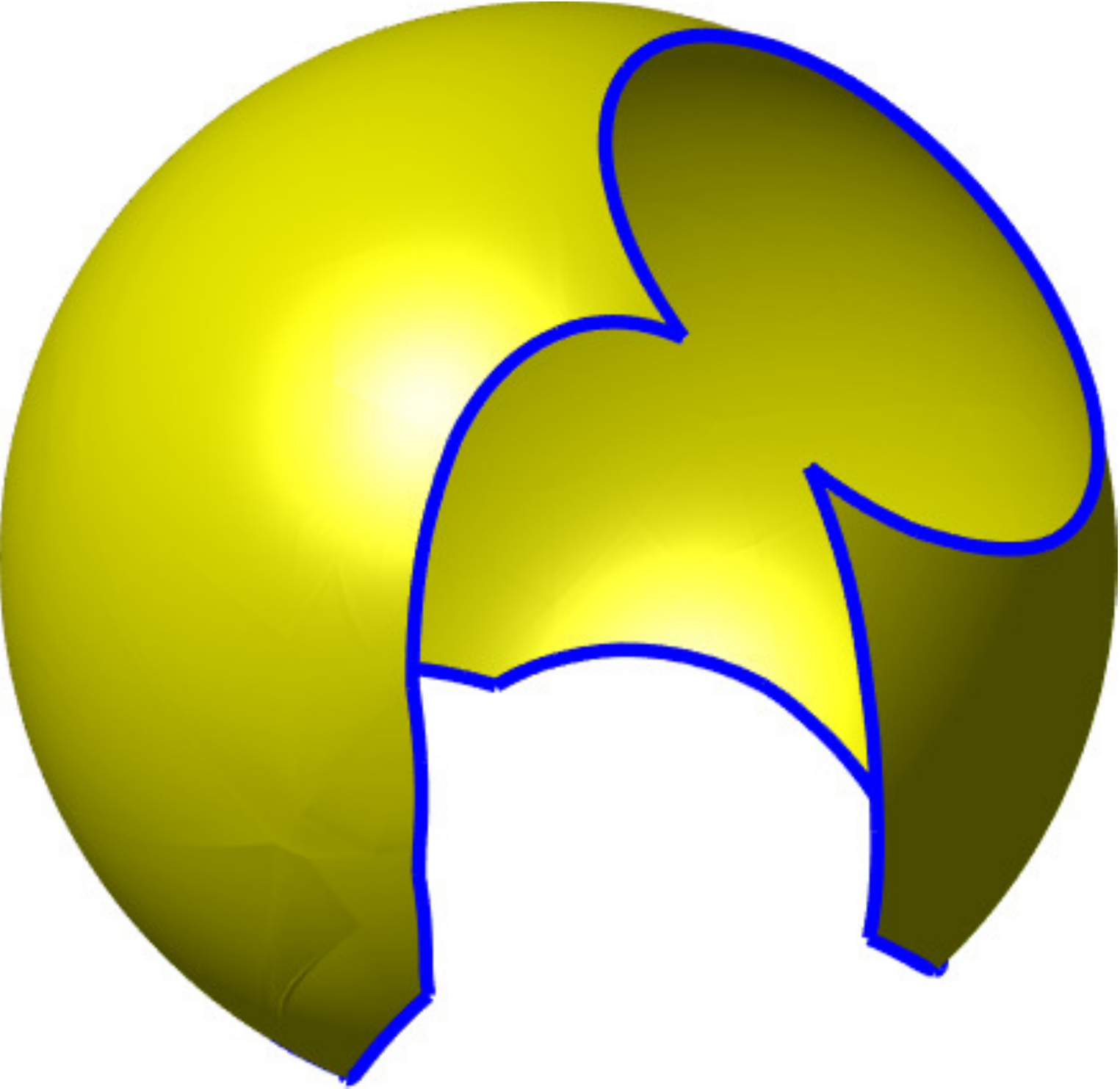}}\qquad\subfigure[ex.~3, manifold]{\includegraphics[height=3cm]{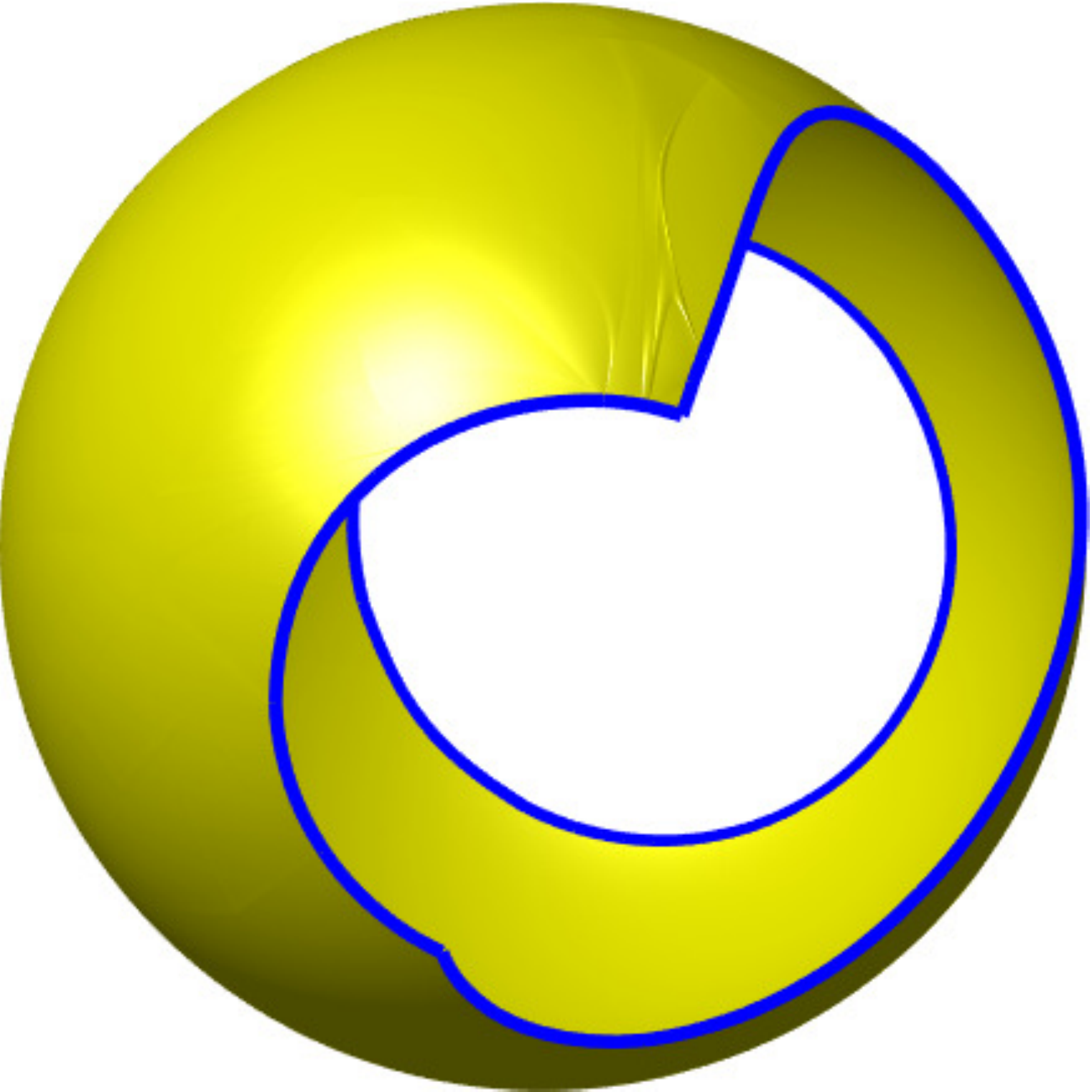}}\qquad\subfigure[ex.~4, manifold]{\includegraphics[width=3cm]{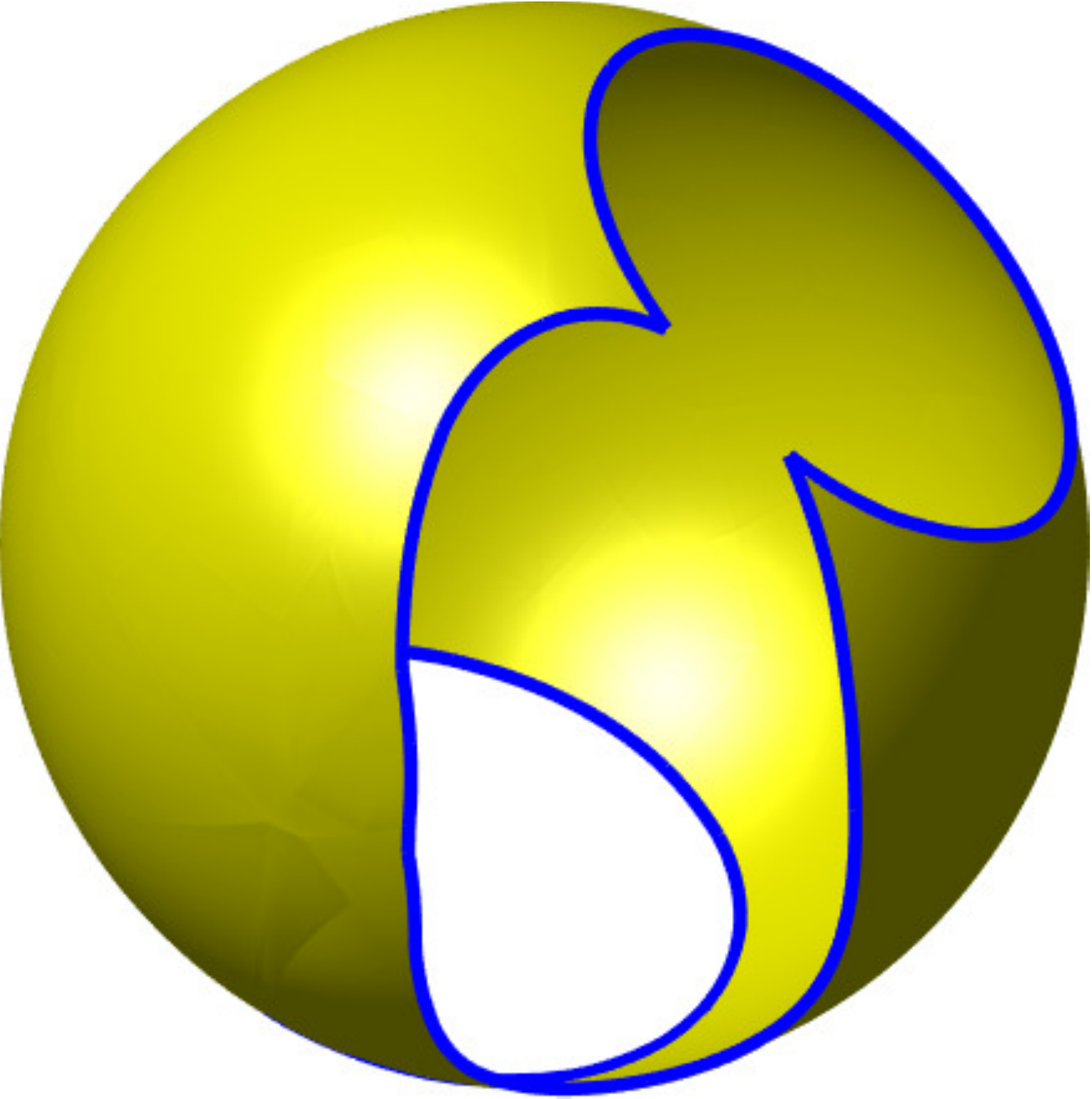}}

\caption{\label{fig:ExampleManifolds}(a) to (d) show some isosurfaces, yellow
of $\phi$ and blue of $\psi^{i}$, (e) to (h) show the implied bounded
manifolds (from different viewpoints), respectively. Several continuous
$\psi^{i}$ are able to define corners in the boundaries of the manifold.}
\end{figure}

Further examples are shown in Fig.~\ref{fig:ExampleManifolds}. In
Fig.~\ref{fig:ExampleManifolds}(a) to (d), zero-isosurfaces are
shown, in yellow of $\phi\left(\vek x\right)$ and in blue of $\psi^{i}\left(\vek x\right)$.
It is seen that complex, bounded manifolds are thereby implied, see
Fig.~\ref{fig:ExampleManifolds}(e) to (h). It is clear, that it
would be very difficult to define such complex manifolds by using
only $\phi\left(\vek x\right)$ and defining subregions $\Omega$.
Therefore, we prefer the second alternative of using multiple level-set
functions which are evaluated in $\mathbb{R}^{3}$. Of course, a combination
of these two alternatives is straightforward: This is in particular
useful when the outer boundary of the manifold is straight and follows
from a box-like subdomain, e.g.~$\Omega=[-1,1]^{3}$, but features
internal boundaries which are rather described by $\psi^{i}$.

An important consequence of using several level-set functions $\psi^{i}\left(\vek x\right)$
to restrict $\Gamma_{\phi}$ is that \emph{corners} in the boundaries
of the manifold may be defined straightforwardly (although each $\psi^{i}$
is smooth). This seen in Fig.~\ref{fig:ExampleManifolds}(f) to (h)
where more than one $\psi^{i}$ is present. It is noted that the employed
level-set functions $\phi$ and $\psi^{i}$ do \emph{not} have to
be signed distance functions.

\begin{figure}
\centering

\subfigure[ex.~1, mesh]{\includegraphics[width=3cm]{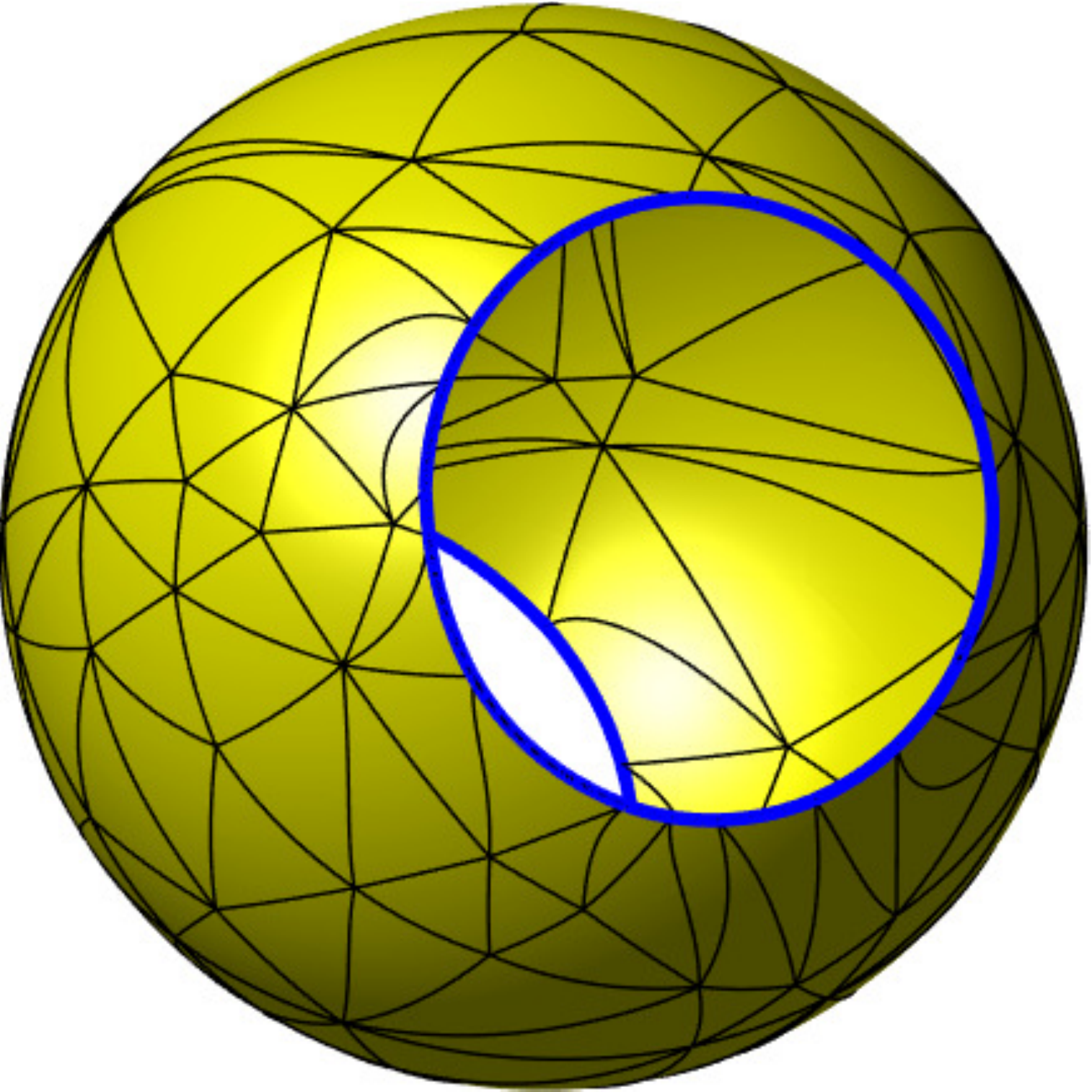}}\qquad\subfigure[ex.~2, mesh]{\includegraphics[width=3cm]{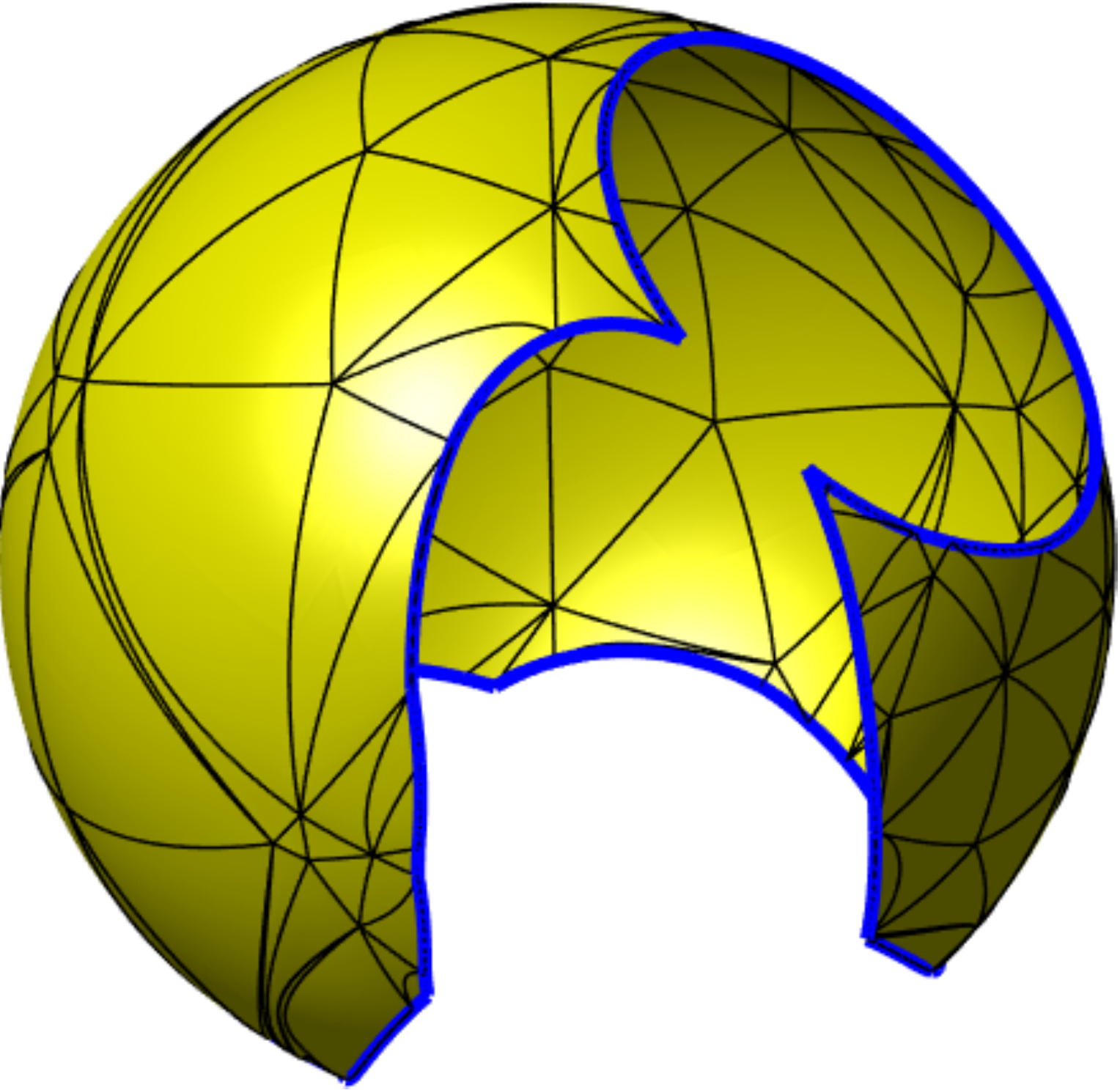}}\qquad\subfigure[ex.~3, mesh]{\includegraphics[height=3cm]{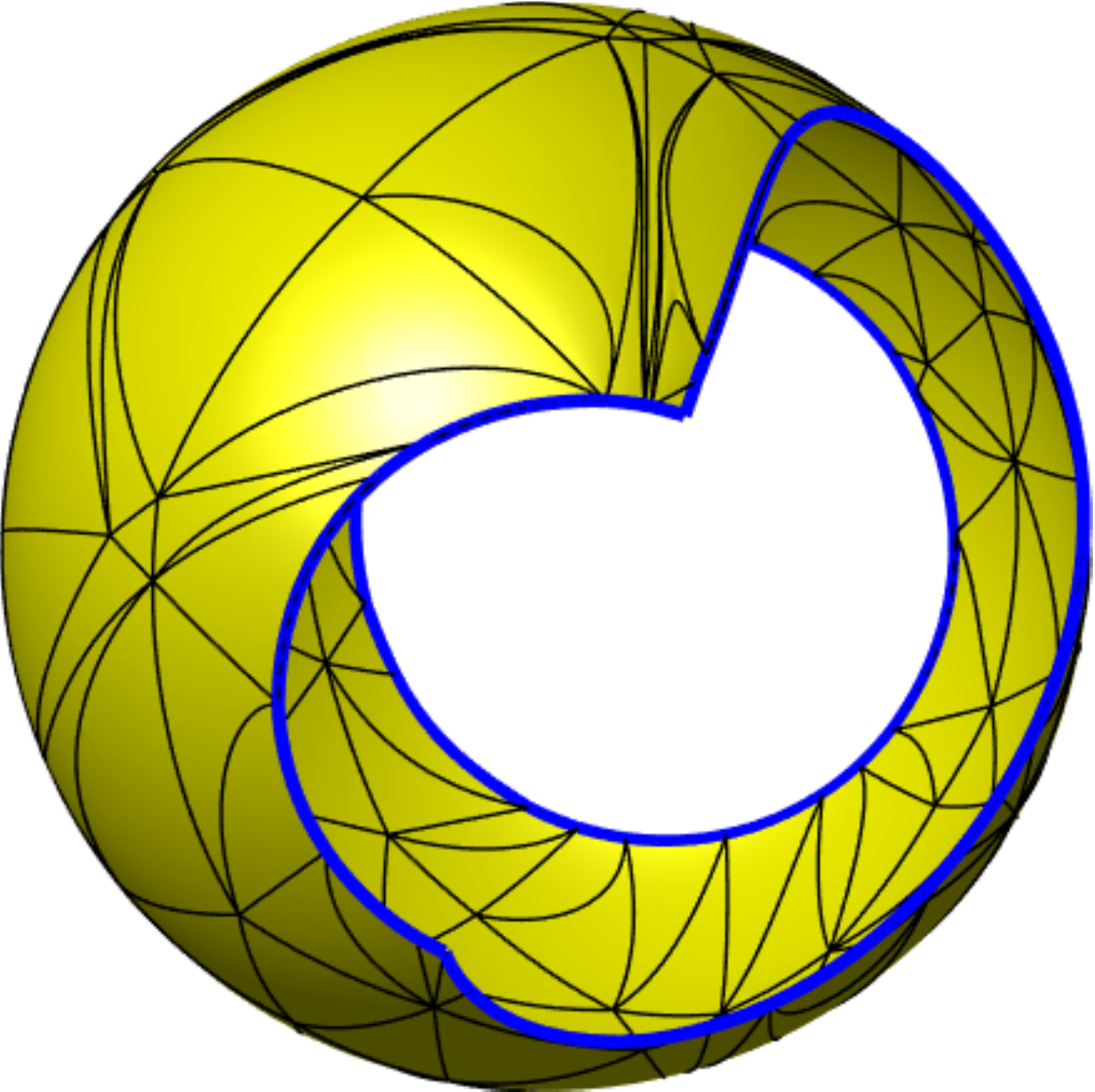}}\qquad\subfigure[ex.~4, mesh]{\includegraphics[width=3cm]{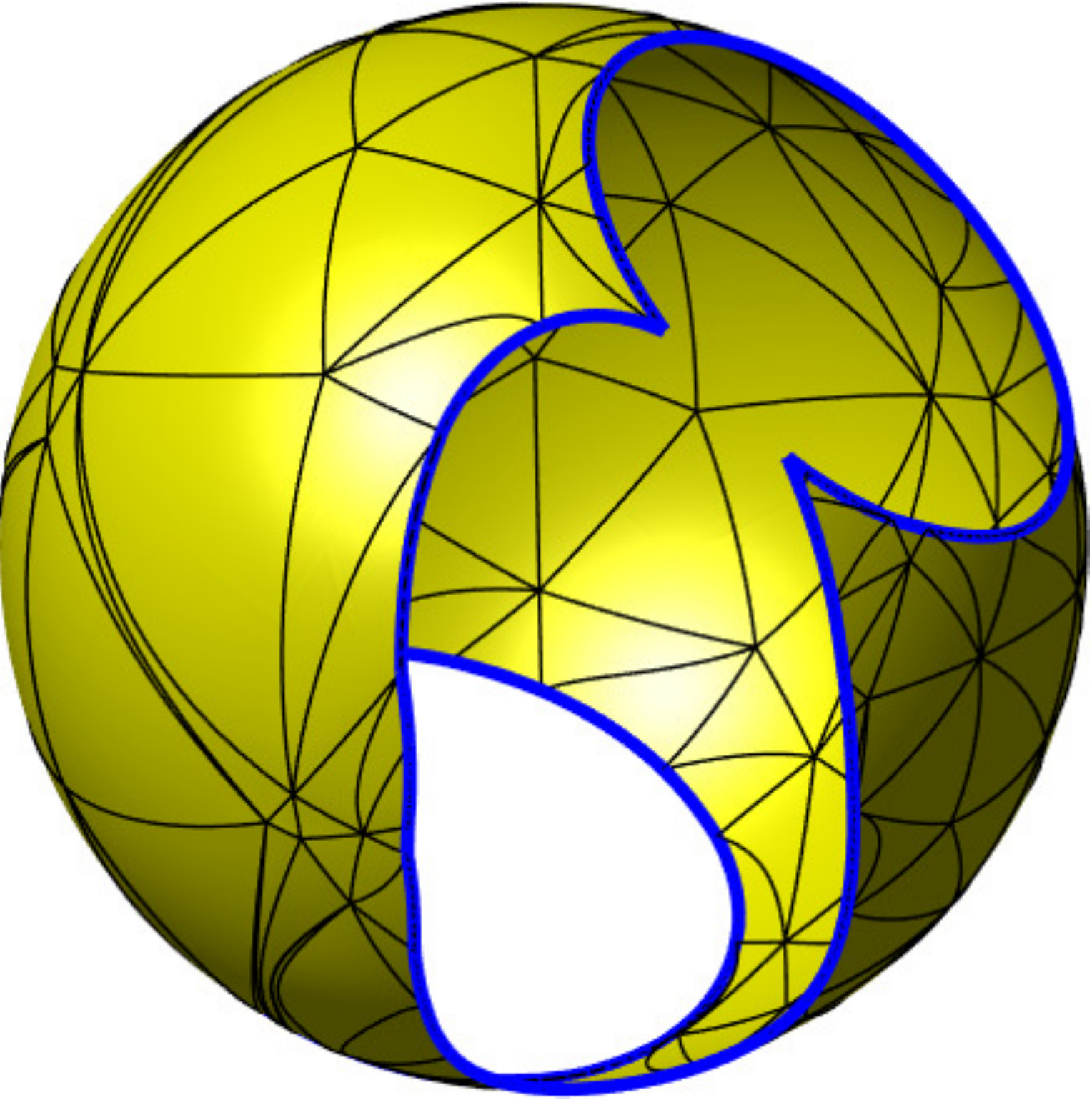}}

\caption{\label{fig:MeshManifolds}Some automatically generated higher-order
meshes of the bounded manifolds shown in Fig.~\ref{fig:ExampleManifolds}.}
\end{figure}

The task is now to automatically generate higher-order meshes of the
bounded manifolds. See Fig.~\ref{fig:MeshManifolds} for some examples
of meshes referring to the bounded isosurfaces shown in Fig.~\ref{fig:ExampleManifolds}.

\section{Mesh generation of the manifold\label{X_MeshGeneration}}

The procedure is outlined as follows: First, a background mesh is
introduced such that the manifold of interest is completely immersed.
The level-set data is evaluated at the nodes of the higher-order elements
of this background mesh. Second, the isosurface $\Gamma_{\phi}$ is
identified and meshed which is called reconstruction. As mentioned
before, this isosurface is either closed or is bounded by the boundary
of the background mesh. Third, the meshed isosurface $\Gamma_{\phi}$
is restricted by the additonal level-set functions $\psi^{i}$.

\subsection{Background mesh\label{XX_BackgroundMesh}}

A background mesh is introduced such that the manifold of interest
is completely immersed. The mesh is composed by higher-order background
elements of the Lagrange class and we restrict ourselves to tetrahedral
elements for simplicity. Hexahedral elements could be decomposed into
tetrahedra and treated similarly, see \cite{Fries_2015a}. There are
no restrictions on the background mesh, i.e.~it may be unstructured
and the element faces may be curved.

\begin{figure}
\centering

\subfigure[building block]{\includegraphics[width=2.5cm]{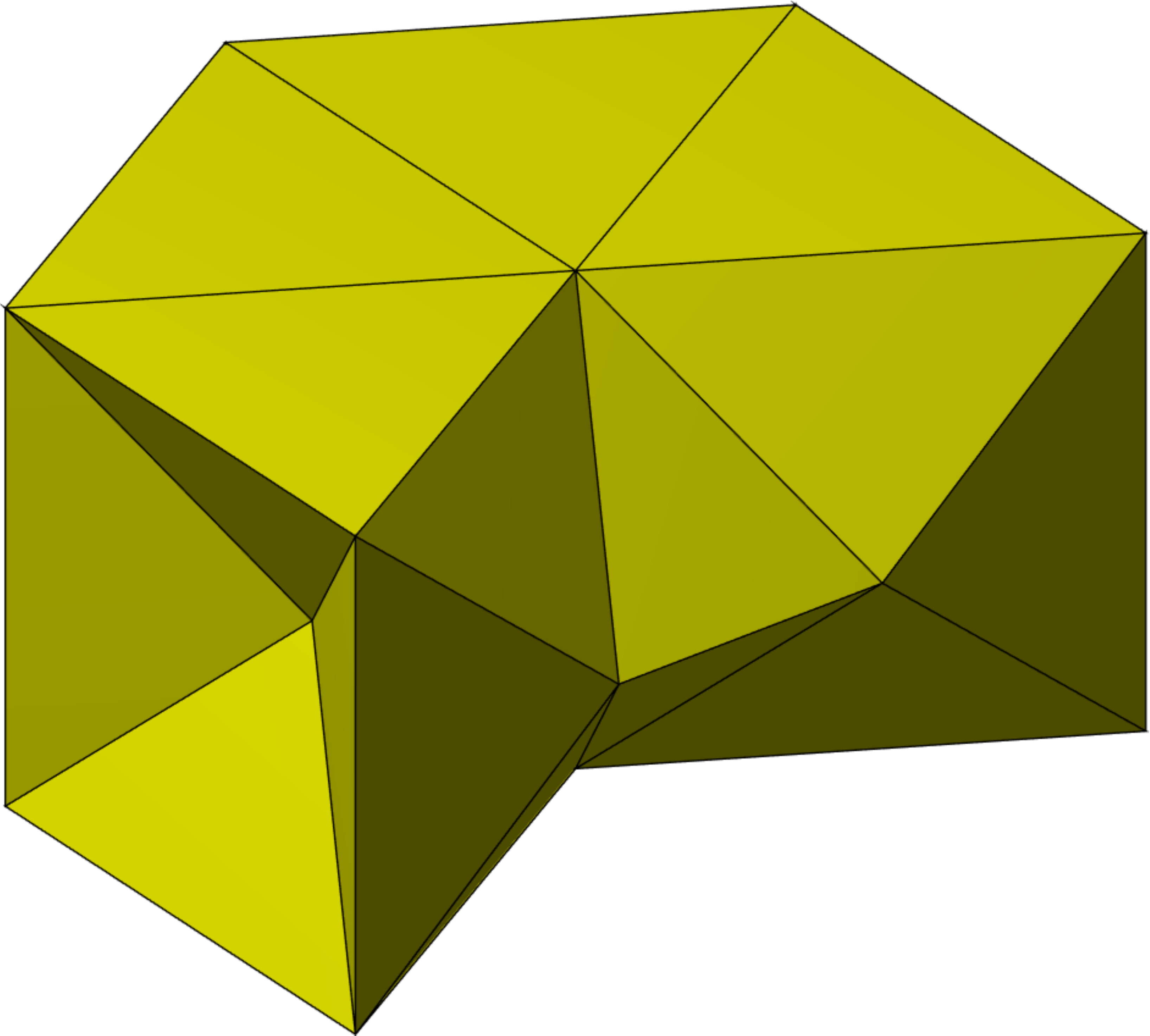}}\qquad\subfigure[universal mesh]{\includegraphics[width=3.5cm]{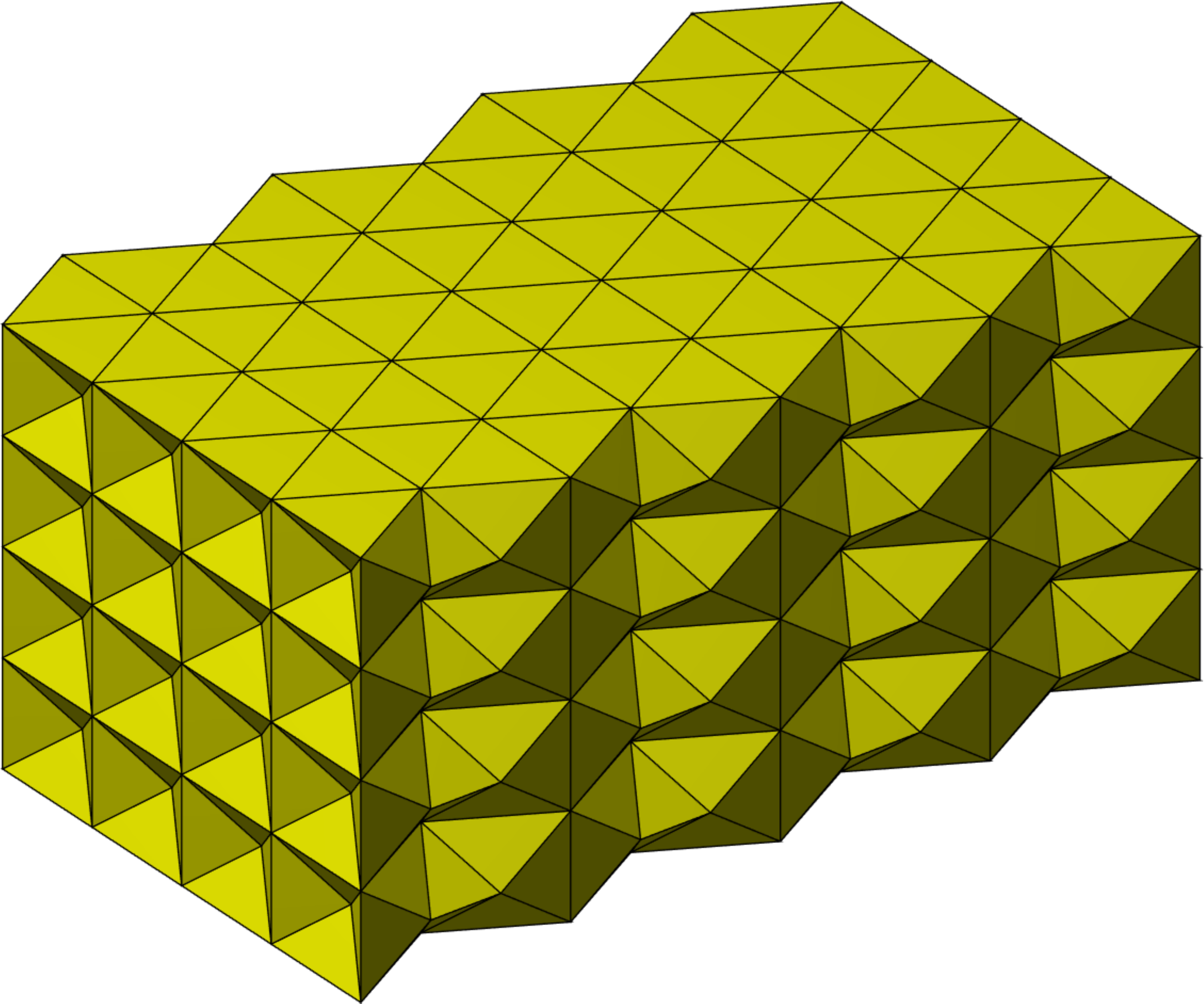}}\qquad\subfigure[building block]{\includegraphics[width=2.5cm]{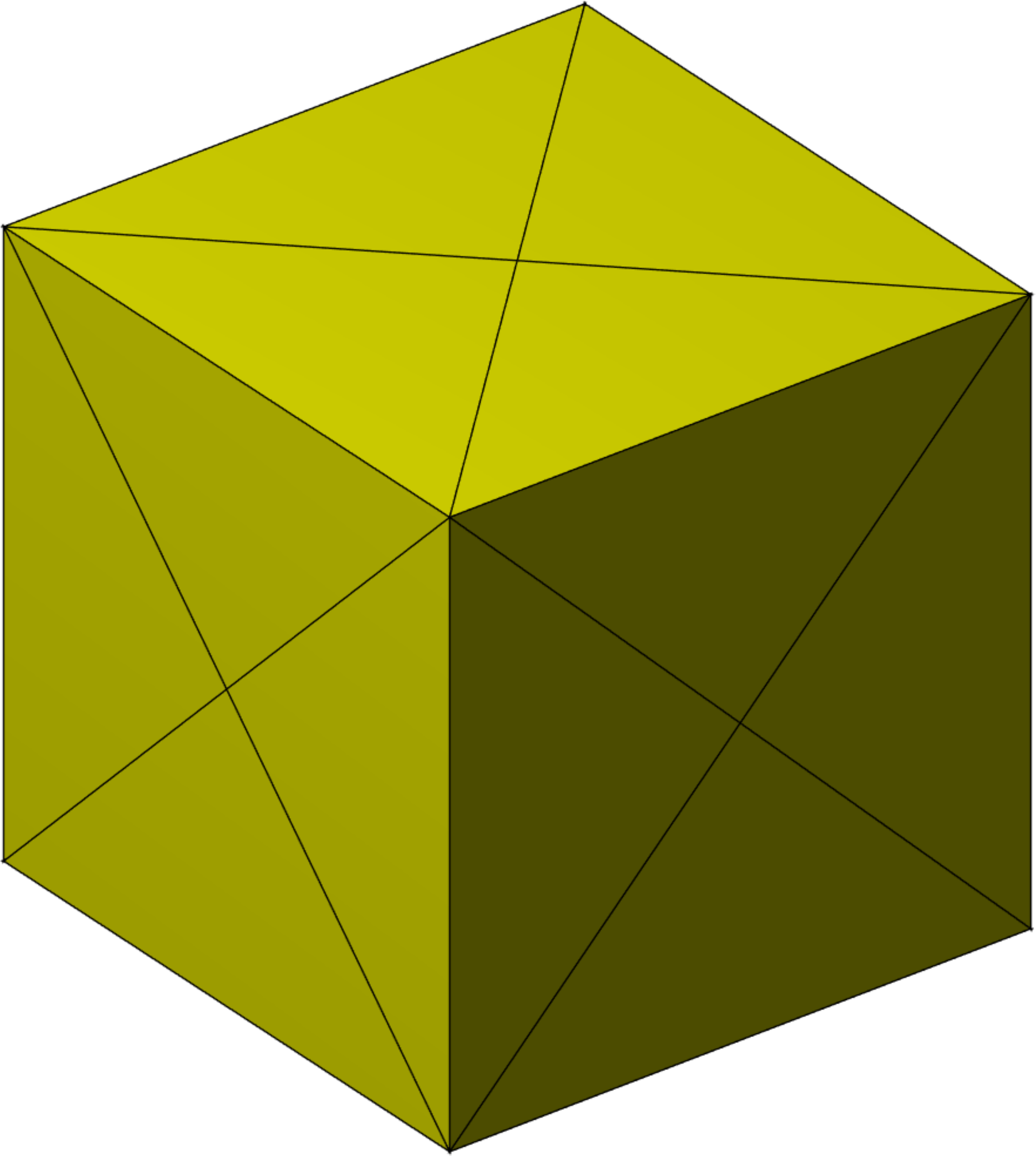}}\qquad\subfigure[box mesh]{\includegraphics[width=3.5cm]{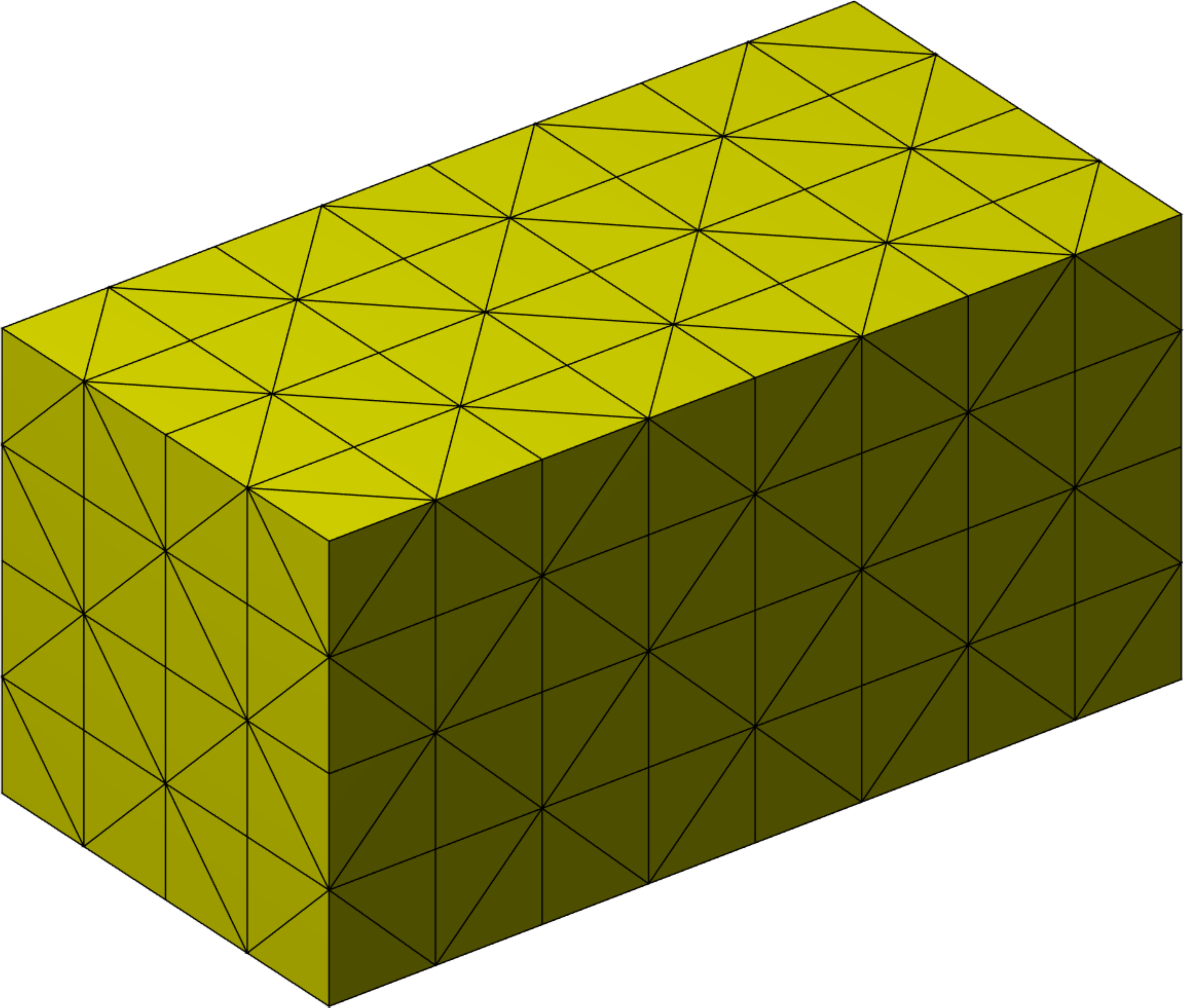}}

\caption{\label{fig:BackgroundMeshes}Different types of typical background
meshes. (a) shows the building block of the universal mesh in (b),
(c) shows the building block of the ``box mesh'' in (d).}
\end{figure}

We shall later discuss the issue of manipulating the background mesh
by moving nodes in order to improve the resulting surface elements
of the manifold mesh. We find that the concept of universal meshes
\cite{Rangarajan_2012a}, extended to three dimensions as shown in
Fig.~\ref{fig:BackgroundMeshes}(b), allows for a maximum of flexibility
in the node movements. But natural choices are also the block-type
meshes shown in Fig.~\ref{fig:BackgroundMeshes}(d). This is in particular
useful when the manifold of interest features straight edges and follows
from Eq.~(\ref{eq:BoundedIsosurfaceOmega}) with a block-type subdomain,
e.g. $\Omega=[-1,1]^{3}$.

\subsection{Reconstruction of $\Gamma_{\phi}$\label{XX_Reconstruction}}

\begin{figure}
\centering

\subfigure[top.~1, ref.~element]{\includegraphics[height=3cm]{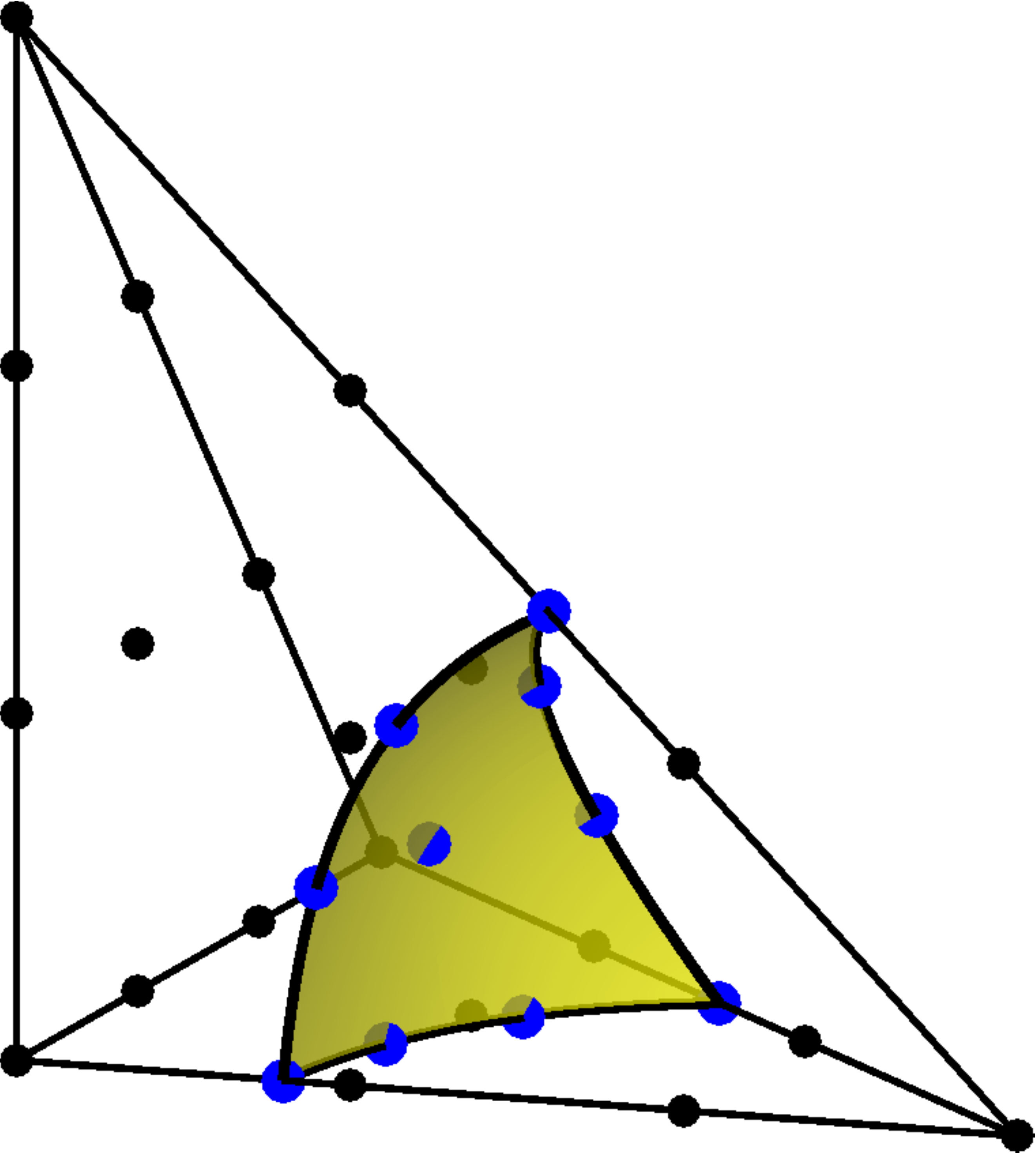}}\qquad\subfigure[top.~2, ref.~element]{\includegraphics[height=3cm]{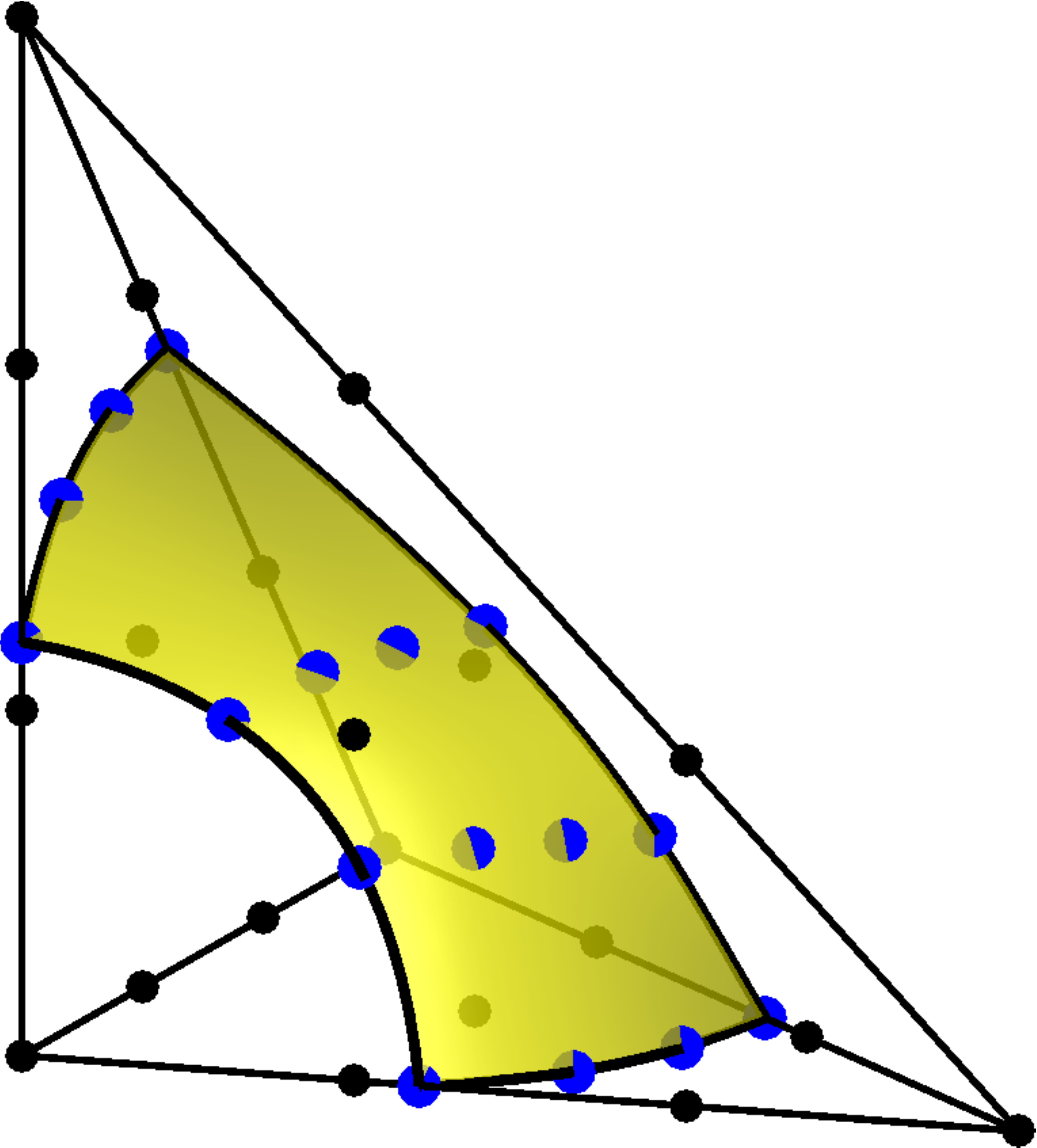}}\qquad\subfigure[top.~1, real element]{\includegraphics[height=3cm]{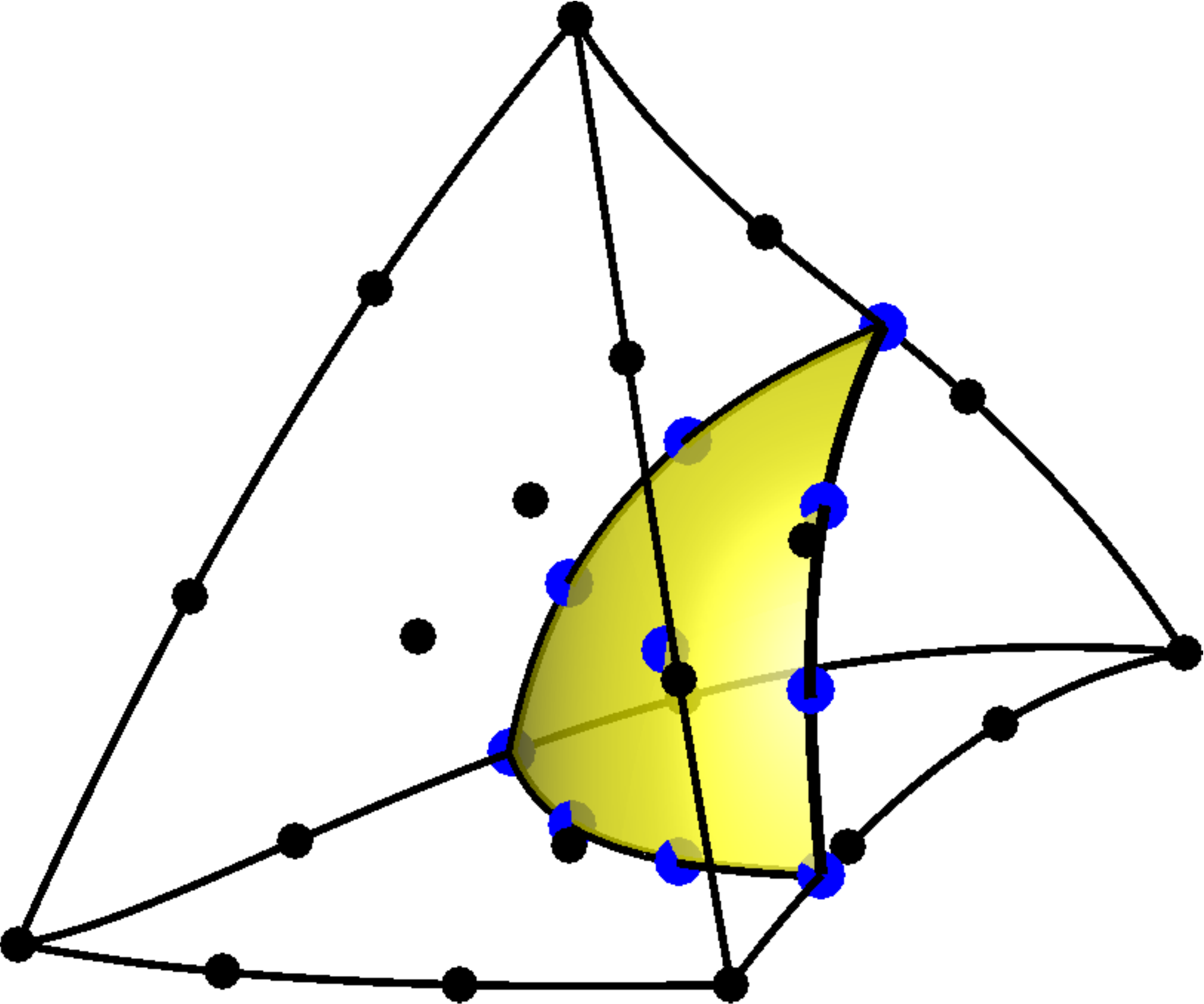}}\qquad\subfigure[top.~2, real element]{\includegraphics[height=3cm]{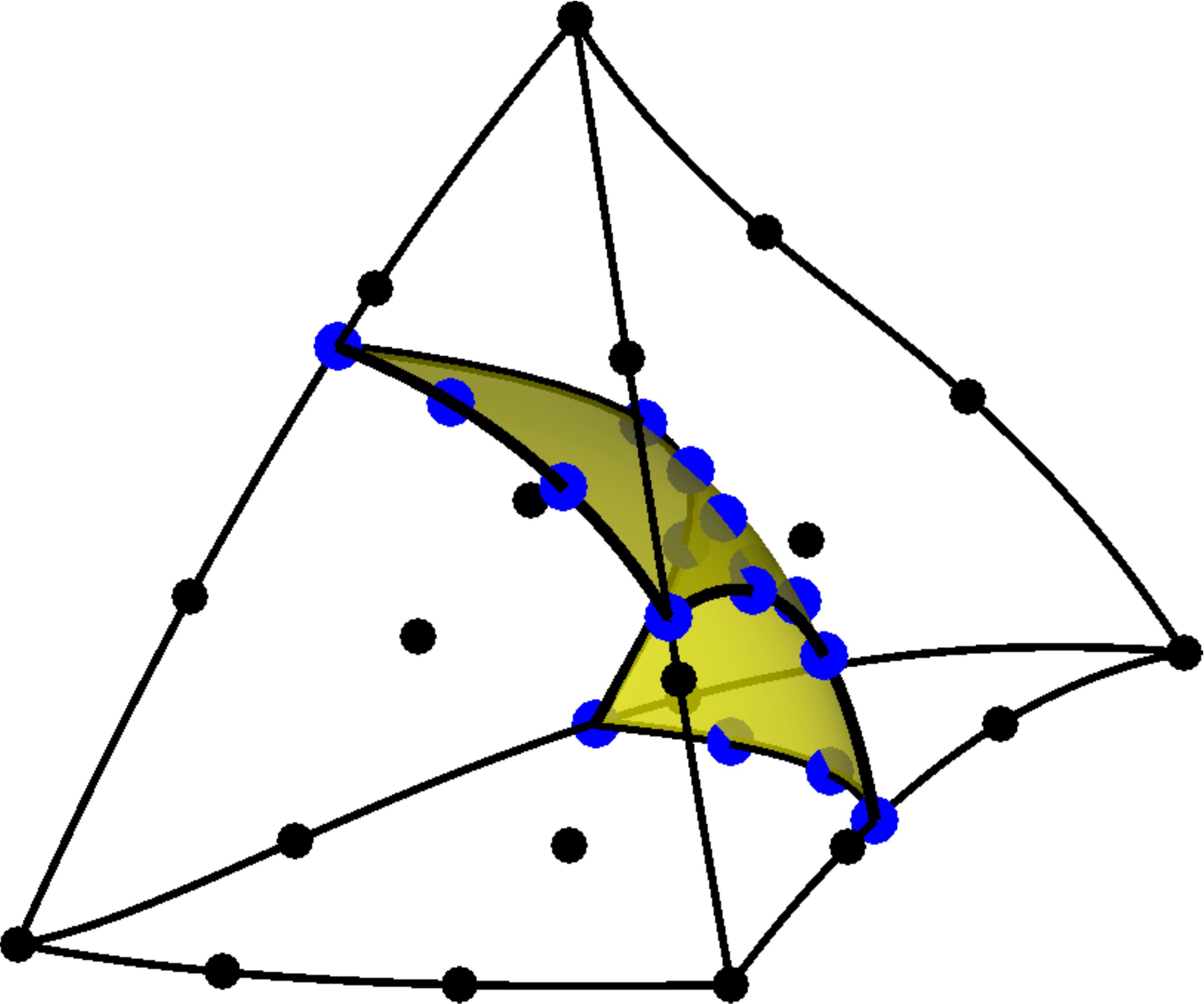}}

\caption{\label{fig:VisOverviewRecon}(a) and (b) show a triangular or quadrilateral
isosurface in the reference background element, respectively, which
is approximated by a surface element, (c) and (d) show the mapped
situation in a physical element of the background mesh.}
\end{figure}

In each element of the background mesh which is cut by the zero-isosurface
of $\phi$, the task is to identify element nodes of a higher-order
surface element on the isosurface. This has been described in detail
in the first part of this series of publications, see \cite{Fries_2016b}
and also \cite{Fries_2015a}. The procedure is shortly summarized
as: In each tetrahedral background element, it is confirmed that valid
level-set data is present which implies only two possible cut scenarios
of the isosurface: The isosurface is either of triangular shape, see
Fig.~\ref{fig:VisOverviewRecon}(a), then three edges are cut. Otherwise,
four edges are cut and it is of quadrilateral shape, see Fig.~\ref{fig:VisOverviewRecon}(b).
Depending on the topology, customized search algorithms are evaluated
in the reference background element. 

First, element nodes on the faces of the tetrahedron are identified
leading to the edge nodes of the sought surface element approximating
the isosurface. This is, in fact, realized by reconstructing a higher-order
line element in a two-dimensional, triangular reference element with
the level-set data from the respective face nodes of the tetrahedron.
This line is mapped to the corresponding face of the tetrahedron.
The inner nodes of the sought surface element are then determined
in the reference tetrahedron wherefore the search algorithm also makes
use of the position of the edge nodes determined beforehand. Finally,
the determined surface element may be mapped to the physical background
element using the isoparametric concept. An example is seen in Fig.~\ref{fig:VisOverviewRecon}(c)
and (d) for the two different topological cases of a triangular or
quadrilateral surface element, respectively.

Obviously, applying the described procedure in \emph{all} cut background
elements leads to a set of higher-order surface elements representing
an accurate approximation of $\Gamma_{\phi}$. Lateron, a mesh suitable
for a finite element analysis of the BVP on the manifold will be automatically
generated from this element set, see Section \ref{XX_MeshGeneration}.
An important property of such a mesh is that it is at least $C_{0}$-continuous
so that there are no gaps in the representation of the manifold. It
is interesting to note that the reconstruction described above does,
in fact, \emph{not }assure this property without further considerations.

\begin{figure}
\centering

\subfigure[reference element]{\includegraphics[height=5cm]{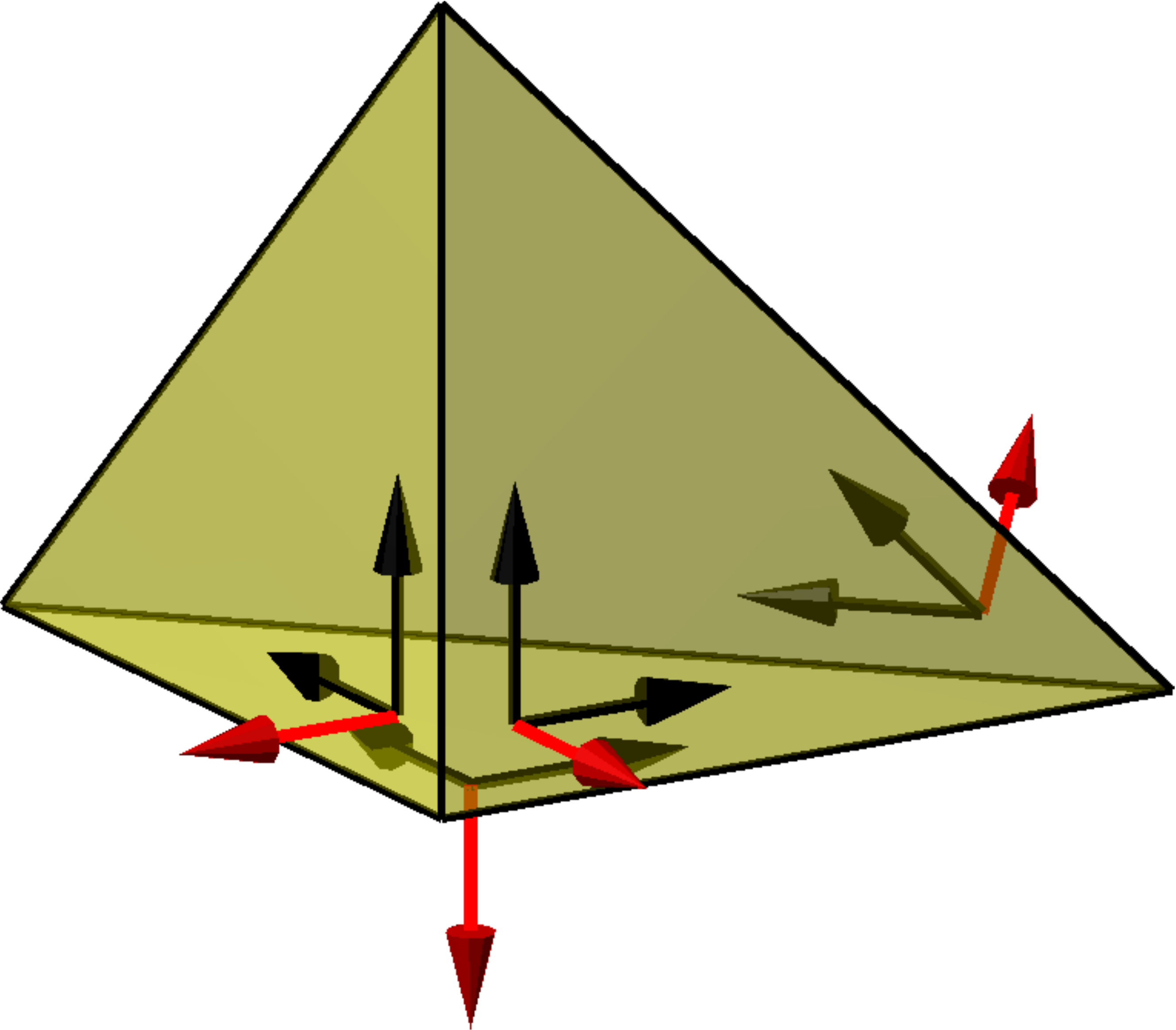}}\qquad\qquad\subfigure[two real elements]{\includegraphics[height=5cm]{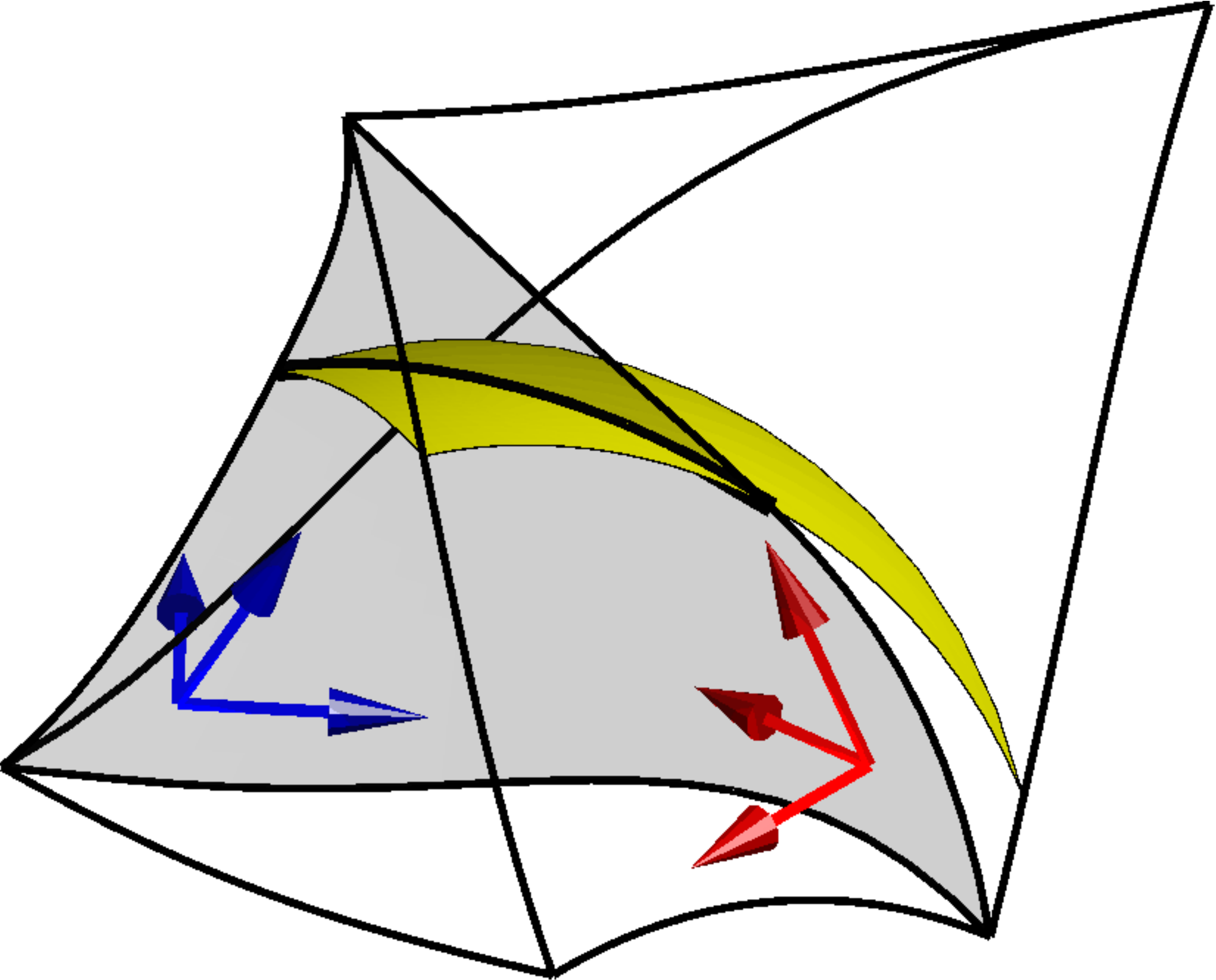}}

\caption{\label{fig:PlotFaceCoordSystems}(a) shows the convention for the
face coordinate systems, normal vectors (red) point outwards, (b)
shows two physical elements where the node numberings imply non-matching
coordinate systems on the shared face. This will lead to slightly
different nodes on the shared edge and, hence, to a small gap in the
(yellow) surface elements. Therefore, for the reconstruction on the
faces, a permutation may be required to ensure that the face coordinate
systems match.}
\end{figure}

The reason for this is found in the orientation of the face coordinate
systems, see Fig.~\ref{fig:PlotFaceCoordSystems}: Only when the
same coordinate system is used for both tetrahedra sharing a face,
this will lead to exactly the same edge nodes of the surface elements.
Otherwise, only the corner nodes of the surface elements perfectly
match. Fortunately, there is a simple solution to the problem: The
face nodes of a tetrahedron which are extracted for the reconstruction
of the edge nodes of the sought surface element should be permuted
such that they imply the same coordinate systems on the shared face.
The resulting set of surface elements is then, as desired, a globally
$C_{0}$-continuous, higher-order approximation of $\Gamma_{\phi}$.
This is an important difference to the reconstructions discussed in
\cite{Fries_2015a,Fries_2016b}.

\subsection{Restriction of $\Gamma_{\phi}$ by additional level-set functions
$\psi^{i}$\label{XX_Restriction}}

\begin{figure}
\centering

\subfigure[ex.~1, reconstruction]{\includegraphics[height=4cm]{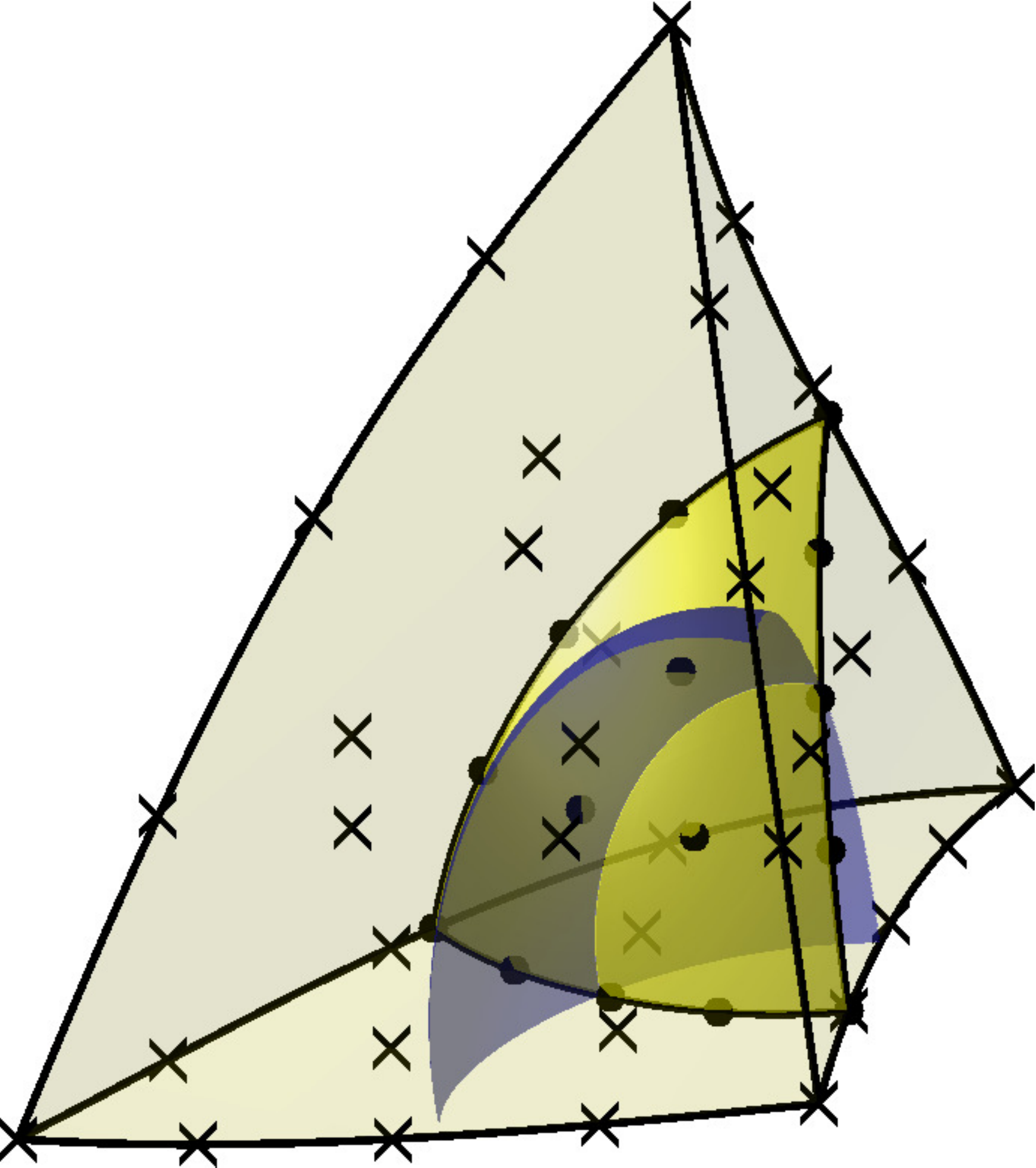}}\qquad\subfigure[ex.~2, reconstruction]{\includegraphics[height=4cm]{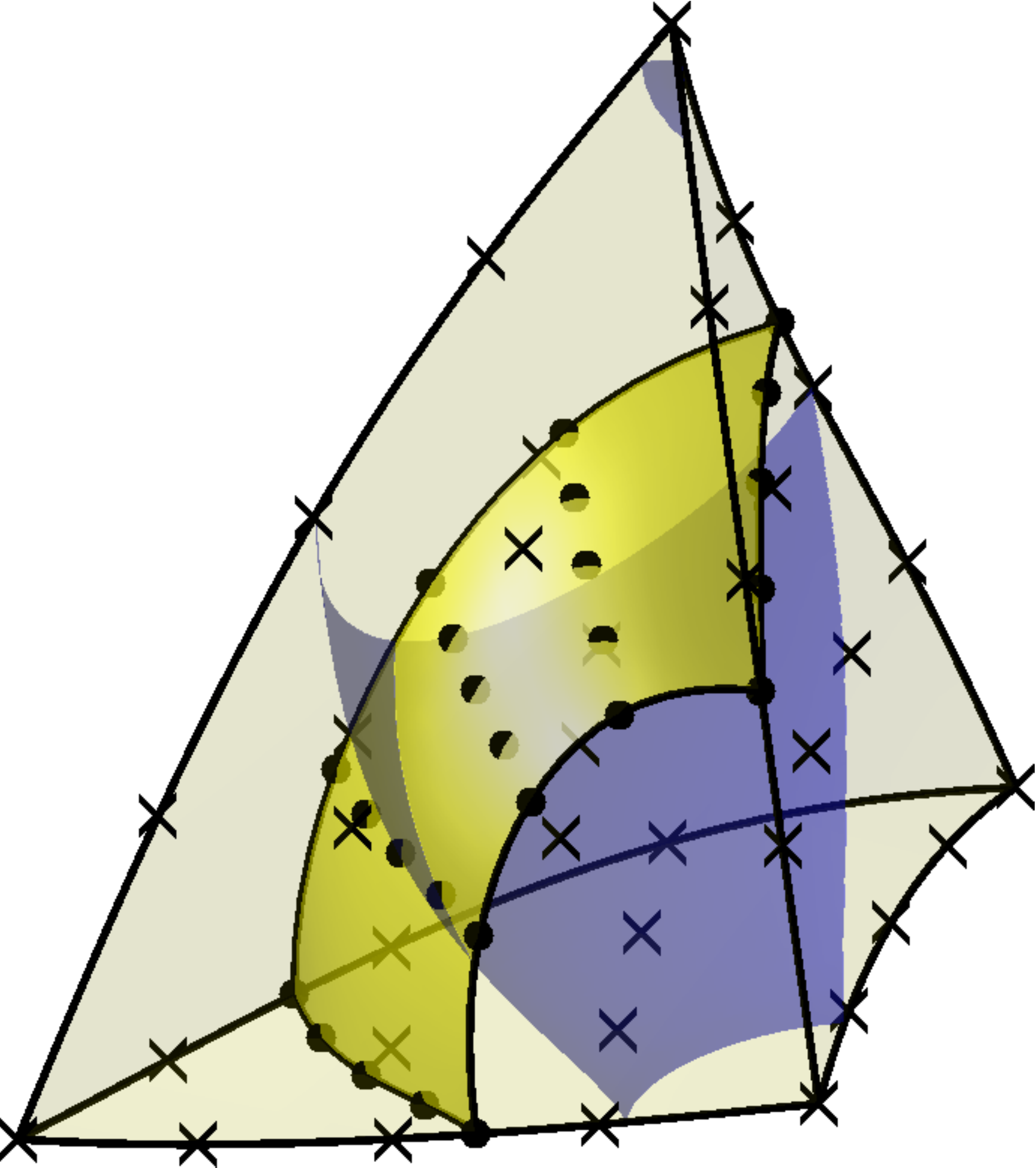}}\qquad\subfigure[ex.~3, reconstruction]{\includegraphics[height=4cm]{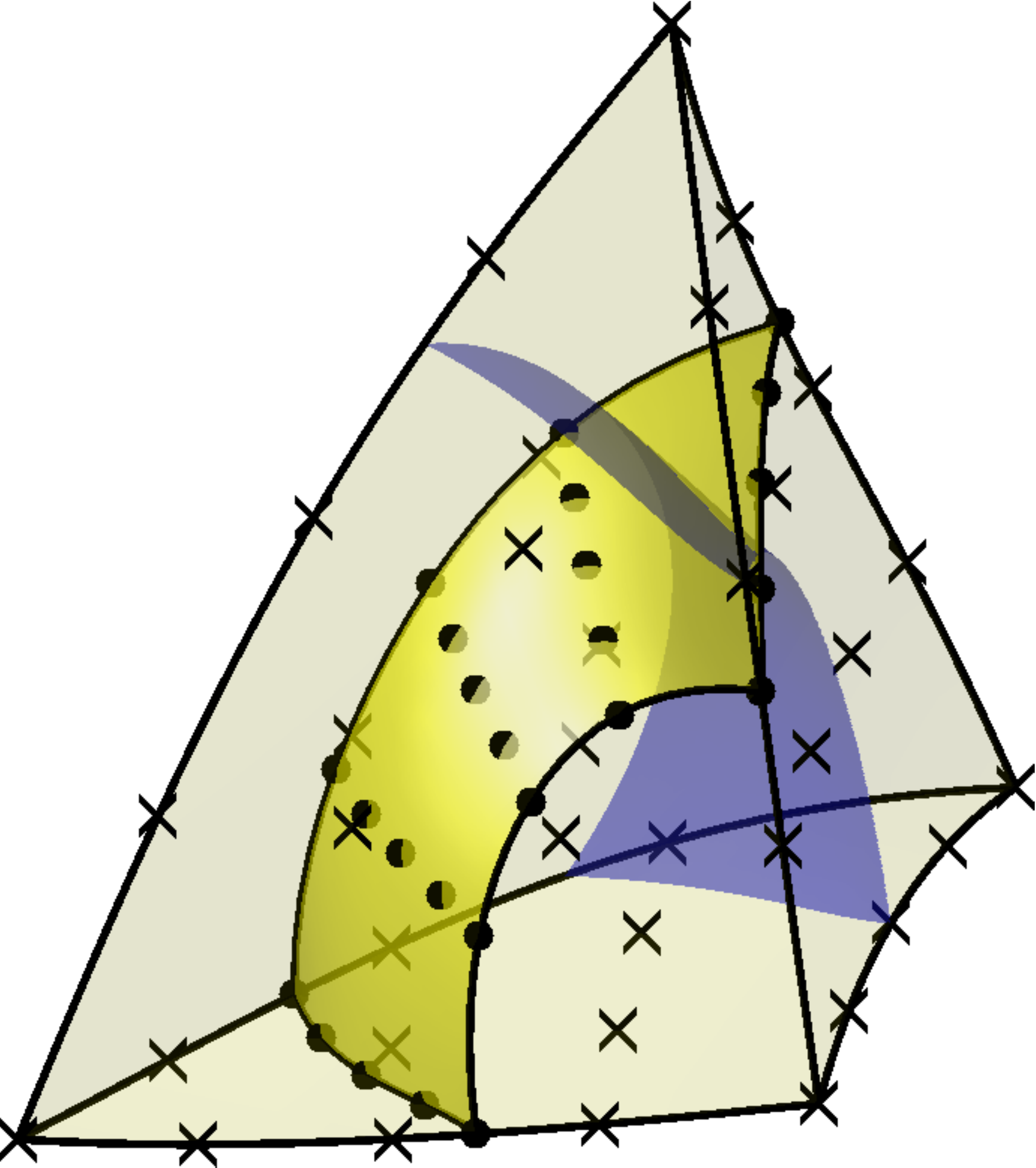}}

\subfigure[ex.~1, decomposition]{\includegraphics[height=4cm]{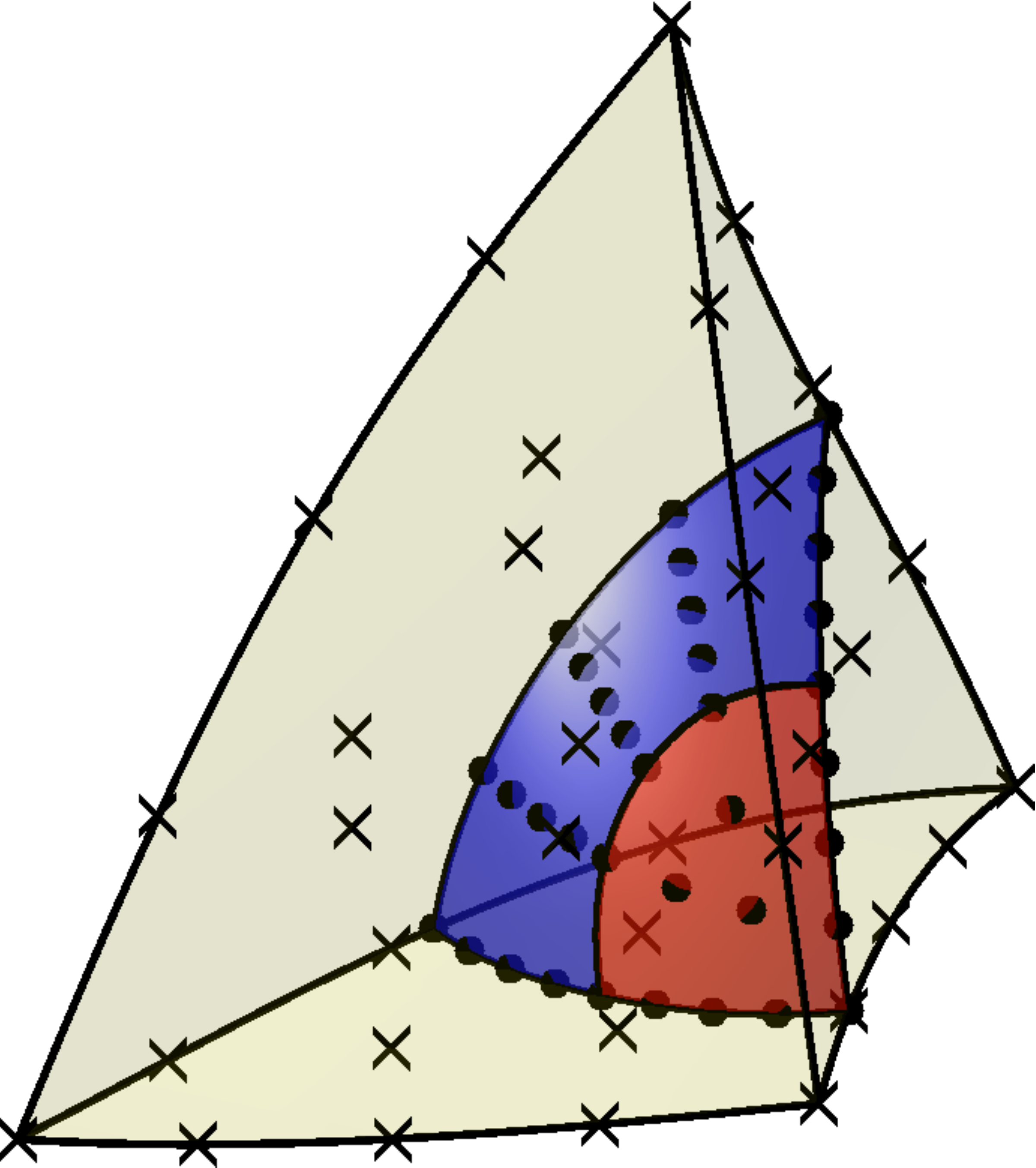}}\qquad\subfigure[ex.~2, decomposition]{\includegraphics[height=4cm]{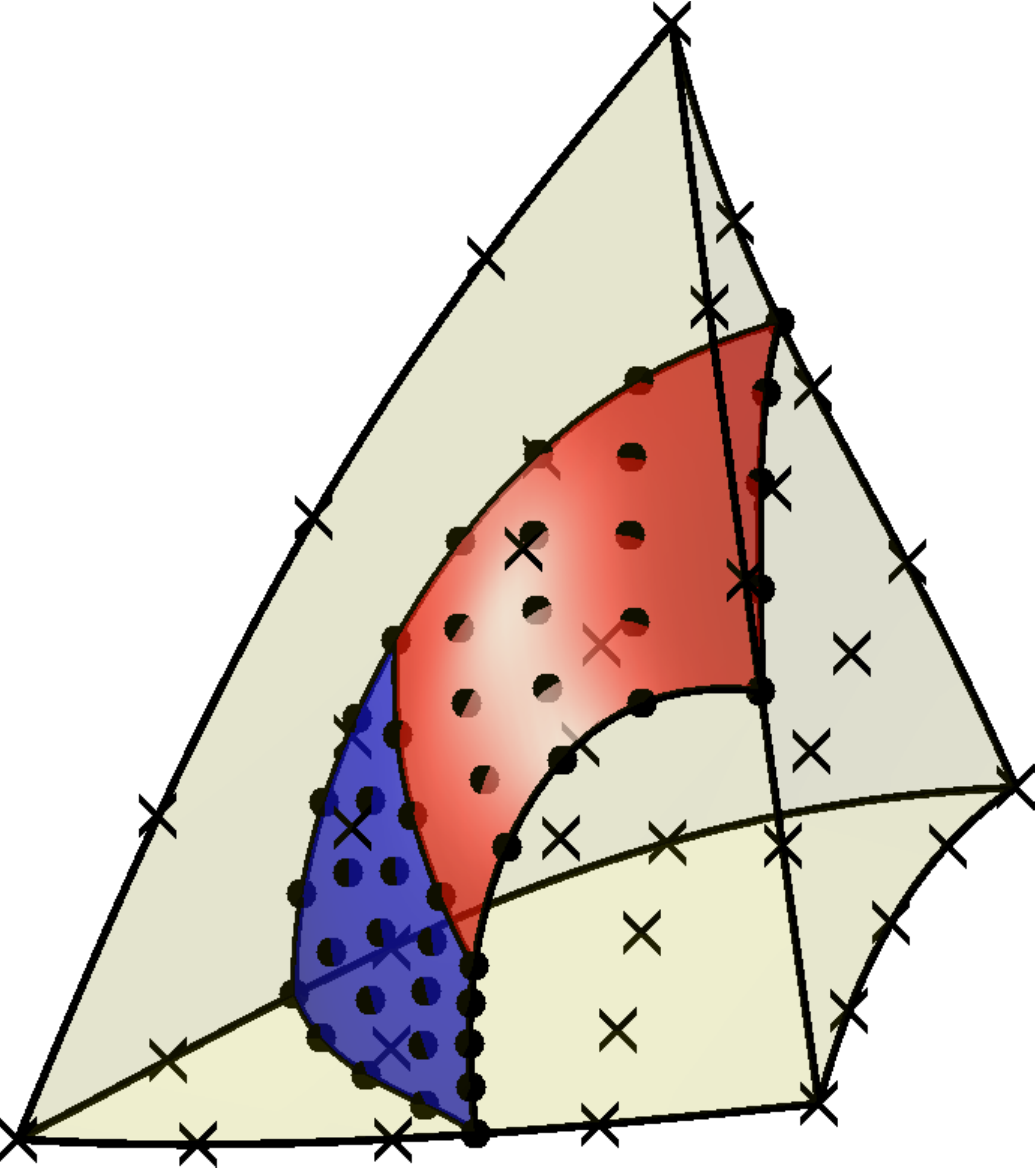}}\qquad\subfigure[ex.~3, decomposition]{\includegraphics[height=4cm]{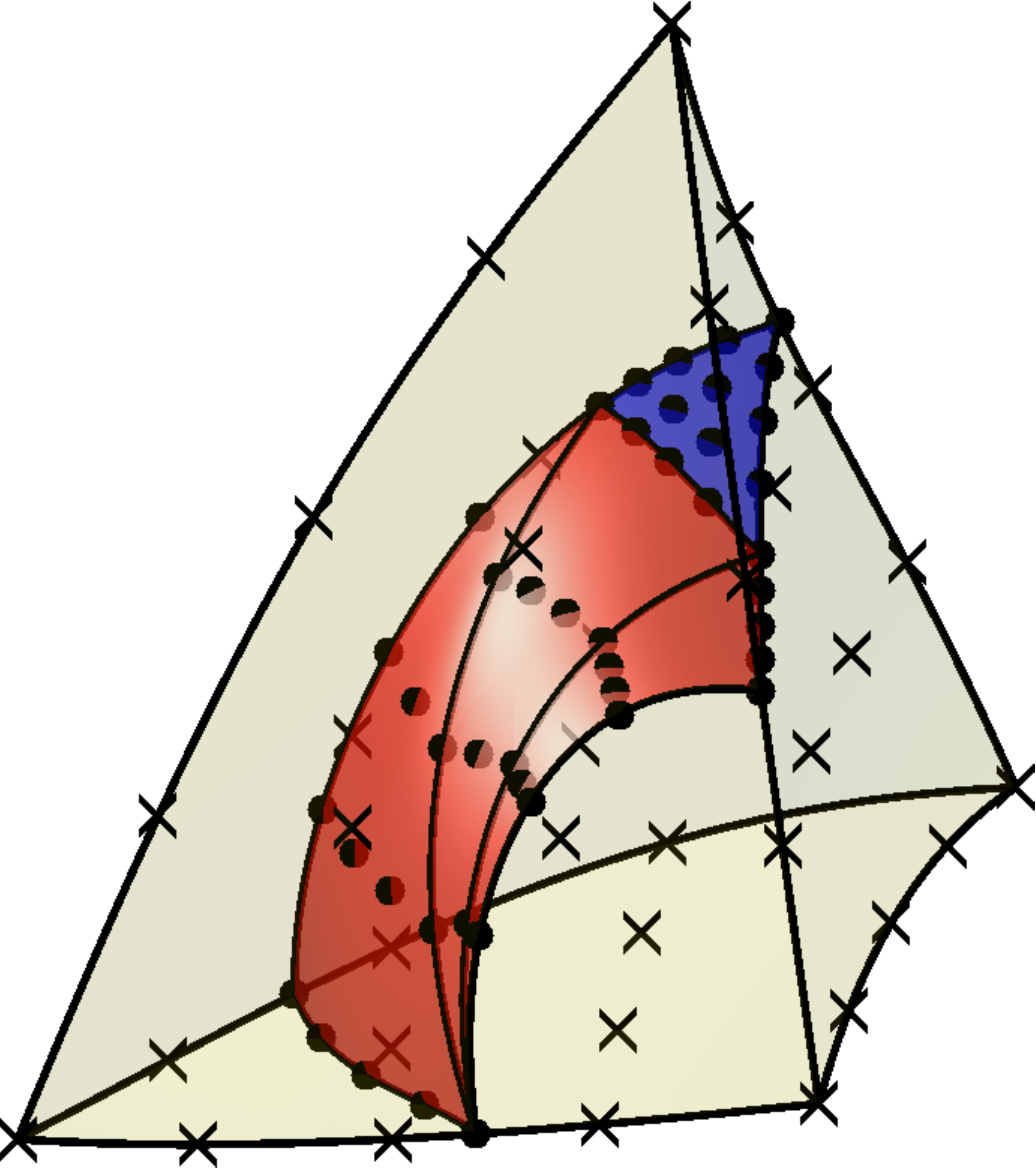}}

\subfigure[ex.~1, 2D ref.~element]{\includegraphics[width=3cm]{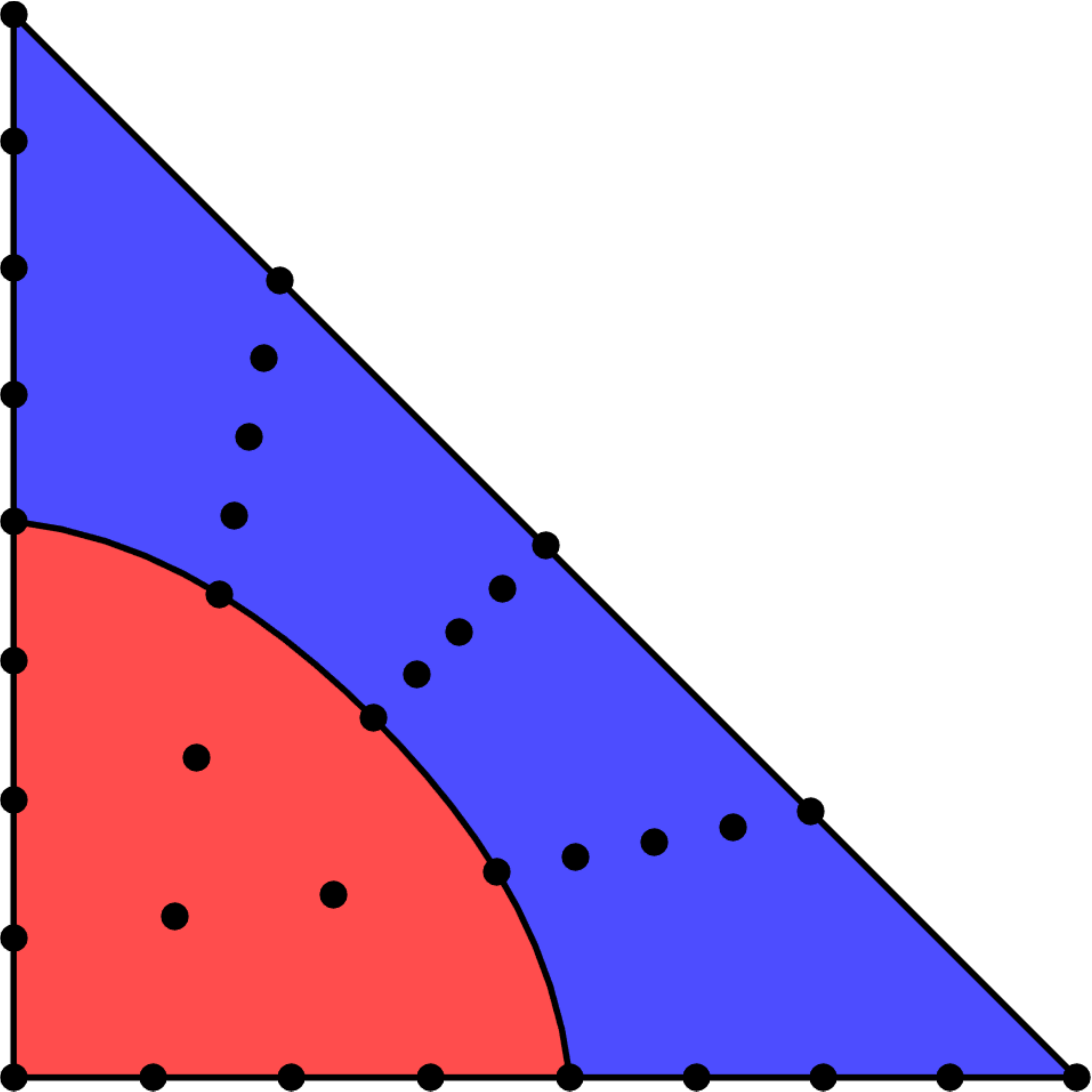}}\qquad\qquad\subfigure[ex.~2, 2D ref.~element]{\includegraphics[width=3cm]{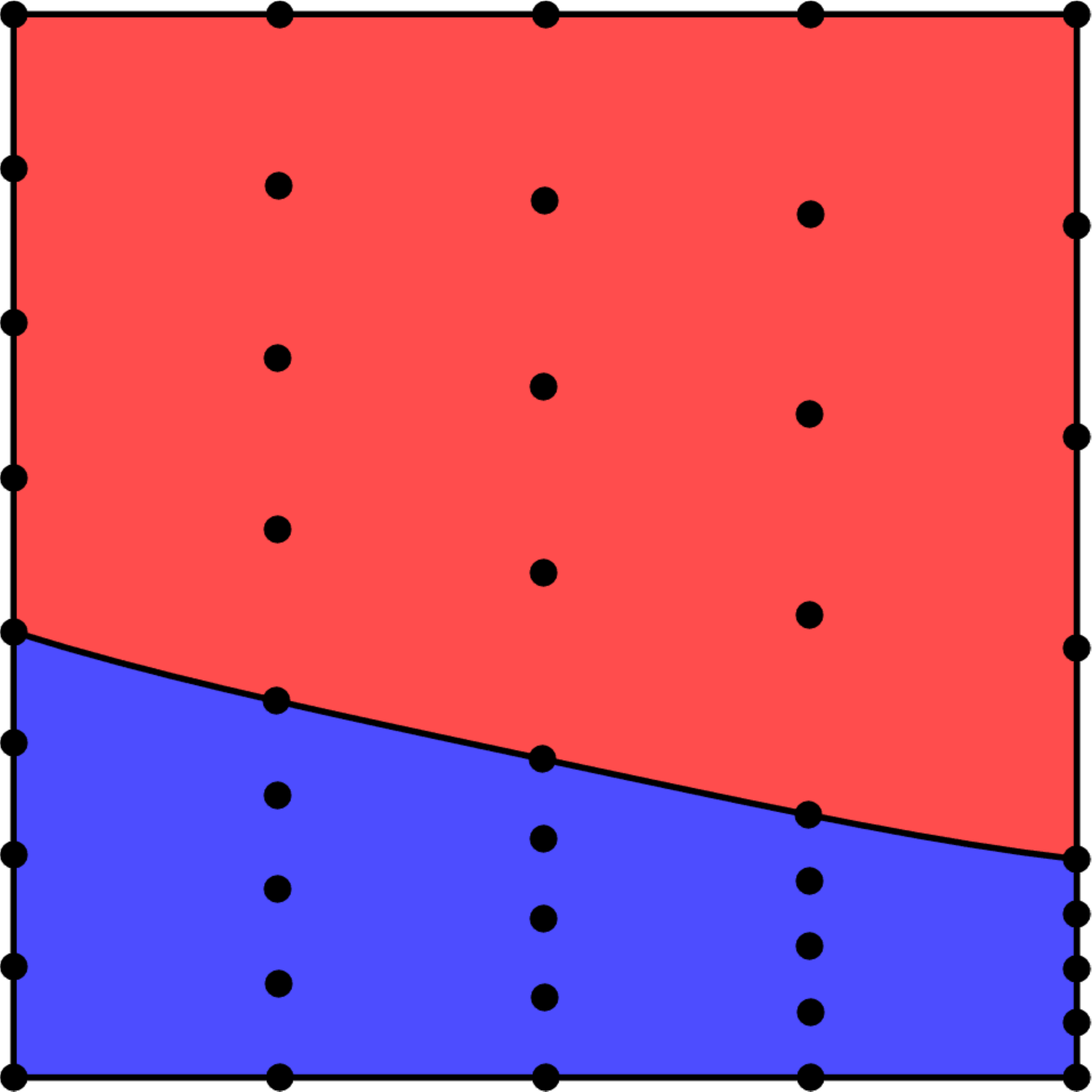}}\qquad\qquad\subfigure[ex.~3, 2D ref.~element]{\includegraphics[width=3cm]{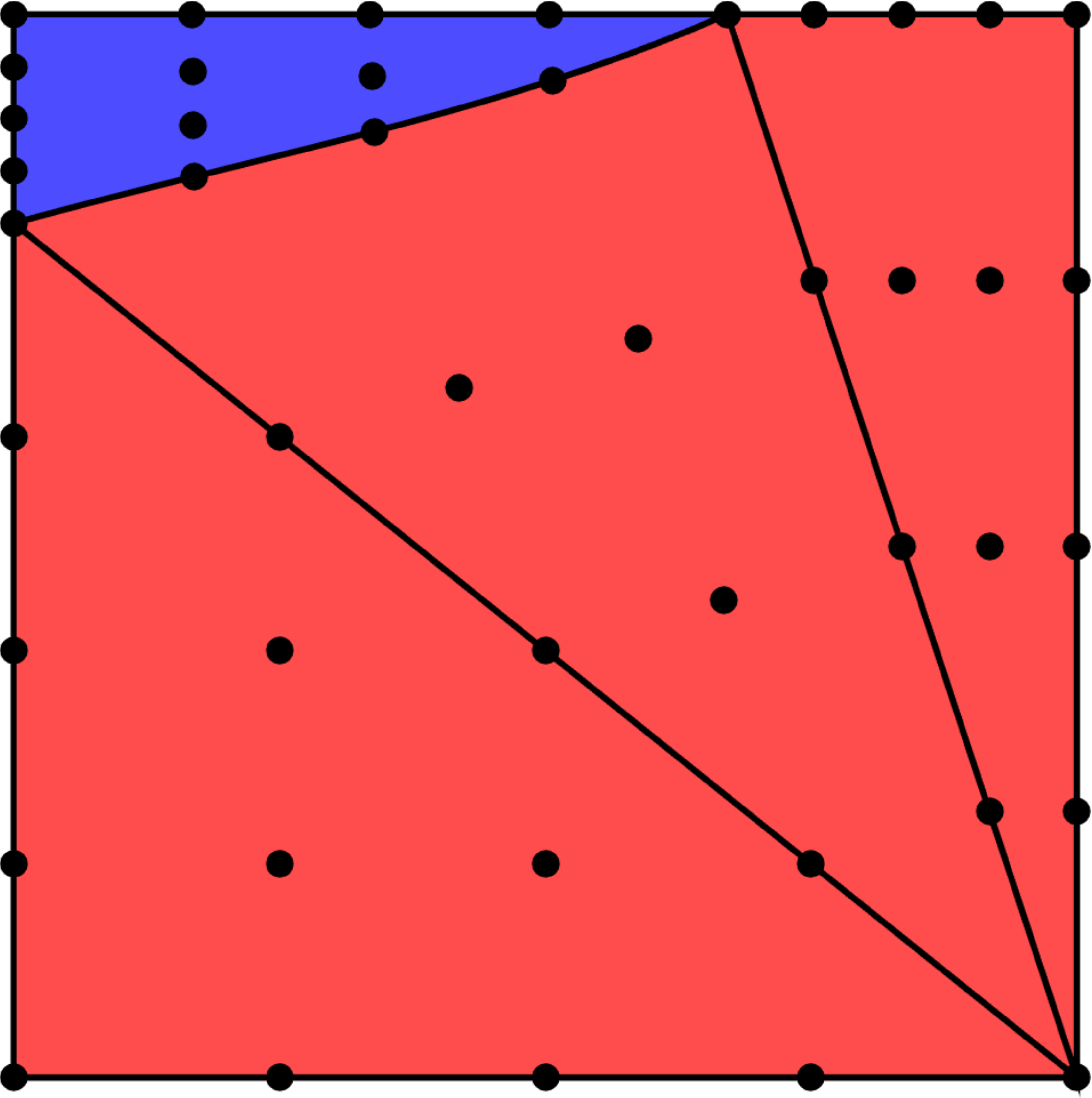}}

\caption{\label{fig:PlotRestriction}Example 1 features a triangular surface
element, examples 2 and 3 a quadrilateral surface element, (a) to
(c) show the reconstructed surface element (yellow) and the zero isosurface
of $\psi^{k}$ (blue) cutting through the elements, (d) to (f) show
the decomposed surface elements with respect to $\psi^{k}$ (red for
$\psi^{k}<0$, blue for $\psi^{k}>0$), (g) to (i) show that the decomposition
is realized in 2D reference elements.}
\end{figure}

The focus is on one background element cut by the zero-level set of
$\phi\left(\vek x\right)$. A triangular or quadrilateral surface
element, approximating the zero isosurface, has been succesfully reconstructed.
Assume that the background element and the surface element are also
cut by another level set function $\psi^{k}\left(\vek x\right)$,
see Fig.~\ref{fig:PlotRestriction}(a) to (c) for different examples.
A decomposition of the surface element is then sought into surface
elements on the two sides, see Fig.~\ref{fig:PlotRestriction}(d)
to (f) for the same examples, respectively. 

For this decomposition, it is necessary to first interpolate the nodal
values of $\psi^{k}\left(\vek x\right)$ in the background element
(crosses) at the nodes of the reconstructed surface element (dots).
Then, the decomposition is carried out in two-dimensional reference
elements with these nodal values, see Fig.~\ref{fig:PlotRestriction}(g)
to (i). As discussed in detail in \cite{Fries_2015a,Fries_2016b},
different cut scenarios may be faced: In triangles, always one of
the corners is on the other side than the other two and one sub-triangle
and one sub-quadrilateral are obtained, see Fig.~\ref{fig:PlotRestriction}(g).
In quadrilaterals, either two neighboring nodes are on the other side
leading to two sub-quadrilaterals, see Fig.~\ref{fig:PlotRestriction}(h),
or only one node, leading to four triangles, see Fig.~\ref{fig:PlotRestriction}(i).
The decomposed two-dimensional elements are mapped to the cut tetrahedral
element using an isoparametric map implied by the original surface
element (from the reconstruction of $\phi$). Typically, considering
Eq.~(\ref{eq:BoundedIsosurfaceMult}), only the subelements on the
negative side of $\psi^{k}$ are needed and the others are neglected,
this may be called restriction.

\begin{figure}
\centering

\subfigure[]{\includegraphics[height=4cm]{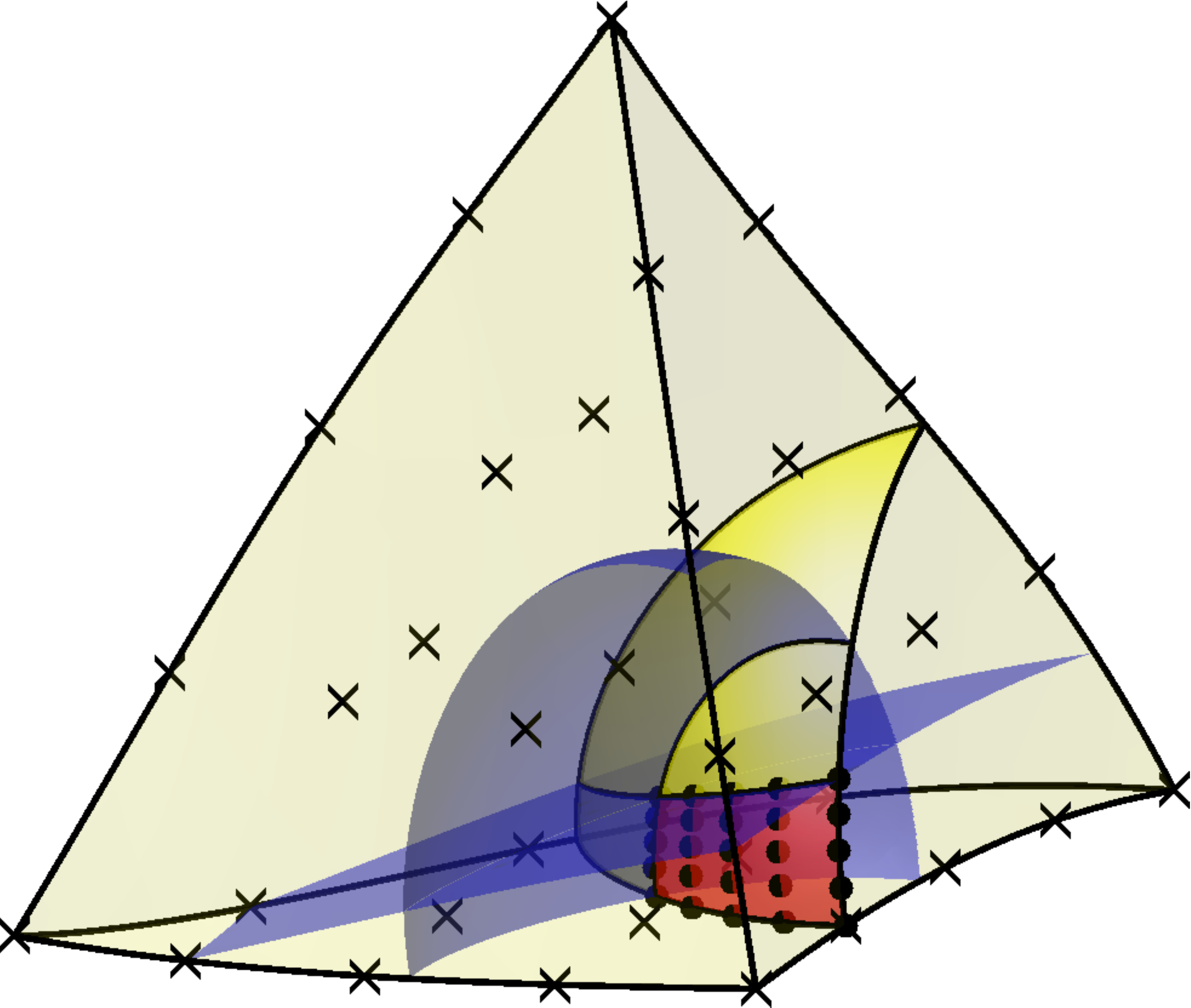}}\qquad\subfigure[]{\includegraphics[height=4cm]{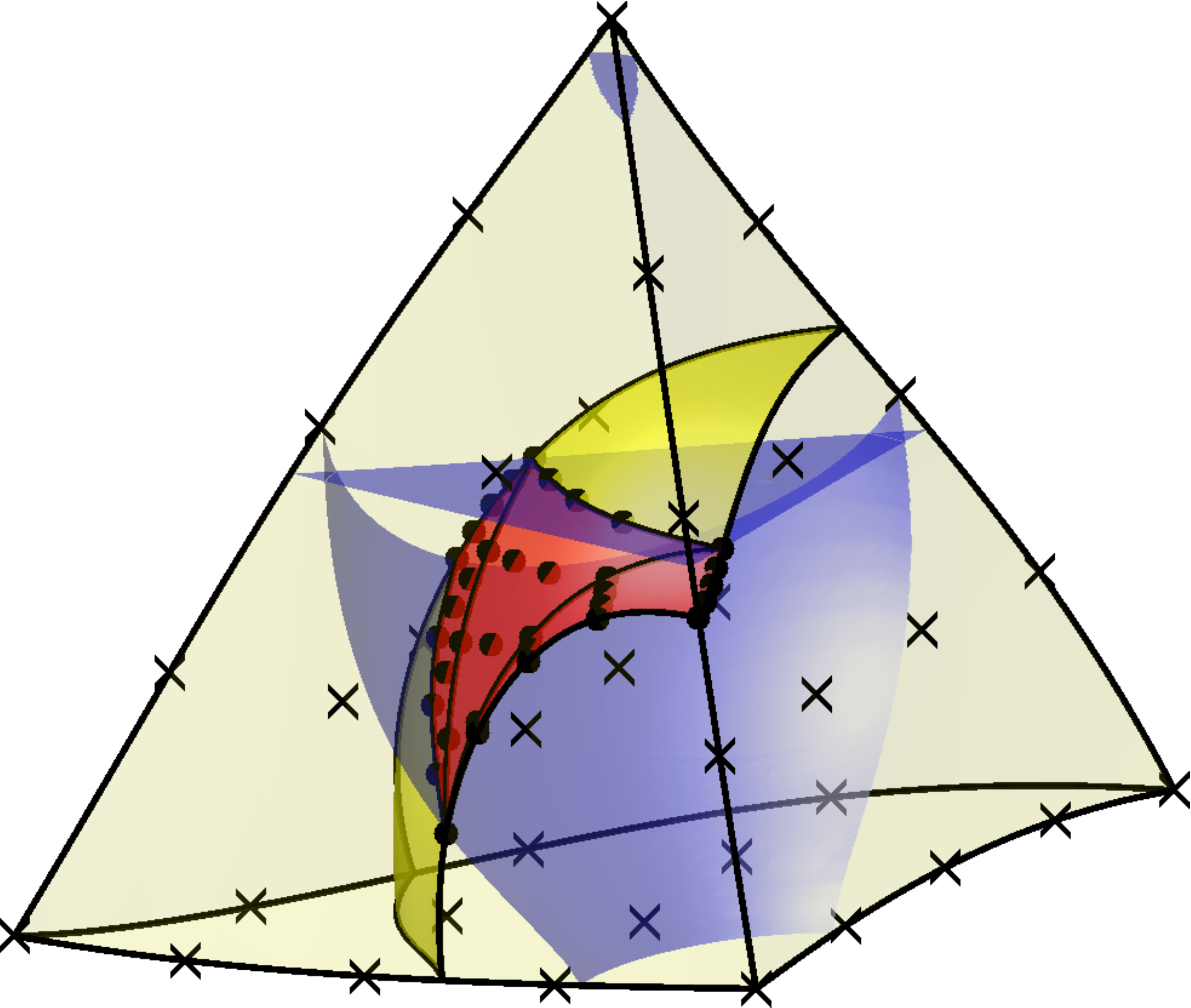}}\qquad\subfigure[]{\includegraphics[height=4cm]{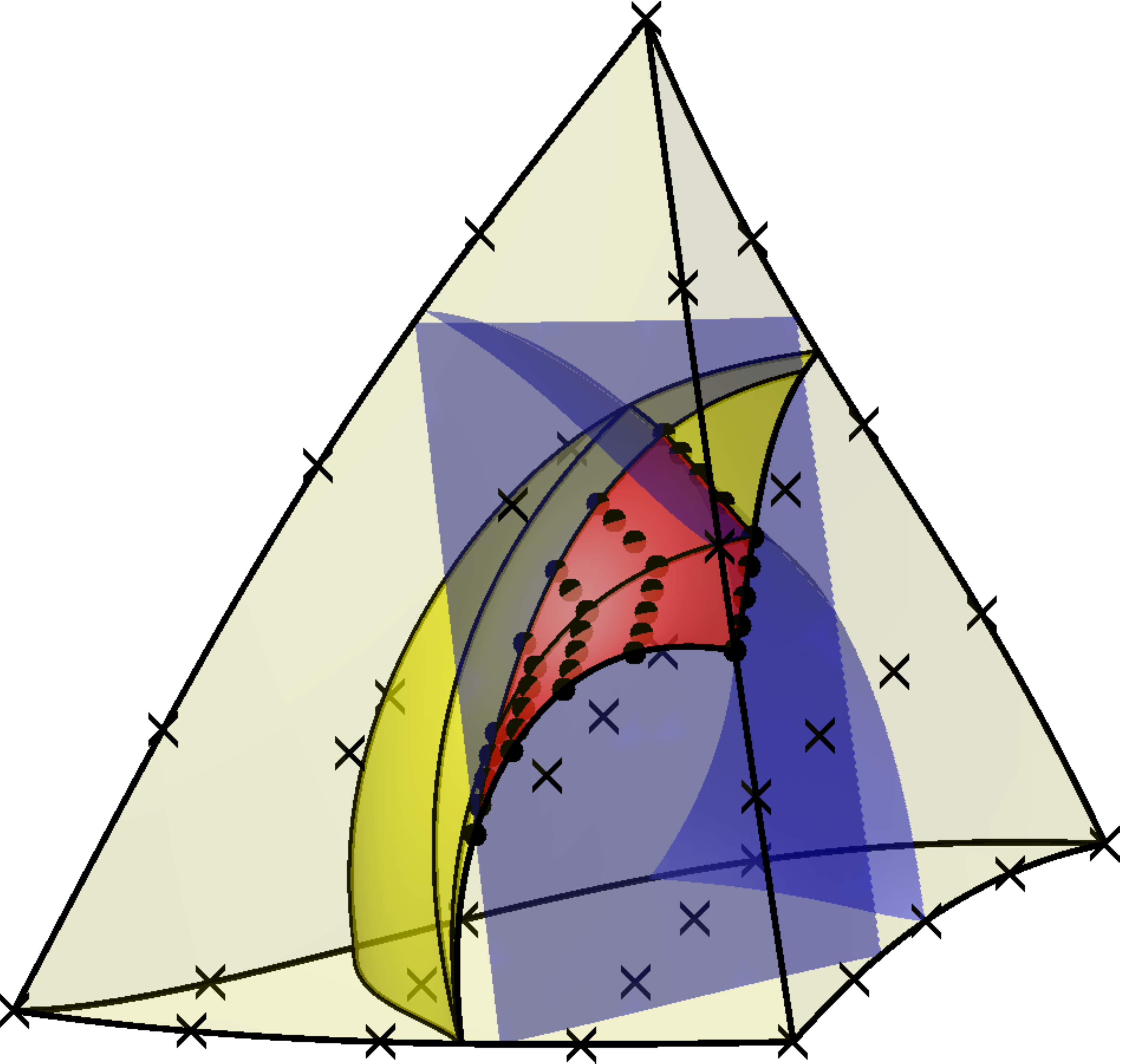}}

\caption{\label{fig:PlotRestrictionMult}Some examples for level-set functions
$\psi^{k}$ further restricting the reconstructed surface element,
the corresponding isosurfaces are shown in blue. The final resulting
reconstructed elements which approximate $\Gamma_{\phi}$ with $\psi^{k}<0$
are shown in red.}
\end{figure}

The procedure is very similar when several level-set functions $\psi^{i}\left(\vek x\right)$
cut the same element, possibly implying one or more corners in the
boundary of the manifold. The starting point is then a set of surface
elements resulting from a previous restriction with respect to other
level-set functions. Then, at each surface element, the procedure
from above is repeated: The nodal values of the current level-set
function at the current surface element are interpolated from the
tetrahedral nodes. Then, the decomposition is carried out in the two-dimensional
reference elements and mapped to the current surface element. Finally,
the resulting sub-elements are restricted to those with $\psi^{i}\left(\vek x\right)<0$.
See Fig.~\ref{fig:PlotRestrictionMult} for some examples. It is
emphasized that the reconstructions are always realized in \emph{reference
}elements (before mapping them to the situation in global coordinates):
For the level-set function $\phi$, these are 3D reference tetrahedra
and for the level-set functions $\psi^{i}$, these are reference triangular
or quadrilateral elements. This enables a very efficient and robust
implementation as the reference elements always feature straight edges
and planar faces, respectively. On the contrary, it seems impossible
to reconstruct \emph{directly} in general, possibly curved, physical
elements.

\subsection{Mesh generation\label{XX_MeshGeneration}}

Starting point is a set of higher-order surface elements which has
been reconstructed from $\Gamma_{\phi}$ and restricted to the manifold
of interest $\Gamma$ using the additional level-set functions $\psi^{k}$.
In order to approximate a BVP on the manifold, a finite element mesh
has to be generated from the element set. The aim is to extract a
set of unique nodes, introduce a global node numbering and set up
a connectivity matrix defining the elements. It is noted that the
generated element set from above contains triangular and quadrilateral
elements, which leads to a mixed finite element mesh. Of course, the
decompositon of quadrilateral elements to triangular elements (or
\emph{vice versa}) is possible if meshes of one element type are preferred.

There are two obvious alternatives for the automatic mesh generation
from the element set: In alternative 1, the nodes are first sorted
with respect to their coordinates and then repeated nodes, such as
they appear on the element edges and corners, are unified and associated
with a node number. This may be implemented very efficiently, however,
it requires some user-defined threshold characterizing a distance
within which nodes are considered equal. Because the overall set of
elements may contain extremely small elements (e.g. when the zero-isosurface
cuts very close to a corner node of a background element) this may
interfere with the threshold leading to invalid data. For example,
all element nodes may be unified to one node in an extreme situation.

Therefore, we prefer alternative 2 where the connectivity information
of the background mesh is employed. One may easily determine the neigboring
elements of each tetrahedron. For the surface elements in the reconstructed
set of elements, the information of the corresponding background elements
is stored, respectively. One may then directly compute which nodes
on a shared face of two tetrahedra must match. The same holds for
the nodes on a shared edge of several tetrahedra. No distance computations
are needed in this approach. It is clear that an advanced implementation
is possible where the computation of nodes with the same coordinates
(at the edges or corners of the surface elements belonging to different
background elements) is completely avoided.

\subsection{Manipulation of the background mesh\label{XX_MeshManipulation}}

Some example meshes for complex manifolds are seen in Fig.~\ref{fig:MeshManifolds}.
It is obvious that awkward element shapes result from the automatic
meshing. Very small inner angles may occur and the size of neighboring
elements can vary strongly so that the shape regularity is an issue.
This may lead to the assumption that the elements are not suitable
for the analysis, even less in a higher-order context. However, as
shown in \cite{Olshanskii_2009a,Olshanskii_2012a} in a low-order
context, only a few manipulations ensure that the ``maximal angle
condition'' is fulfilled and the spectral condition number is bounded
uniformly. More generally, it has recently been shown in several publications
by the groups around Hansbo, Burman and Reusken, e.g.~\cite{Burman_2015a,Olshanskii_2012a,Reusken_2014a},
that optimal results for BVPs on implicitly defined manifolds are
possible with such meshes. They use finite element shape functions
based on the background mesh which are only evaluated on the manifold;
these approaches are labelled TraceFEM or CutFEM. They introduce tailored
stabilization terms which avoid negative effects resulting from the
ill-shaped elements. We shall report on the use of such methods in
this context in a future publication. Herein, we wish to avoid the
detailed (mathematical) discussion of these approaches.

A simple and intuitive attempt is to manipulate the background mesh
such that the resulting reconstructed elements remain shape regular.
In particular, the area of the surface elements shall be bounded from
below. When the condition numbers in a finite element analysis obtained
by a tailored, manufactured surface mesh, ideally featuring the same
element shapes and areas in the whole mesh, are compared to those
of the automatically generated meshes from above, it is seen that
there is typically an increase by the factor of the largest to the
smallest element area $A_{\mathrm{max}}/A_{\mathrm{min}}$. Because
$A_{\mathrm{max}}$ typically scales with the size of the (very regular)
background elements, this clearly motivates the need to bound $A_{\mathrm{min}}$
from below.

For simplicity, we focus on the set of elements resulting from a reconstruction
with respect to $\phi$, i.e.~the set of elements approximating $\Gamma_{\phi}$
on a finite background mesh. An algorithm is proposed which moves
the nodes of the background mesh away from the zero-level set of $\phi$.
Basically, this ensures that the zero-isosurface is not too close
to the corner nodes of the background mesh and, consequently, that
the surface elements are not too small. Only the nodes in a close
band around the zero isosurface are moved.

Although the elements of the background mesh may be arbitrarily curved,
it is quite useful to enforce straight edges for simplicity of the
node manipulation. The algorithm for the node movement is described
for \emph{one} concrete node at position $\vek x^{'}$ in the background
mesh; the same holds for all other nodes. It may be seen as a fix-point
iteration and the following steps are repeated until all nodes have
been moved sufficiently away from the zero-level set.
\begin{enumerate}
\item Approximate the distance of the node at $\vek x'$ to the zero-level
set. When the level-set function $\phi$ is given in analytical form,
this is realized by a Newton-Raphson procedure, 
\[
\vek x_{i+1}=\vek x_{i}-\frac{\phi\left(\vek x_{i}\right)}{\left\Vert \nabla\phi\left(\vek x_{i}\right)\right\Vert }\cdot\nabla\phi\left(\vek x_{i}\right),
\]
using $\vek x_{0}=\vek x'$ as the starting point. This algorithm
converges to some position $\vek x^{\star}$ on the zero-level set
which is \emph{not} necessarily the closest point. However, it is
found that $\vek x^{\star}$ is approaching the closest point quickly
the closer $\vek x^{'}$ is to the zero-level set. Because the node
manipulation is only realized for nodes which are close to the zero-level
set, it turned out to be sufficient to continue with the computed
position $\vek x^{\star}$. The estimated signed distance is 
\[
D\left(\vek x'\right)=\left\Vert \vek x'-\vek x^{\star}\right\Vert \cdot\mathrm{sign}\left(\phi\left(\vek x'\right)\right)
\]
and the direction 
\[
\vek d\left(\vek x'\right)=\frac{\left(\vek x'-\vek x^{\star}\right)}{D\left(\vek x^{'}\right)}
\]
is computed as well; there holds $\left\Vert \vek d\left(\vek x'\right)\right\Vert =1$.
In the case that $\phi$ is only given at the nodes of a higher-order
background mesh (and no analytic information is available), one may
estimate the distance through a linear reconstruction of the zero
isosurface which is much cheaper than the desired higher-order reconstruction
performed after the node manipulation.
\item Move the node provided that it is sufficiently close to the zero-level
set, i.e.~when $\left|D\left(\vek x'\right)\right|<D_{\mathrm{crit}}$.
$D_{\mathrm{crit}}$ scales with the resolution of the background
mesh and we typically set $D_{\mathrm{crit}}=3\cdot h$ with $h$
being some characteristic element length. The node is moved as
\[
\vek x''=\vek x'+q\left(D,D_{\mathrm{crit}},D_{\mathrm{step}}\right)\cdot\vek d\left(\vek x'\right)\cdot\mathrm{sign}\left(D\left(\vek x'\right)\right)
\]
with
\[
q\left(D,D_{\mathrm{crit}},D_{\mathrm{step}}\right)=D_{\mathrm{step}}\cdot\left(1-\frac{\left|D\right|}{D_{\mathrm{crit}}}\right).
\]
That is, the maximum distance a node is moved in one step is $D_{\mathrm{step}}$
which is often set as $0.1\cdot h$. Obviously, $q$ scales linearly
with the distance of a node from the zero-level set.
\item Set $\vek x'=\vek x''$ and repeat the two steps from above until
$\left|D\left(\vek x'\right)\right|>D_{\mathrm{min}}$ where, typically,
$D_{\mathrm{min}}=0.25\cdot h$.
\end{enumerate}
This algorithm works very fast and scales linearly with the number
of nodes in the background mesh. It is controlled by the three parameters
which are summarized as: $D_{\mathrm{crit}}$ and $D_{\mathrm{min}}$
define local regions around the zero-level sets (``bands''): Nodes
within the band controlled by $D_{\mathrm{crit}}$ are considered
for movement and $D_{\mathrm{min}}$ defines the band within which
no nodes shall be present. $D_{\mathrm{step}}$ defines the maximum
moving length in one iteration which is scaled linearly in the band
controlled by $D_{\mathrm{crit}}$. Obviously, we need $D_{\mathrm{step}}<D_{\mathrm{min}}<D_{\mathrm{crit}}$.

It is mentioned that in our first approximations of BVPs on manifolds
with the automatically generated surface meshes, we were surprised
by the accuracy of the results even without any mesh manipulations.
Actually, it seems that as long as direct solvers are employed (in
our case Matlab's backslash solver), there is a large range of condition
numbers still leading to accurate results even in the frame of systematic
convergence studies. Nevertheless, the node manipulations are used
by default in the numerical studies presented in Section \ref{X_NumericalResults}.
It is an important advantage that they also render recursive refinements
of the background elements unnecessary as long as the resolution of
the background elements is reasonably adjusted to the complexity of
the zero-level sets. The effect of the proposed node manipulation
is now demonstrated in one, two an three dimensions.

\subsubsection{Node manipulation in 1D}

\begin{figure}
\centering

\subfigure[$\phi ,D,d$]{\includegraphics[width=3.5cm]{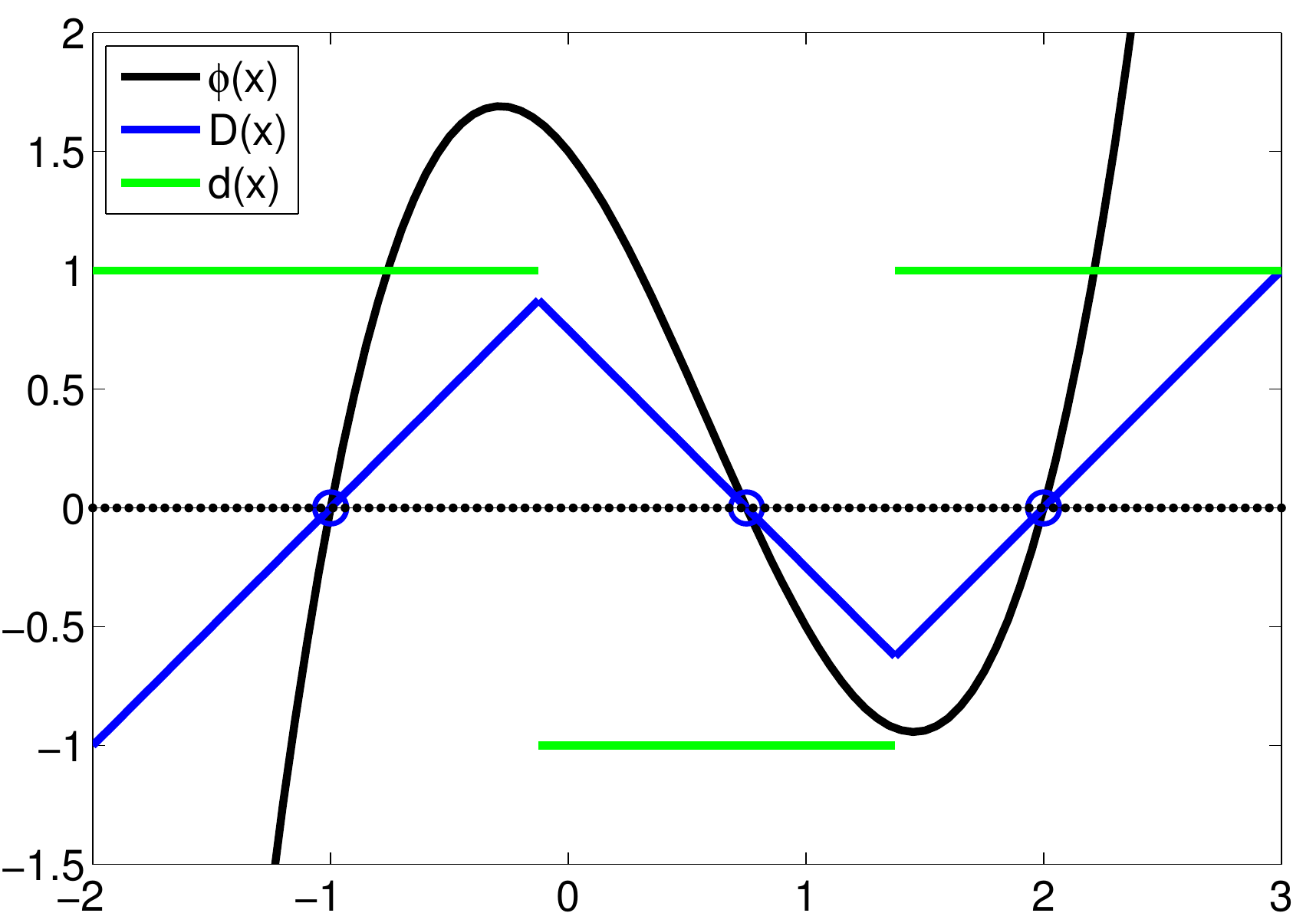}}\quad\subfigure[Example 1]{\includegraphics[width=3.5cm]{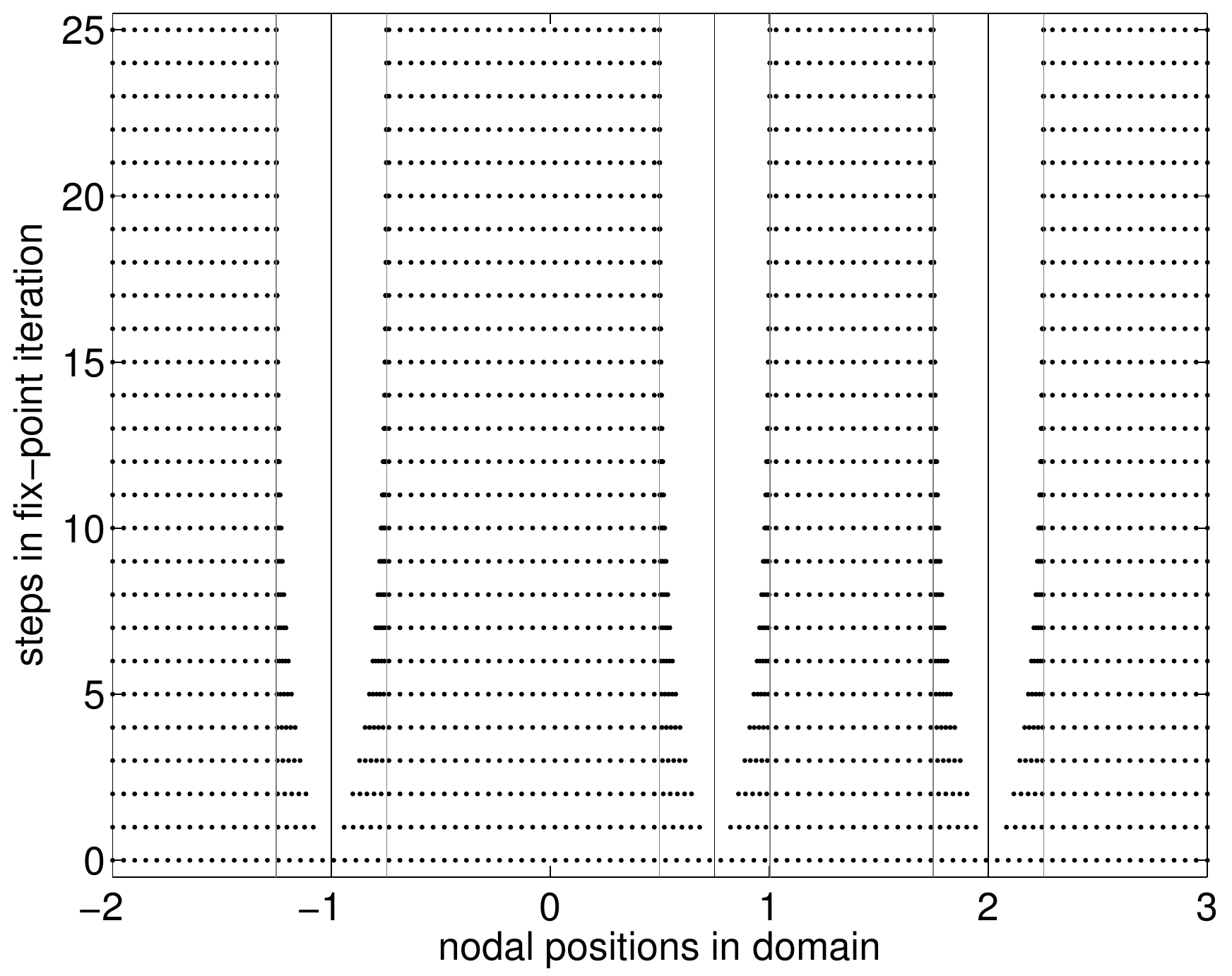}}\quad\subfigure[Example 2]{\includegraphics[width=3.5cm]{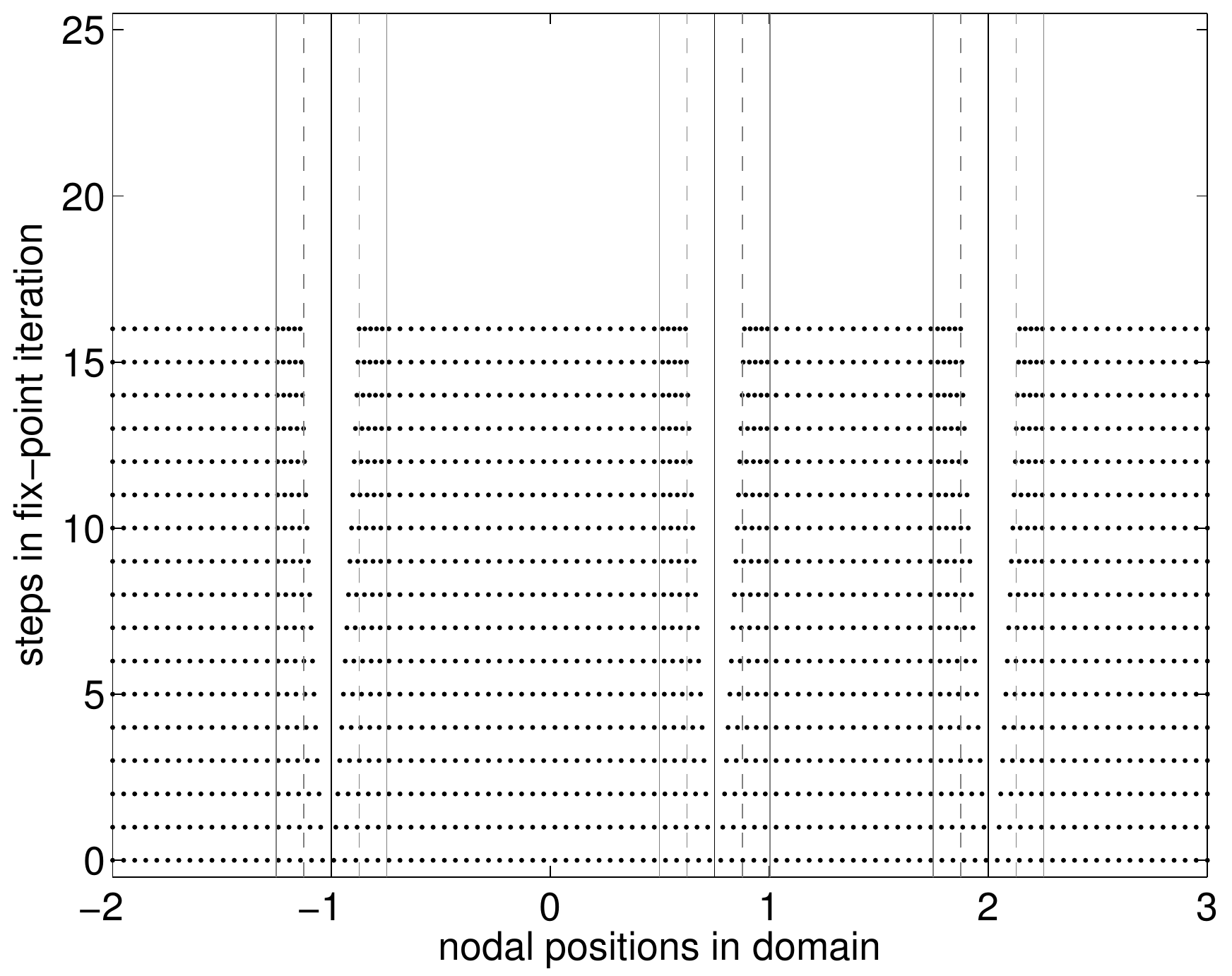}}\quad\subfigure[Example 3]{\includegraphics[width=3.5cm]{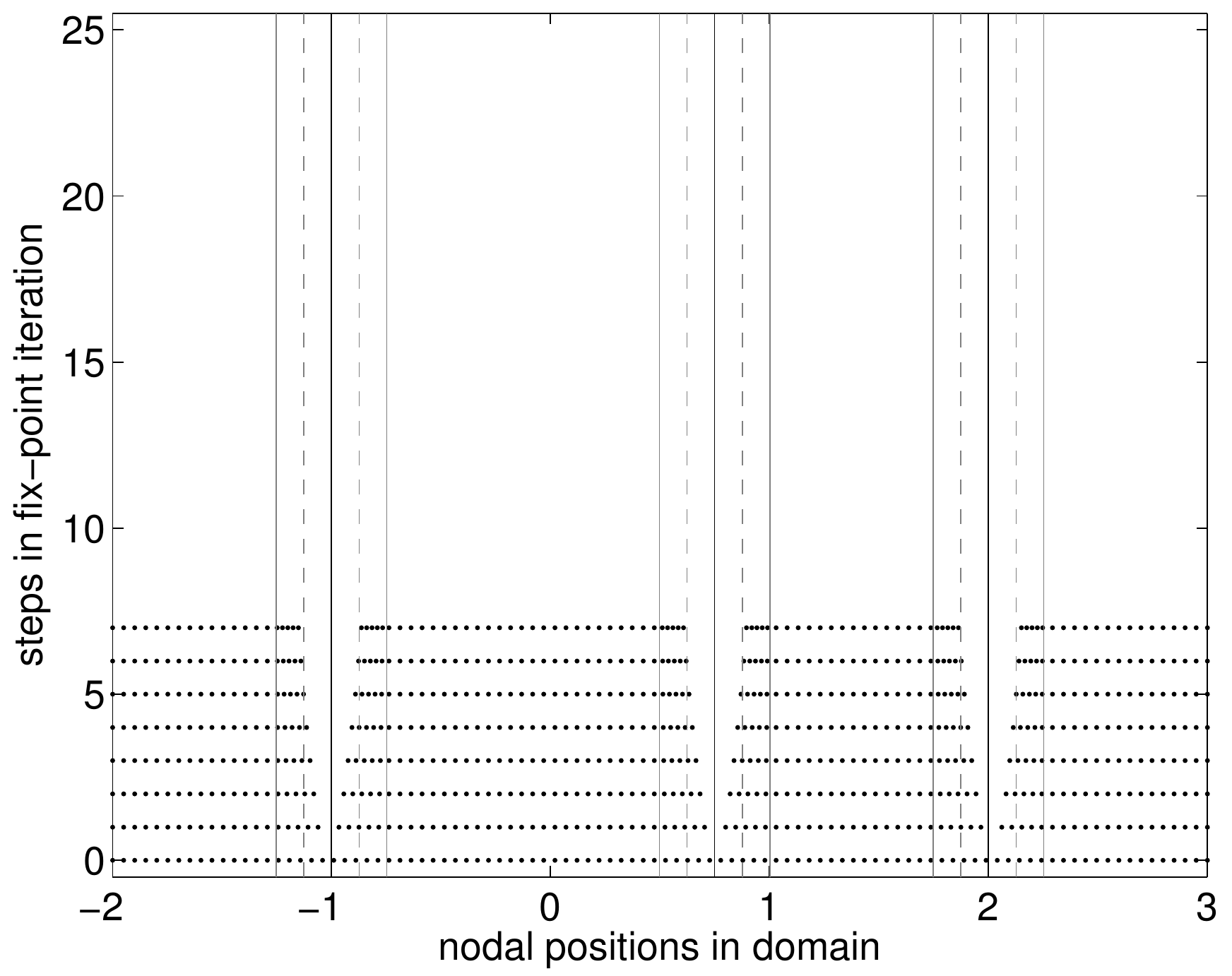}}

\caption{\label{fig:PlotMeshManip1d}Examples for node movements in 1D, (a)
shows a level-set function $\phi$ (black), the corresponding signed-distance
function $D$ (blue) and its gradient $d$ (green), (b) to (d) show
how the nodes within $D_{\mathrm{crit}}$ are moved away from the
roots during the fix-point iteration.}
\end{figure}

Consider the one-dimensional example shown in Fig.~\ref{fig:PlotMeshManip1d}.
In Fig.~\ref{fig:PlotMeshManip1d}(a), a level-set function $\phi$
is plotted in black and the corresponding signed-distance function
$D$, which stores the (signed) distance to the closest root of $\phi$,
is shown in blue. On the $x$-axis, there are regularly distributed
nodes marked as blacked dots, the distance between two nodes is $h$.
The aim is to move these nodes away from the roots (blue circles).
In Fig.~\ref{fig:PlotMeshManip1d}(b), we set $D_{\mathrm{crit}}=D_{\mathrm{min}}=5\cdot h$
and $D_{\mathrm{step}}=h$. It is then seen that after about 10 steps
of the fix-point iteration, the nodes in the region $D_{\mathrm{crit}}$
are squeezed to the boundary between the band scaled by $D_{\mathrm{crit}}$
and the rest of the domain. The resulting nodal distribution leads
to undesired, locally concentrated background elements. Therefore,
we set $D_{\mathrm{crit}}=5\cdot h$ and $D_{\mathrm{min}}=\nicefrac{1}{2}D_{\mathrm{crit}}$
next. Results with $D_{\mathrm{step}}=\nicefrac{1}{5}h$ and $\nicefrac{1}{2}h$
are shown in Figs.~\ref{fig:PlotMeshManip1d}(c) and (d), respectively.
Obviously, the larger the value for $D_{\mathrm{step}}$, the faster
the nodes are moved out of the band controlled by $D_{\mathrm{min}}$.
Nevertheless, in more than one dimension, setting $D_{\mathrm{step}}$
considerably smaller than $D_{\mathrm{min}}$ allows for curved node
paths throughout the nodal movements leading to smoother element distributions
in the background meshes.

\subsubsection{Node manipulation in 2D and 3D}

\begin{figure}
\centering

\subfigure[$\phi_1$, coarse]{\includegraphics[width=5cm]{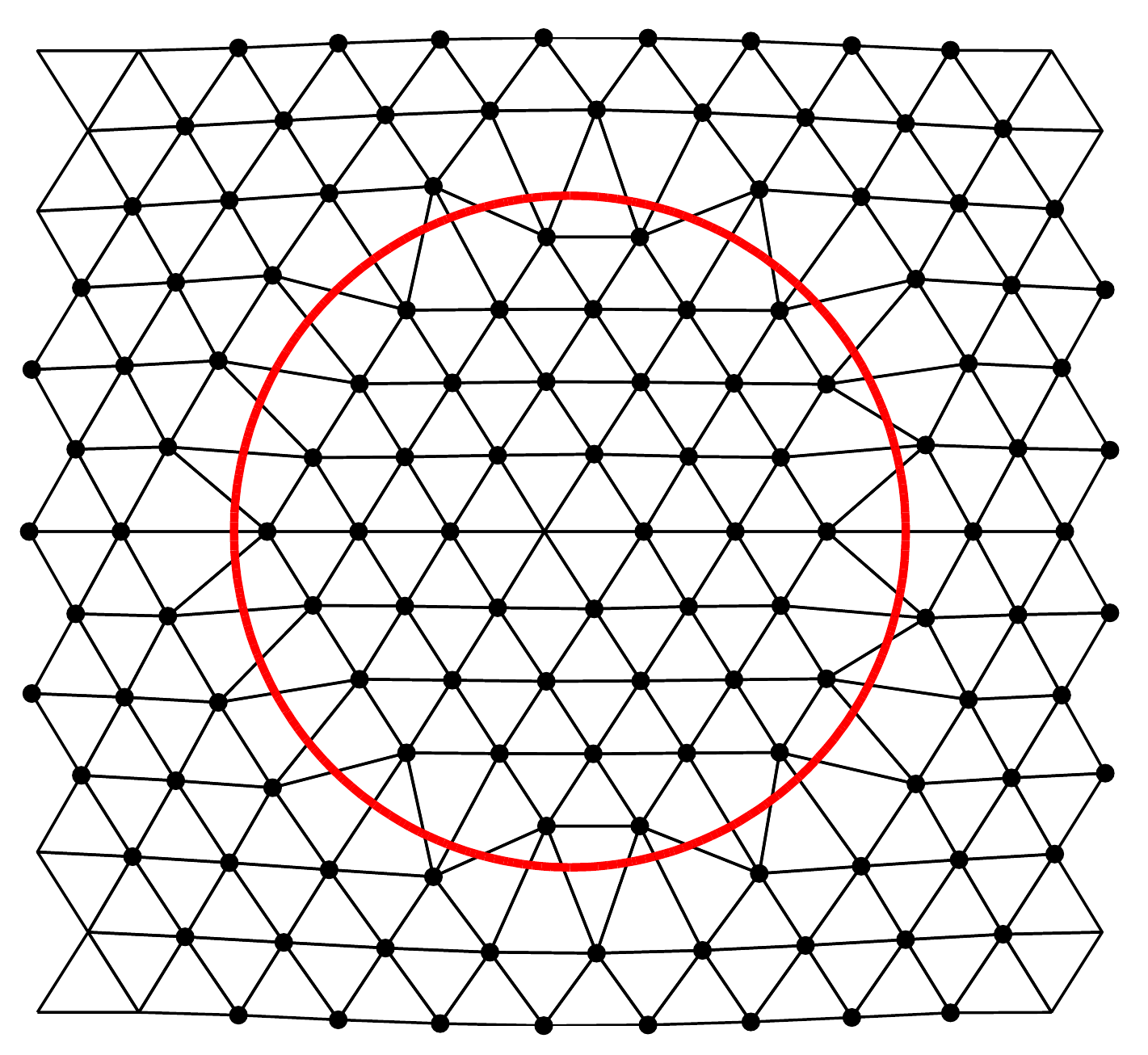}}\quad\subfigure[$\phi_1$, fine]{\includegraphics[width=5cm]{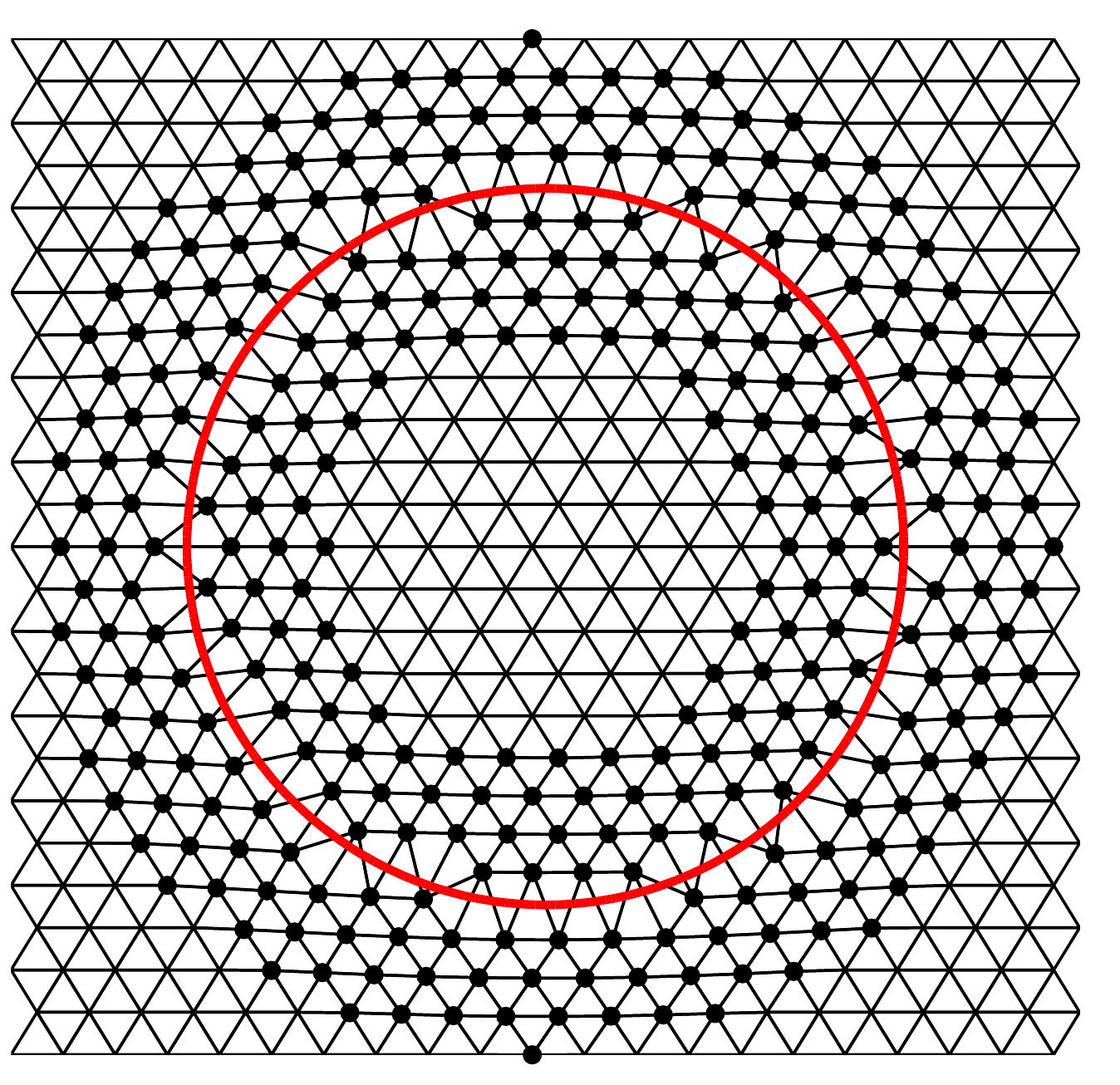}}\quad\\\subfigure[$\phi_2$, coarse]{\includegraphics[width=5cm]{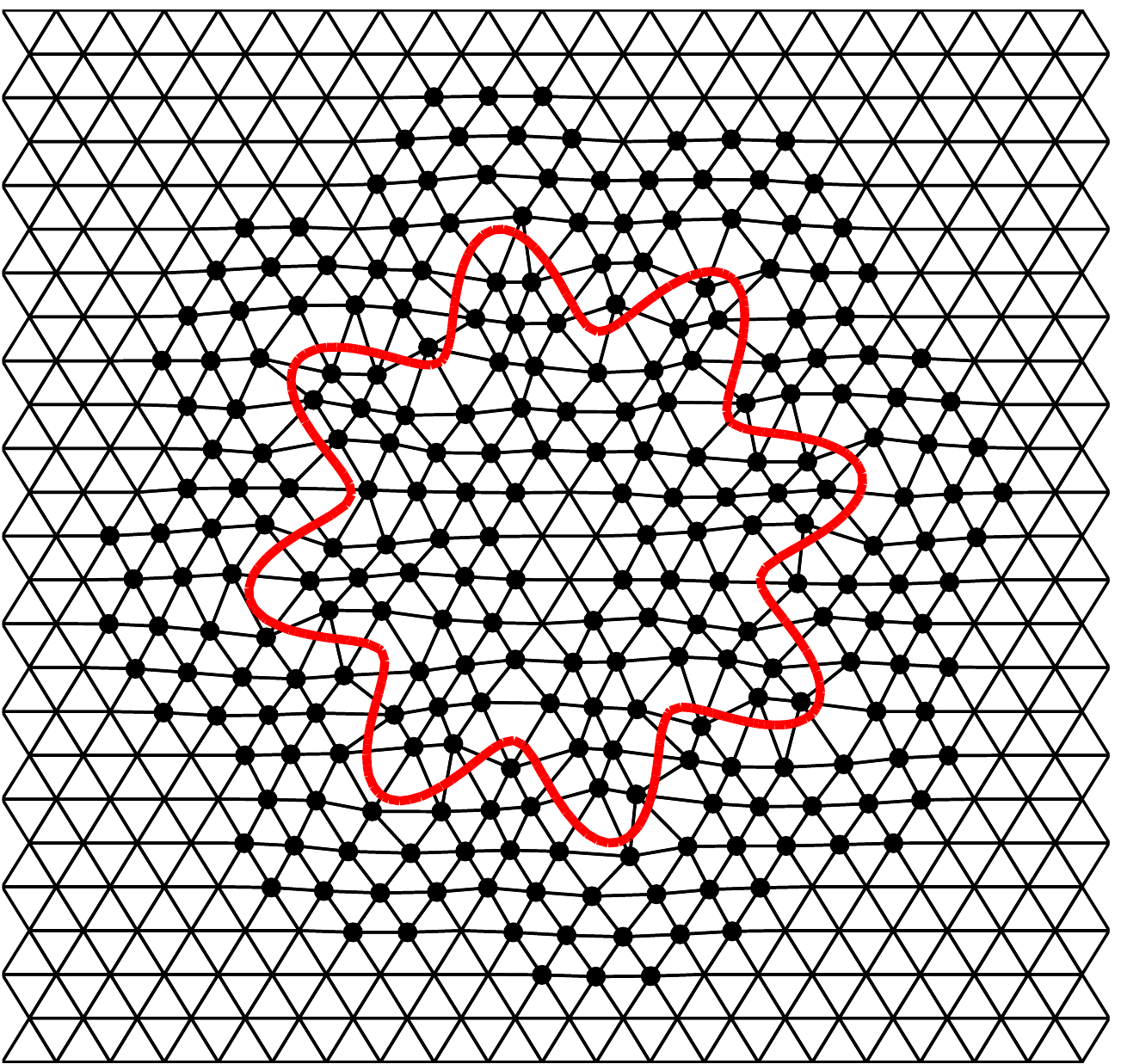}}\quad\subfigure[$\phi_2$, fine]{\includegraphics[width=5cm]{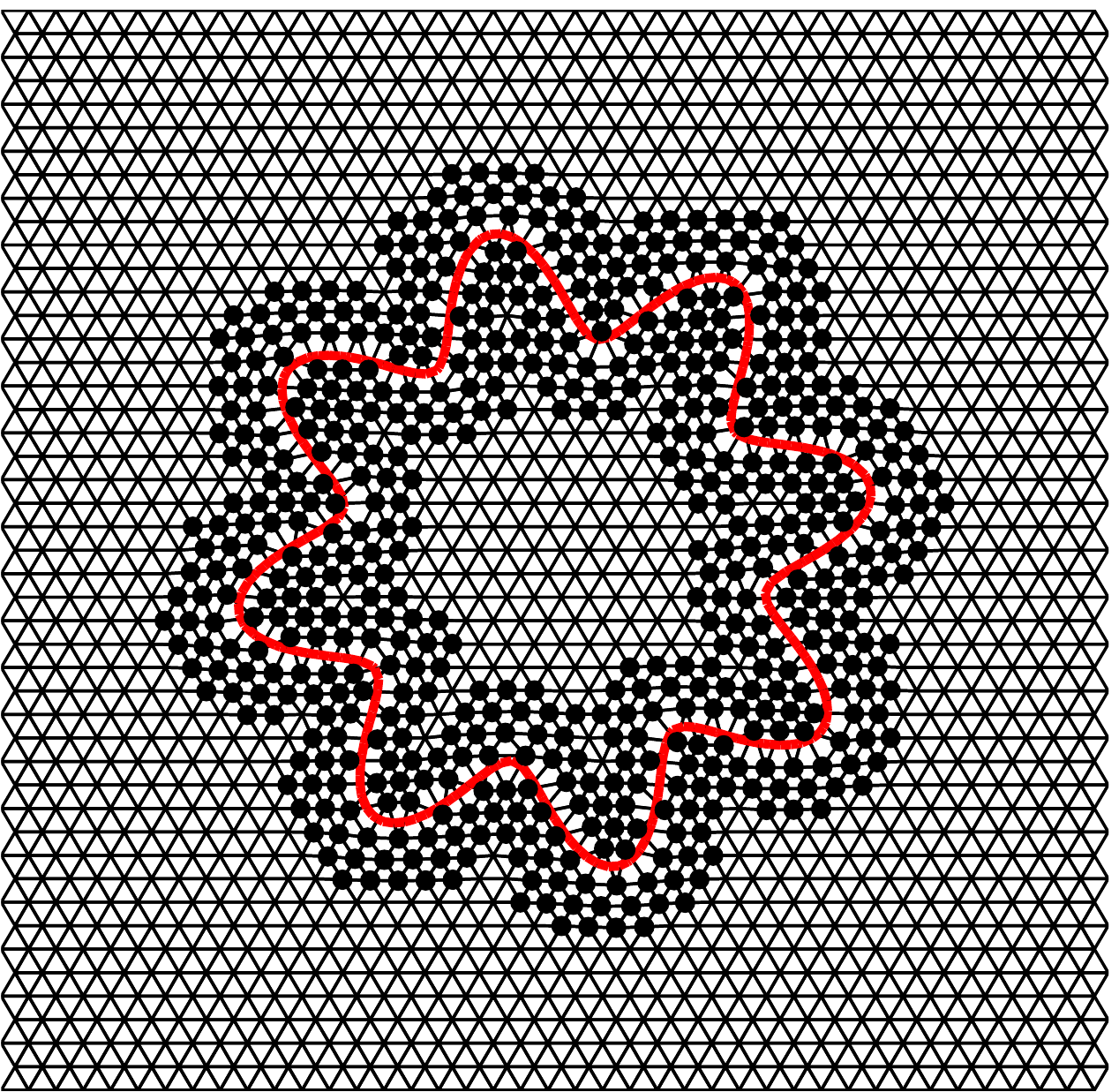}}

\caption{\label{fig:PlotMeshManip2d}Examples for node movements in 2D based
on different level-set functions $\phi$ and background meshes. The
red lines are the zero isolines $\Gamma_{\phi_{i}}$ which are meshed
by higher-order line elements.}
\end{figure}

Examples for manipulations of background meshes in 2D are shown in
Fig.~\ref{fig:PlotMeshManip2d}. Two different level-set functions
$\phi$ on a coarse and a fine mesh are considered, respectively.
The red lines indicate the zero-level sets which are meshed by higher-order
line elements. It is seen that due to the nodal movements, the lengths
of the line elements in the background mesh are quite regular, hence,
the ratio $l_{\mathrm{max}}/l_{\mathrm{min}}$ is bounded. The results
are achieved with $D_{\mathrm{crit}}=3\cdot h$, $D_{\mathrm{min}}=0.25\cdot h$,
and $D_{\mathrm{step}}=0.1\cdot h$.

\begin{figure}
\centering

\subfigure[$\Gamma_{\phi_1}$, coarse]{\includegraphics[width=4cm]{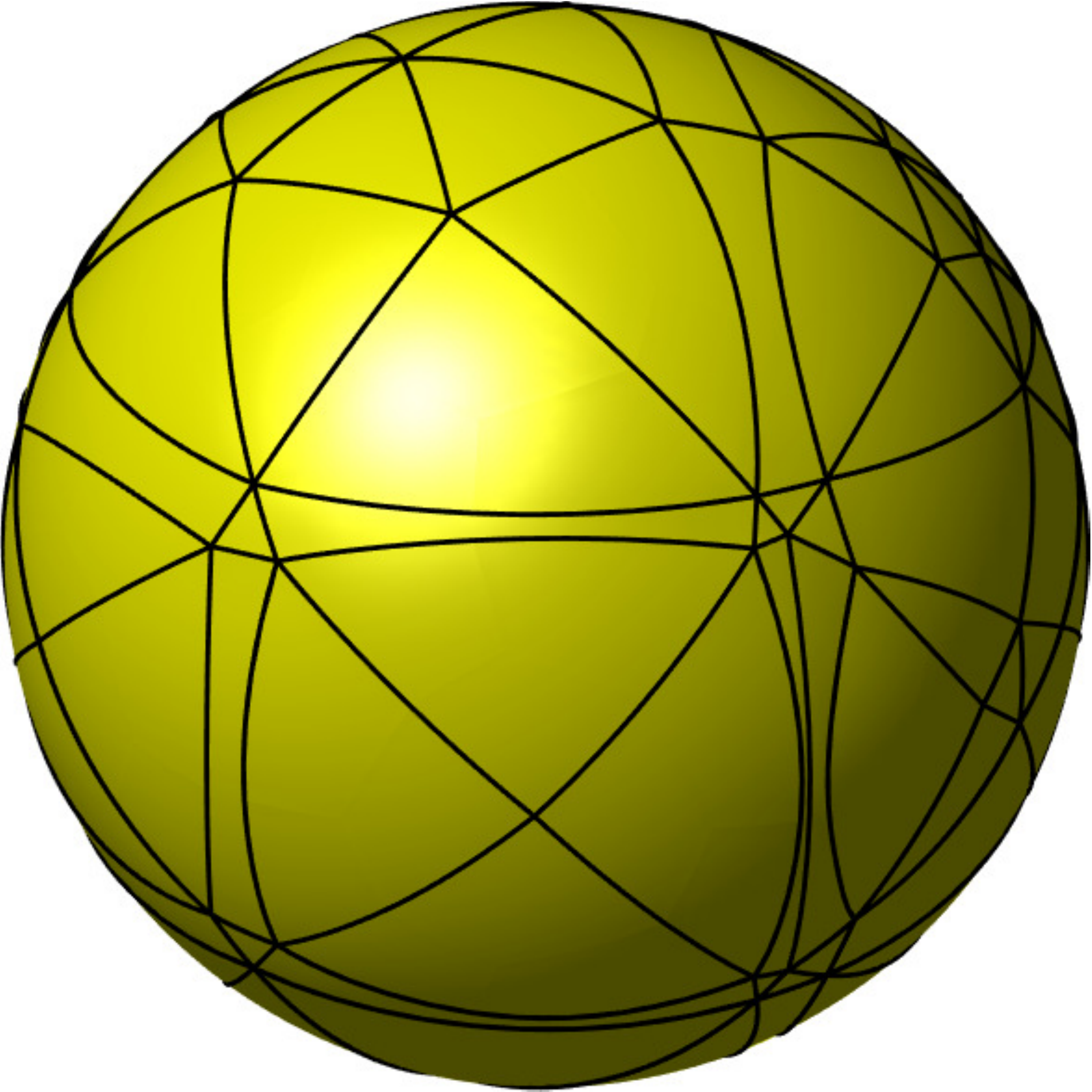}}\quad\subfigure[$\Gamma_{\phi_1}$, medium]{\includegraphics[width=4cm]{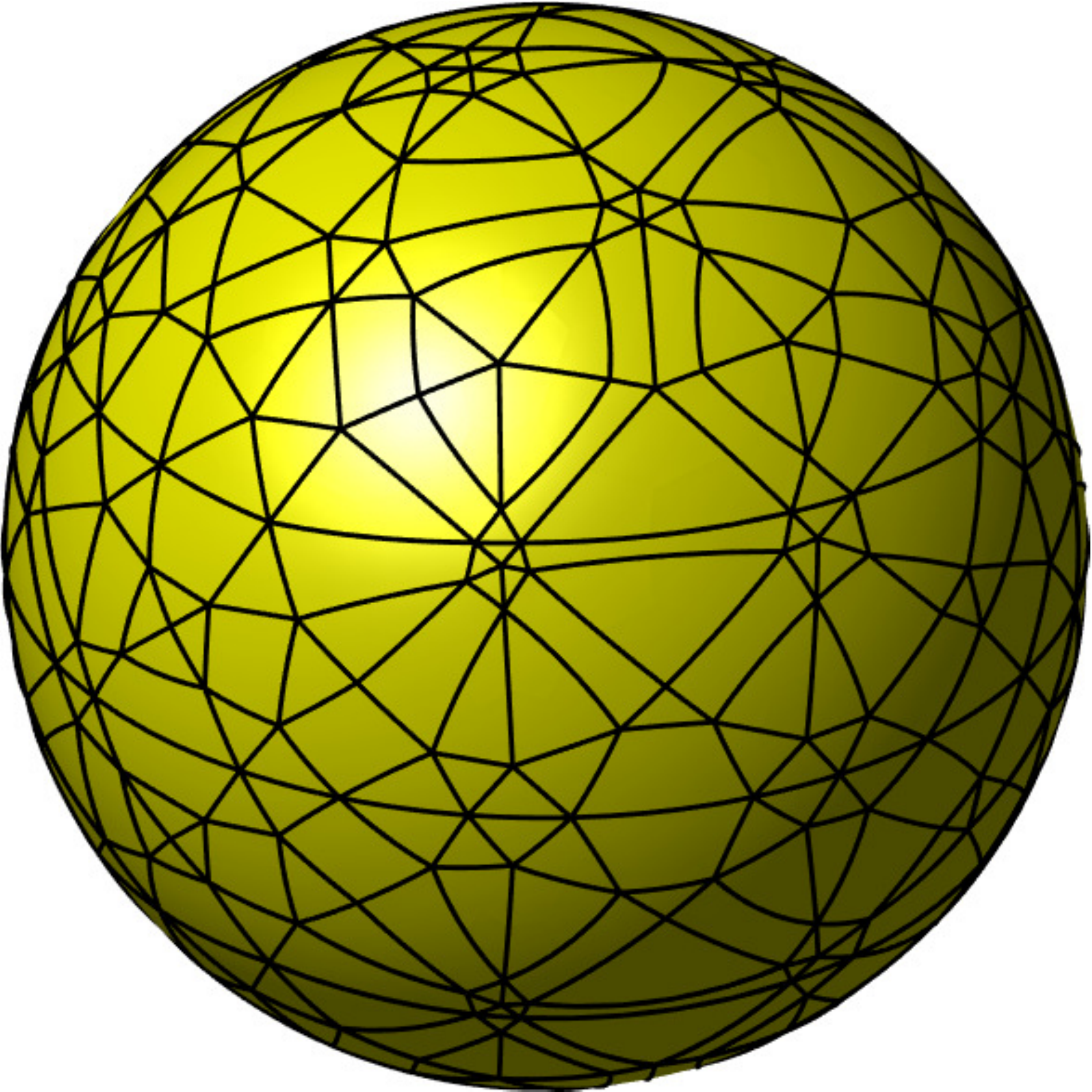}}\quad\subfigure[$\Gamma_{\phi_1}$, fine]{\includegraphics[width=4cm]{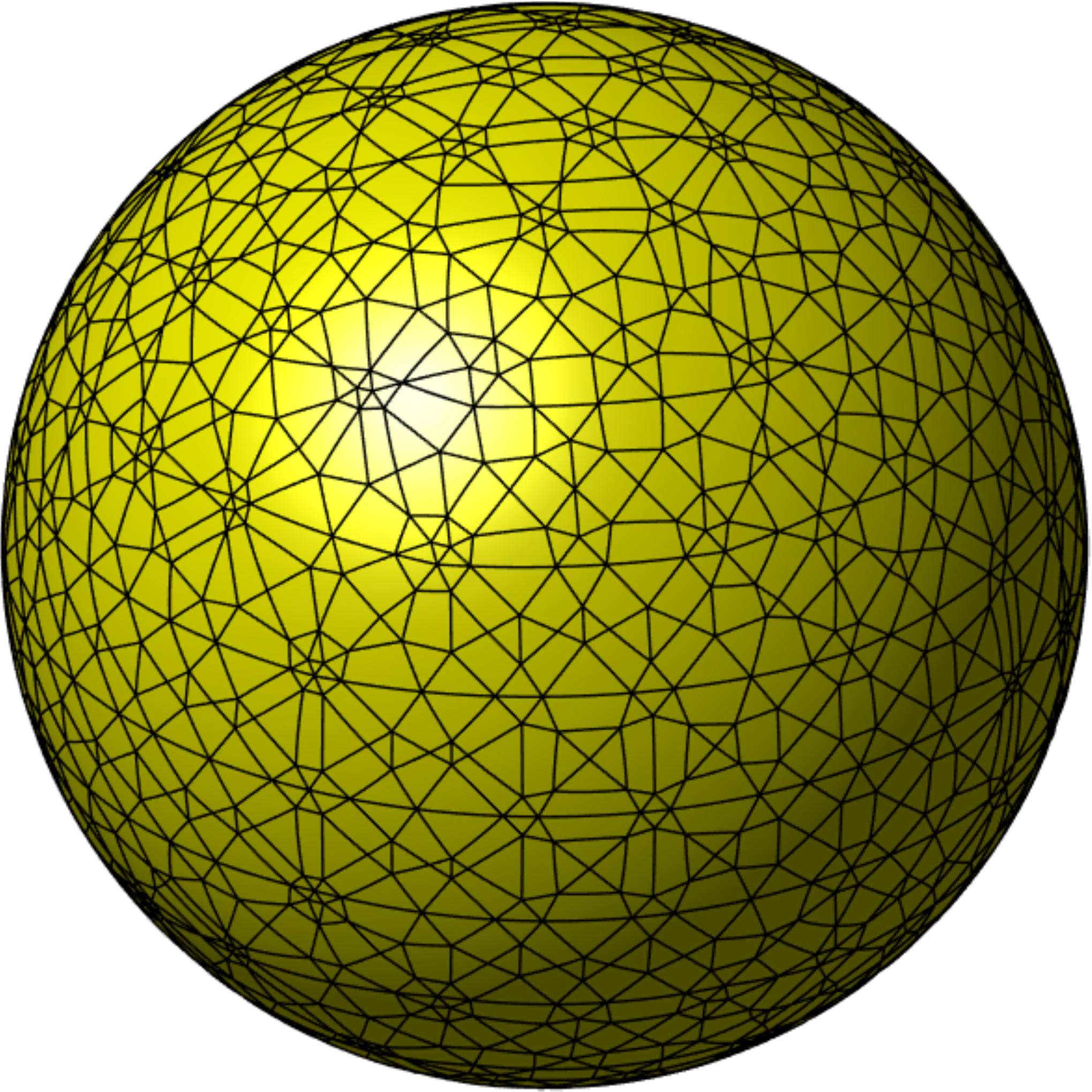}}

\subfigure[$\Gamma_{\phi_2}$, coarse]{\includegraphics[width=4cm]{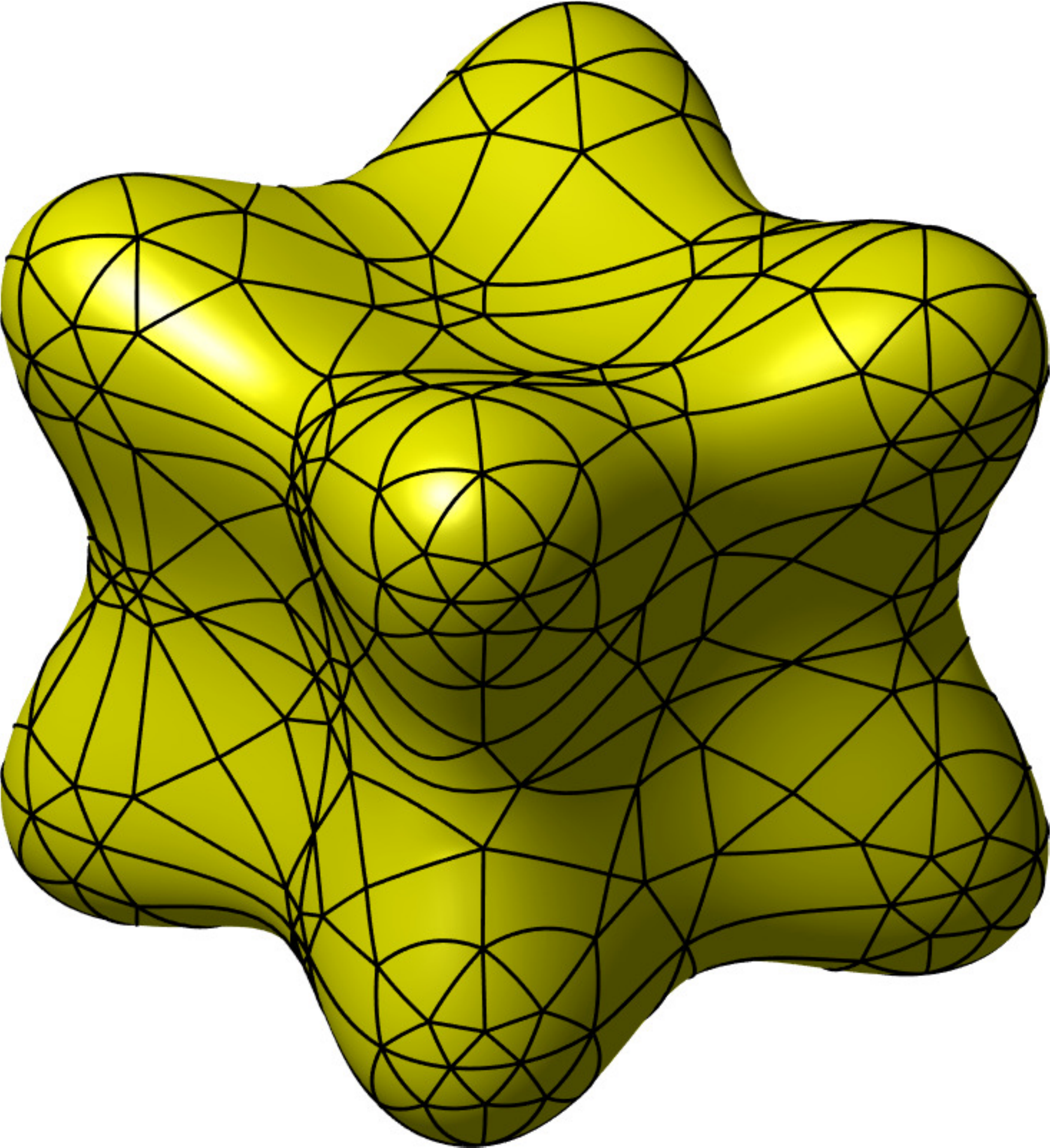}}\quad\subfigure[$\Gamma_{\phi_2}$, medium]{\includegraphics[width=4cm]{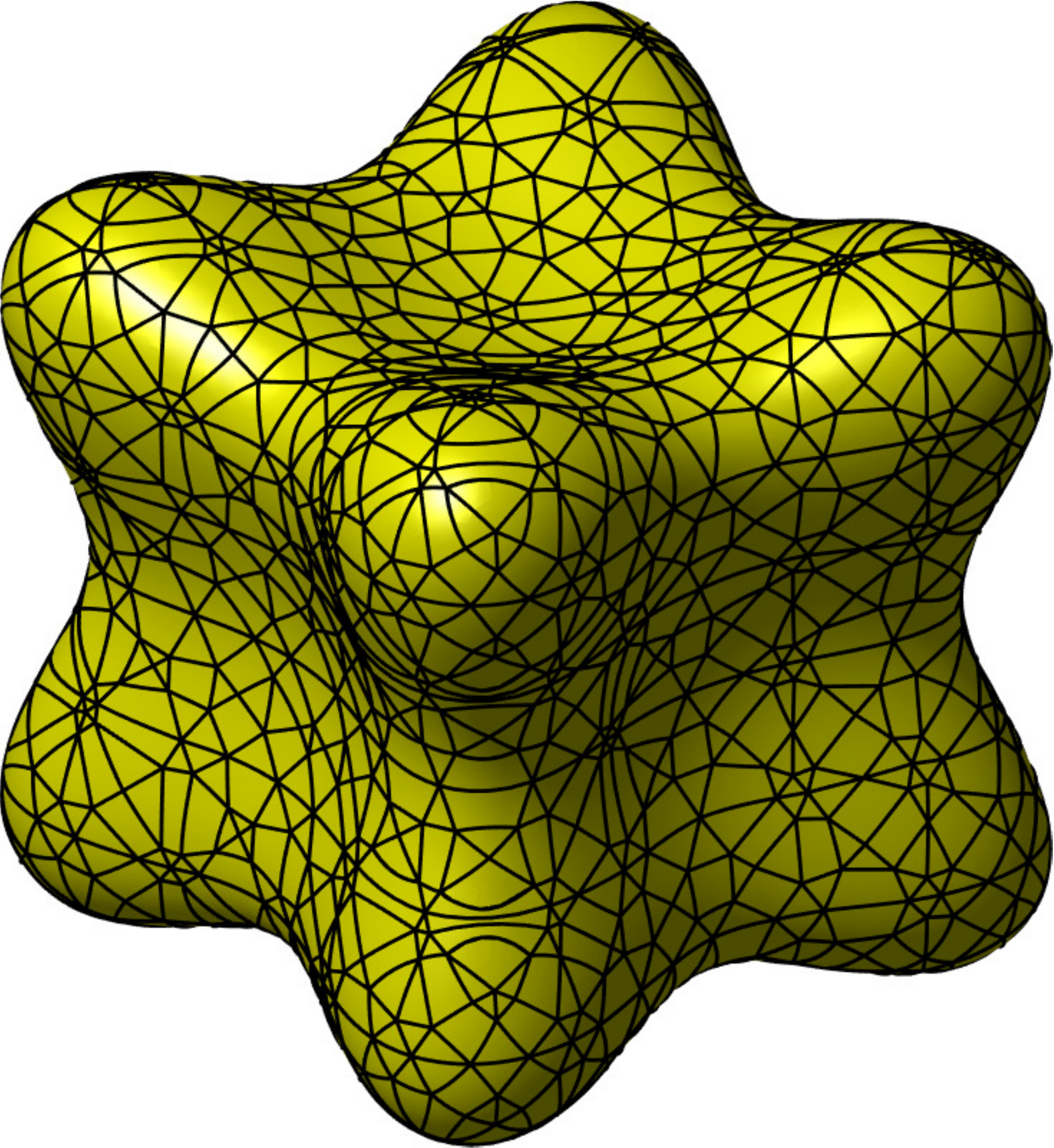}}\quad\subfigure[$\Gamma_{\phi_2}$, fine]{\includegraphics[width=4cm]{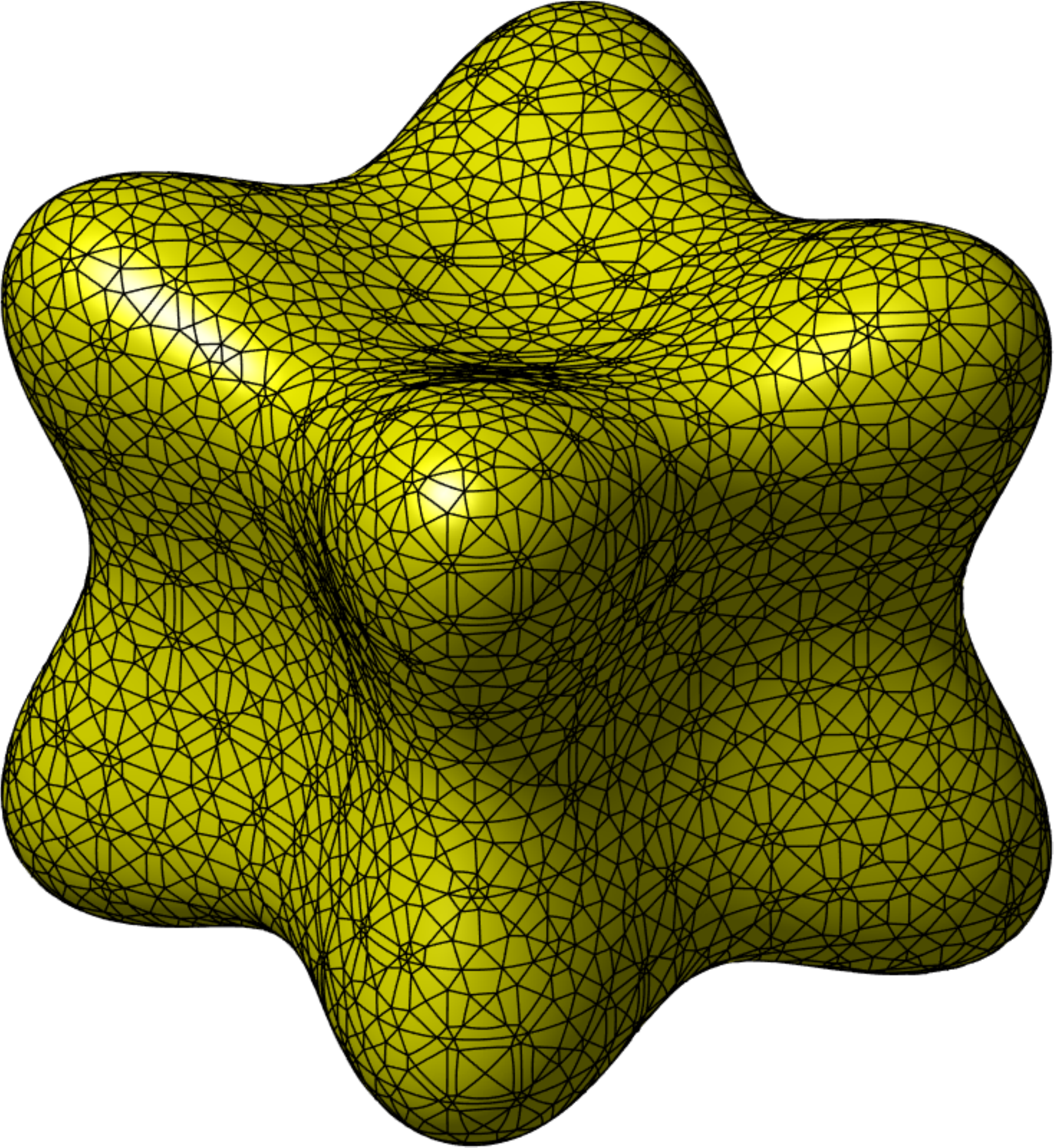}}

\caption{\label{fig:PlotMeshManip3d}Examples for node movements in 3D based
on different level-set functions $\phi$ and background meshes (not
shown).}
\end{figure}

A similar study is repeated in 3D. In Fig.~\ref{fig:PlotMeshManip3d},
the resulting surface meshes are shown for two different level-set
functions on different background meshes, respectively. The manipulated
background meshes are omitted for clarity and the focus is on the
shape and size of the surface elements. The element areas of the largest
and smallest elements are computed and given in Tables \ref{tab:ElementAreasPhi1}
and \ref{tab:ElementAreasPhi2}. As can be seen, the ratio $A_{\mathrm{max}}/A_{\mathrm{min}}$
is well bounded; it is typically less than $100$. The increase in
the condition number when using the automatically generated meshes
as proposed herein (compared to manually constructed meshes) is in
the same range. This is not a problem as long as direct solvers are
employed, however, will increase the effort in the context of iterative
solvers.

\begin{table}
\centering

\begin{tabular}{|c||c|c|c|}
\hline 
$\phi_{1}$ & $A_{\mathrm{max}}$ & $A_{\mathrm{min}}$ & $A_{\mathrm{max}}/A_{\mathrm{min}}$\tabularnewline
\hline 
\hline 
coarse mesh & $1.51\cdot10^{-1}$ & $1.63\cdot10^{-3}$ & $93.1$\tabularnewline
\hline 
medium mesh & $3.69\cdot10^{-2}$ & $6.23\cdot10^{-4}$ & $59.2$\tabularnewline
\hline 
fine mesh & $6.72\cdot10^{-3}$ & $1.84\cdot10^{-4}$ & $36.5$\tabularnewline
\hline 
\end{tabular}\caption{\label{tab:ElementAreasPhi1}Element areas of the surface meshes resulting
from $\phi_{1}$ as shown in Fig.~\ref{fig:PlotMeshManip3d} (a)
to (c).}
\end{table}

\begin{table}
\centering

\begin{tabular}{|c||c|c|c|}
\hline 
$\phi_{2}$ & $A_{\mathrm{max}}$ & $A_{\mathrm{min}}$ & $A_{\mathrm{max}}/A_{\mathrm{min}}$\tabularnewline
\hline 
\hline 
coarse mesh & $3.20\cdot10^{-2}$ & $7.35\cdot10^{-4}$ & $43.6$\tabularnewline
\hline 
medium mesh & $1.07\cdot10^{-2}$ & $1.33\cdot10^{-4}$ & $81.0$\tabularnewline
\hline 
fine mesh & $1.88\cdot10^{-3}$ & $4.31\cdot10^{-5}$ & $43.7$\tabularnewline
\hline 
\end{tabular}\caption{\label{tab:ElementAreasPhi2}Element areas of the surface meshes resulting
from $\phi_{2}$ as shown in Fig.~\ref{fig:PlotMeshManip3d} (d)
to (f).}
\end{table}

\section{PDEs on manifolds\label{X_FEMonManifolds}}

Herein, we consider the stationary Poisson equation and the instationary
advection-diffusion equation on curved manifolds. Both involve the
Laplace-Beltrami operator. The PDEs are conveniently defined based
on the tangential differential calculus, i.e.~on surface gradients
(sometimes refered to as intrinsic gradients). The manifold $\Gamma$
is implicitly defined by the zero-level set of a sufficiently smooth
level-set function. The tangential gradient operator $\nabla_{\Gamma}$
of a scalar function $u:\Gamma\to\mathbb{R}$ on the manifold is defined
as
\begin{align*}
\nabla_{\Gamma}u=\mat P\cdot\nabla\tilde{u}=\left(\mathbb{I}-\vek n_{\Gamma}\otimes\vek n_{\Gamma}\right)\cdot\nabla\tilde{u} & .
\end{align*}
$\mat P$ is the discrete tangential projector, $\vek n_{\Gamma}$
is the unit normal vektor on the manifold $\Gamma$ and $\tilde{u}$
is a smooth extension in a neighborhood $\mathcal{U}$ of the manifold
$\Gamma$. With the tangential divergence of a vector $\vek v$, $\nabla_{\Gamma}\cdot\vek v=\text{tr}(\nabla_{\Gamma}\vek v)$,
the Laplace-Beltrami operator is \cite{Olshanskii_2009a,Dziuk_2013a}
\begin{align}
\Delta_{\Gamma}u=\nabla_{\Gamma}\cdot\nabla_{\Gamma}u\qquad\text{on}\;\Gamma.
\end{align}

Assume a sucessfull reconstruction of the manifold $\Gamma$, the
discrete surface $\Gamma_{h}$ is the union of the elements $\Gamma_{h}=\underset{T\in\tau_{h}}{\bigcup}T$,
where $\tau_{h}$ is the set of automatically generated Lagrange elements
as described above. The standard surface FEM as described in \cite{Demlow_2009a,Dziuk_2013a}
is used to solve the PDEs on the discrete manifold. In the context
of the FEM the function $u$ is often a shape function $N\left(\vek r\right)$
living in a reference element $\bar{T}$ with coordinates $\vek r$.
Using the standard isoparametric mapping from the reference element
to the physical element, $\vek\chi:=\sum_{i}N_{i}\left(\vek r\right)\cdot\vek x_{i},\ \vek r\in\bar{T}$,
and $\vek x_{i}\in T$ being nodal coordinates, the finite element
space on $\Gamma_{h}$ with the polynomial degree $k$ is then defined
by
\begin{align}
\mathcal{S}_{h}^{k}=\left\{ u_{h}\in C\text{\textsuperscript{0}}\left(\Gamma_{h}^{k}\right)\mid u_{h}\circ\vek\chi^{-1}\vert_{T}\in\mathbb{P}_{k}\left(\bar{T}\right)\;\text{or}\;\mathbb{Q}_{k}\left(\bar{T}\right),T\in\tau_{h}\right\} .
\end{align}
$\mathbb{P}_{k}\left(\bar{T}\right)$ is the polynomial basis in a
triangular element and $\mathbb{Q}_{k}\left(\bar{T}\right)$ in a
quadrilateral element, both being complete of order $k$. This space
is spanned by the nodal basis $M_{1},\ldots,\, M_{n}$, with $n$
being the total number of nodes, which is given by
\begin{align*}
M_{i}\left(\vek x\right)\in S_{h}^{k},\qquad M_{i}(\vek x_{j})=\delta_{ij}\ .
\end{align*}
A base function is defined by $M_{i}=N_{i}\circ\vek\chi\vert_{T}$.
Note that the function space $\mathcal{S}_{h}^{k}$ is a subspace
of the continuous space $\mathcal{H}^{1}(\Gamma_{h})$. For more information
we refer to \cite{Demlow_2009a,Dziuk_2013a}. 

In general, it is desirable to compute the surface gradient without
explicitly computing an extension $\tilde{u}$. Assume that the surface
element $T$ in the physical space $\vek x\in\mathbb{R}^{3}$ follows
from the standard isoparametric mapping $\vek\chi$, then the surface
gradient is directly obtained from
\begin{equation}
\nabla_{\Gamma_{h}}u_{h}\left(\vek x\left(\vek r\right)\right)=\mat J\cdot\mat G^{-1}\cdot\nabla_{\vek r}u_{h}\left(\vek r\right)\label{eq:SurfGradShapeFct}
\end{equation}
with $\mat J=\nicefrac{\partial\vek\chi}{\partial\vek r}$ being the
($3\times2$)-Jacobi matrix and $\mat G=\mat J^{T}\cdot\mat J$ being
the first fundamental form.

\subsection{Stationary Poisson equation on manifolds\label{XX_LaplaceBeltrami}}

The complete boundary value problem in strong form for the Poisson
equation on manifolds states that we seek $u:\Gamma\to\mathbb{R}$
such that: \begin{subequations} 
\begin{alignat}{2}
-\Delta_{\Gamma}u & =f &  & \text{on}\quad\Gamma,\label{eq:lap}\\
u & =g_{\mathrm{D}}\quad &  & \text{on}\quad\partial\Gamma_{\text{D}},\\
\vek n_{\partial\Gamma}\cdot\nabla_{\Gamma}u & =g_{\mathrm{N}} &  & \text{on}\quad\partial\Gamma_{\text{N}},
\end{alignat}
\end{subequations} where $f\in L^{2}(\Gamma)$ is a source function
and $\vek n_{\partial\Gamma}$ is the co-normal vector being normal
to the boundary yet in the tangent plane of the manifold at each point
on the boundary.

For a closed (compact) manifold, where no boundary exists, one needs
an additional condition for the problem to be well-posed. Therefore,
typically the zero mean constraint is imposed,
\begin{align}
\int\limits _{\Gamma}u\,\text{d}s=0.\label{eq:const}
\end{align}
The boundary value problem from above is converted to a weak form
and discretized using the surface FEM. Therefore, we introduce the
following trial and test function spaces
\begin{eqnarray}
\mathcal{S}_{u} & = & \left\lbrace u_{h}\in\mathcal{S}_{h}^{k}(\Gamma_{h})\ :\ u\rvert_{\partial\Gamma_{h,D}}=g_{\mathrm{D}}\right\rbrace ,\label{eq:ApproxSpacesA}\\
\mathcal{V}_{u} & = & \left\lbrace v_{h}\in\mathcal{S}_{h}^{k}(\Gamma_{h})\ :\ v\rvert_{\partial\Gamma_{h,D}}=0\right\rbrace \qquad\mathrm{for\; manifolds\; with\; boundary},\label{eq:ApproxSpacesB}\\
\mathcal{V}_{u} & = & \big\{ v_{h}\in\mathcal{S}_{h}^{k}(\Gamma_{h})\ :\ \smallint_{\Gamma_{h}}v_{h}\,\text{\text{\text{d}}}s=0\big\}\qquad\mathrm{for\; compact\; manifolds}.\label{eq:ApproxSpacesC}
\end{eqnarray}
Then, the discrete weak form for the Poisson equation on a manifold
$\Gamma_{h}$ is stated as: Find $\ensuremath{u_{h}\in\mathcal{S}_{u}}$
such that
\begin{equation}
\int\limits _{\Gamma_{h}}\nabla_{\Gamma_{h}}v_{h}\cdot\nabla_{\Gamma_{h}}u_{h}\;\text{\text{d}}s=\int\limits _{\Gamma_{h}}v_{h}\cdot f\;\text{\text{d}}s+\int_{\partial\Gamma_{h}}v_{h}\cdot g_{\mathrm{N}}^{h}\;\text{\text{d}}q\qquad\forall v_{h}\in\mathcal{V}_{u}.\label{eq:varLap}
\end{equation}
The existence and uniqueness of this BVP is shown in \cite{Dziuk_2013a,Burman_2015c}.
The zero mean constraint for compact manifolds is enforced by a Lagrange
multiplier.

\subsection{Instationary advection-diffusion equation on manifolds\label{XX_AdvectionDiffusion}}

The strong form of the instationary advection-diffusion equation in
the time interval from $0$ to $T$ reads \begin{subequations} 
\begin{alignat}{2}
\dot{u}+\vek c_{\Gamma}\cdot\nabla_{\Gamma}u-\lambda\Delta_{\Gamma}u & =f & \text{on} & \quad\Gamma\times\left(0,T\right),\label{eq:instat}\\
u & =g_{\mathrm{D}} & \text{on} & \quad\partial\Gamma_{\text{D}}\times\left(0,T\right),\\
\vek n_{\partial\Gamma}\cdot\nabla_{\Gamma}u & =g_{\mathrm{N}} & \text{on} & \quad\partial\Gamma_{\mathrm{N}}\times\left(0,T\right),\\
u(\vek x,0) & =u_{0}(\vek x)\quad & \text{on} & \quad\Gamma\;\mathrm{at}\; t=0,
\end{alignat}
\end{subequations}with $\dot{u}=\nicefrac{\partial u}{\partial t}$,
$\vek c_{\Gamma}$ being the advection velocity in the tangent space
of the manifold and $\lambda$ the diffusion coefficient. Using the
same function spaces from above, the semi-discrete weak form becomes:
Find $u_{h}\in\mathcal{S}_{u}\times\left(0,T\right)$ with $\dot{u}_{h}\in L^{2}(\Gamma_{h})$
such that $\forall v_{h}\in\mathcal{V}_{u}$
\begin{equation}
\int\limits _{\Gamma_{h}}\left[v_{h}\cdot\dot{u}+v_{h}\cdot\vek c_{\Gamma}\cdot\nabla_{\Gamma_{h}}u_{h}+\lambda\cdot\nabla_{\Gamma_{h}}v_{h}\cdot\nabla_{\Gamma_{h}}u_{h}\right]\text{\text{d}}s=\int\limits _{\Gamma_{h}}v_{h}\cdot f(\vek x,t)\;\text{\text{d}}s+\int_{\partial\Gamma_{h}}v_{h}\cdot g_{\mathrm{N}}^{h}\;\text{\text{d}}q.\label{eq:WeakFormAdvDiff}
\end{equation}
The zero mean constraint from (\ref{eq:ApproxSpacesC}) for compact
manifolds is not needed here. The existence and uniqueness of the
solution is guaranteed by the Lax-Milgramm lemma and is proven for
the stationary problem in \cite{Olshanskii_2014a}.

After the FEM discretization in space, a first-order system of ordinary
differential equations in time is obtained as usual. We use a $6$th-order
accurate, implicit, $3$-step Runge-Kutta method based on Gauss-Legendre
points to discretize in time. The corresponding Butcher tableau for
this Runge-Kutta scheme is given as

\begin{equation}
\left[\begin{array}{c|c}
\vek c & \mat A\\
\hline  & \vek b^{T}
\end{array}\right]=\left[\begin{array}{c|ccc}
\frac{1}{2}-\frac{1}{10}\sqrt{15} & \frac{5}{36} & \frac{2}{9}-\frac{1}{15}\sqrt{15} & \frac{5}{36}-\frac{1}{30}\sqrt{15}\\
\frac{1}{2} & \frac{5}{36}+\frac{1}{24}\sqrt{15} & \frac{2}{9} & \frac{5}{36}-\frac{1}{24}\sqrt{15}\\
\frac{1}{2}+\frac{1}{10}\sqrt{15} & \frac{5}{36}+\frac{1}{30}\sqrt{15} & \frac{2}{9}-\frac{1}{15}\sqrt{15} & \frac{5}{36}\\
\hline  & \frac{5}{18} & \frac{4}{9} & \frac{5}{18}
\end{array}\right]\label{eq:ButcherTableau}
\end{equation}

It is useful to employ this highly accurate scheme in time to virtually
eliminate the temporal error in the numerical studies and focus on
the spatial approximation error of the surface finite elements (with
different orders).

\section{Numerical results\label{X_NumericalResults}}

As defined in the previous section, the stationary Poisson equation
and the instationary advection-diffusion equation are considered here.
The focus above is on curved surfaces in three domensions, however,
results are also shown for curved lines in two dimensions. Convergence
studies for different element orders up to $6$ are conducted where
the analytic solutions are generated by the method of manufactured
solutions. In each test case, we compare the performance of handcrafted
meshes composed by very regular elements with the automatically generated
meshes as described in Section \ref{X_MeshGeneration}. The manipulation
of the background mesh from Section~\ref{XX_MeshManipulation} is
enabled in all examples. 

According to the method of manufactured solutions in the context of
the Poisson equation, the source function $f$ on the right hand side
is computed such that Eq.~(\ref{eq:lap}) is fulfilled for a chosen
solution $u(\vek x)$. That is, the Laplace-Beltrami operator (in
local coordinates)
\begin{align}
-\Delta_{\Gamma}u=-\dfrac{1}{\sqrt{\vert\det\mat G\vert}}\dfrac{\partial}{\partial x_{i}}\left(g^{ij}\sqrt{\vert\det\mat G\vert}\dfrac{\partial u}{\partial x_{j}}\right)=f,\label{eq:lapLoc}
\end{align}
is applied to $u(\vek x)$ to achieve $f(\vek x)$. There, $\mat G=\mat J^{\text{T}}\cdot\mat J$
with entries $g_{ij}$ is the first fundamental form and $\mat G^{-1}$
with entries $g^{ij}$ its inverse.

For the convergence studies, the error is shown over the element lengths
$h$ of the corresponding meshes. For the handcrafted meshes, $h$
is the element length of a typical element of the surface mesh (where
elements are quite uniformly shaped). For an automatically reconstructed
surface mesh, the element lengths vary largely so that $h$ is the
characteristic element length of the \emph{background} mesh (from
which the surface mesh was reconstructed).

\subsection{Poisson equation on curved lines in 2D}

\subsubsection{Circular manifold\label{sec:lap1dEx1}}

The Poisson equation is first solved on a circle with radius $r=1$.
The circle is defined by the implicit level-set function $\phi(\vek x)=\Vert x\Vert-1$.
Using polar coordinates $(r,\,\varphi)$, the solution $u$ which
satisfies the zero-mean condition from Eq.~(\ref{eq:const}) is chosen
as \cite{Kamilis_2013} 
\begin{align}
u(r,\varphi)=12\sin3\varphi
\end{align}
with the corresponding right hand side resulting as 
\begin{align}
f(r,\varphi)=\dfrac{108}{r^{2}}\sin3\varphi\ .
\end{align}

\begin{figure}
\centering

\subfigure[meshes]{\includegraphics[width=0.3\textwidth]{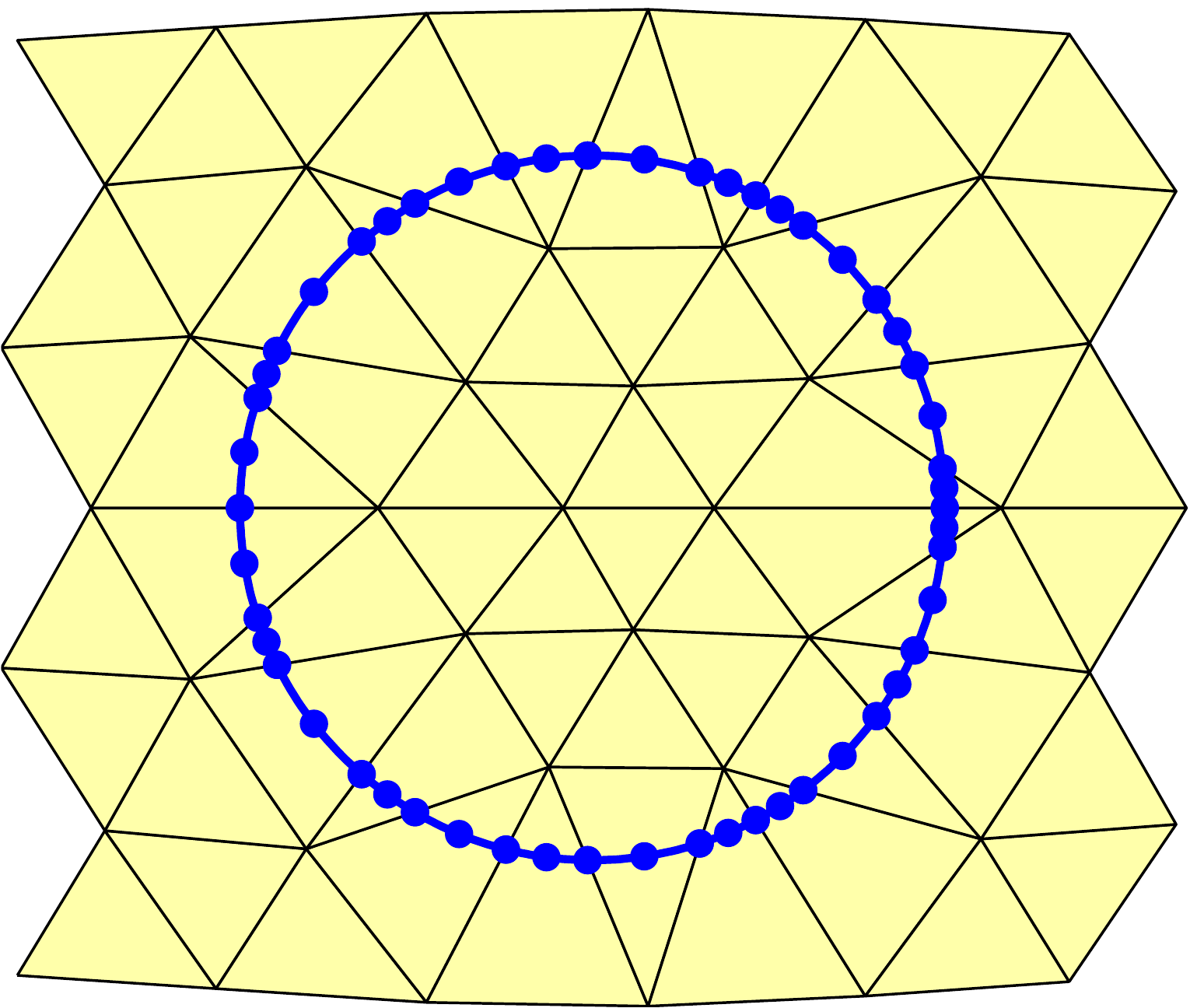}}\subfigure[exact solution]{\includegraphics[width=0.3\textwidth]{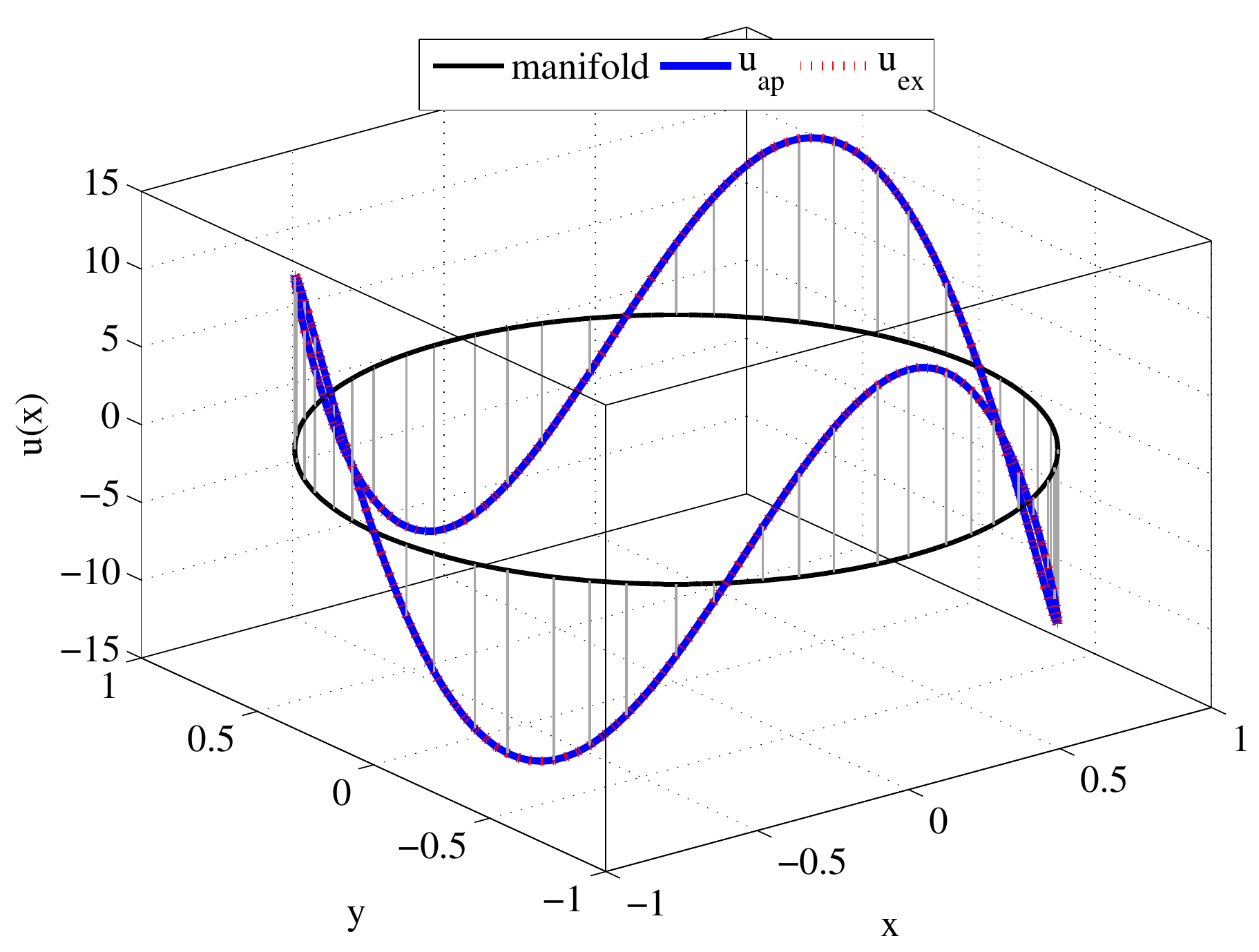}}\subfigure[ratio $h_{\mathrm{max}}/h_{\mathrm{min}}$]{\includegraphics[width=0.3\textwidth]{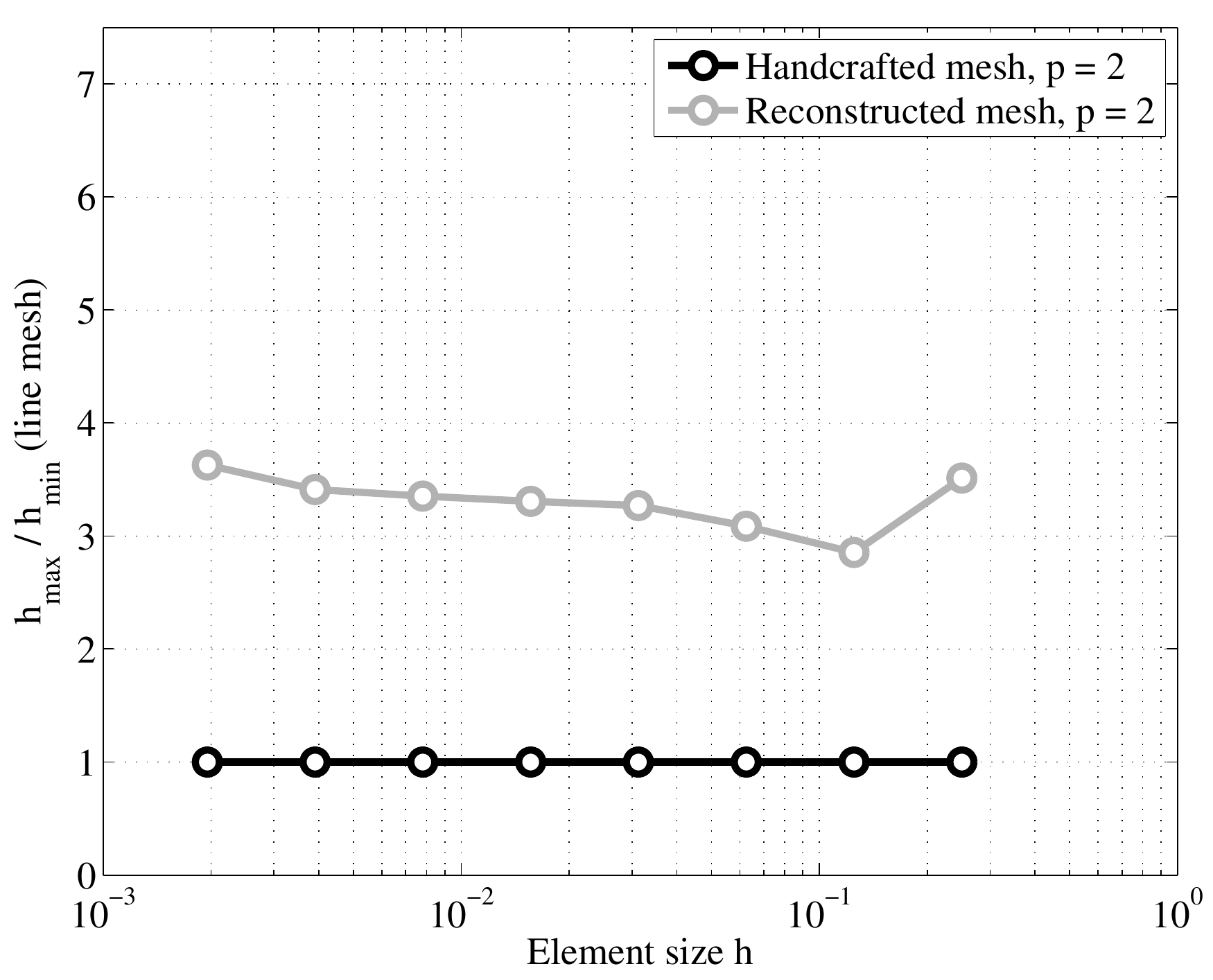}}

\caption{(a) Background and reconstructed line mesh for the circular manifold,
(b) exact solution of the Poisson equation, (c) the largest ratio
$h_{\mathrm{max}}/h_{\mathrm{min}}$ of the reconstructed line elements.}

\label{fig:PoissonCircleSetup} 
\end{figure}

The element sizes of the background meshes are systematically varied
as $h=r/\left\{ 4,\,8,\,16,\,32,\,64,\,128,\,256,\,512\right\} $
for the convergence studies. One example of a coarse background mesh
and a reconstruction by quadratic elements is shown in Fig.~\ref{fig:PoissonCircleSetup}(a).
The corresponding numerical solution is shown in Fig.~\ref{fig:PoissonCircleSetup}(b)
together with the analytical solution. Fig.~\ref{fig:PoissonCircleSetup}(c)
shows the ratio of the largest element length of the line mesh $h_{\mathrm{max}}$
over the smallest $h_{\mathrm{min}}$ for the different background
meshes. These numbers are very similar for the different element orders,
so that we only show one curve being representative for all orders.
It is seen that the worst ratio is less than $4$ in this example
which shows the success of the manipulation of the background mesh.
Of course, for handcrafted meshes on the circle, it is not a problem
to achieve the optimal ratio of $1$ in this test case.

\begin{figure}
\centering

\subfigure[handcr.~mesh, $L_2$-norm]{\includegraphics[width=0.35\textwidth]{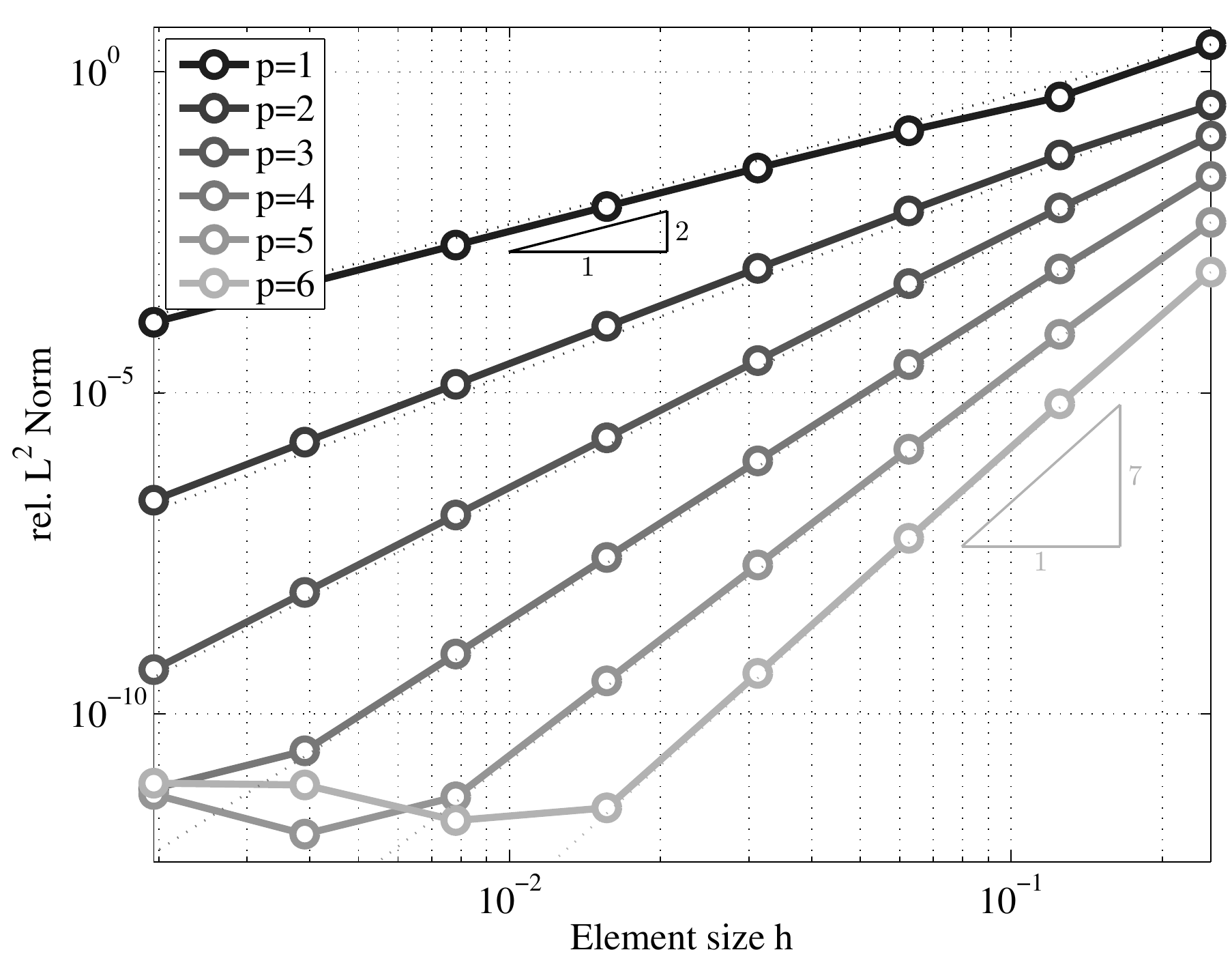}}\quad\subfigure[recon.~mesh, $L_2$-norm]{\includegraphics[width=0.35\textwidth]{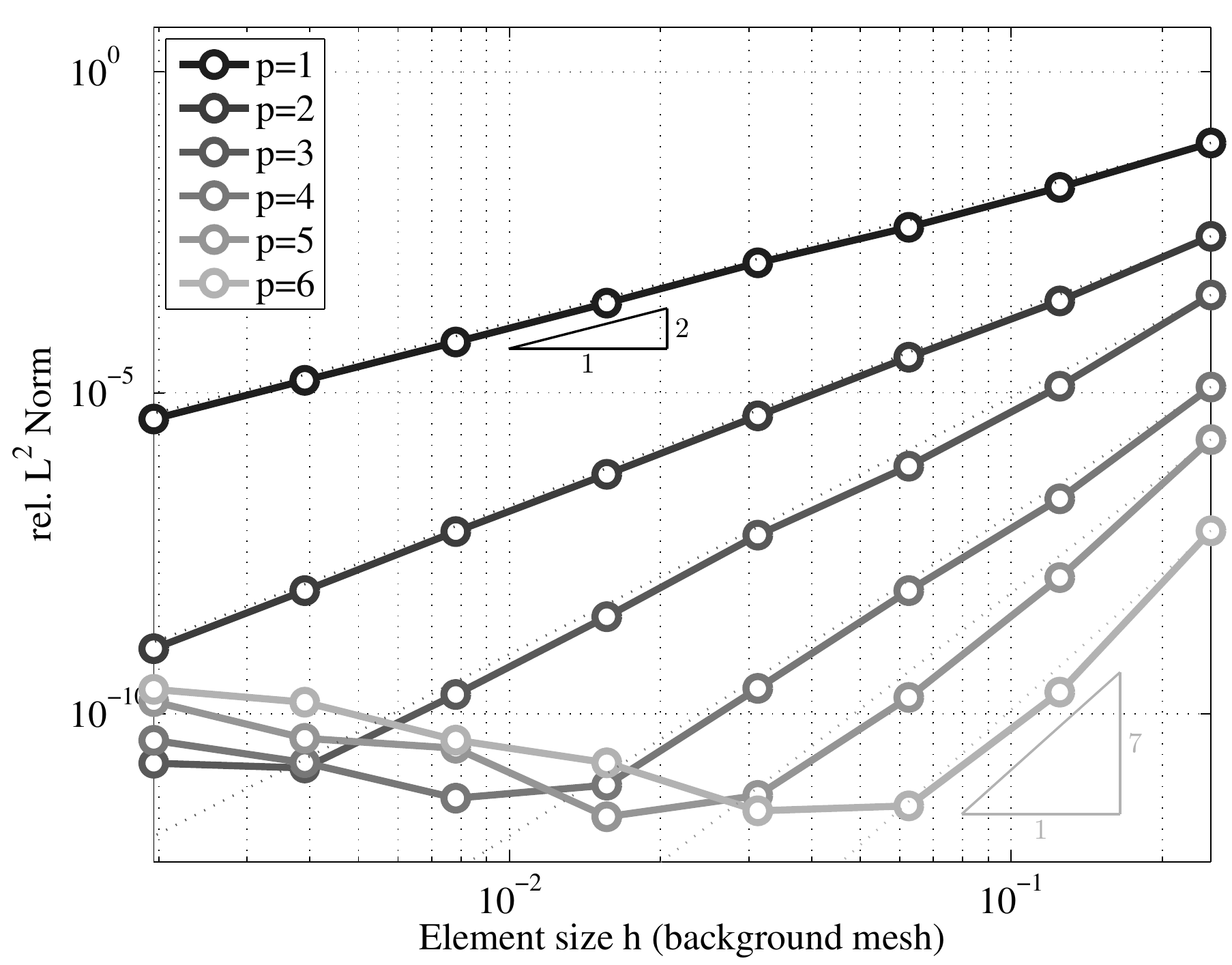}}

\subfigure[handcr.~mesh, cond.]{\includegraphics[width=0.35\textwidth]{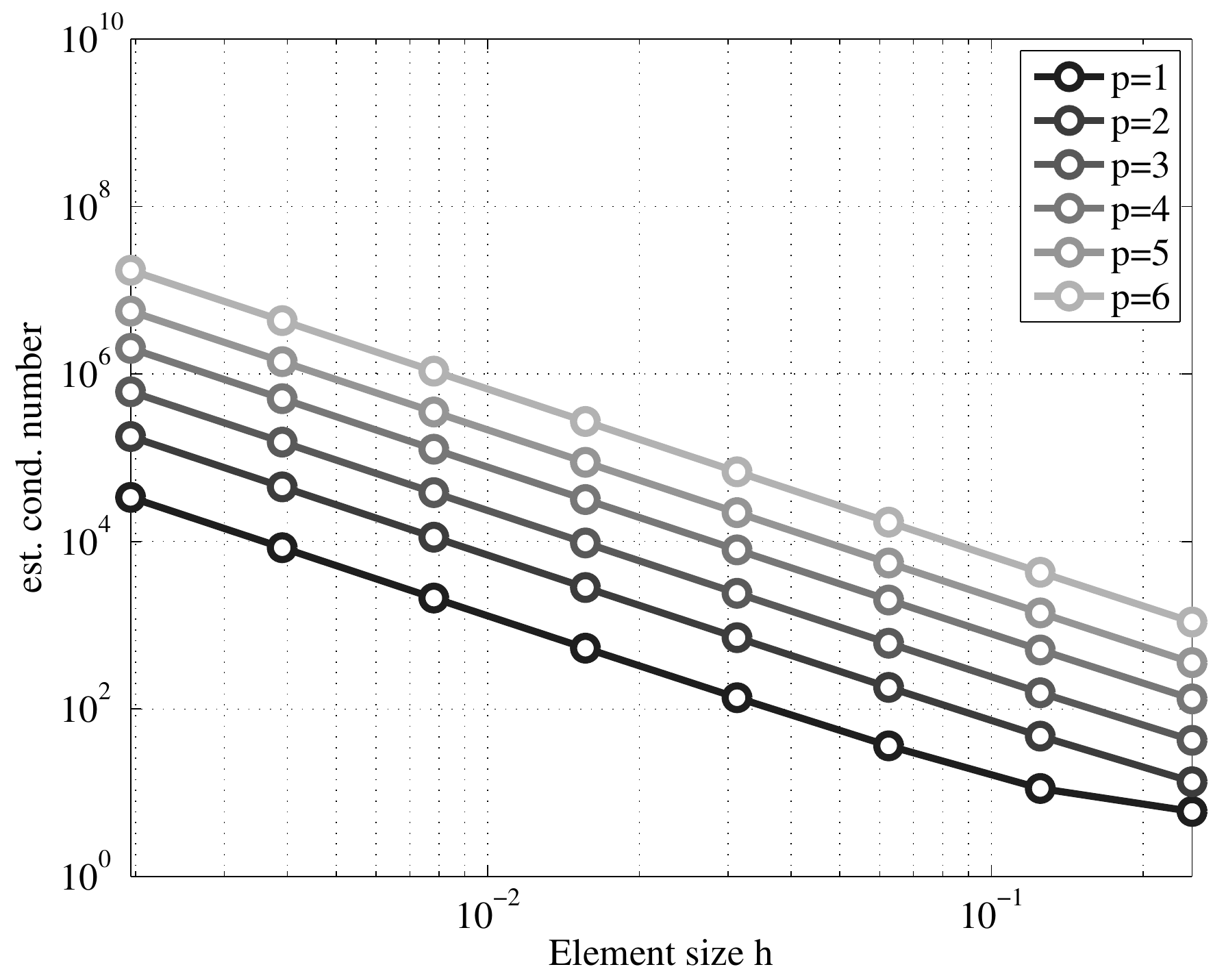}}\quad\subfigure[recon.~mesh, cond.]{\includegraphics[width=0.35\textwidth]{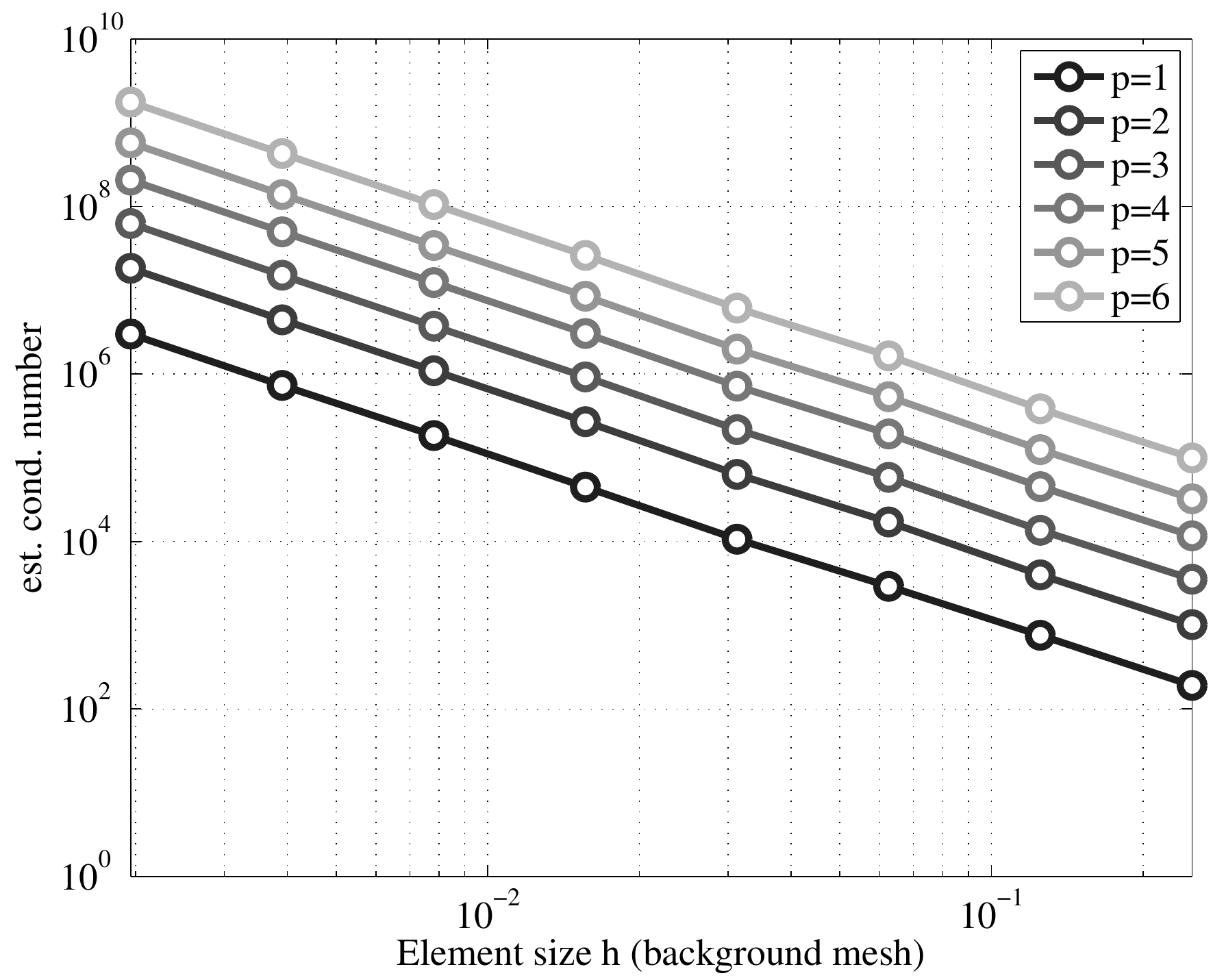}}

\caption{\label{fig:PoissonCircleRes}Convergence results and condition numbers
for the Poisson equation on a circular manifold.}
\end{figure}

Fig.~\ref{fig:PoissonCircleRes} shows the convergence rates and
condition numbers. Fig.~\ref{fig:PoissonCircleRes}(a) and (c) are
achieved on handcrafted meshes with equi-distant line elements. As
expected, optimal convergence rates are seen in the $L_{2}$-norm
and the condition number scales with $O\left(h^{2}\right)$. Fig.~\ref{fig:PoissonCircleRes}(b)
and (d) show the corresponding results on the automatically generated
meshes. The convergence rates are again optimal (the dotted lines
show the optimal slope) and the condition numbers also scale with
$O\left(h^{2}\right)$. Obviously, compared to the handcrafted meshes,
the error level is modestly larger and the condition number is larger
by about two order of magnitude.

\subsubsection{Flower-shaped manifold\label{sec:lap1dEx2}}

\begin{figure}
\centering

\subfigure[meshes]{\includegraphics[width=0.3\textwidth]{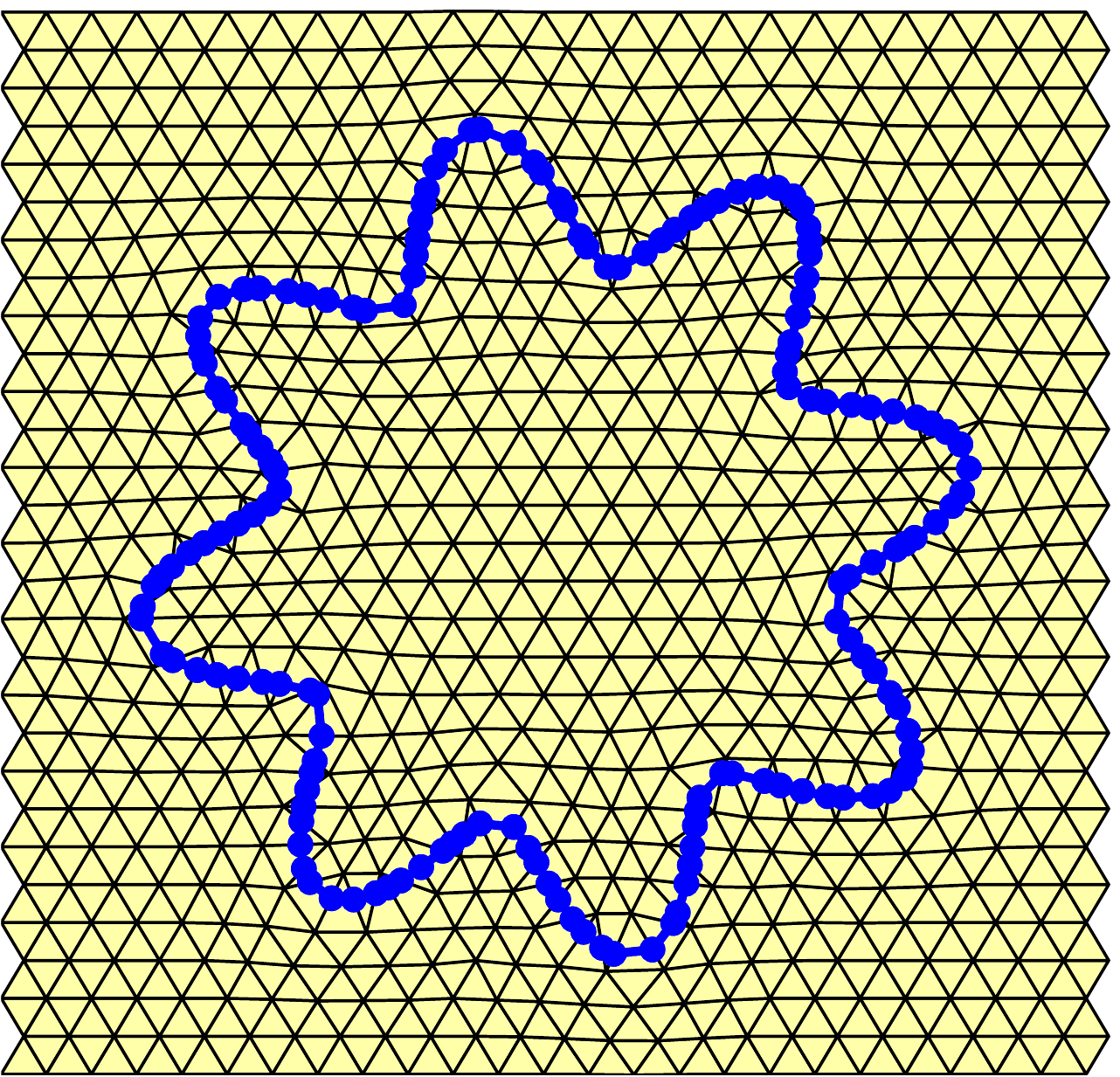}}\subfigure[exact solution]{\includegraphics[width=0.3\textwidth]{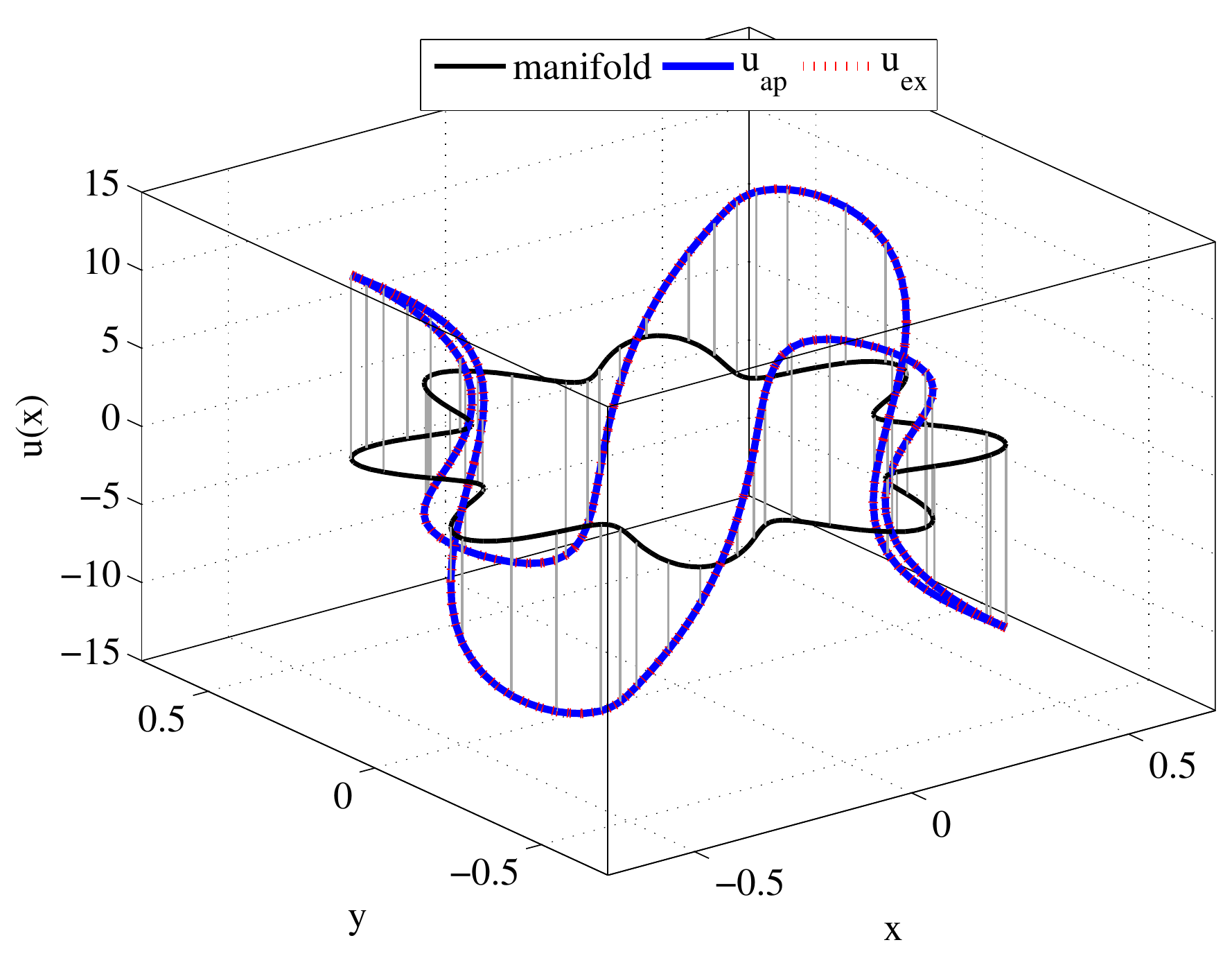}}\subfigure[ratio $h_{\mathrm{max}}/h_{\mathrm{min}}$]{\includegraphics[width=0.3\textwidth]{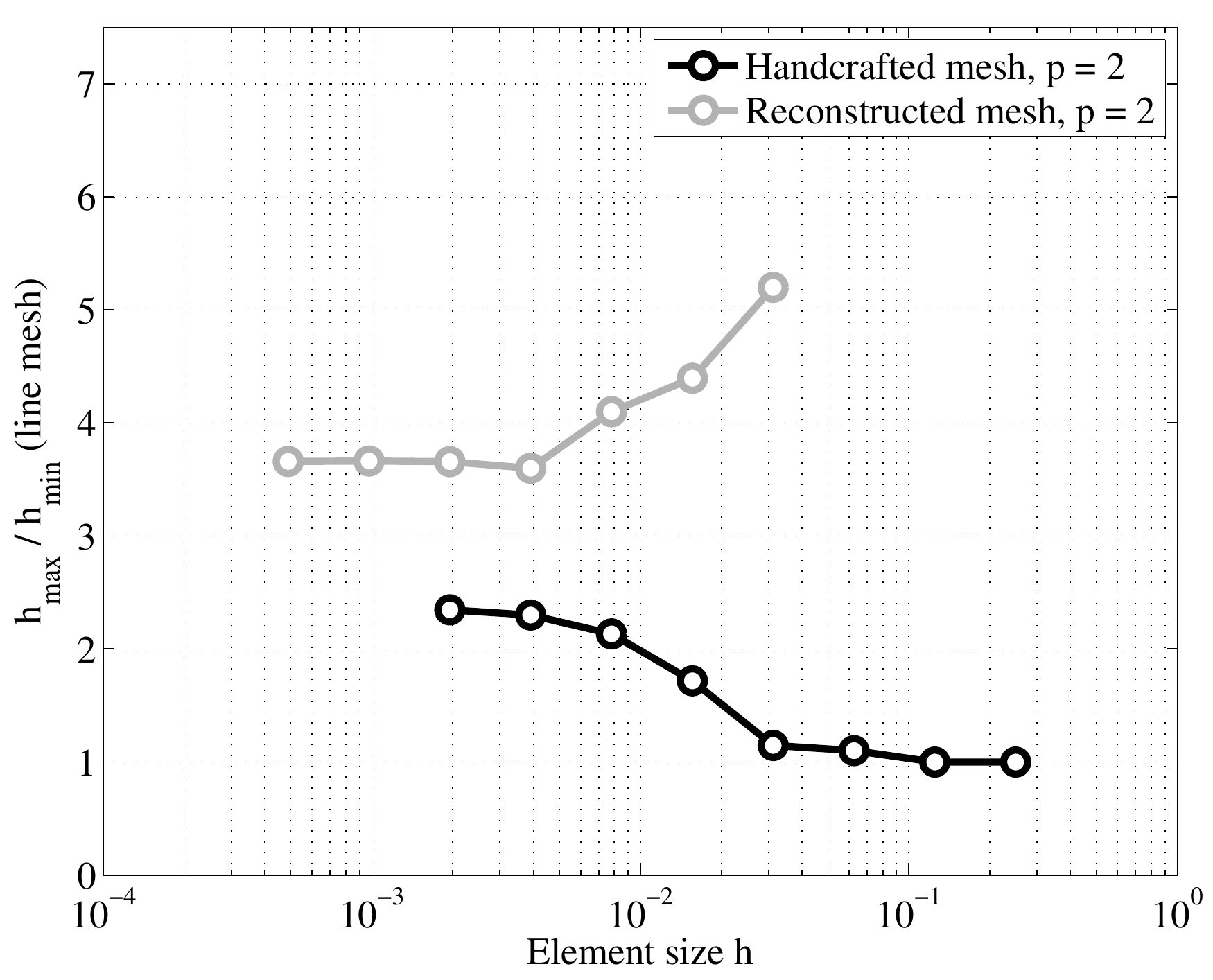}}

\caption{(a) Background and reconstructed line mesh for the flower-shaped manifold,
(b) exact solution of the Poisson equation, (c) the largest ratio
$h_{\mathrm{max}}/h_{\mathrm{min}}$ of the reconstructed line elements.}

\label{fig:PoissonFlowerSetup} 
\end{figure}

Next, the Laplace-Beltrami operator is solved on a closed, flower-shaped
manifold implied by the level-set function
\[
\phi\left(\vek x\right)=\sqrt{x^{2}+y^{2}}-R\left(\theta\right),
\]
with $R\left(\theta\right)=0.5+0.1\cdot\sin\left(8\theta\right)$
and $\theta\left(\vek x\right)=\mathrm{atan}\left(y/x\right)$. The
function $u$ is defined as $u(\theta)=12\sin3\theta$ and the right
hand side is computed accordingly. See Fig.~\ref{fig:PoissonFlowerSetup}(a)
for a coarse background mesh and a reconstructed line mesh of order
$1$. The exact solution is seen in Fig.~\ref{fig:PoissonFlowerSetup}(b)
and the largest ratio $h_{\mathrm{max}}/h_{\mathrm{min}}$ of the
reconstructed line elements resulting from the different background
meshes in Fig.~\ref{fig:PoissonFlowerSetup}(c).

\begin{figure}
\centering

\subfigure[handcr.~mesh, $L_2$-norm]{\includegraphics[width=0.35\textwidth]{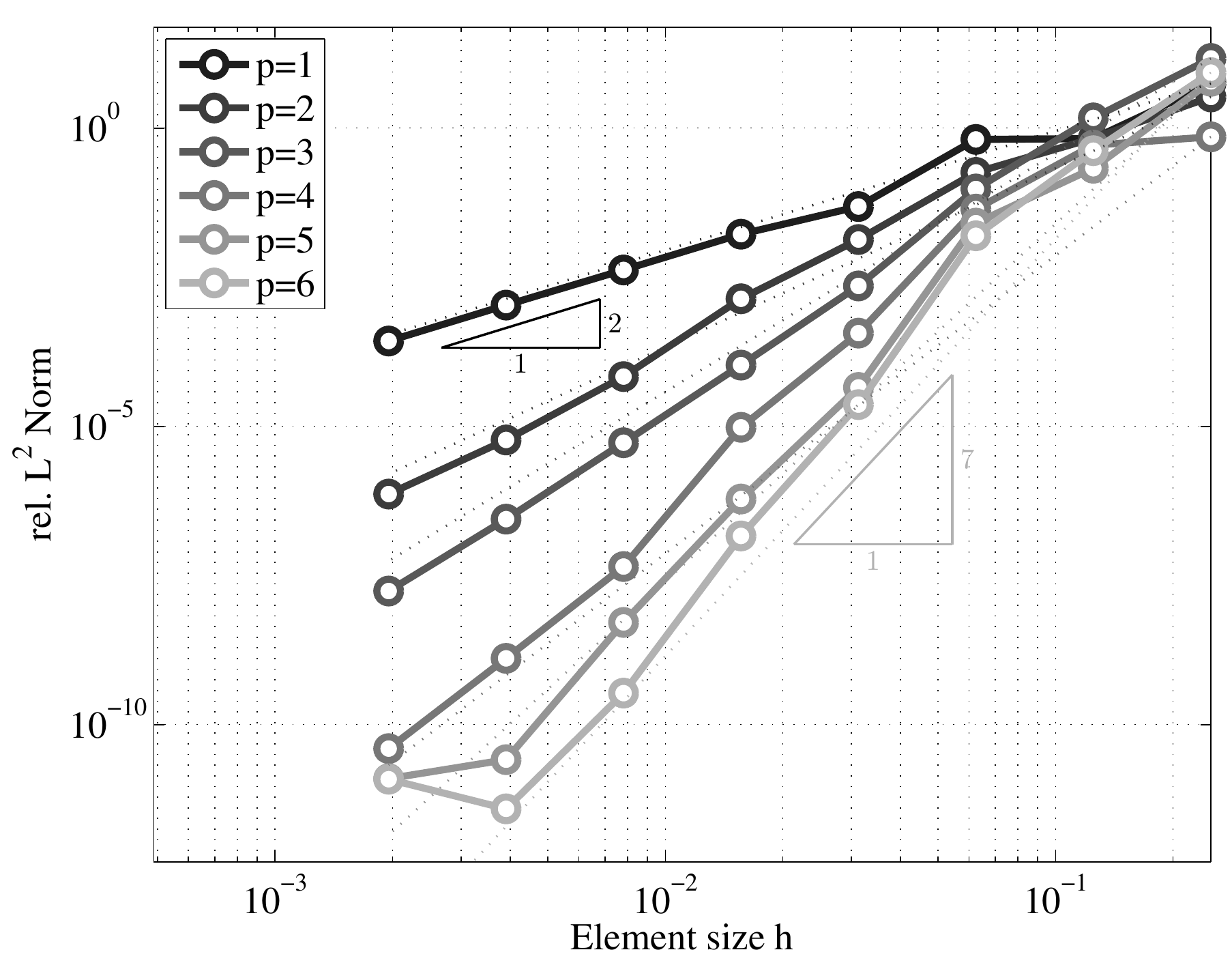}}\quad\subfigure[recon.~mesh, $L_2$-norm]{\includegraphics[width=0.35\textwidth]{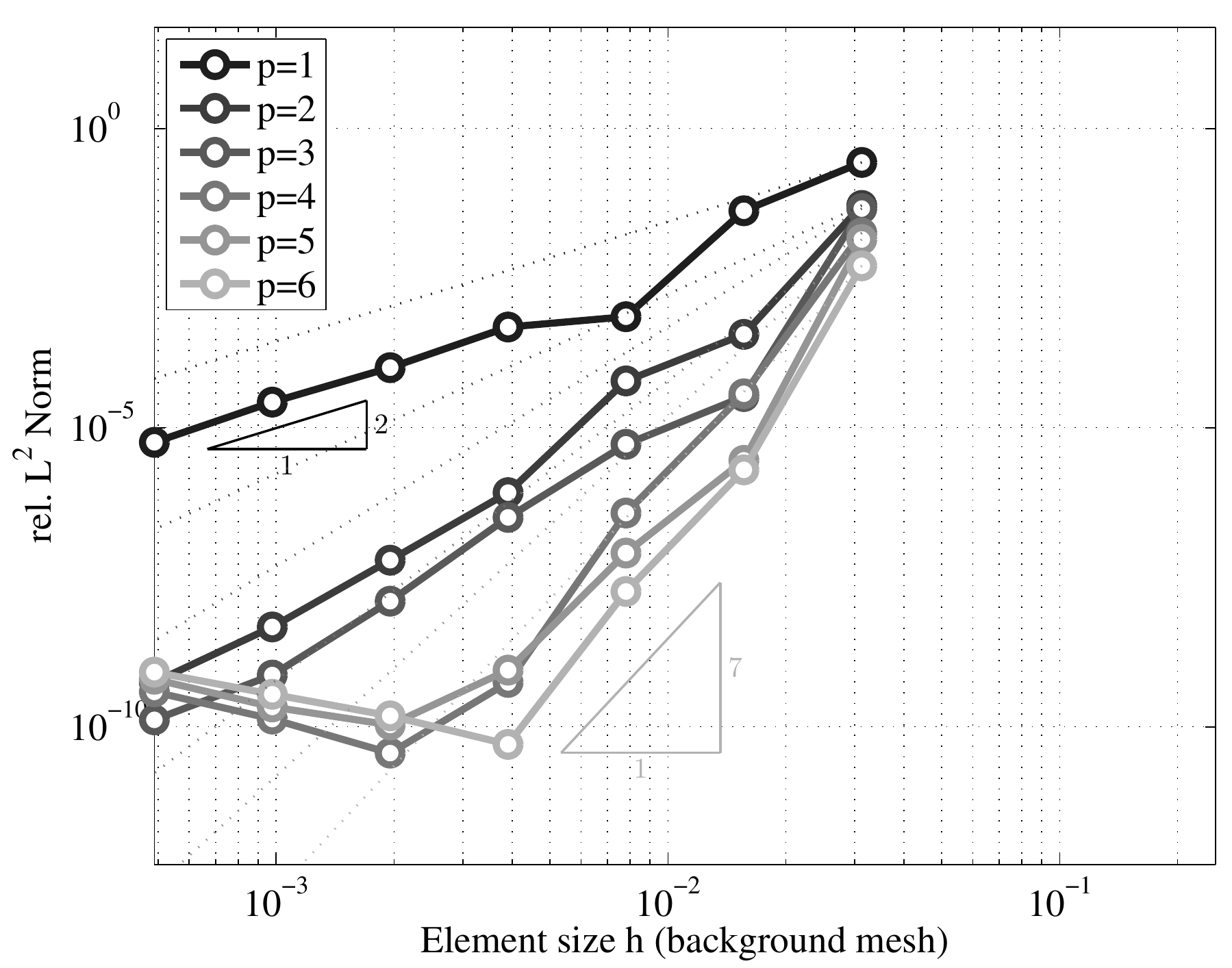}}

\subfigure[handcr.~mesh, cond.]{\includegraphics[width=0.35\textwidth]{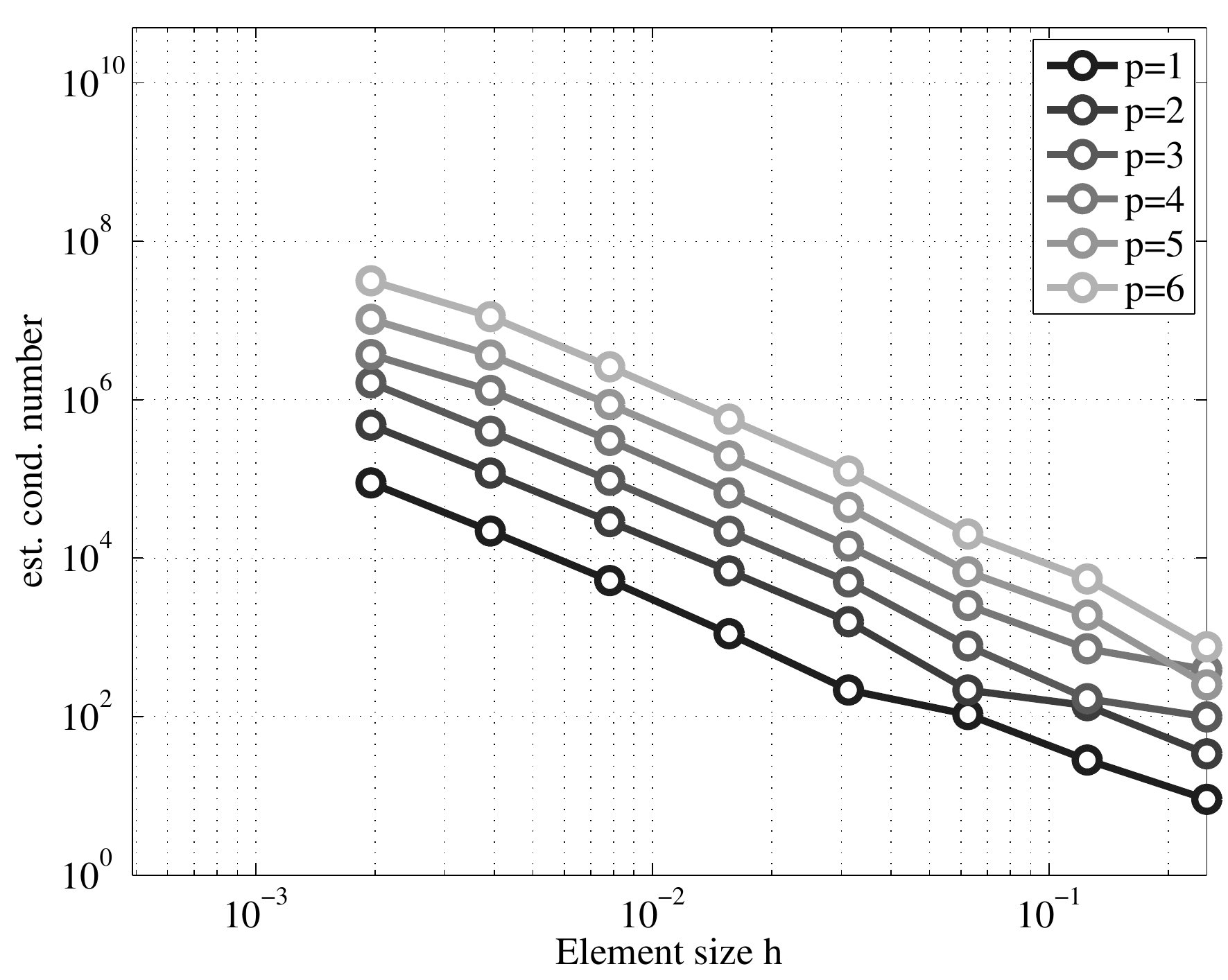}}\quad\subfigure[recon.~mesh, cond.]{\includegraphics[width=0.35\textwidth]{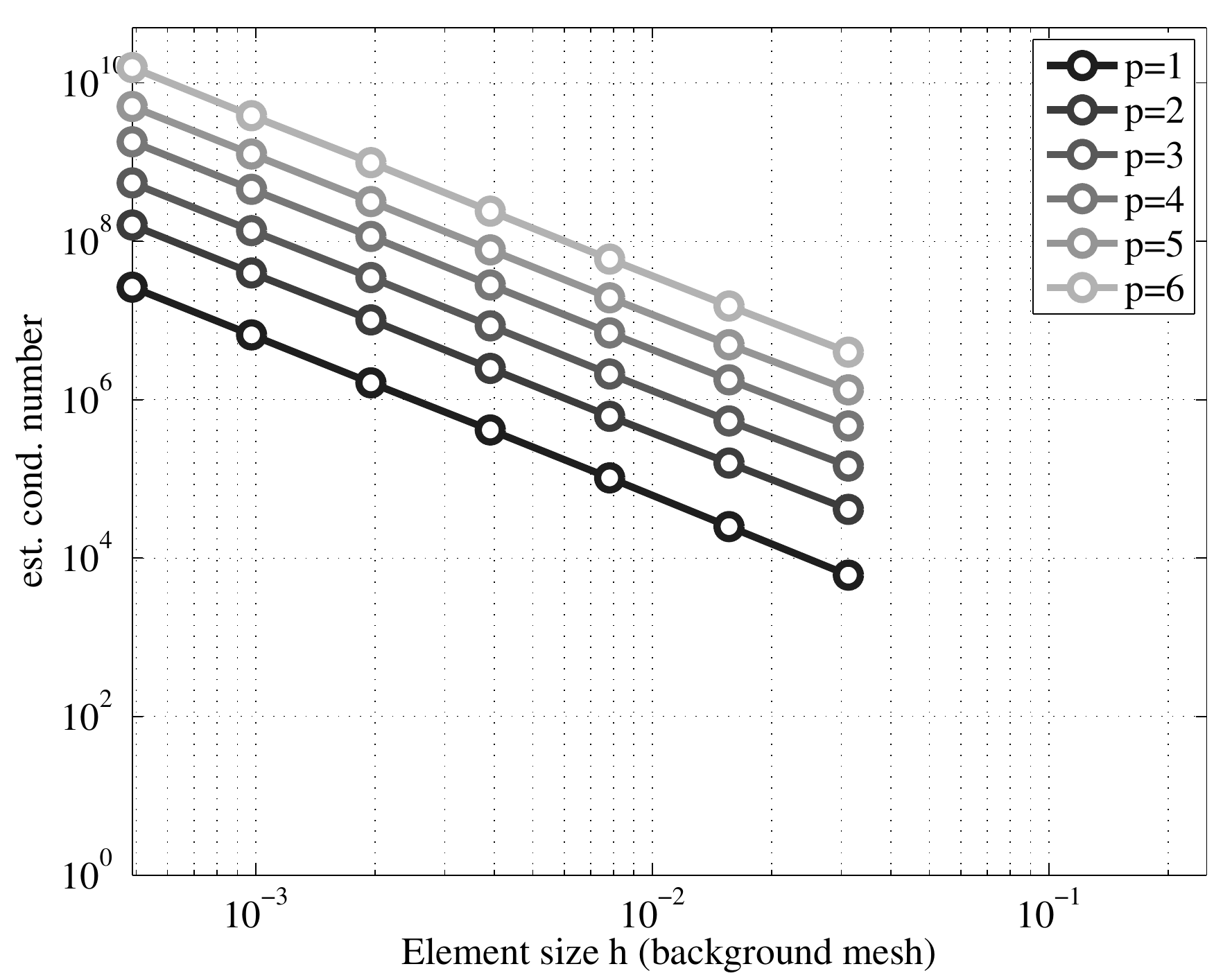}}

\caption{\label{fig:PoissonFlowerRes}Convergence results and condition numbers
for the Poisson equation on a flower-shaped manifold.}
\end{figure}

In the convergence studies, we use background meshes with $h=R/\left\{ 32,\,64,\,128,\,256,\,512,\,1024,\,2048\right\} $
which is fine enough to ensure that the reconstruction is successfully
achieved without any (recursive/adaptive) refinements of the background
elements. For the handcrafted meshes, significantly larger element
lengths $h$ may be used. The results presented in Fig.~\ref{fig:PoissonFlowerRes}
follow the same style than in the previous test case. It is confirmed
that optimal convergence rates are achieved on handcrafted and automatically
reconstructed meshes and the condition numbers behave as expected.

\subsubsection{S-shaped manifold\label{sec:lap1dEx3}}

\begin{figure}
\centering

\subfigure[meshes]{\includegraphics[width=0.3\textwidth]{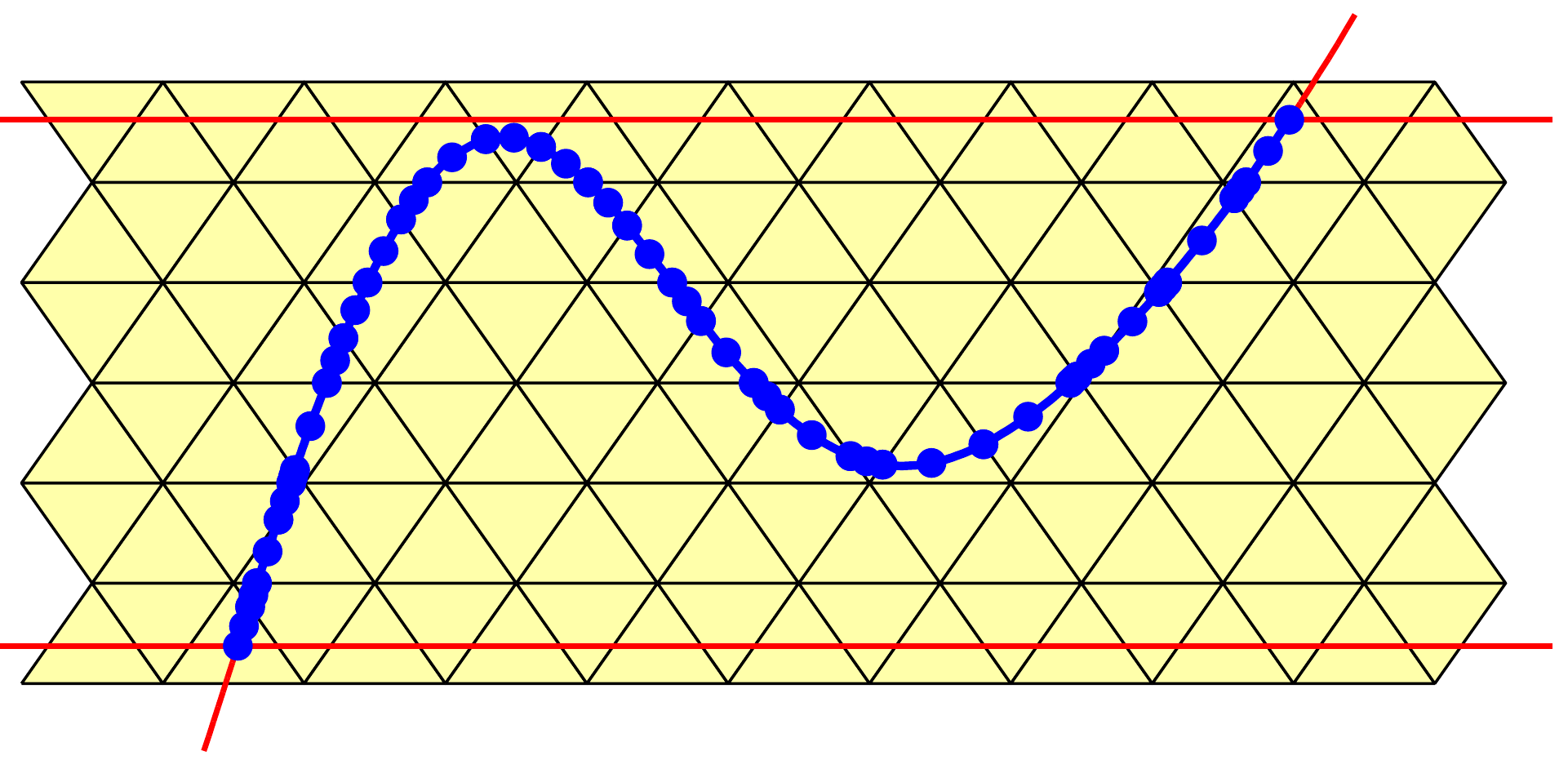}}\subfigure[exact solution]{\includegraphics[width=0.3\textwidth]{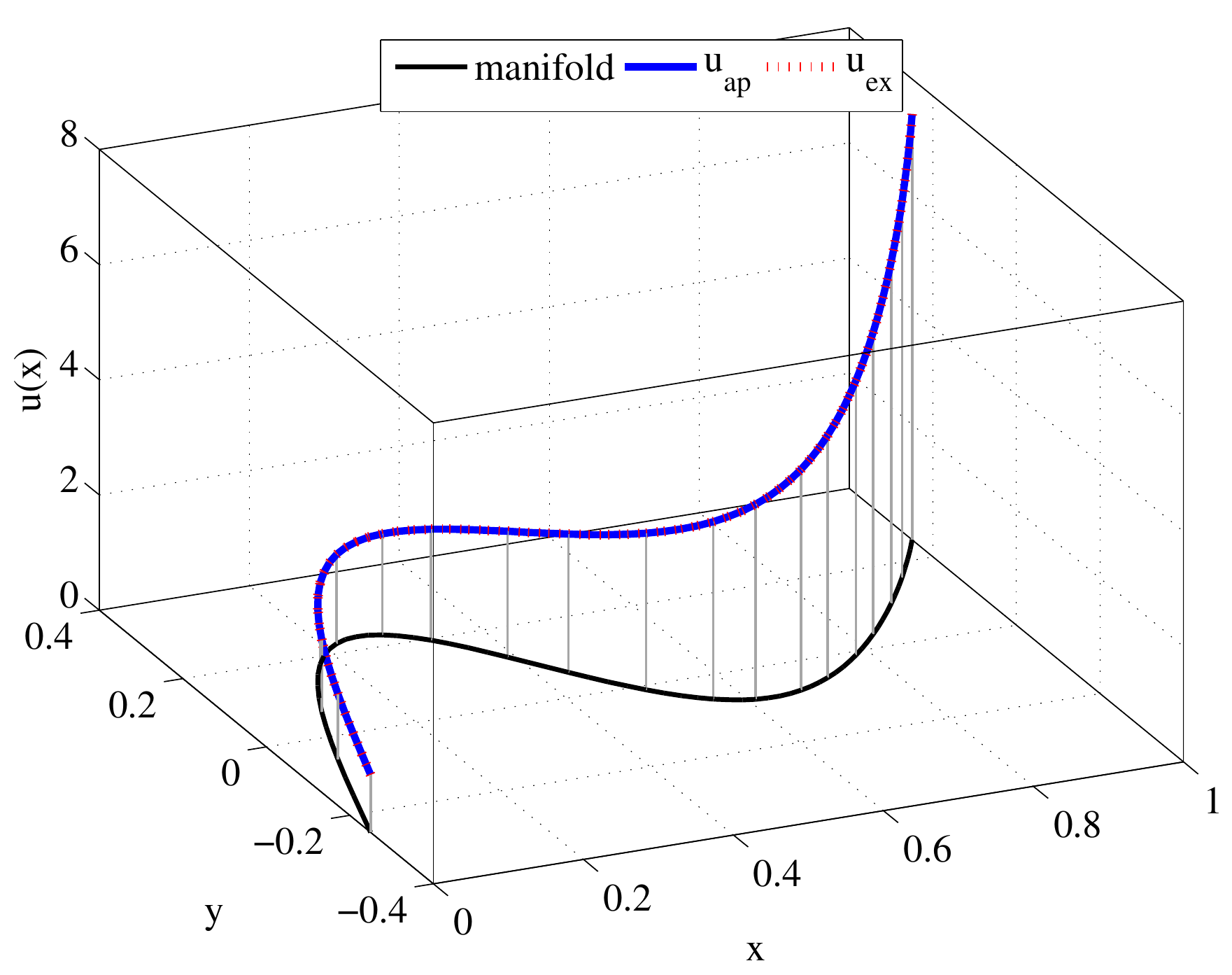}}\subfigure[ratio $h_{\mathrm{max}}/h_{\mathrm{min}}$]{\includegraphics[width=0.3\textwidth]{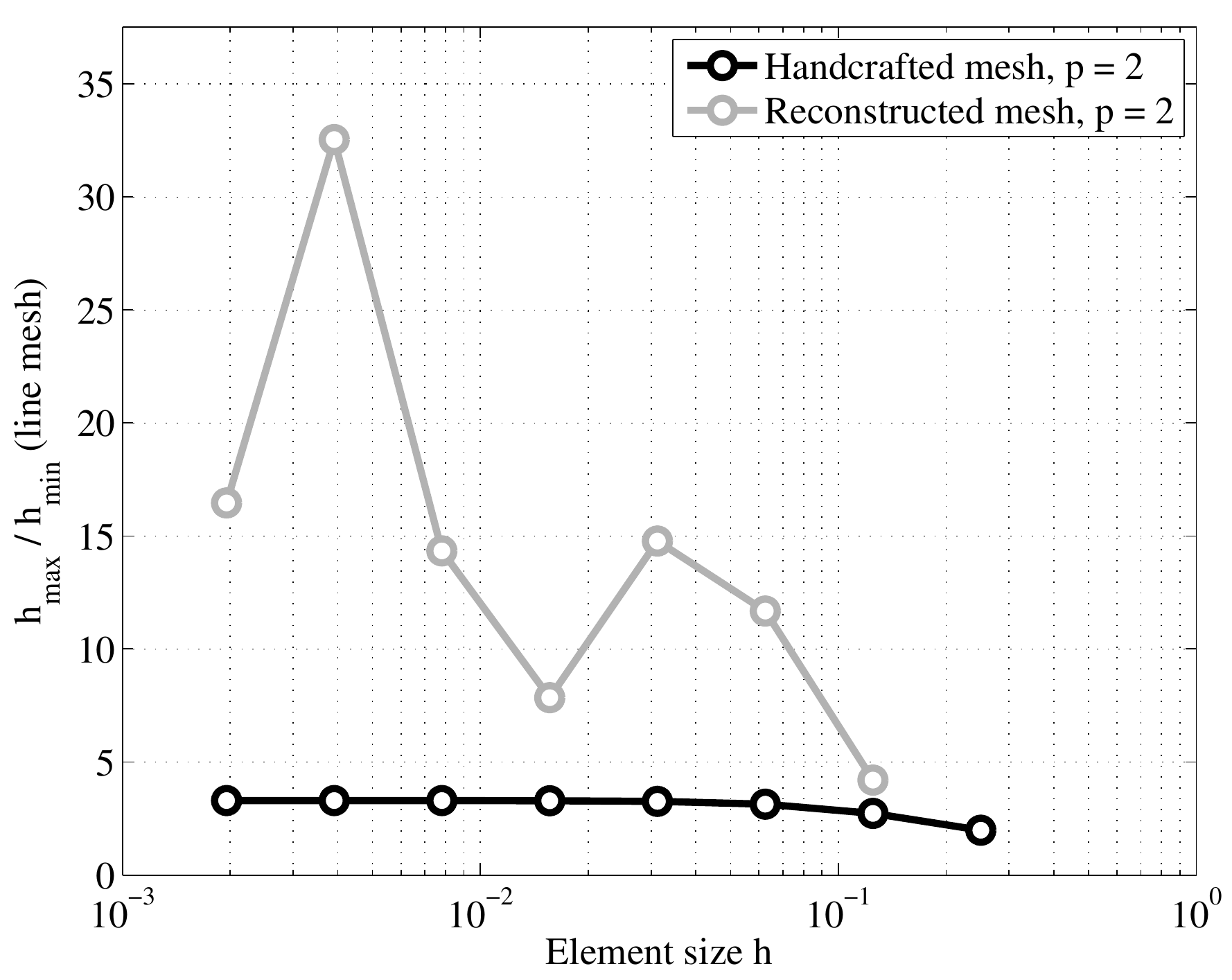}}

\caption{(a) Background and reconstructed line mesh for the S-shaped manifold,
(b) exact solution of the Poisson equation, (c) the largest ratio
$h_{\mathrm{max}}/h_{\mathrm{min}}$ of the reconstructed line elements.}

\label{fig:PoissonSLineSetup} 
\end{figure}

In this example, we consider an \emph{open} manifold where the boundary
is defined by additional level-set functions $\psi^{i}\left(\vek x\right)$.
The manifold in parametric form is given as 
\begin{align}
\Gamma=\begin{bmatrix}x\\
f(x)
\end{bmatrix}=\begin{bmatrix}x\\
\frac{x^{3}}{2}+\sin(\pi(1-x))\sin^{5}\left(\frac{\pi}{2}(1-x)\right)-\frac{1}{4}
\end{bmatrix}, & \quad\mathrm{with}\; x\in\left[0,\,1\right].\label{eq:sline}
\end{align}
The corresponding level-set function may be defined as 
\begin{align}
\phi(\vek x)=f(x)-y\ .
\end{align}
The boundaries are defined by the two additional level-set functions
$\psi^{1}(\vek x)=\nicefrac{1}{4}-y$ and $\psi^{2}(\vek x)=y-\nicefrac{1}{4}$.
An example for a reconstructed line mesh is seen in Fig.~\ref{fig:PoissonSLineSetup}(a)
where also the zero-isolines of all three level-set functions are
shown. The exact solution $u(\vek x)=\exp\left(2x\right)$ is plotted
in Fig.~\ref{fig:PoissonSLineSetup}(b). Again, the ratio of the
largest element length of the automatically reconstructed line mesh
to the smallest is shown in Fig.~\ref{fig:PoissonSLineSetup}(c).
This ratio may be worse than in the previous examples due to the presence
of several level-set functions which renders the manipulation of the
background mesh more complicated.

\begin{figure}
\centering

\subfigure[handcr.~mesh, $L_2$-norm]{\includegraphics[width=0.35\textwidth]{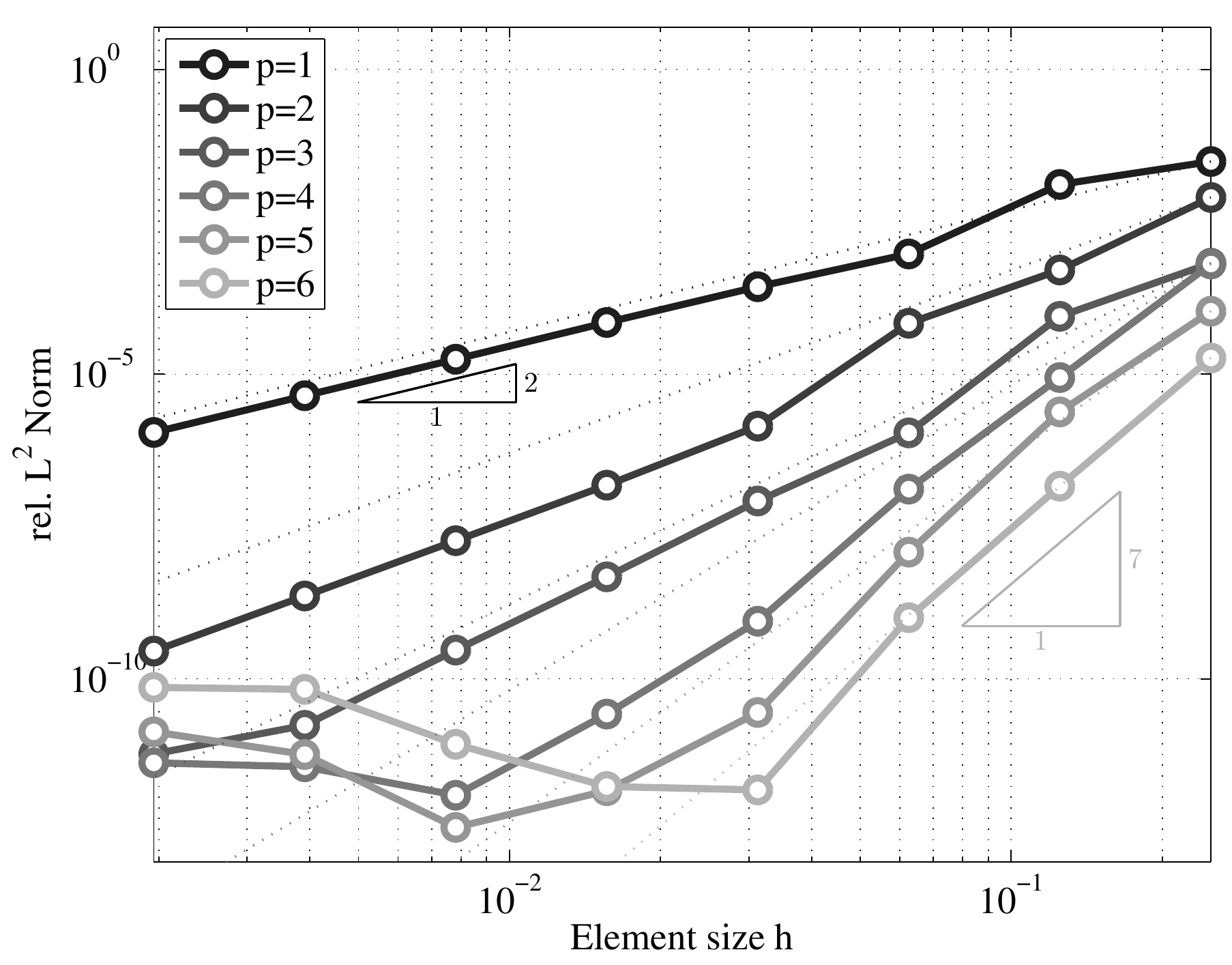}}\quad\subfigure[recon.~mesh, $L_2$-norm]{\includegraphics[width=0.35\textwidth]{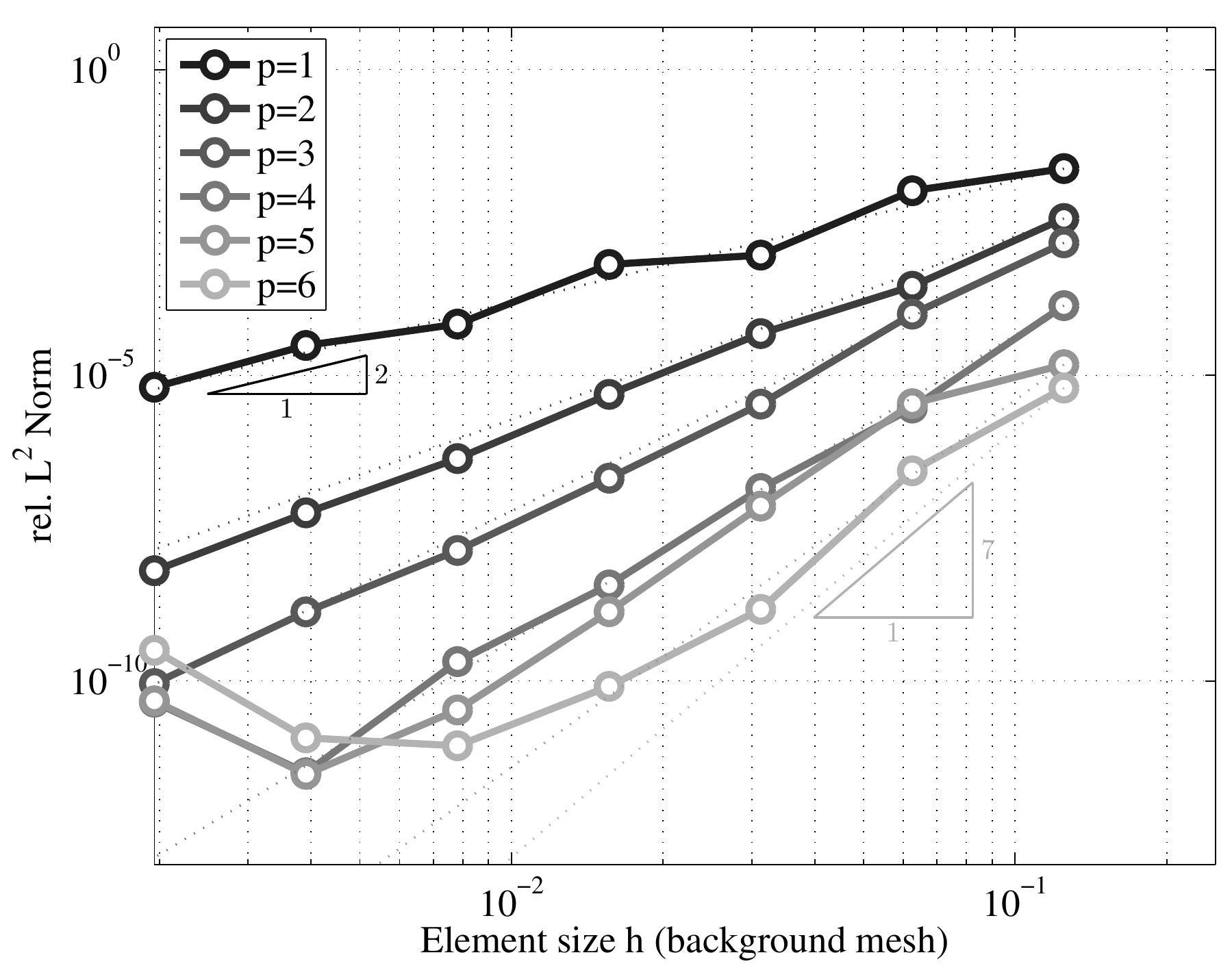}}

\subfigure[handcr.~mesh, cond.]{\includegraphics[width=0.35\textwidth]{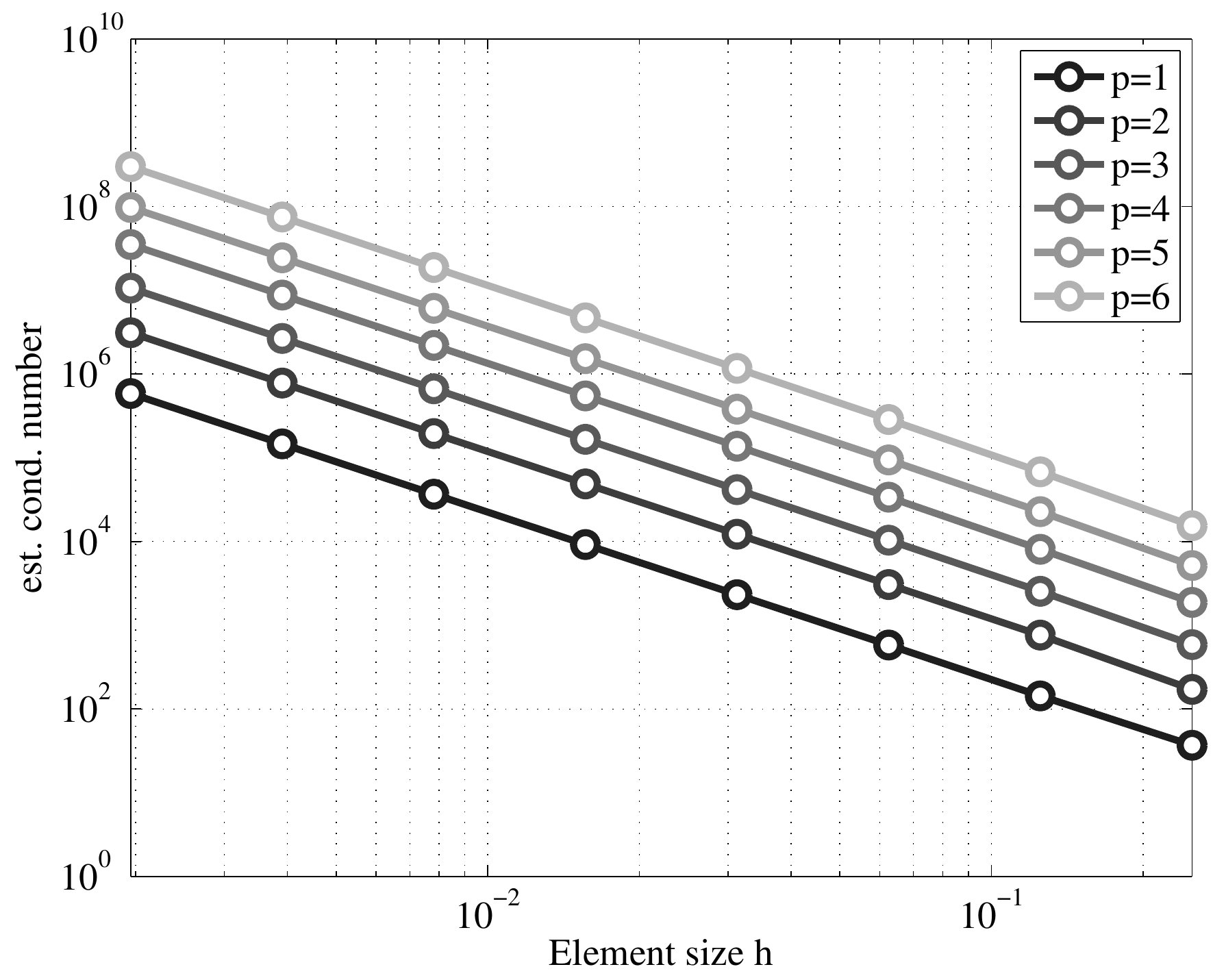}}\quad\subfigure[recon.~mesh, cond.]{\includegraphics[width=0.35\textwidth]{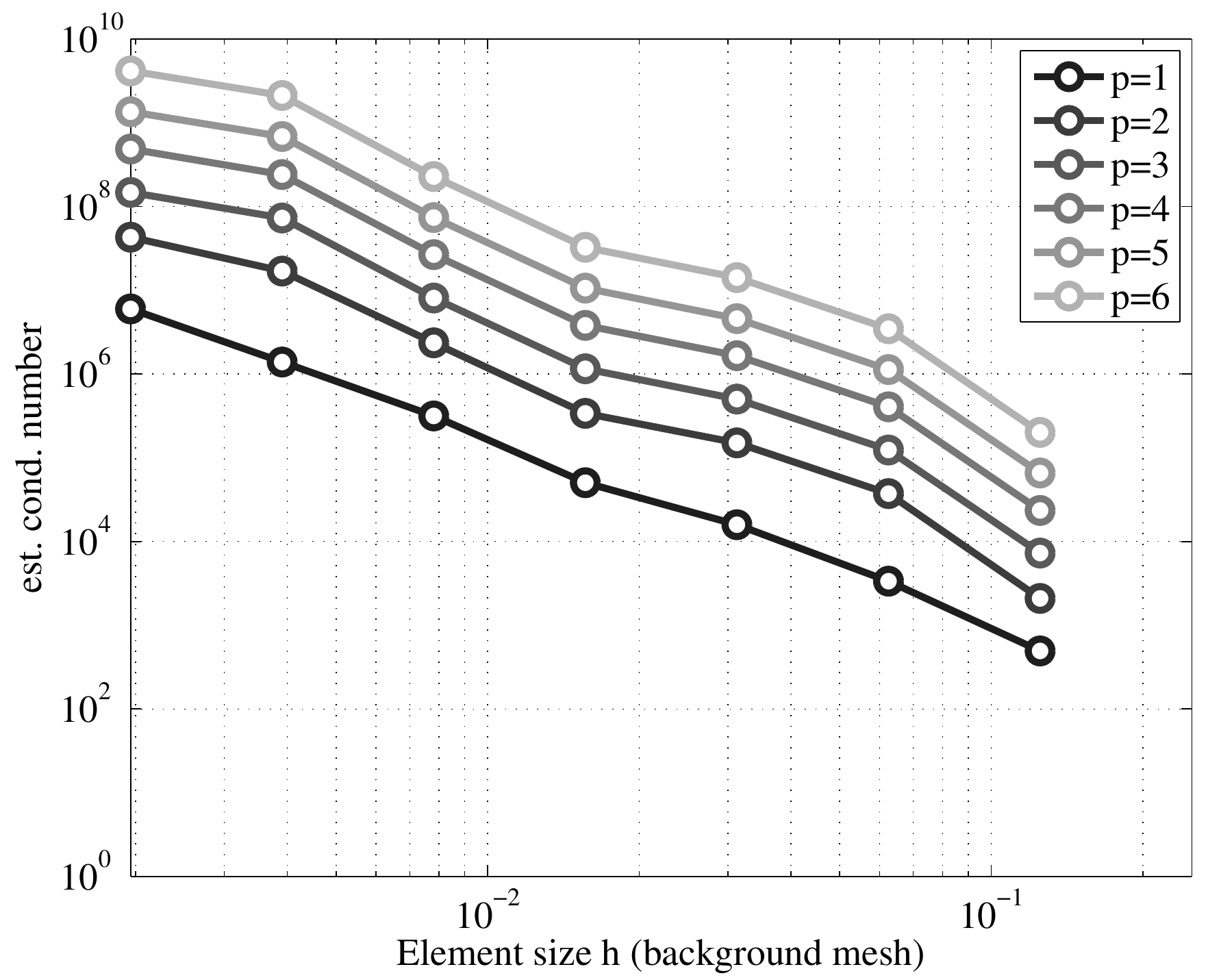}}

\caption{\label{fig:PoissonSLineRes}Convergence results and condition numbers
for the Poisson equation on a S-shaped manifold.}
\end{figure}

In Fig.~\ref{fig:PoissonSLineRes}, the results of the convergence
analysis is presented in the same style than above with the same conclusions
to be drawn. Hence, we conclude that approximations of the Poisson
equation with optimal convergence rates are possible for the automatically
generated meshes discretizing closed and open curved lines in 2D.

\subsection{Poisson equation on curved surfaces in 3D}

\subsubsection{Quarter Cylinder\label{sec:lap2dEx1}}

Here, the Laplace-Beltrami operator is solved on the surface of a
quarter cylinder with the radius $r=1$ and height $L=4$. The analytical
solution in cylindrical coordinates $(r,\,\varphi,\, z)$ is given
in \cite{Kamilis_2013} based on the following parameters and functions:
\begin{align}
g_{\varphi,1}(\varphi) & =\left(1-\cos\varphi\right)\left(1-\sin\varphi\right),\\
g_{\varphi,2}(\varphi) & =\left(\cos\varphi+\sin\varphi-4\sin\varphi\cos\varphi\right),\\
g_{z}(z) & =\sin\dfrac{\alpha\pi z}{L},\\
\alpha & =3\,,\qquad\beta=\dfrac{1}{\frac{3}{2}-\sqrt{2}}\ .
\end{align}
The analytic solution is then defined as $u(r,\,\varphi,\, z)=\beta g_{\varphi,1}(\varphi)g_{z}(z)$.
Applying the Laplace-Beltrami operator yields the source function
\[
f=\beta\cdot g_{z}(z)\left[\left(\dfrac{\alpha\pi}{L}\right)^{2}g_{\varphi,1}(\varphi)-g_{\varphi,2}(\varphi)\right].
\]
The Dirichlet boundary condition is $u\rvert_{\partial\Gamma_{\text{D}}}=0$.
In this special case, it is particularly simple to use a tensor-product
background mesh in the volume $\Omega=(0,r)\times(0,r)\times(0,L)$
to reconstruct the higher-order surface mesh with the desired boundaries.
Hence, the use of addional level-set functions $\psi^{i}$ would unnecessarily
complicate the situation here.

\begin{figure}
\centering

\subfigure[Surface mesh]{\includegraphics[height=5cm]{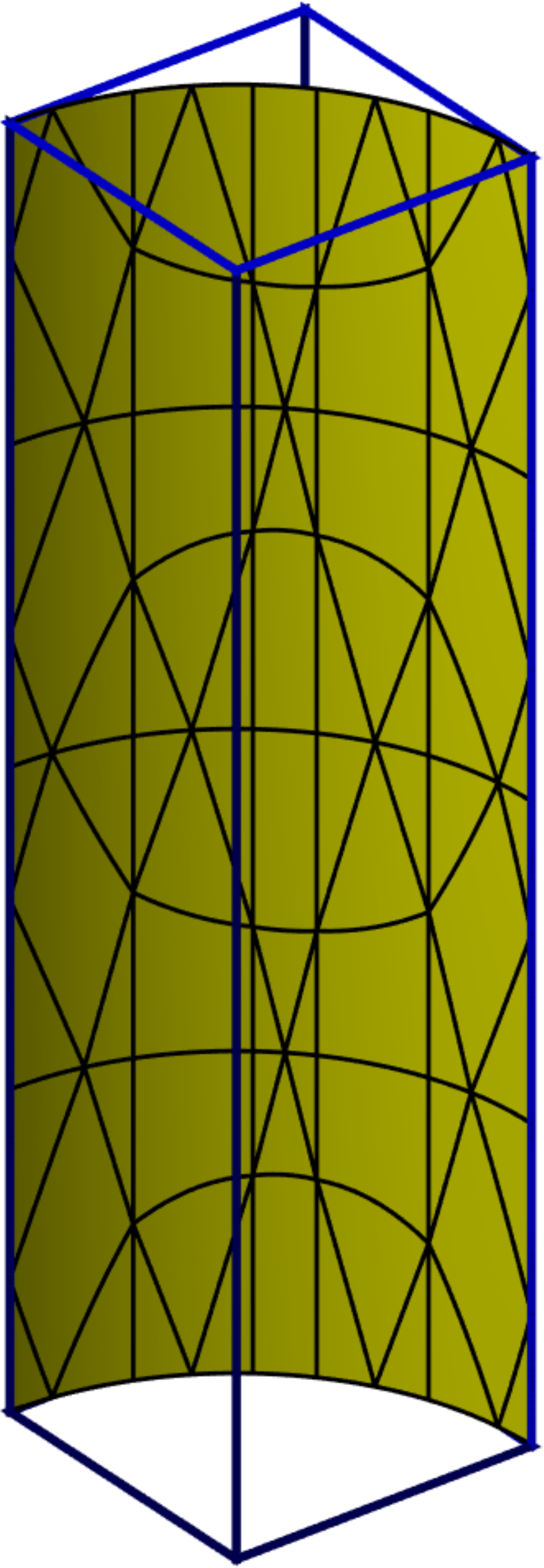}}\subfigure[exact solution]{\includegraphics[height=5cm]{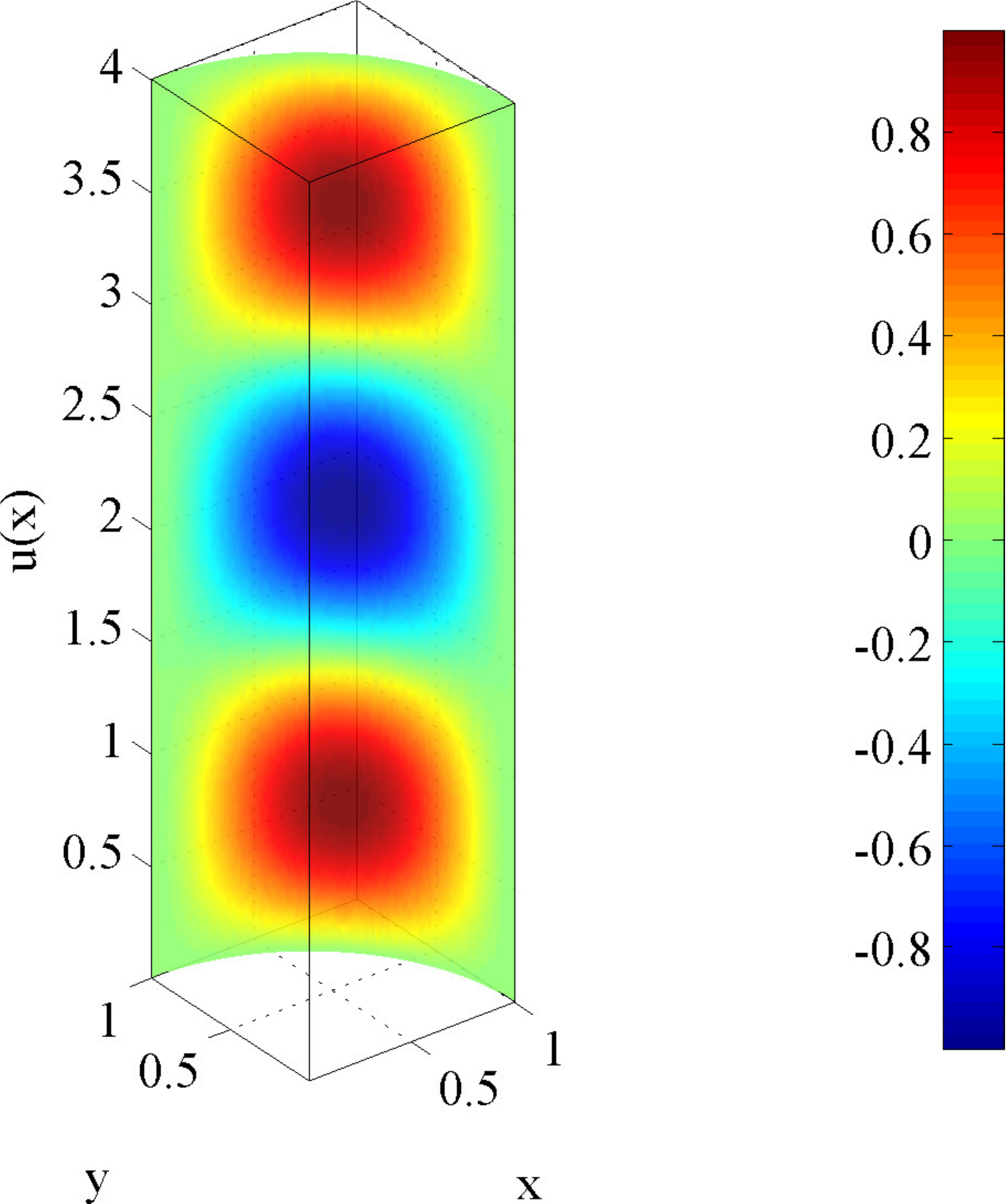}}\subfigure[ratio $A_{\mathrm{max}}/A_{\mathrm{min}}$]{\includegraphics[width=0.3\textwidth]{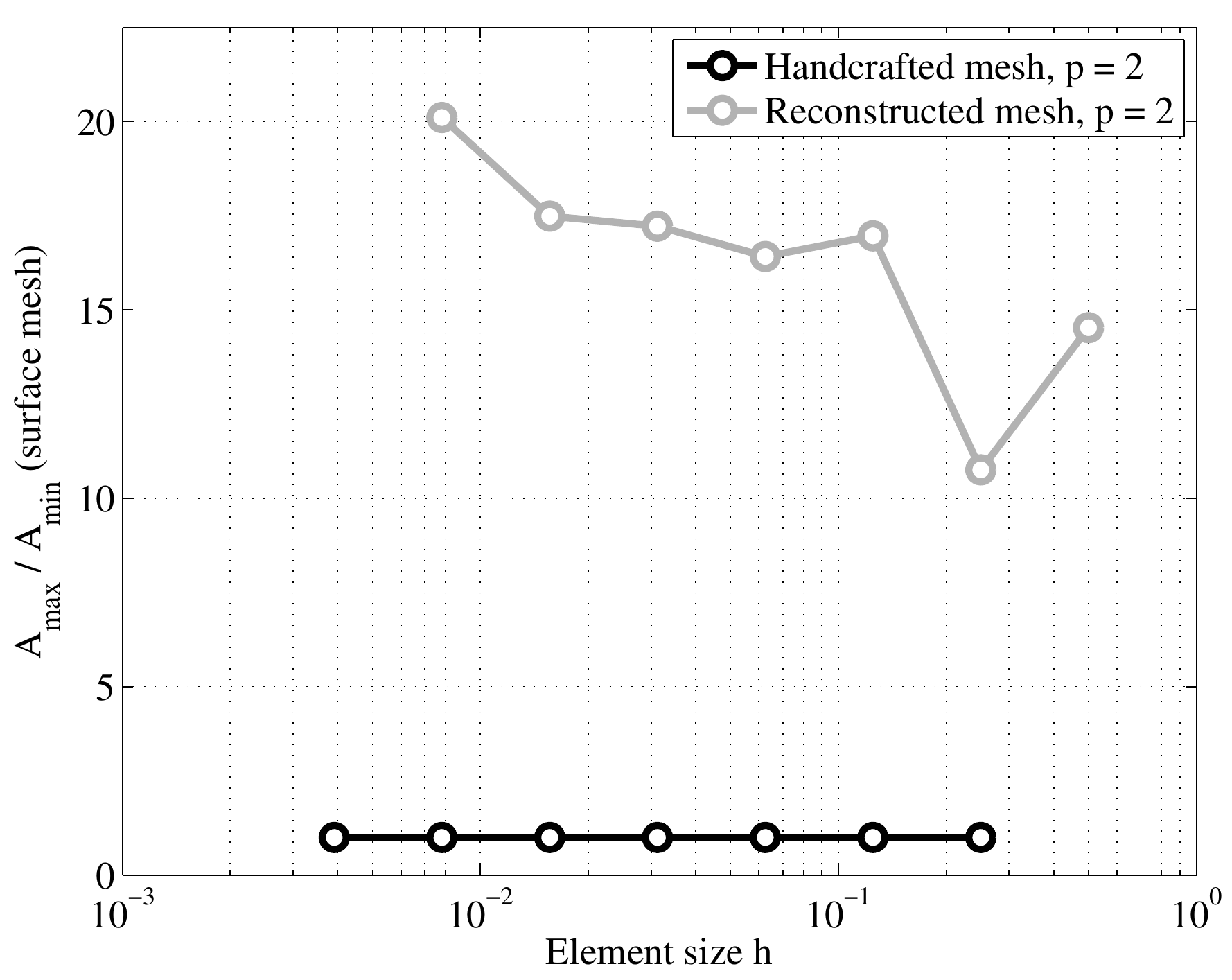}}

\caption{(a) Reconstructed surface mesh for the quarter cylinder, (b) exact
solution of the Poisson equation, (c) the largest ratio $A_{\mathrm{max}}/A_{\mathrm{min}}$
of the reconstructed surface elements.}

\label{fig:PoissonQuartCylSetup} 
\end{figure}

An example of a reconstructed, higher-order surface mesh is shown
in Fig.~\ref{fig:PoissonQuartCylSetup}(a). The exact solution is
given in Fig.~\ref{fig:PoissonQuartCylSetup}(b). The ratio of the
largest surface element to the smallest is given in Fig.~\ref{fig:PoissonQuartCylSetup}(c)
and is, again, clearly bounded (in the range of $<20$ for this example).

For the convergence studies we use background meshes with $h=r/\left\{ 4,\,8,\,16,\,32,\,64,\,128\right\} $
and results are shown in Fig.~\ref{fig:PoissonSLineRes}. As for
the previous convergence studies for curved \emph{lines}, Fig.~\ref{fig:PoissonSLineRes}(a)
and (c) refer to convergence rates and condition numbers of handcrafted
\emph{surface} meshes, respectively. Fig.~\ref{fig:PoissonSLineRes}(b)
and (d) give the results of the automatically generated surface meshes.
The resulting convergence plots show optimal convergence rates which
are surprisingly smooth given the irregularity of the reconstructed
surface elements. The condition numbers are about two orders of magnitude
larger than in the handcrafted meshes.

\begin{figure}
\centering

\subfigure[handcr.~mesh, $L_2$-norm]{\includegraphics[width=0.35\textwidth]{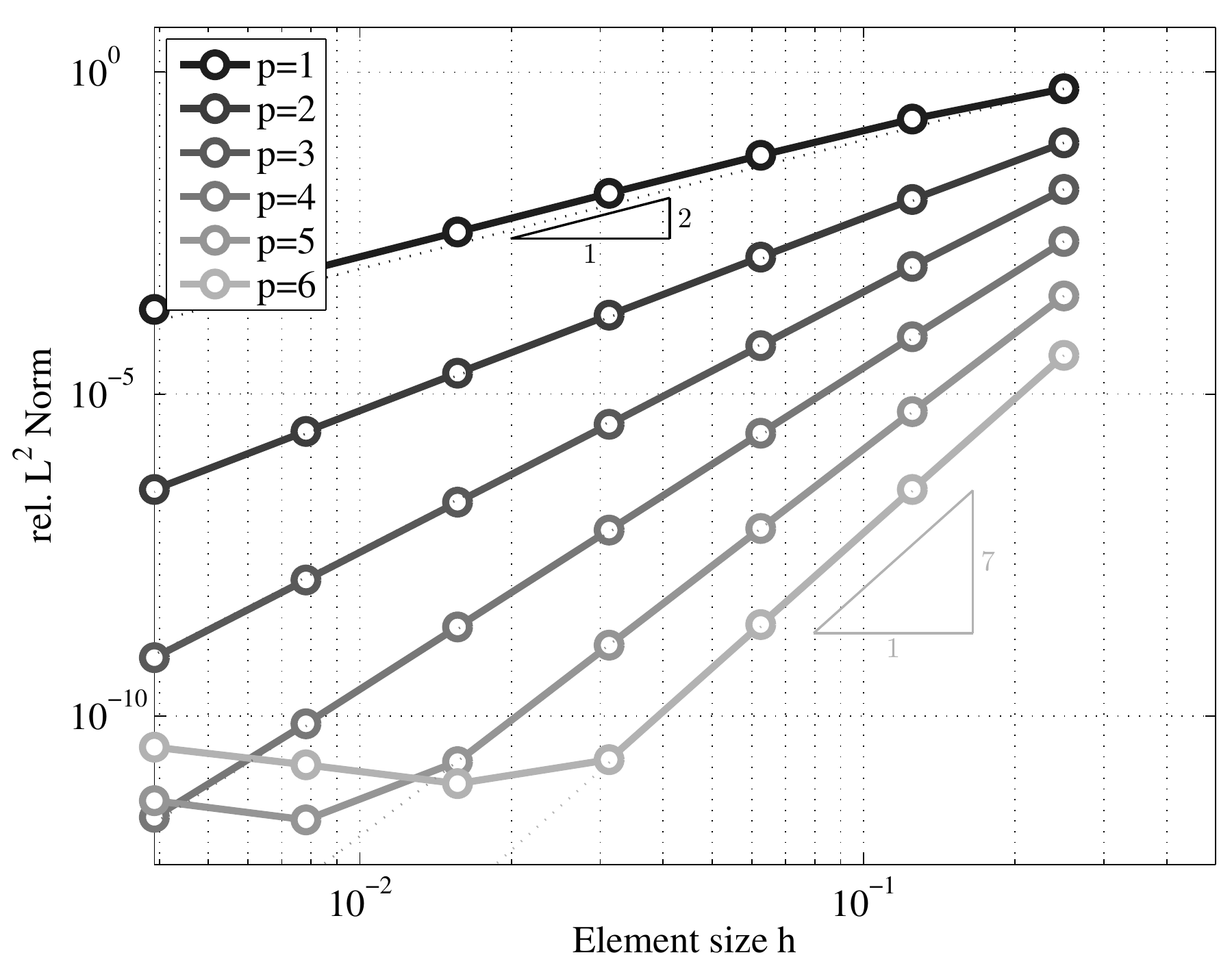}}\quad\subfigure[recon.~mesh, $L_2$-norm]{\includegraphics[width=0.35\textwidth]{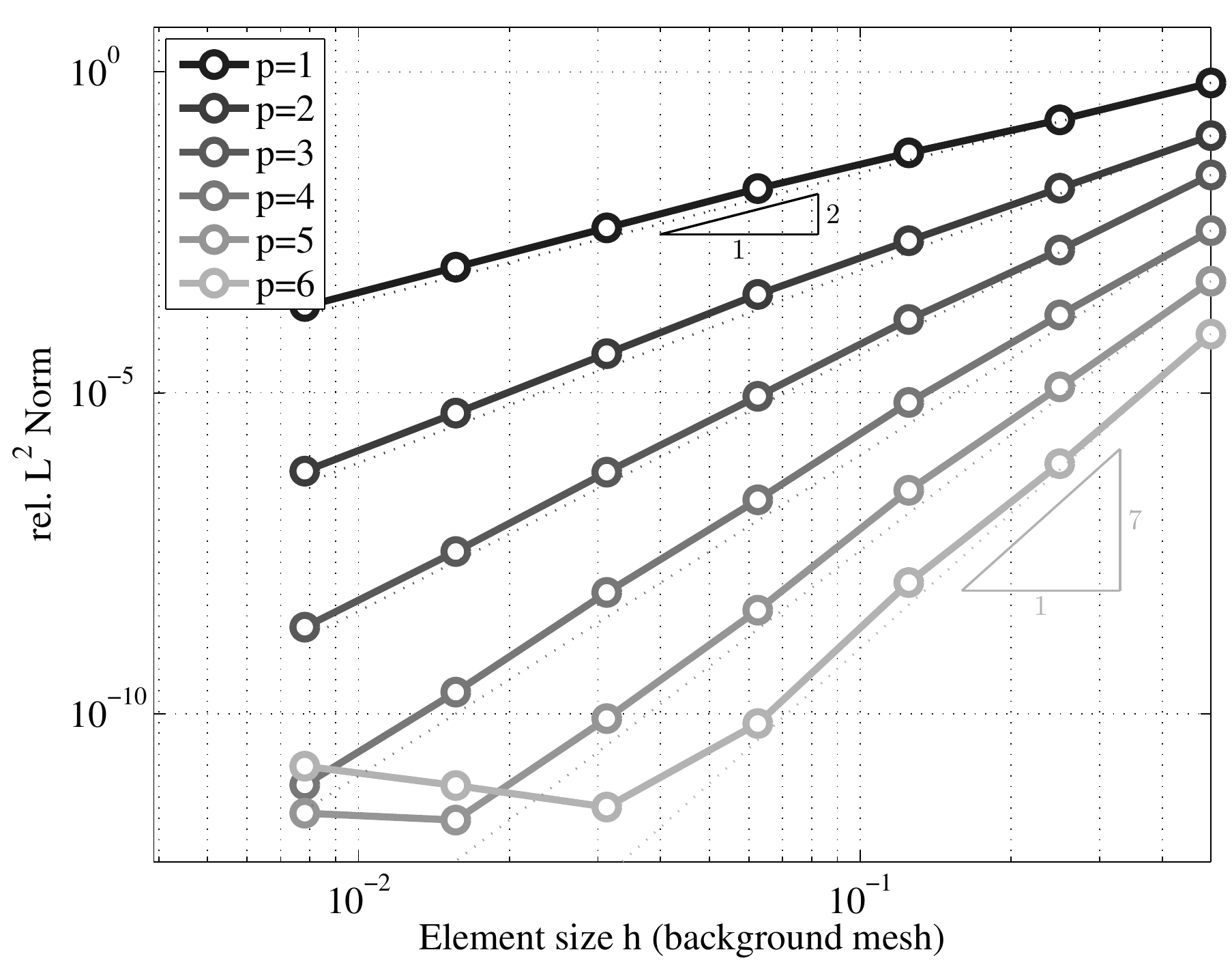}}

\subfigure[handcr.~mesh, cond.]{\includegraphics[width=0.35\textwidth]{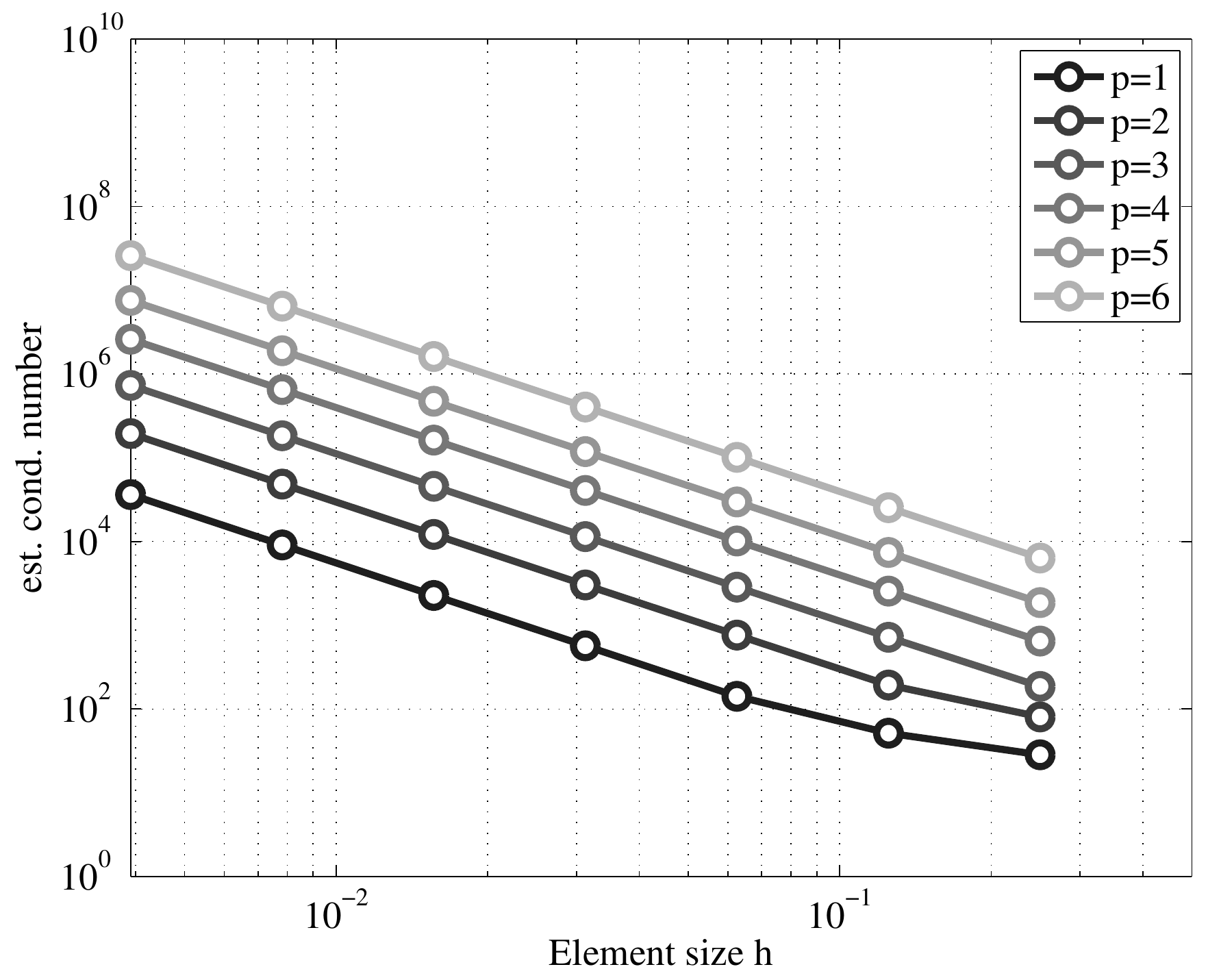}}\quad\subfigure[recon.~mesh, cond.]{\includegraphics[width=0.35\textwidth]{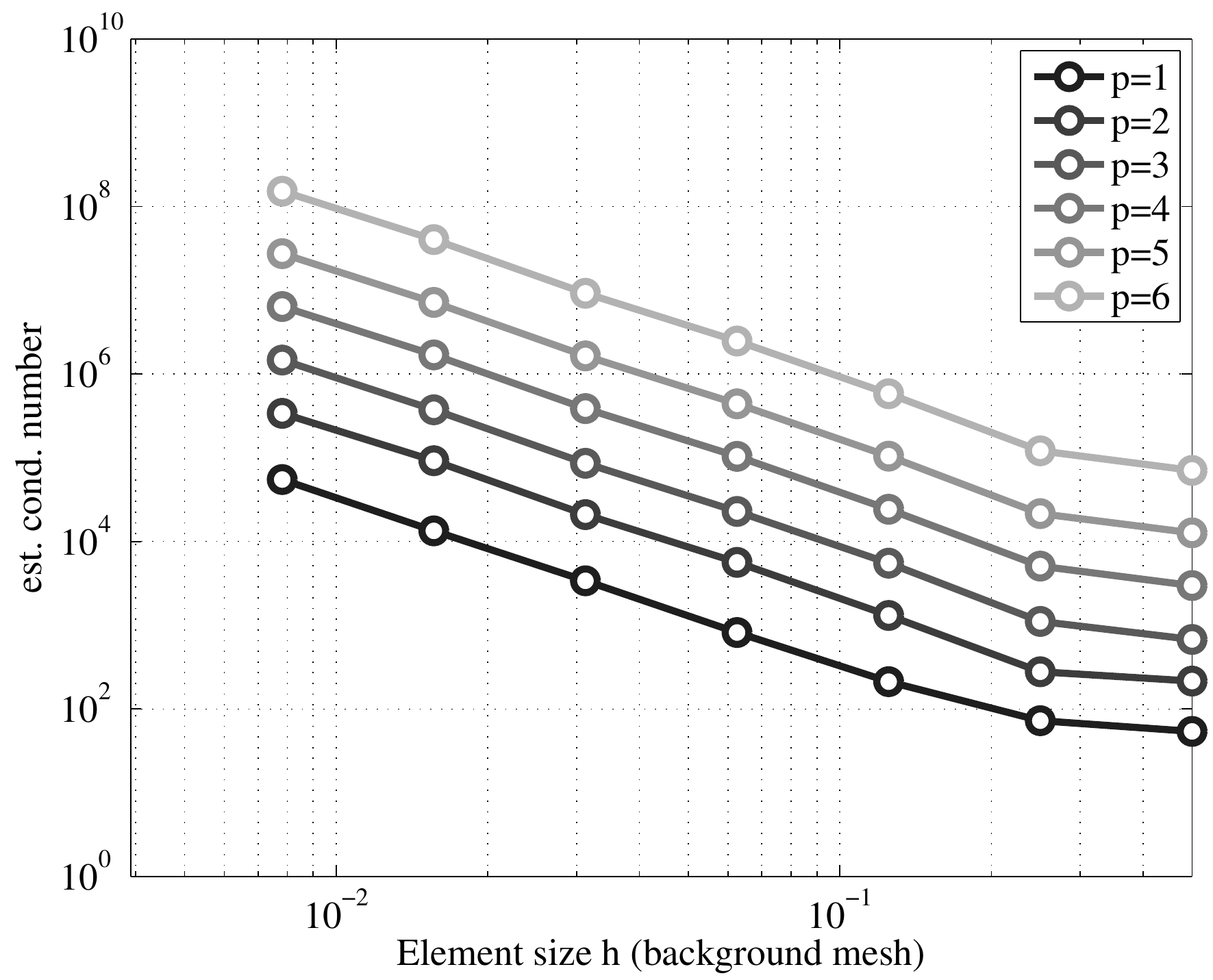}}

\caption{\label{fig:PoissonQuartCylRes}Convergence results and condition numbers
for the Poisson equation on a quarter cylinder.}
\end{figure}

\subsubsection{Sphere\label{sec:lap2dEx2}}

\begin{figure}
\centering

\subfigure[Surface mesh]{\includegraphics[height=4cm]{PlotManifolds3dB_1}}\subfigure[exact solution]{\includegraphics[height=5cm]{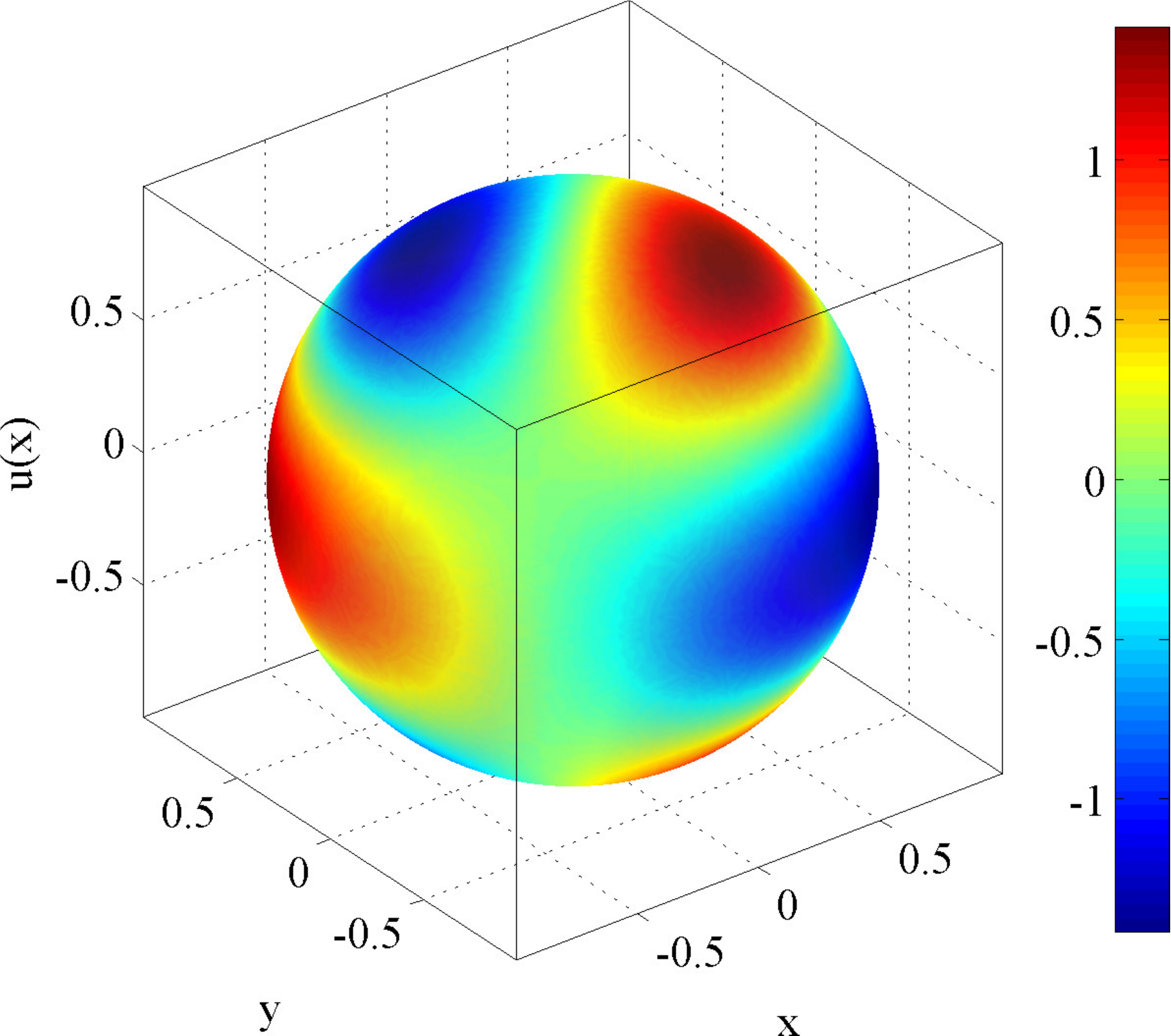}}\subfigure[ratio $A_{\mathrm{max}}/A_{\mathrm{min}}$]{\includegraphics[width=0.3\textwidth]{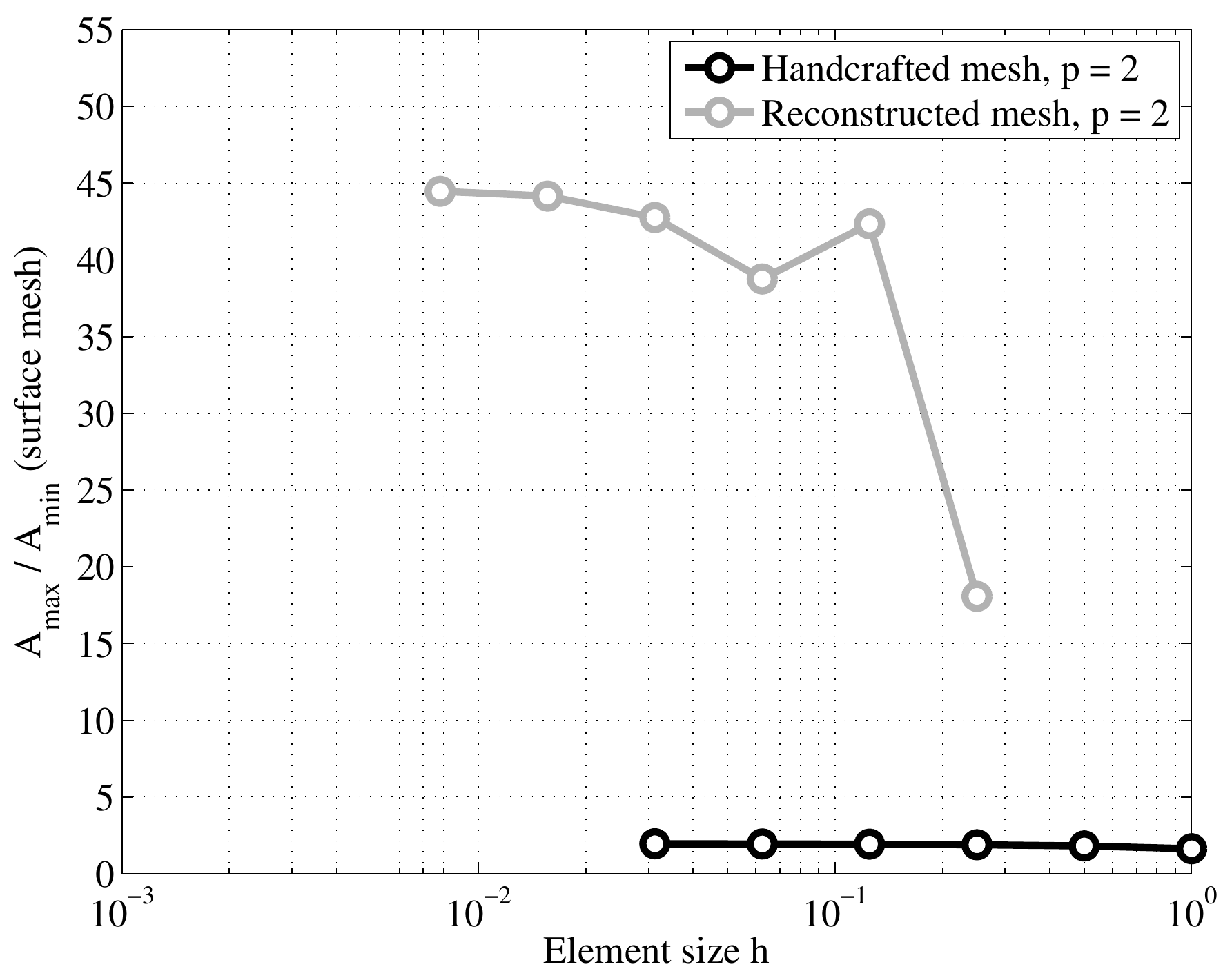}}

\caption{(a) Reconstructed surface mesh for the sphere, (b) exact solution
of the Poisson equation, (c) the largest ratio $A_{\mathrm{max}}/A_{\mathrm{min}}$
of the reconstructed surface elements.}

\label{fig:PoissonSphereSetup} 
\end{figure}

\begin{figure}
\centering

\subfigure[recon.~mesh, $L_2$-norm]{\includegraphics[width=0.35\textwidth]{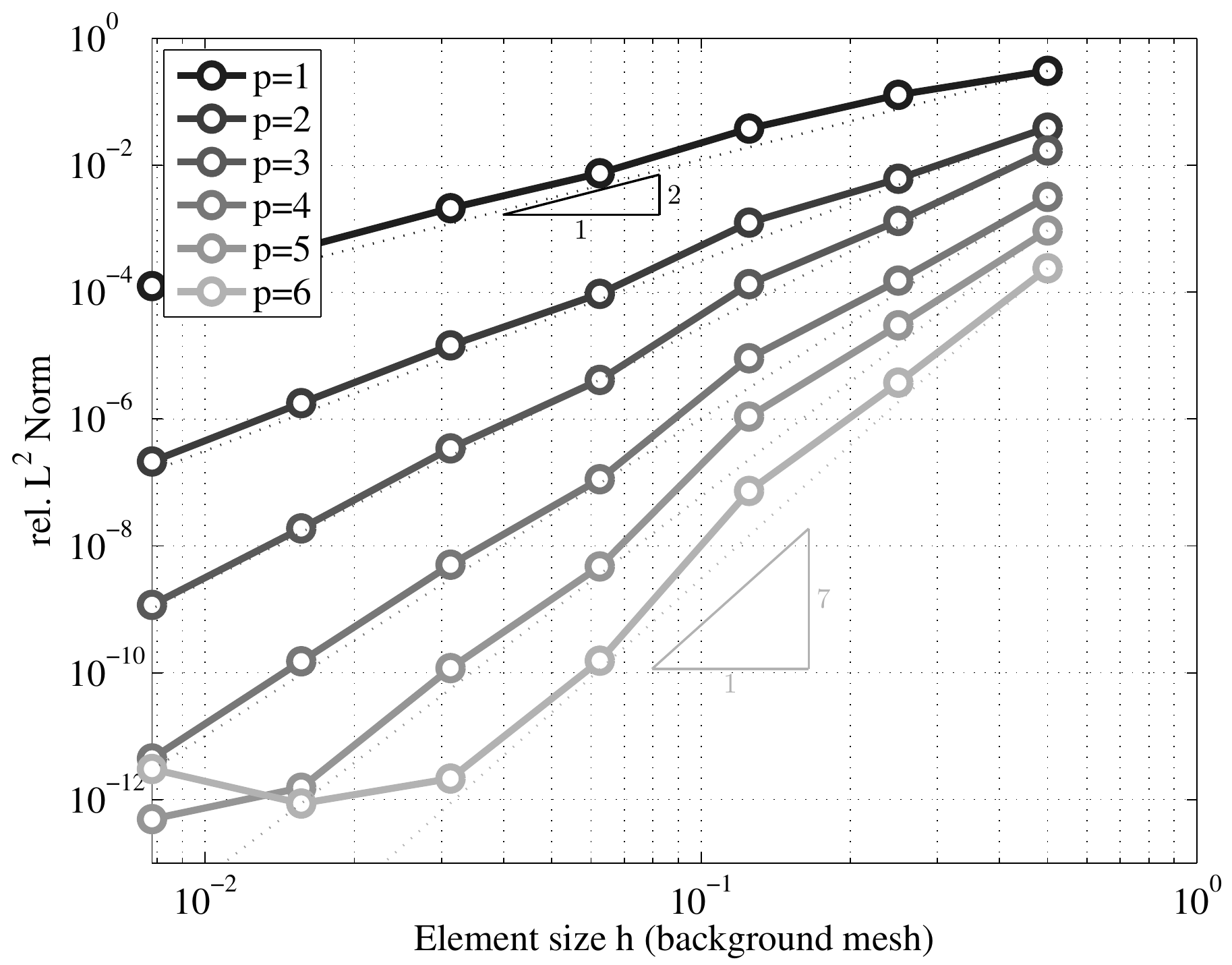}}\quad\subfigure[recon.~mesh, cond.]{\includegraphics[width=0.35\textwidth]{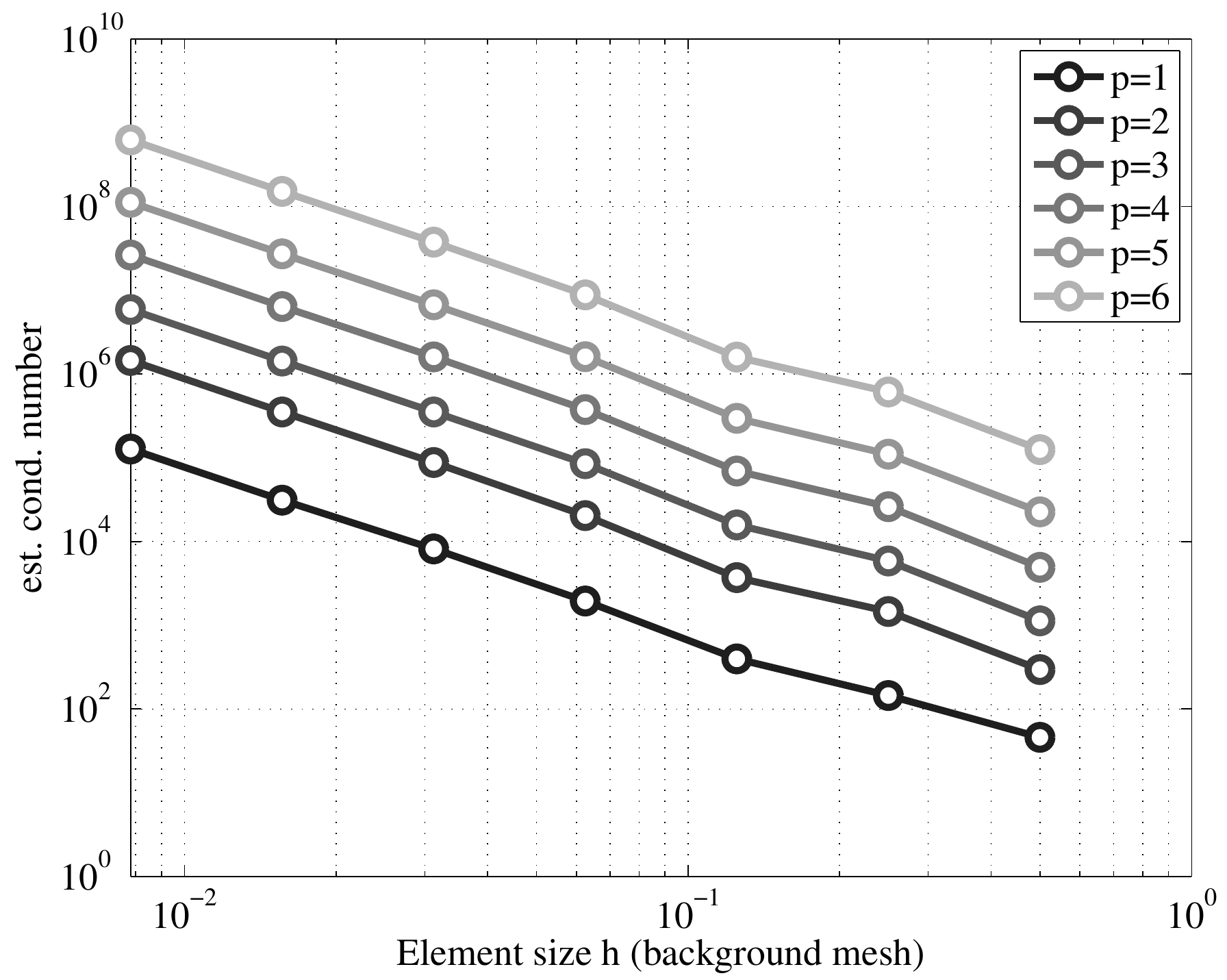}}

\caption{\label{fig:PoissonSphereRes}Convergence results and condition numbers
for the Poisson equation on a sphere.}
\end{figure}

\begin{figure}
\centering

\includegraphics[width=0.35\textwidth]{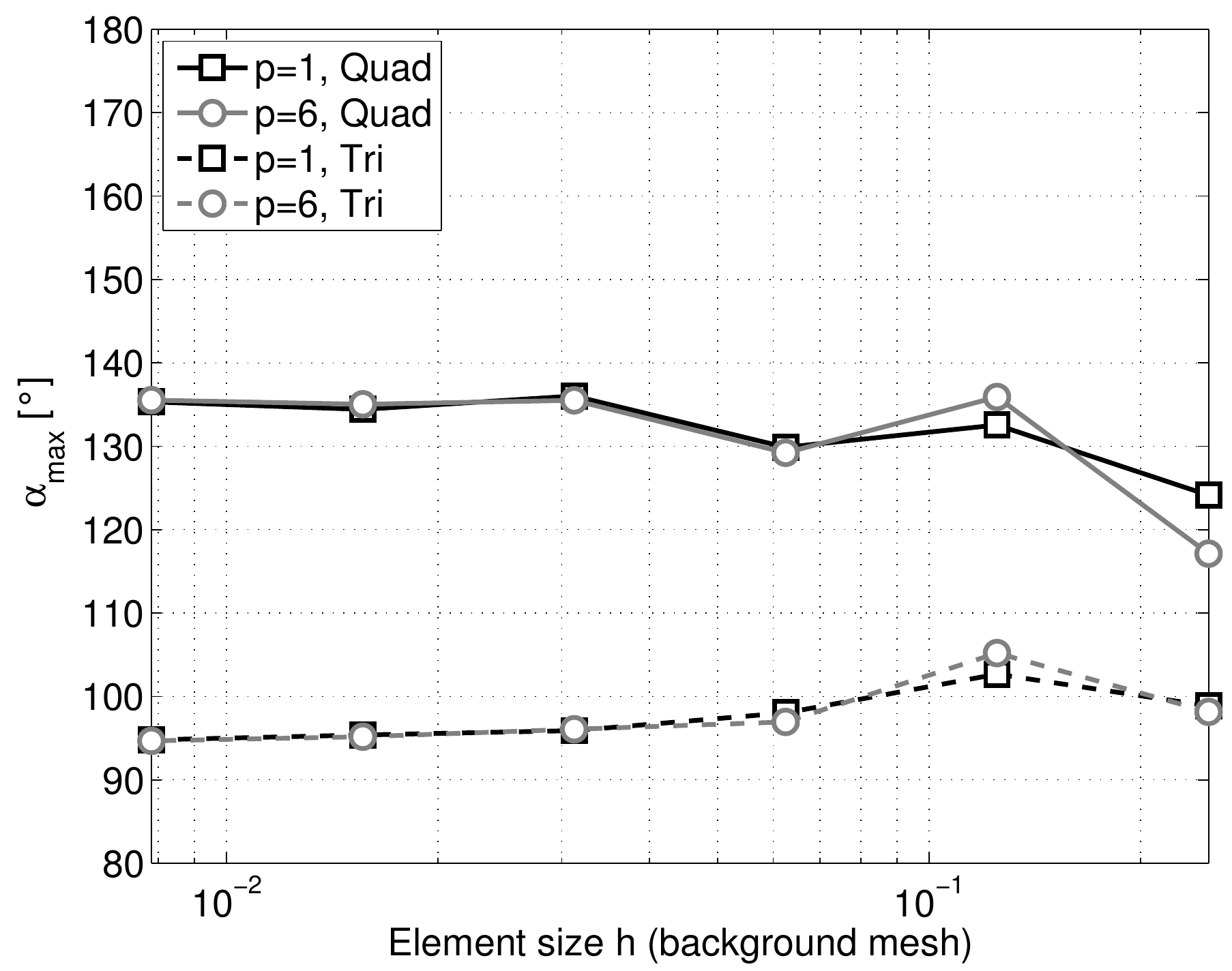}

\caption{\label{fig:PoissonSphereAngles}Maximum inner angles of the automatically
generated meshes of the sphere.}
\end{figure}

Next, we consider a sphere with radius $r=1$. Using spherical coordinates
$(r,\,\varphi,\,\theta)$, an exact solution $u(r,\,\varphi,\,\theta)$
fulfilling the zero-mean constraint from Eq.~(\ref{eq:const}) is
chosen as 
\begin{align}
u(r,\,\varphi,\,\theta)=\sin(3\theta)\left(\cos\varphi-\sin\varphi\right)\ .
\end{align}
Applying the Laplace-Beltrami operator yields the source function
\begin{align}
f=-\dfrac{2\sqrt{2}}{\sin(\theta)}\cos\left(\dfrac{\pi}{4}+\varphi\right)\,\left(24\cos^{4}\theta-29\cos^{2}\theta+5\right).
\end{align}
The test case is sketched in Fig.~\ref{fig:PoissonSphereSetup} showing
(a) a reconstructed example mesh and (b) the exact solution on the
sphere. The area ratios are given in Fig.~\ref{fig:PoissonSphereSetup}(c).
In Fig.~\ref{fig:PoissonSphereRes}, convergence results and condition
numbers are only shown for the automatically reconstructed surface
meshes because the findings coincide with the previous test case.
It is well-known that the inner angles in the elements play an important
role for the quality of a mesh. In particular, large inner angles
may hinder optimal convergence \cite{Deuflhard_2012a}. Therefore,
for all meshes the maximum inner angles have also been computed. The
generated meshes are typically mixed, i.e., composed by quadrilateral
and triangular elements and the maximum angles have been computed
for these two element types individually. Results are shown in Fig.~\ref{fig:PoissonSphereAngles}
for first order and $6$th-order meshes only because all other orders
yield very similar results. As can be seen the maximum inner angles
in quadrilaterals are below $140^{\circ}$ and in triangles below
$110^{\circ}$. This was also confirmed in the other test cases where
one level-set function defines the manifold.

\subsubsection{Hyperbolic paraboloid with bumps\label{sec:lap2dEx3}}

Now we consider the Laplace-Beltrami operator on a hyperbolic paraboloid
with bumps. The surface is defined by the function
\[
f\left(x,y\right)=\dfrac{1}{2}\left(x^{2}-y^{2}\right)+\dfrac{3}{20}\sin\left(2\pi x\right)\sin\left(2\pi y\right)
\]

with $x,y\in\left[-0.5,0.5\right]$. The analytic solution $u$ is
chosen as 
\begin{align}
u(x,y)=\sin\left[\pi\left(x-\dfrac{1}{2}\right)\right]\,\sin\left[\pi\left(y-\dfrac{1}{2}\right)\right]
\end{align}
and the source function $f$ determined accordingly. The level-set
function $\phi\left(\vek x\right)$ implying the manifold is
\begin{equation}
\phi(\vek x)=f(x,y)-z\ .\label{eq:HyperbolSurface}
\end{equation}
If the reconstruction is achieved based on a cube-like background
mesh with suitable dimensions, the manipulation of the background
mesh is restricted in that nodes on the outer contour of the background
mesh cannot be moved freely without changing the boundary $\partial\Gamma$
of the implied manifold. Therefore, we use four additional level-set
functions $\psi^{i}\left(\vek x\right)$ which restrict the manifold
properly. The mesh manipulation of the (universal) background mesh,
see Section \ref{X_Preliminaries}, is then applied with respect to
\emph{all} level-set functions and the element ratios are sufficently
bounded.

\begin{figure}
\centering

\subfigure[Zero-level sets]{\includegraphics[width=0.4\textwidth]{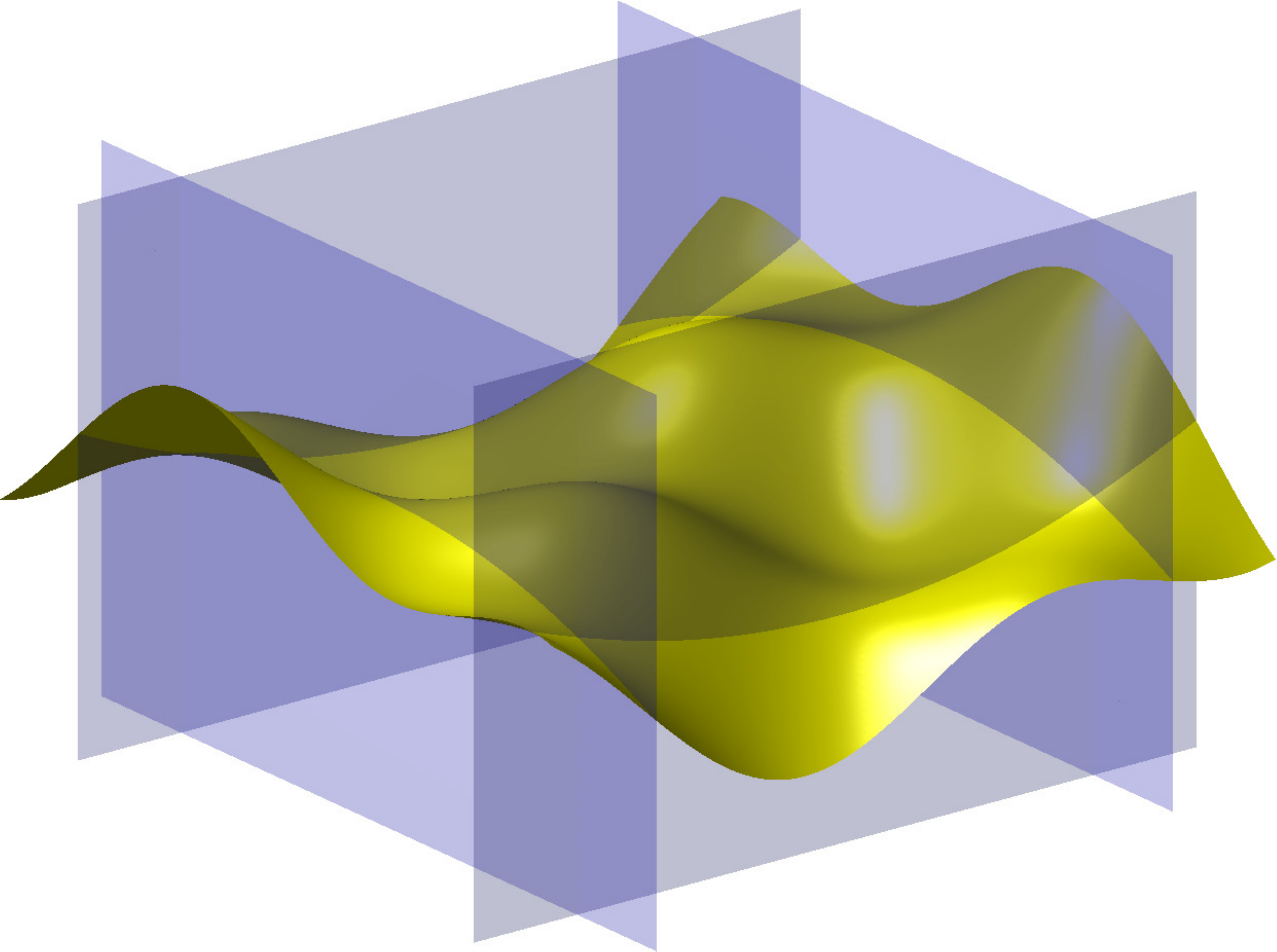}}\qquad\subfigure[Surface mesh]{\includegraphics[width=0.35\textwidth]{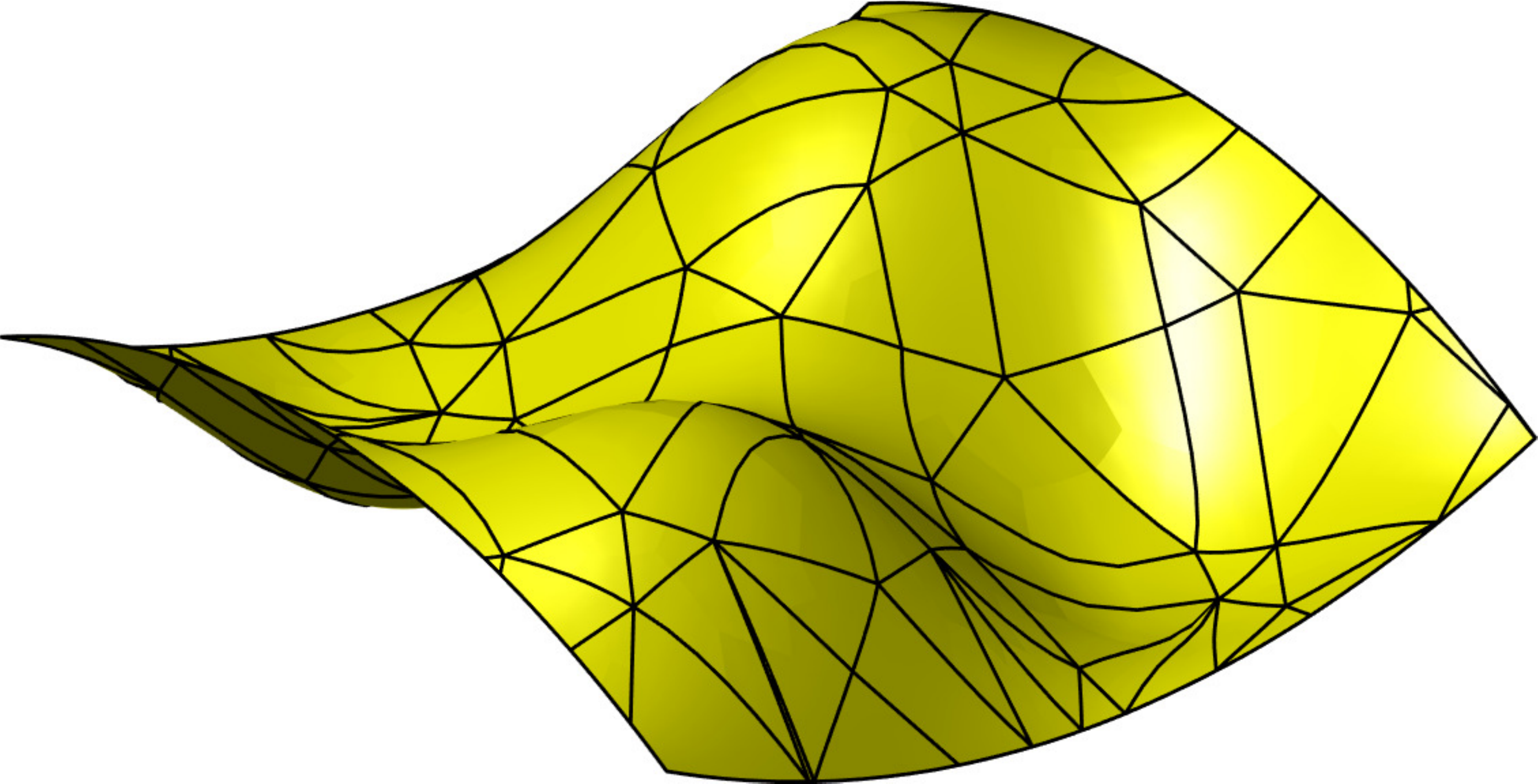}}

\subfigure[exact solution]{\includegraphics[width=0.35\textwidth]{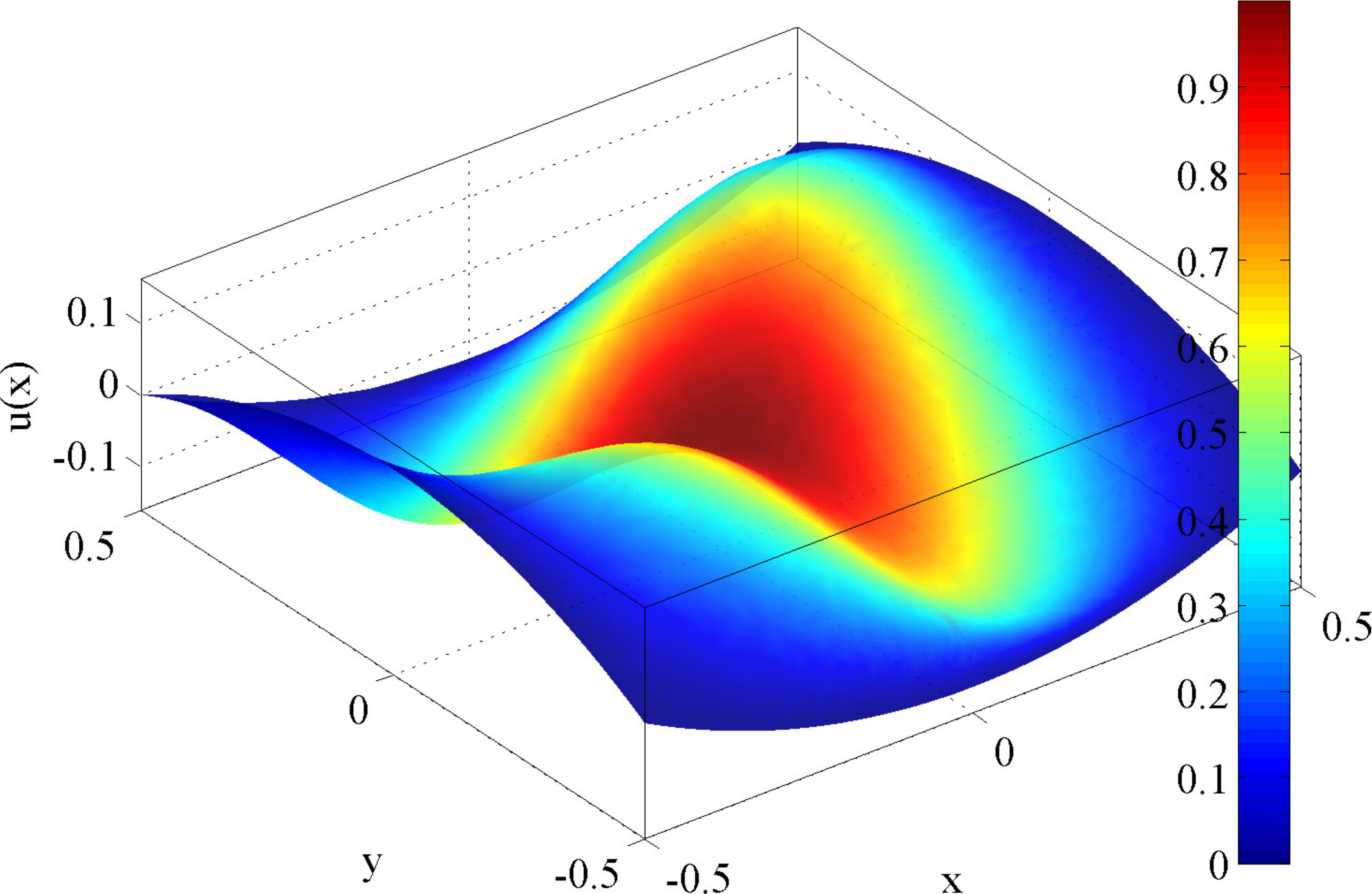}}\qquad\subfigure[ratio $A_{\mathrm{max}}/A_{\mathrm{min}}$]{\includegraphics[width=0.35\textwidth]{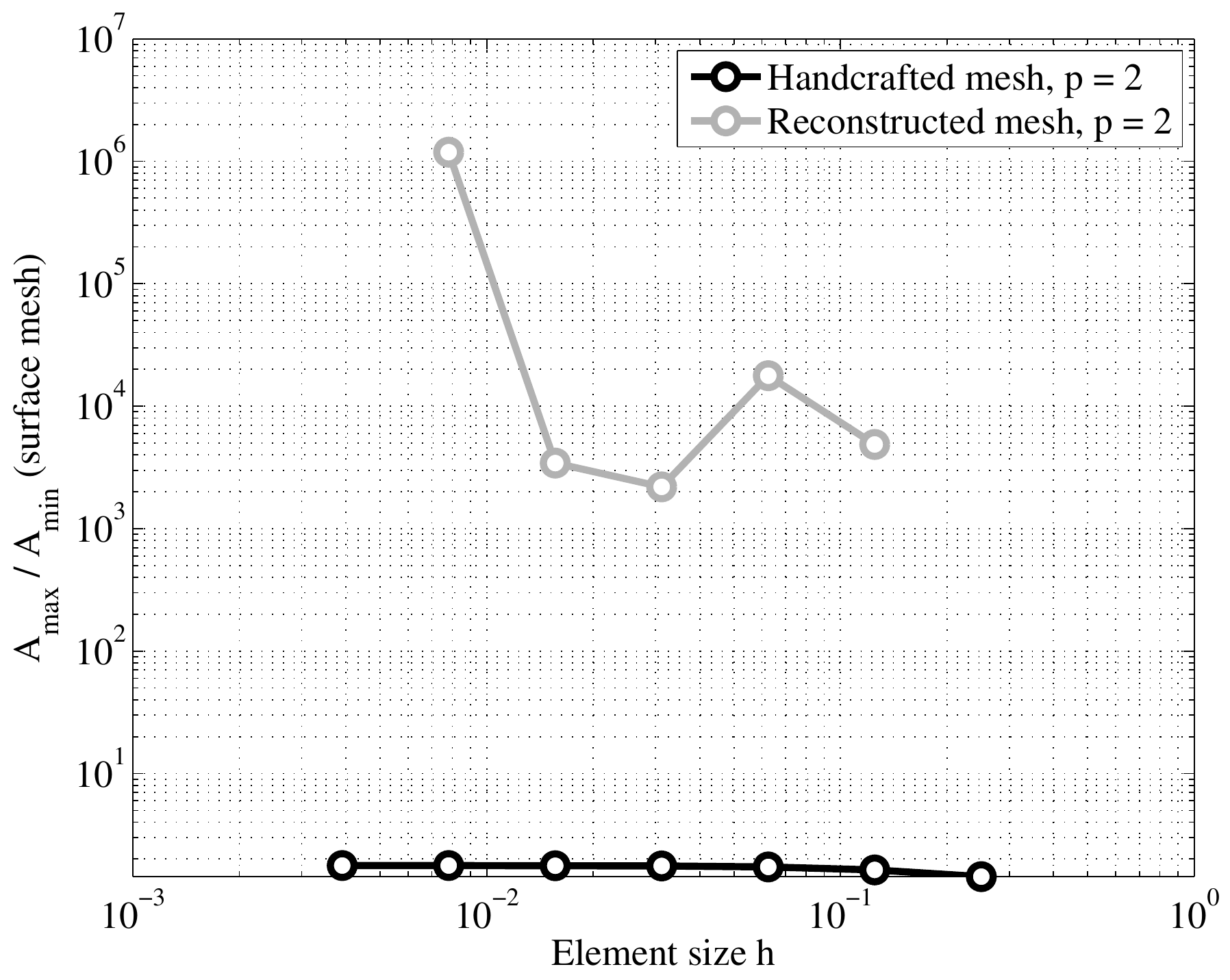}}

\caption{(a) Zero iso-surfaces of $\phi$ and $\psi^{i}$, (b) reconstructed
surface mesh for the hyperbolic paraboloid with bumps, (c) exact solution
of the Poisson equation, (d) the largest ratio $A_{\mathrm{max}}/A_{\mathrm{min}}$
of the reconstructed surface elements.}

\label{fig:PoissonHyperBolSetup} 
\end{figure}

\begin{figure}
\centering

\subfigure[recon.~mesh, $L_2$-norm]{\includegraphics[width=0.35\textwidth]{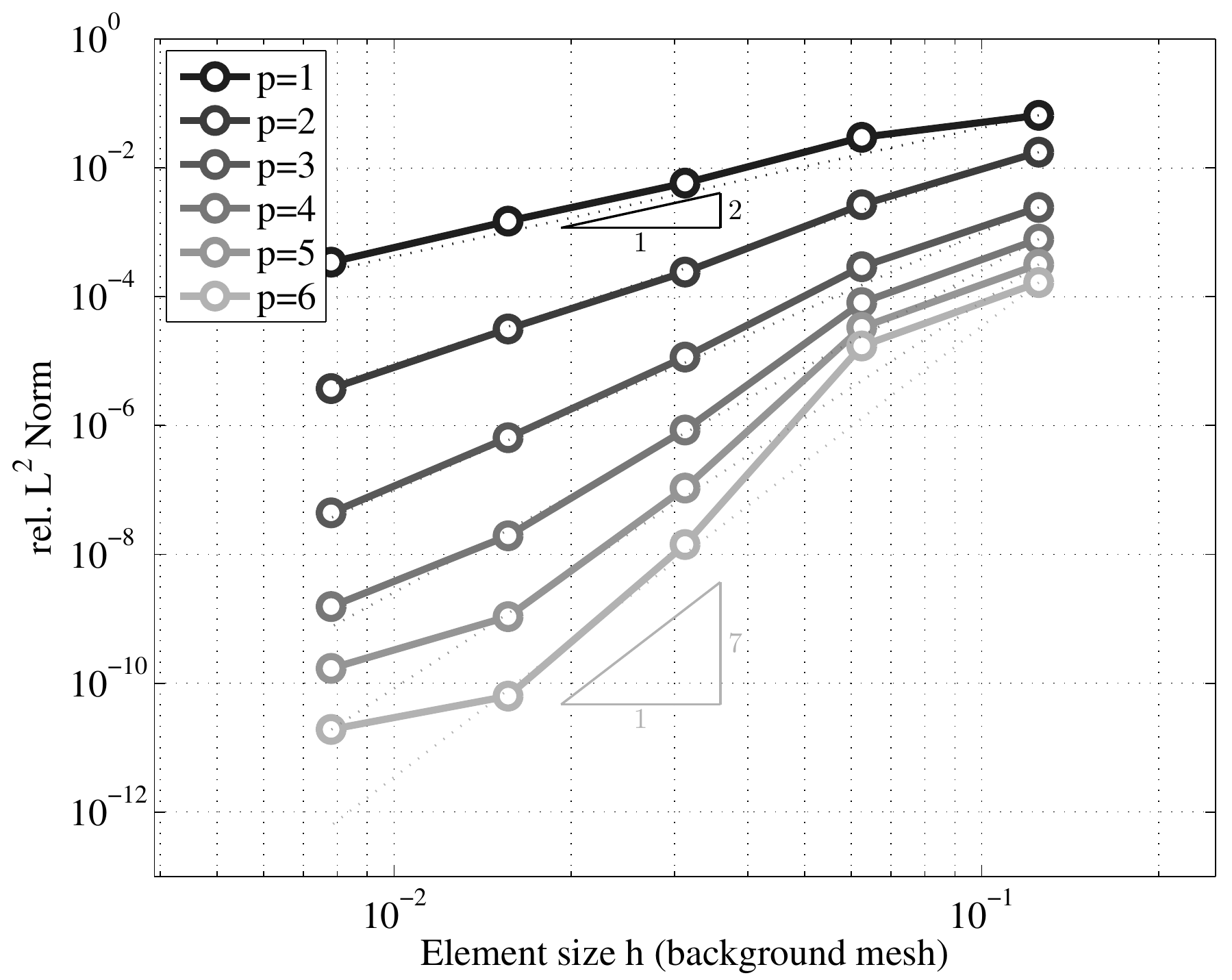}}\quad\subfigure[recon.~mesh, cond.]{\includegraphics[width=0.35\textwidth]{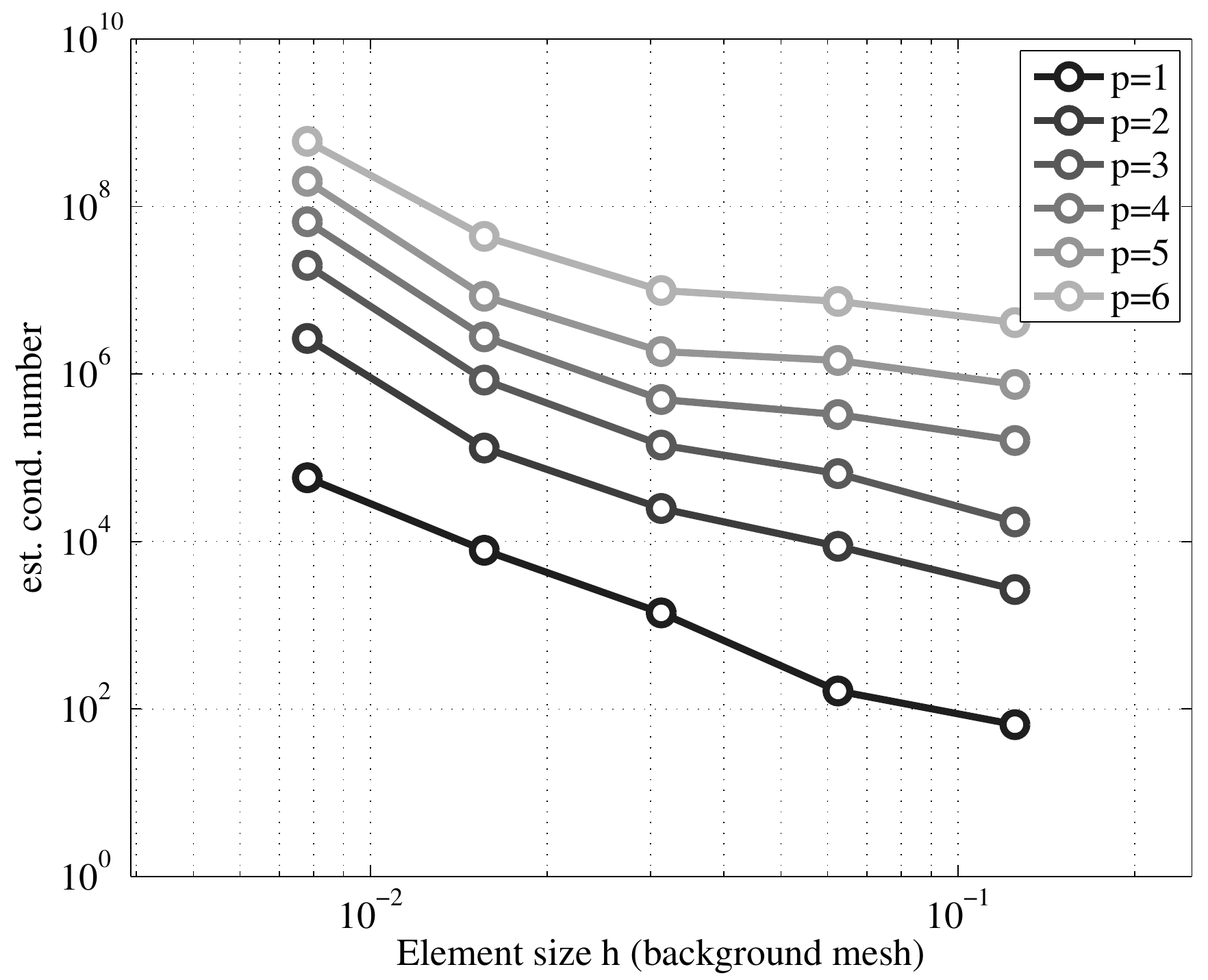}}

\caption{\label{fig:PoissonHyperBolRes}Convergence results and condition numbers
for the Poisson equation on a hyperbolic paraboloid with bumps.}
\end{figure}

The additional level-set functions are 
\begin{align*}
\psi^{i}(\vek x) & =\left\langle \vek x-\vek P_{i},\,\vek n_{i}\right\rangle \qquad & i=1,\,\ldots,\,4
\end{align*}
where $\left\langle \;,\;\right\rangle $ is a scalar product and
\begin{alignat*}{4}
\vek P_{1} & =\left[\begin{array}{c}
0.5\\
0\\
0
\end{array}\right]\quad & \vek P_{2} & =\left[\begin{array}{c}
0\\
0.5\\
0
\end{array}\right]\quad & \vek P_{3} & =\left[\begin{array}{c}
-0.5\\
0\\
0
\end{array}\right]\quad & \vek P_{4} & =\left[\begin{array}{c}
0\\
-0.5\\
0
\end{array}\right]\\
\vek n_{1} & =\left[\begin{array}{c}
1\\
0\\
0
\end{array}\right] & \vek n_{2} & =\left[\begin{array}{c}
0\\
1\\
0
\end{array}\right] & \vek n_{3} & =\left[\begin{array}{c}
-1\\
0\\
0
\end{array}\right] & \vek n_{4} & =\left[\begin{array}{c}
0\\
-1\\
0
\end{array}\right]
\end{alignat*}
The resulting zero-level sets are planes restricting the zero-level
set $\Gamma_{\phi}$ implied by Eq.~(\ref{eq:HyperbolSurface}).
See Fig.~\ref{fig:PoissonHyperBolSetup}(a) for a visualization of
the zero-level sets and (b) the resulting bounded manifold. Fig.~\ref{fig:PoissonHyperBolSetup}(c)
and (d) show the exact solution and the area ratios of the resulting
surface elements, respectively. Convergence results and condition
numbers for the automatically reconstructed meshes are seen in Fig.~\ref{fig:PoissonHyperBolRes}.

\subsection{Instationary transport equation in 2D and 3D\label{sec:instat}}

The instationary advection-diffusion equation on manifolds is considered
here. The semi-discrete weak form is given in Eq.~(\ref{eq:WeakFormAdvDiff}).
The time discretization is carried out using a $6$th-order accurate,
implicit, $3$-step Runge-Kutta method as defined in Eq.~(\ref{eq:ButcherTableau}).
We use a constant number of 4096 time steps in all upcoming examples
which virtually eliminates the discretization error in time. First,
we consider pure advection problems on a circle and a sphere, i.e.~without
any diffusion. Then, two test cases are considered with diffusion
where no analytic solution is available. In the error studies, the
$L_{2}$-error is measured in space at the final time step. Of course,
this error is effected by the discretizations in space \emph{and}
time, wherefore we employ the highly accurate time stepping scheme
mentioned above.

\subsubsection{Pure advection on a circular manifold\label{sec:veri1d}}

\begin{figure}
\centering

\subfigure[Initial state]{\includegraphics[width=0.32\textwidth]{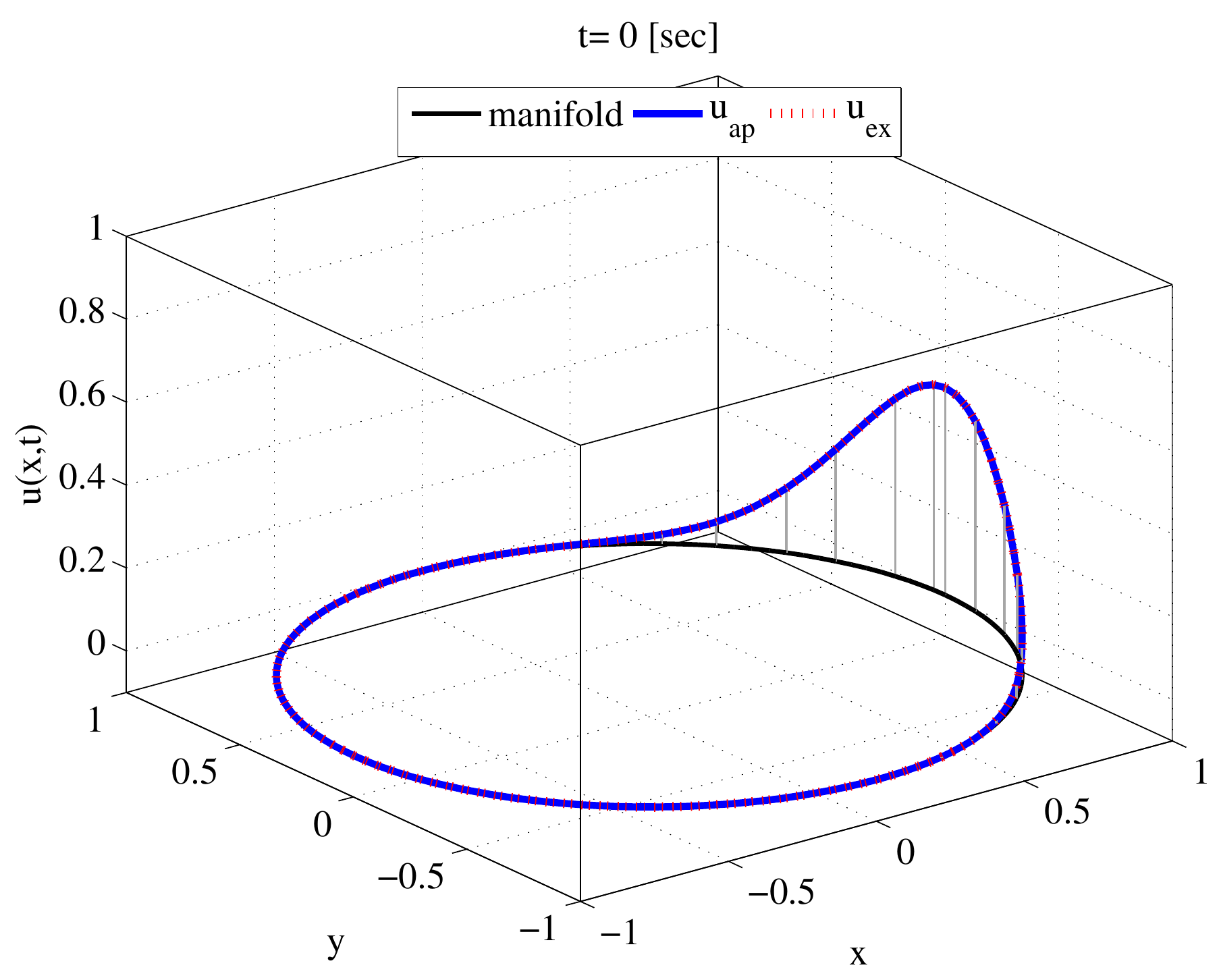}}\hfill\subfigure[State at $0.20$]{\includegraphics[width=0.32\textwidth]{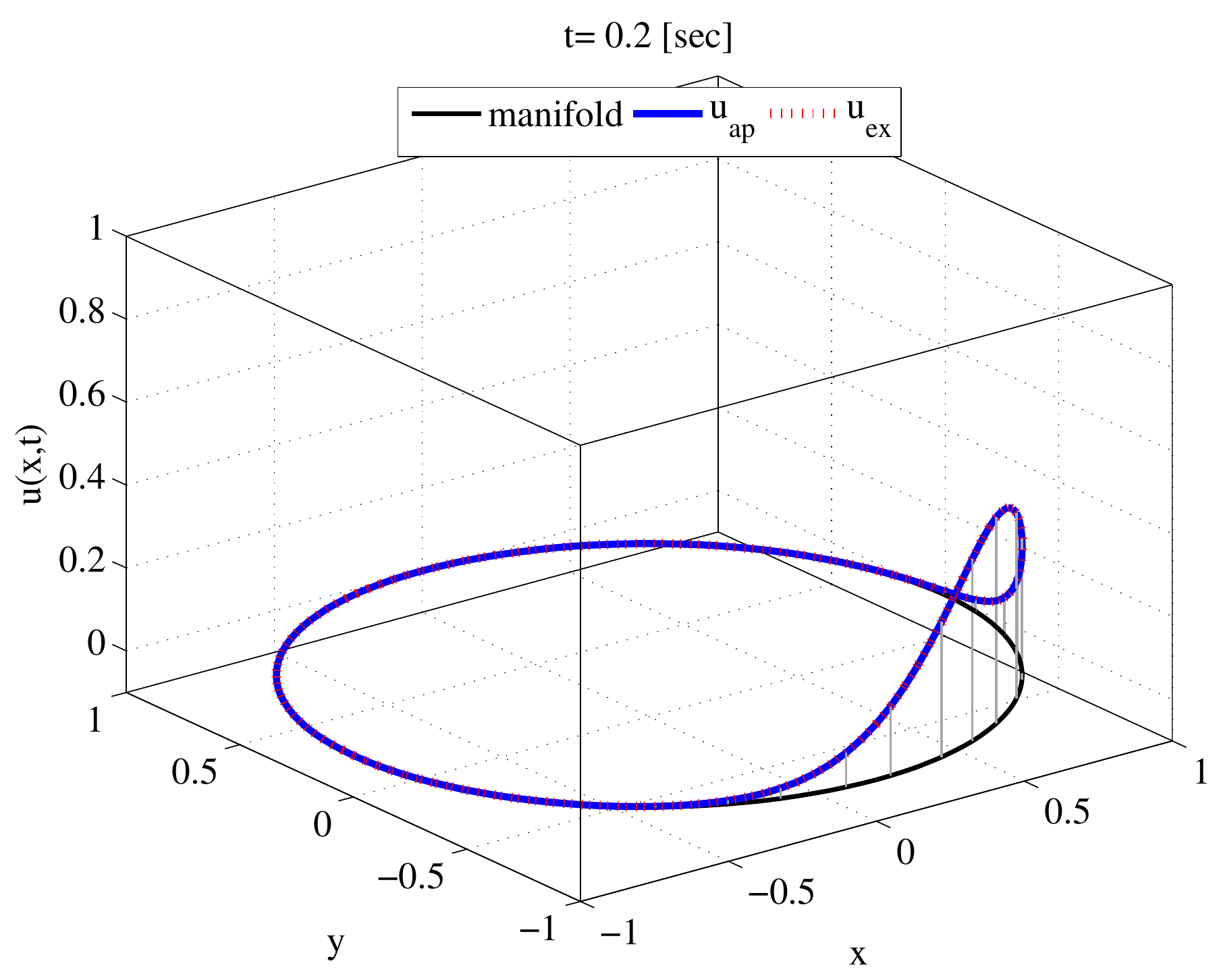}}\hfill\subfigure[State at $1.00$]{\includegraphics[width=0.32\textwidth]{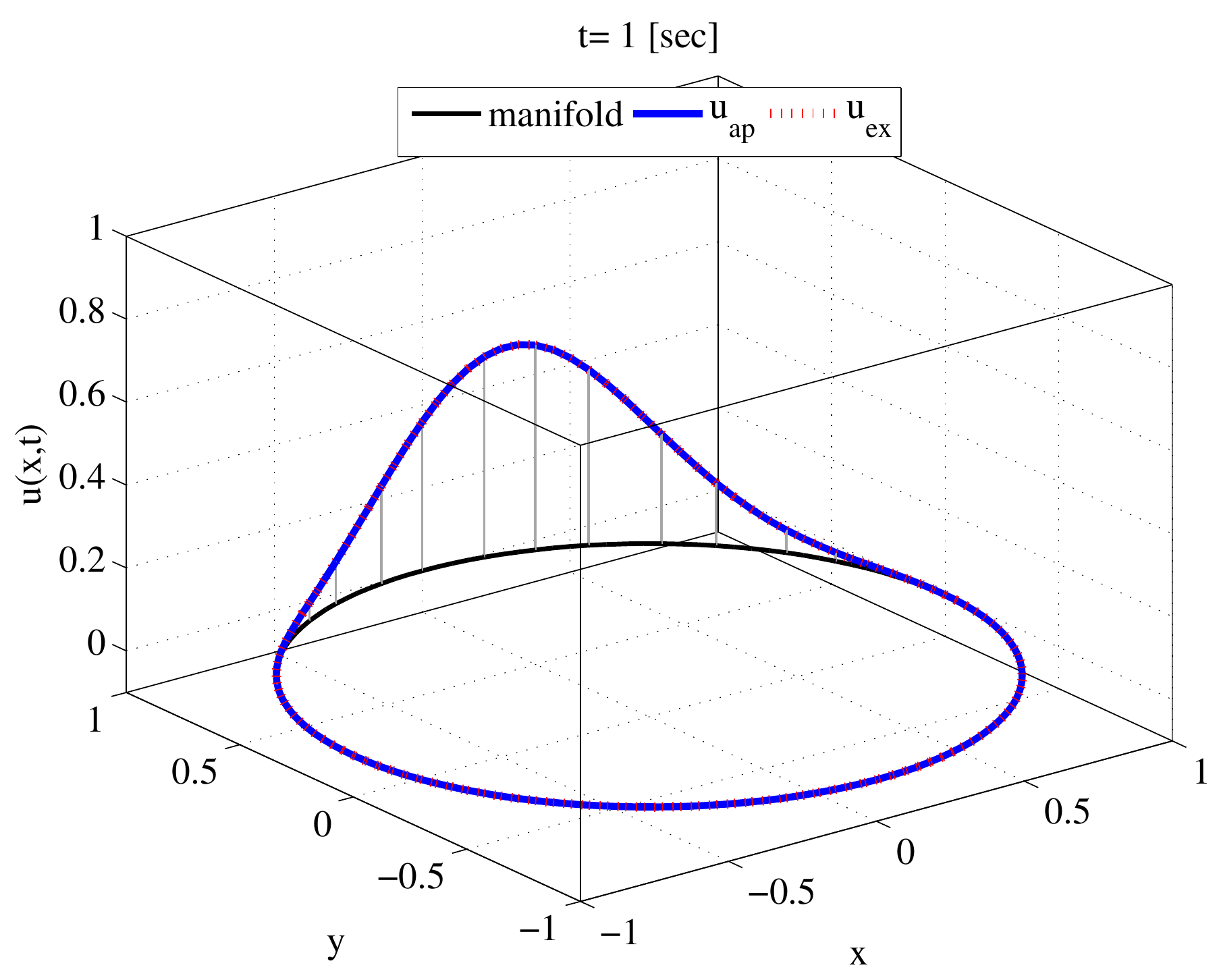}}

\caption{\label{fig:TranspCircleSol} Solution of the instationary transport
problem on a circle without diffusion.}
\end{figure}

The instationary transport problem is solved on a circle with radius
$r=1$ in the time interval $t\in(0,1)$. The advection velocity is
$\left\Vert \vek c_{\Gamma}\right\Vert =5$ in tangential direction
of the circle and the diffusion coefficient $\lambda=0$. The initial
condition is
\begin{equation}
u_{0}(\vek x)=\exp\left(-4\cdot\varphi\left(\vek x\right)^{2}\right)\ \mathrm{with}\;\varphi=\mathrm{atan}\left(y/x\right),\;\varphi\in\left[-\pi,\,\pi\right],\ \vek x\in\Gamma.\label{eq:InitCondTranspCircle}
\end{equation}
This initial distribution is simply advected around the circle without
changing its shape. That is, the exact solution is (\ref{eq:InitCondTranspCircle})
with a shift in the angle of $\Delta\varphi=\left\Vert \vek c_{\Gamma}\right\Vert /r\cdot t$.
See Fig.~\ref{fig:TranspCircleSol} for a graphical representation
of the solution at different instances in time at $t=\left[0.00,\,0.20,\,1.00\right]$.

\begin{figure}
\centering

\subfigure[handcr.~mesh, $L_2$-norm]{\includegraphics[width=0.35\textwidth]{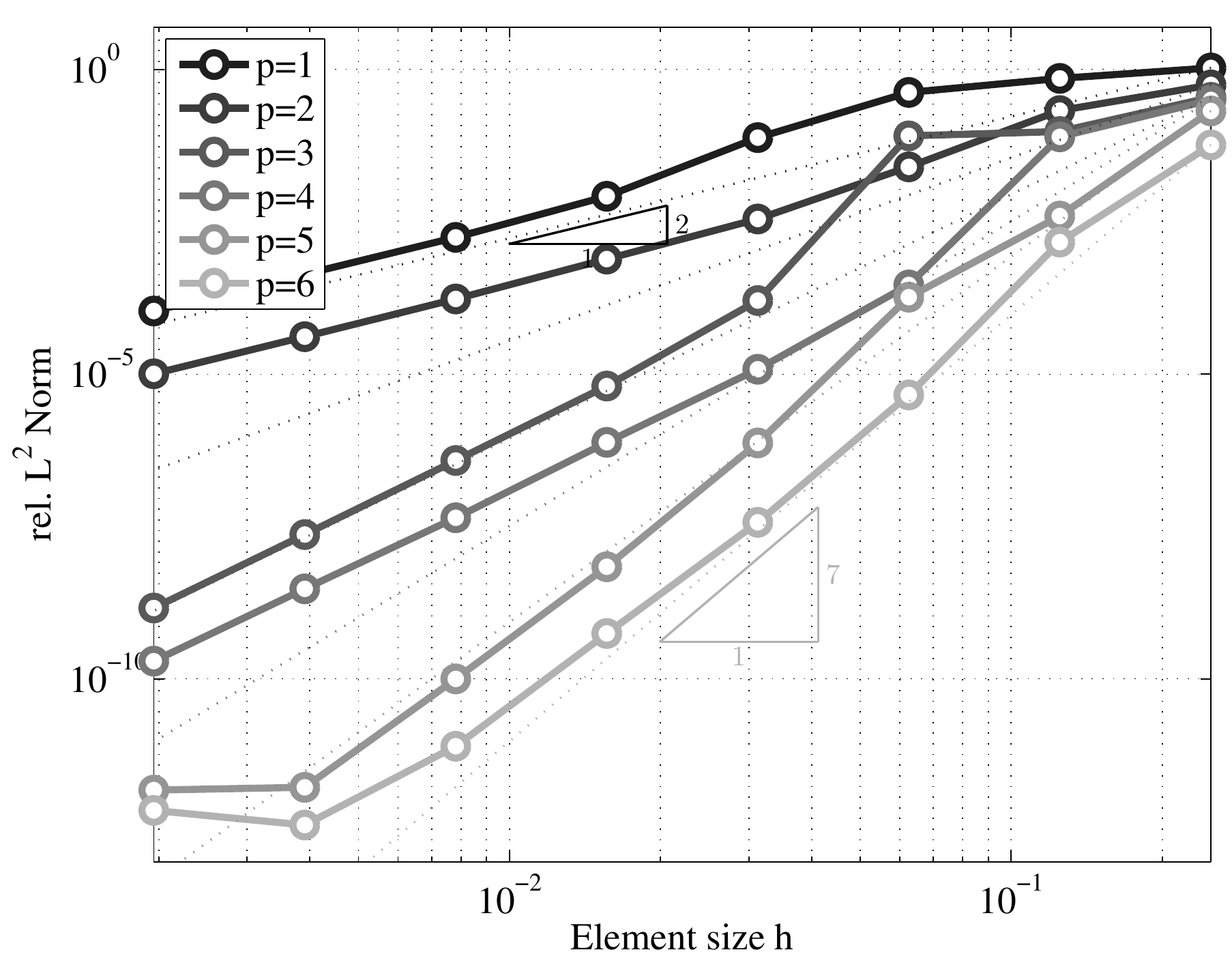}}\quad\subfigure[recon.~mesh, $L_2$-norm]{\includegraphics[width=0.35\textwidth]{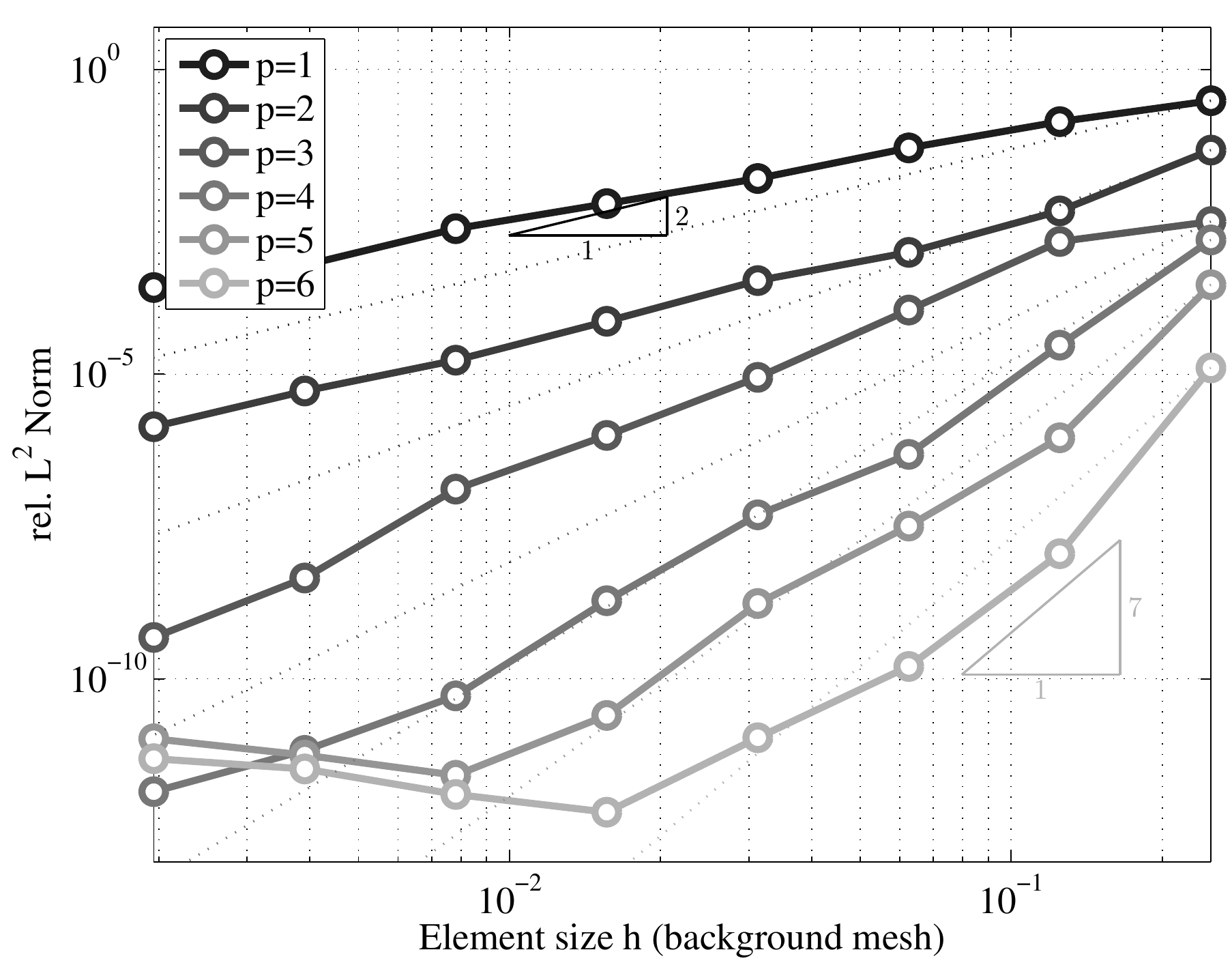}}

\subfigure[handcr.~mesh, cond.]{\includegraphics[width=0.35\textwidth]{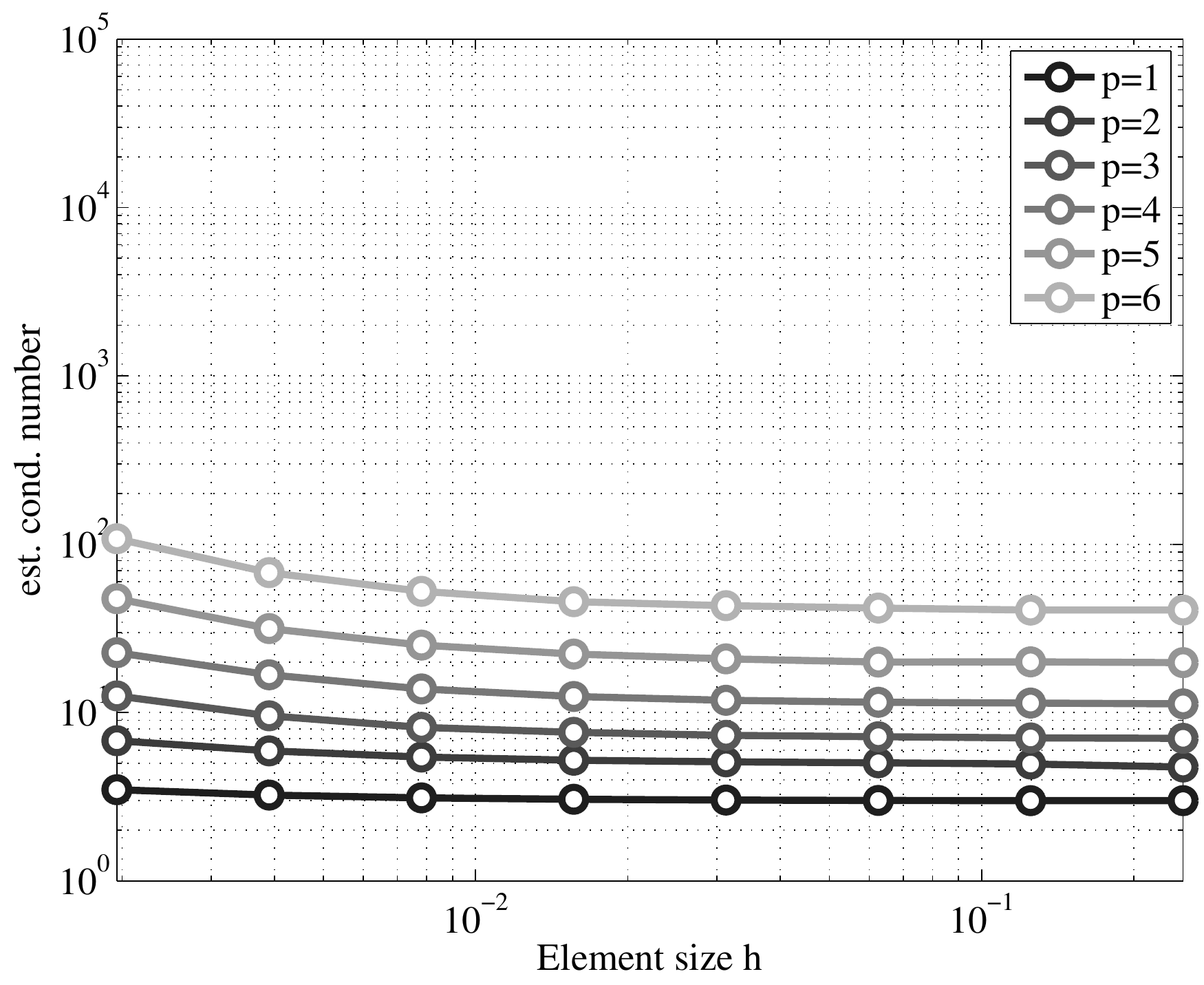}}\quad\subfigure[recon.~mesh, cond.]{\includegraphics[width=0.35\textwidth]{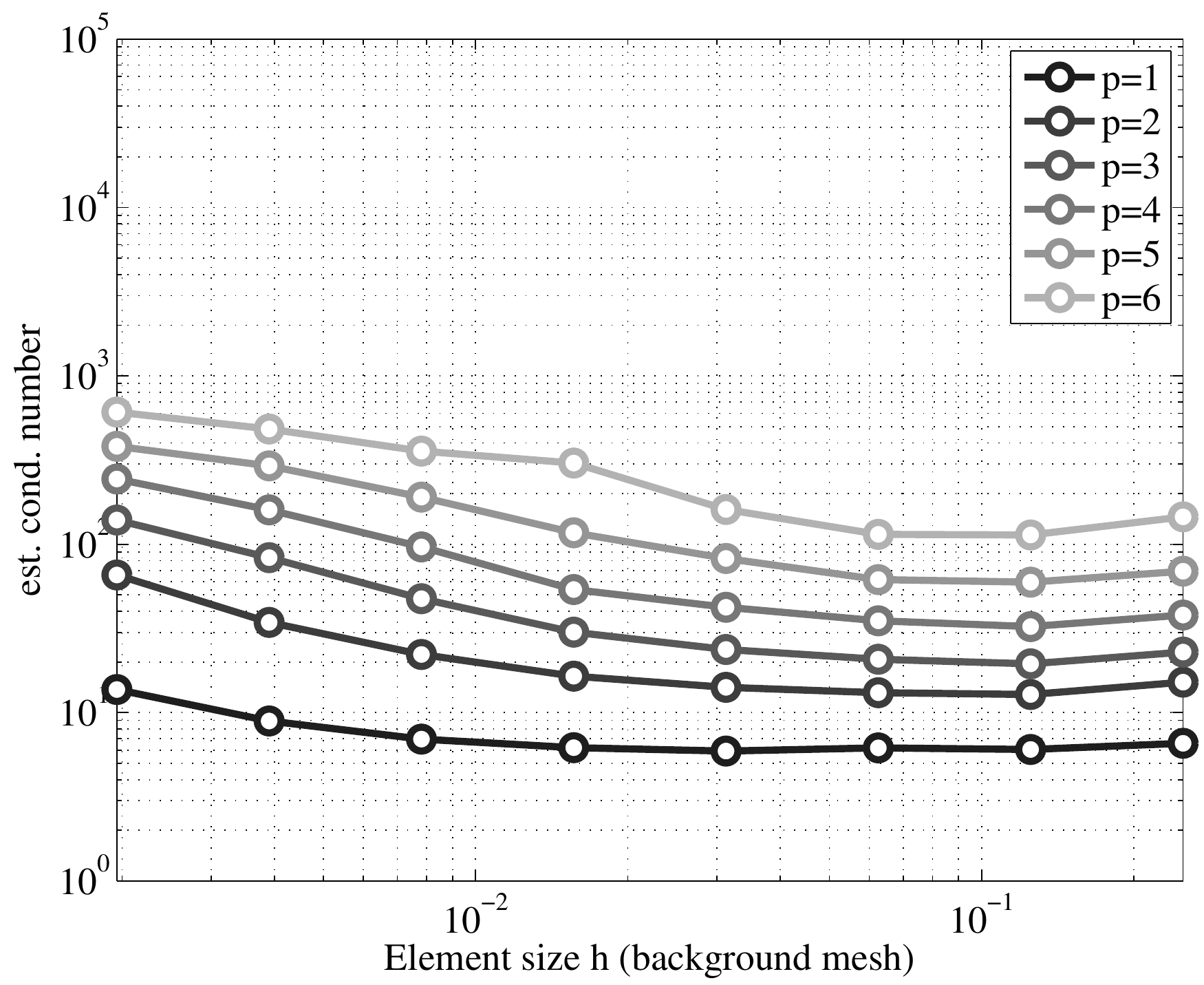}}

\caption{\label{fig:TranspCircleRes}Convergence results and condition numbers
for the transport equation on a circular manifold.}
\end{figure}

For this study, the same meshes from Section \ref{sec:lap1dEx1} are
used. Results are presented in Fig.~\ref{fig:TranspCircleRes} in
the same style than above. It is seen that higher-order convergence
rates are achieved. However, for this pure advection problem, one
order of the optimal convergence rate is lost for \emph{even} orders
of the elements, i.e.~$p=2,4,6,\dots$ This is the same for handcrafted
meshes as well as for the automatically reconstructed meshes. We have
confirmed that the same behavior also occurs for planar advection
problems so this has nothing to do with the fact that the problem
is solved on a curved manifold nor that no boundaries are present.

\subsubsection{Pure advection on a sphere\label{sec:veri2d}}

\begin{figure}
\centering

\subfigure[Initial state]{\includegraphics[width=0.32\textwidth]{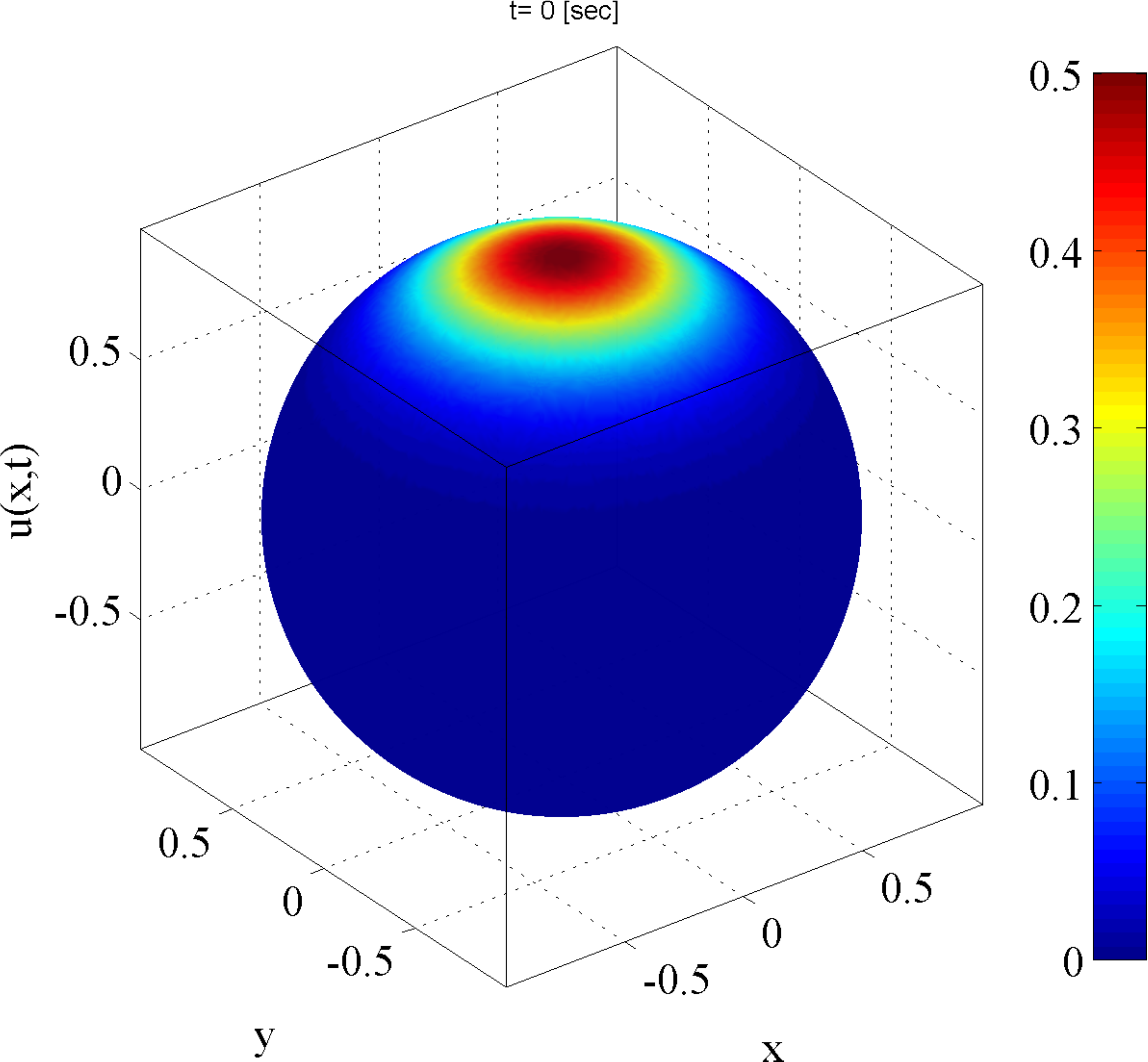}}\hfill\subfigure[State at $0.50\,\text{s}$]{\includegraphics[width=0.32\textwidth]{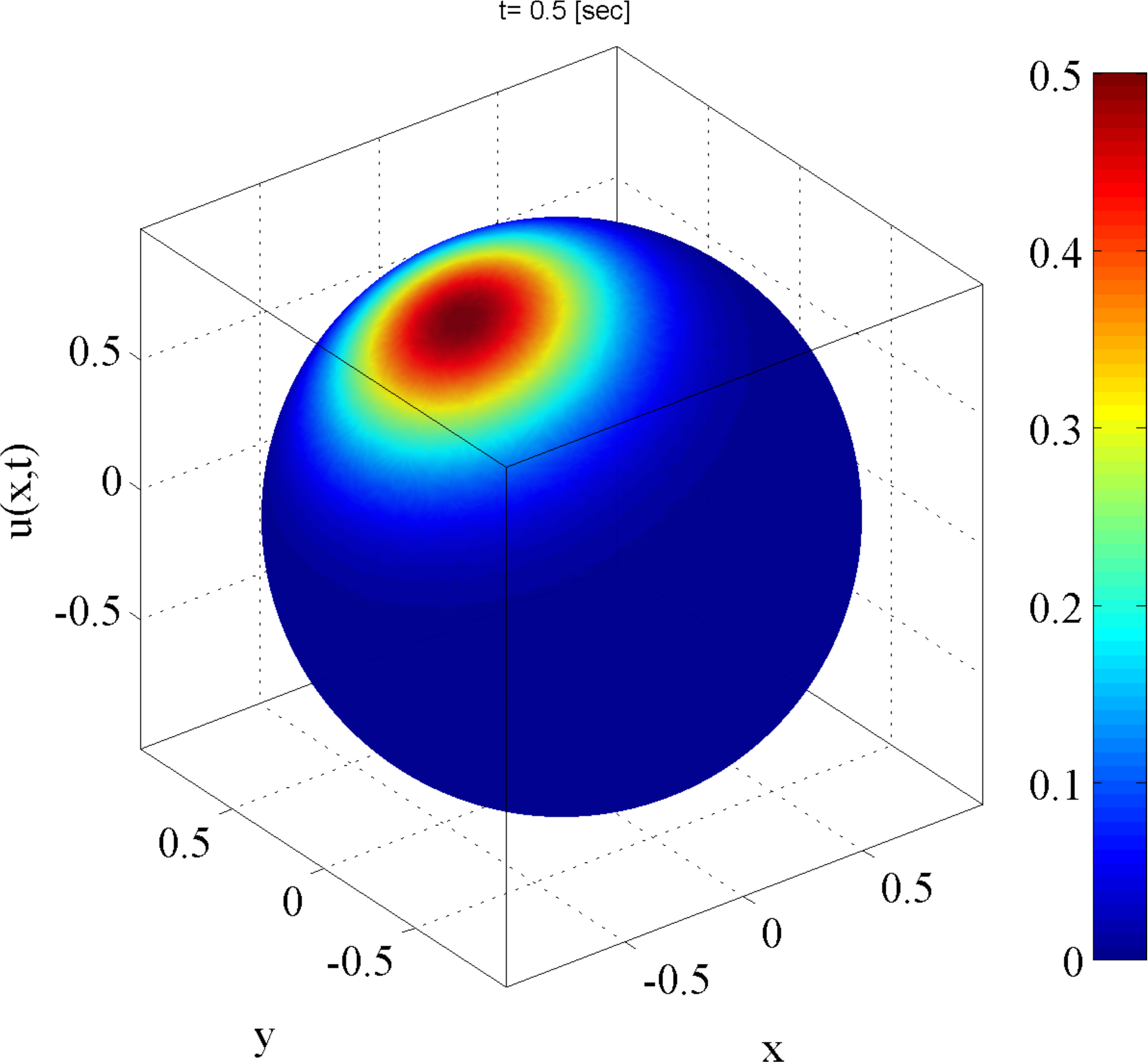}}\hfill\subfigure[State at $1.00\,\text{s}$]{\includegraphics[width=0.32\textwidth]{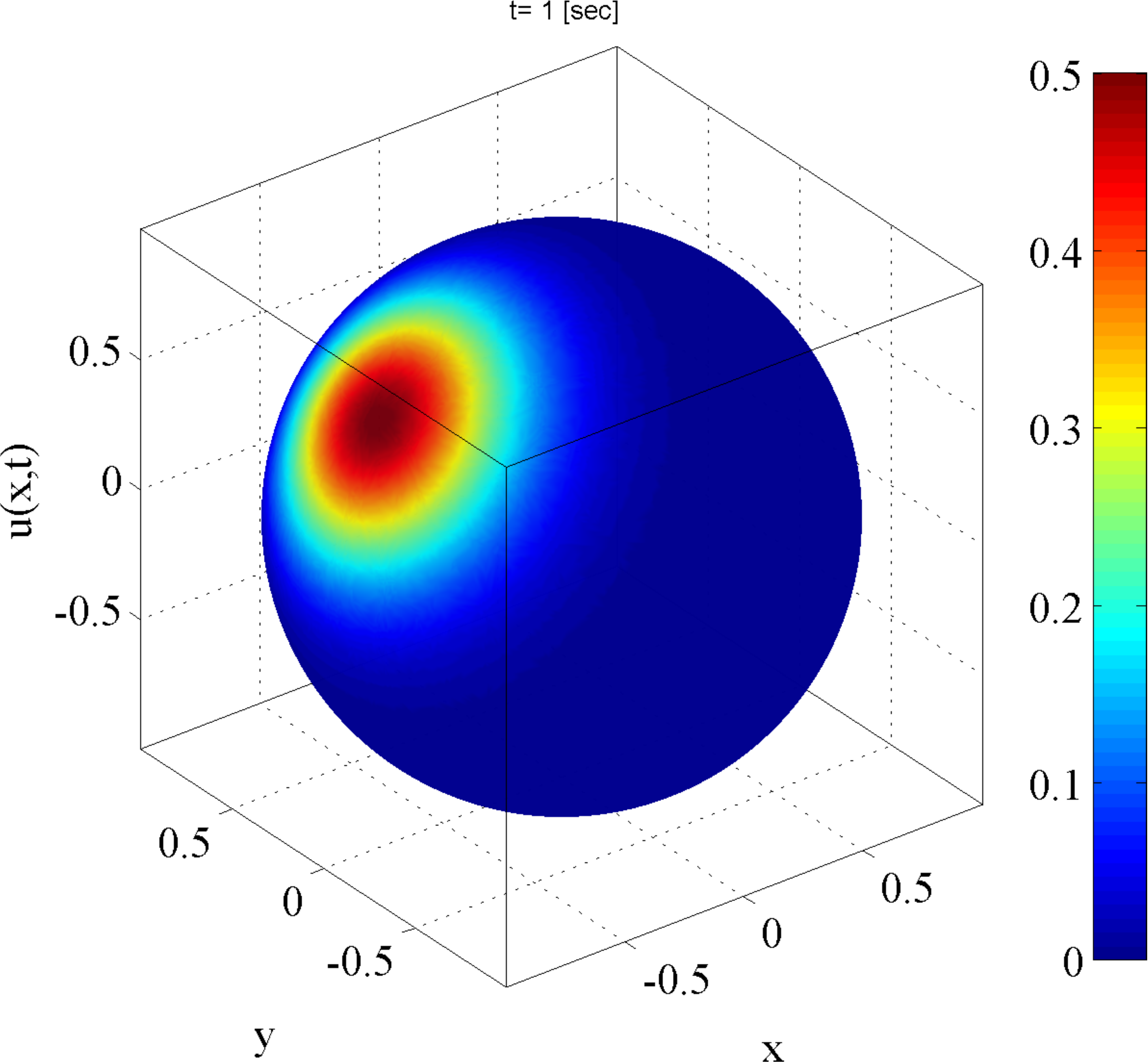}}

\caption{\label{fig:exIn2dex} Solution of the instationary transport problem
on a sphere without diffusion.}
\end{figure}

A similar problem without diffusion is now solved on a sphere with
radius $r=1$ in the time interval $t\in(0,1)$. The advection velocity
is
\[
\vek c_{\Gamma}=\left[\begin{array}{c}
c_{x}\\
0\\
c_{z}
\end{array}\right]=\left[\begin{array}{c}
z\\
0\\
-x
\end{array}\right]\cdot\frac{c}{r},\quad c=-\frac{7}{8}
\]

taking place only in the $xz$-plane tangential to the sphere. The
initial state is set to
\begin{align}
u_{0}(\vek x)=\exp\left(-4\theta(\vek x)^{2}\right)\ \mathrm{with}\;\theta=\mathrm{acos}\left(z/r\right),\ \theta\in\left[0,\,\pi\right],\ \vek x\in\Gamma.
\end{align}
Analogously to the example above, the exact solution is obtained by
the rotation angle in the $xz$-plane. In Fig.~\ref{fig:exIn2dex},
the analytical solution at $t=\left[0.00,\,0.50,\,1.00\right]$ is
illustrated.

\begin{figure}
\centering

\subfigure[handcr.~mesh, $L_2$-norm]{\includegraphics[width=0.35\textwidth]{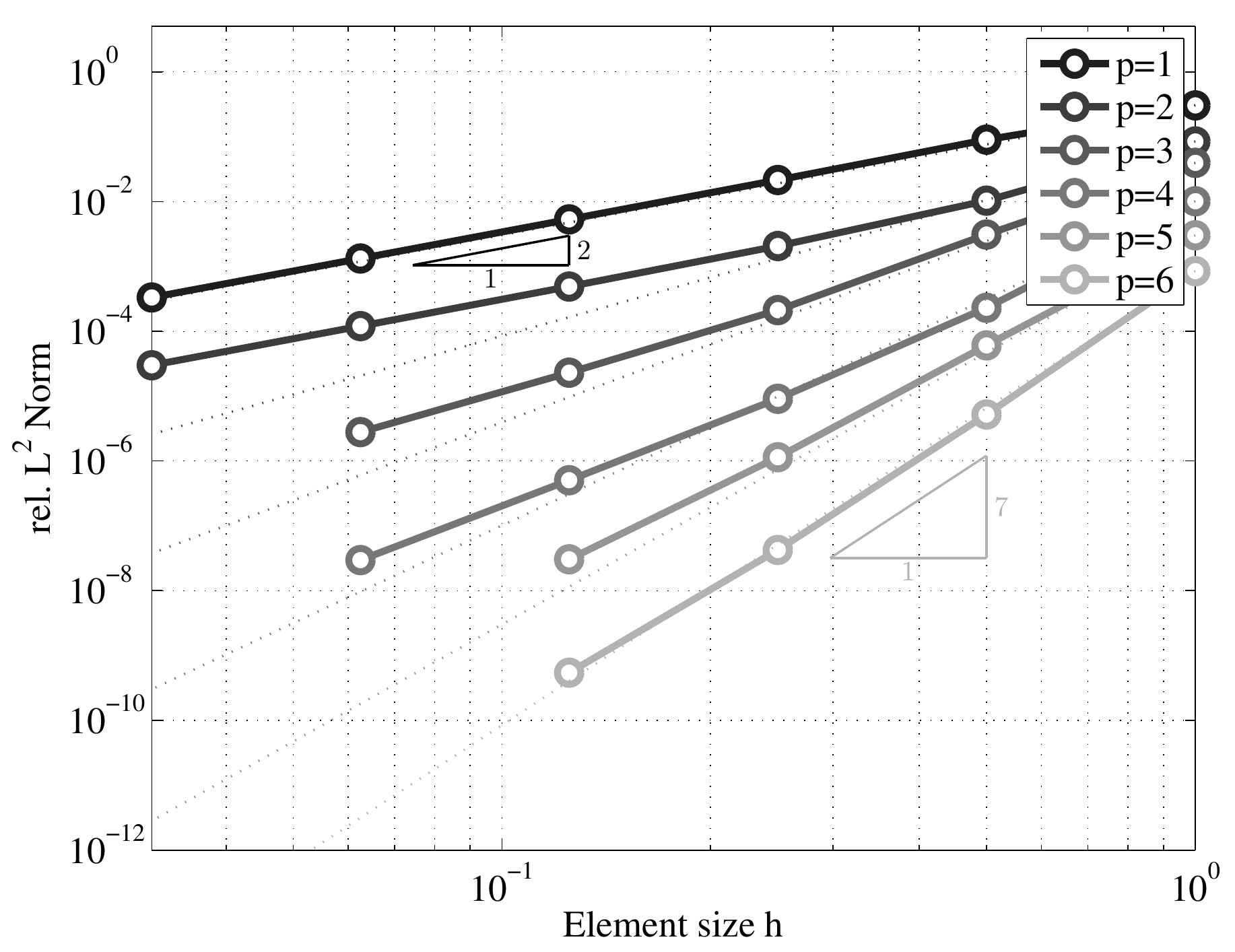}}\quad\subfigure[recon.~mesh, $L_2$-norm]{\includegraphics[width=0.35\textwidth]{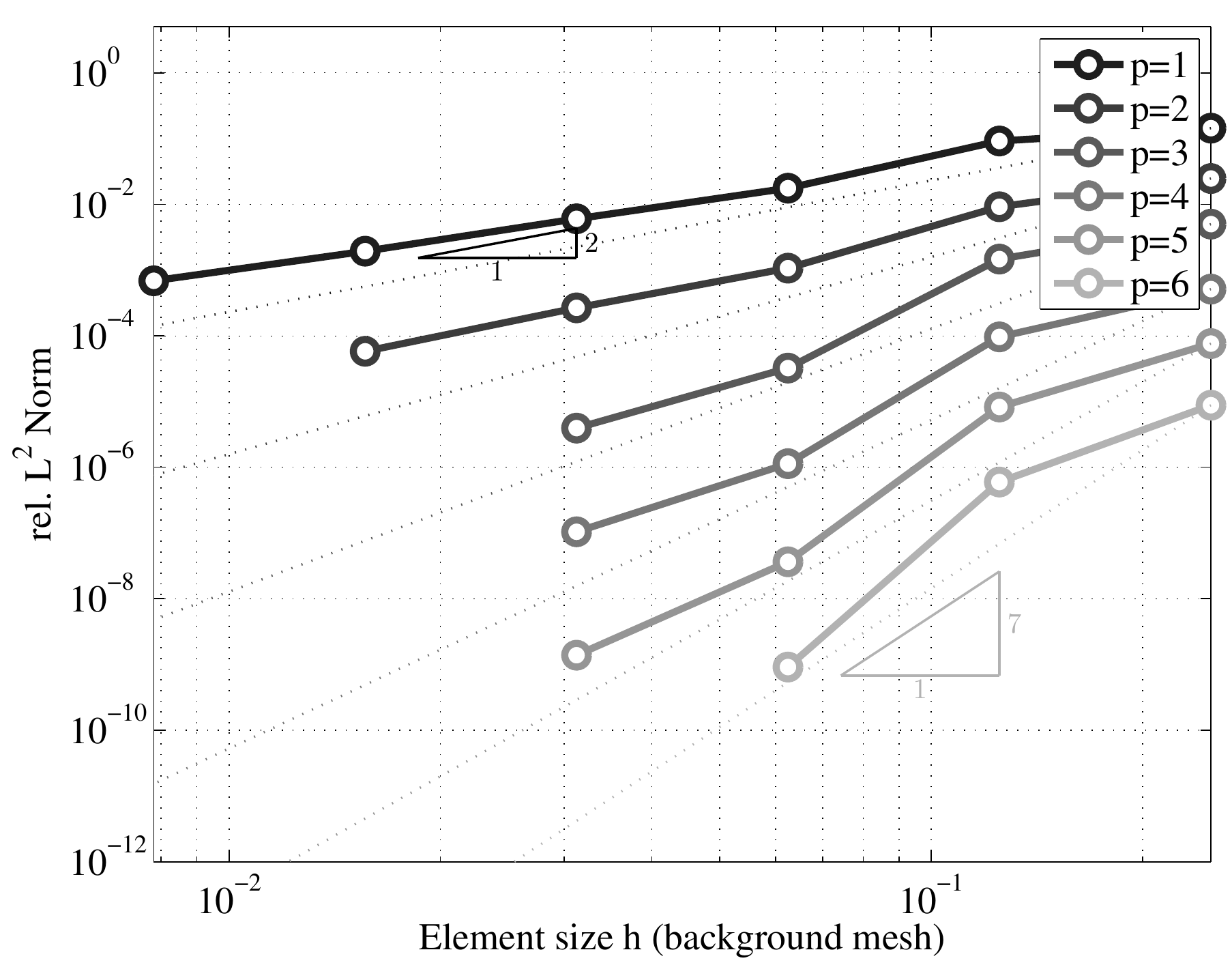}}

\subfigure[handcr.~mesh, cond.]{\includegraphics[width=0.35\textwidth]{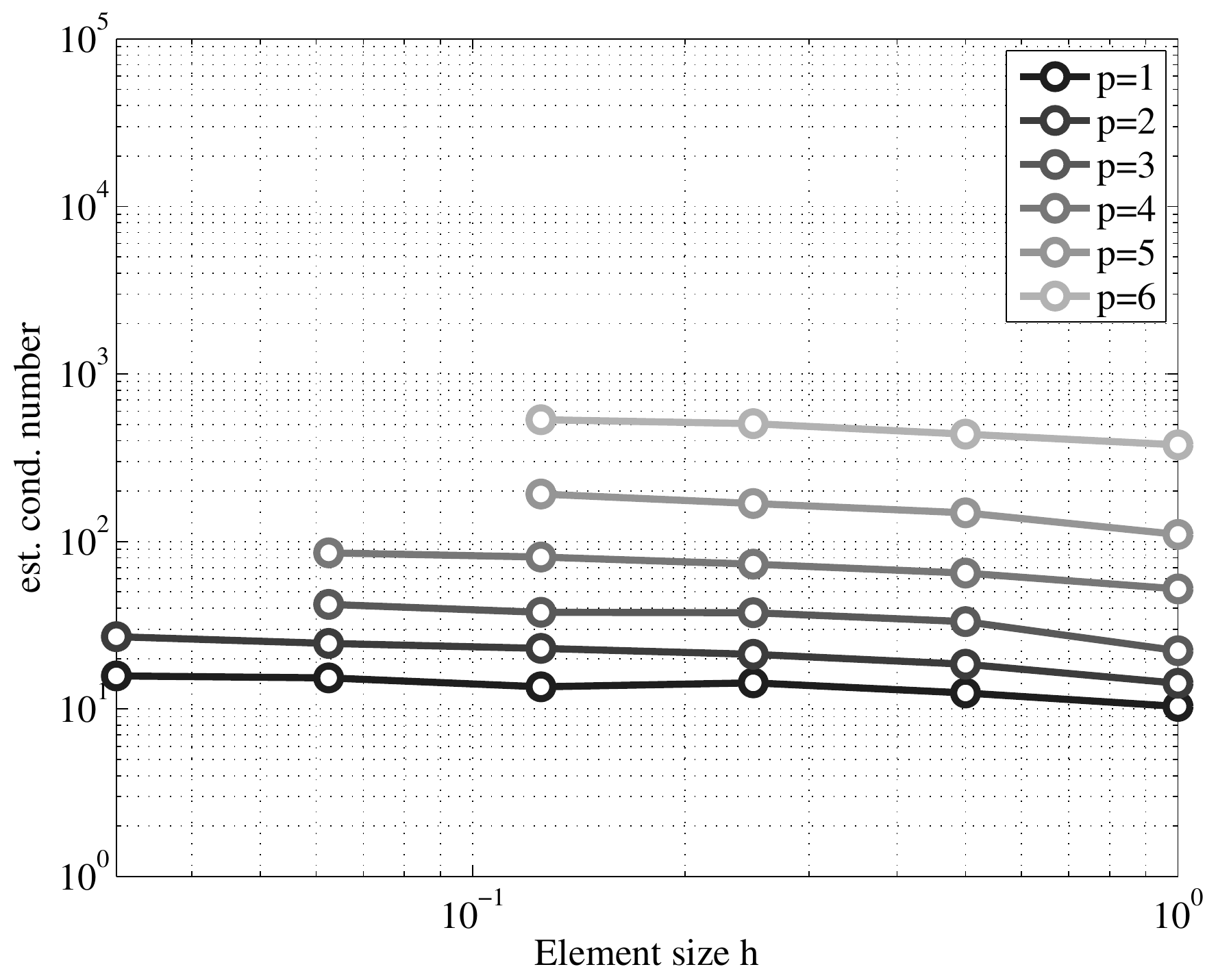}}\quad\subfigure[recon.~mesh, cond.]{\includegraphics[width=0.35\textwidth]{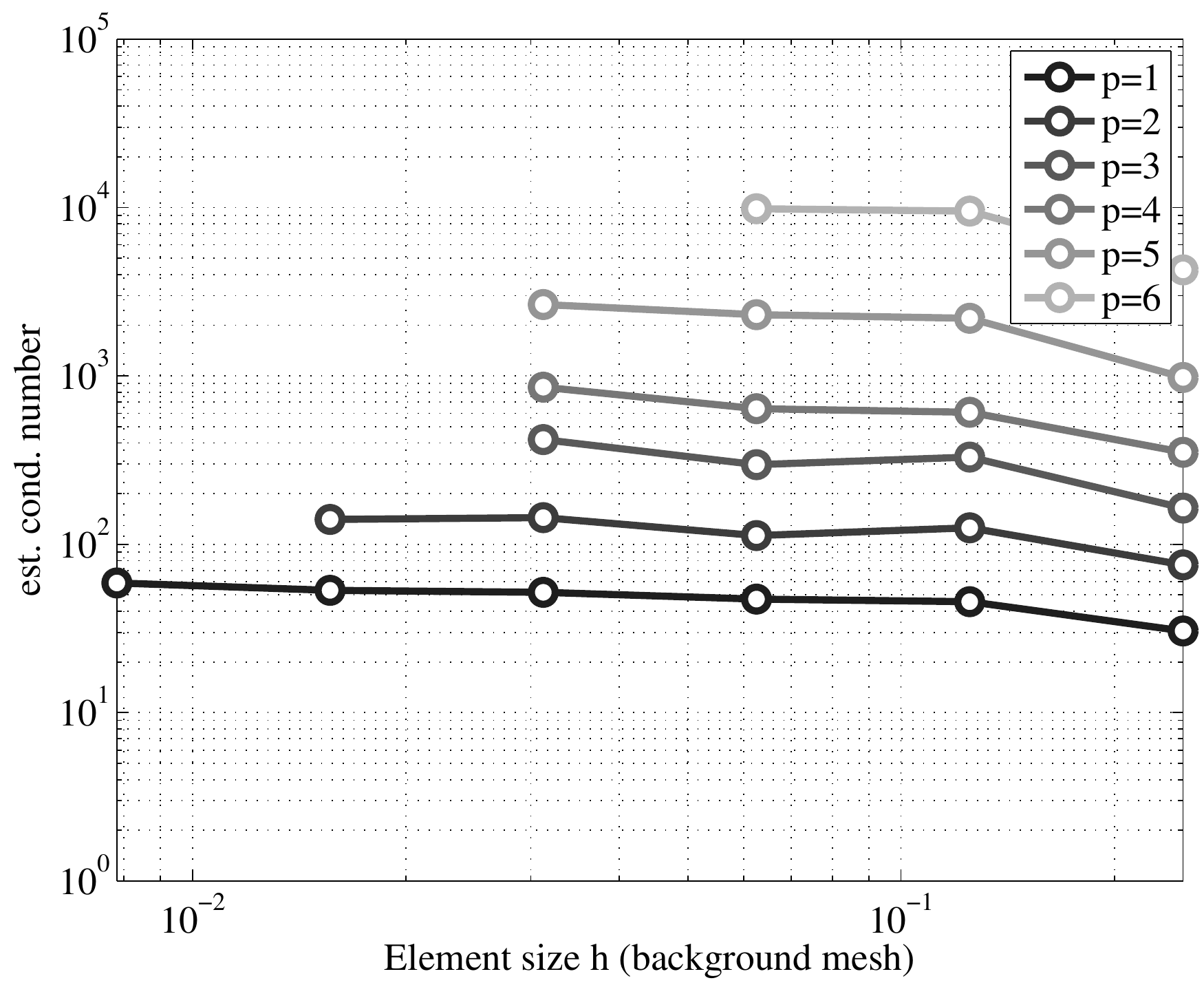}}

\caption{\label{fig:TranspSphereRes}Convergence results and condition numbers
for the transport equation on a sphere.}
\end{figure}

The same meshes from Section \ref{sec:lap2dEx2} are used and results
are seen in Fig.~\ref{fig:TranspSphereRes}. The same conclusions
from above may be drawn, in particular the loss of one order in the
convergence rate for elements with even orders. It is recalled that
this behaviour is typical for \emph{pure} advection and optimal convergence
rates are recovered in the presence of diffusion. This is equivalent
to results of planar problems with handcrafted meshes.

\subsubsection{Advection-diffusion on an S-shaped manifold\label{sec:ex1dIn}}

\begin{figure}
\centering

\subfigure[Initial state]{\includegraphics[width=0.32\textwidth]{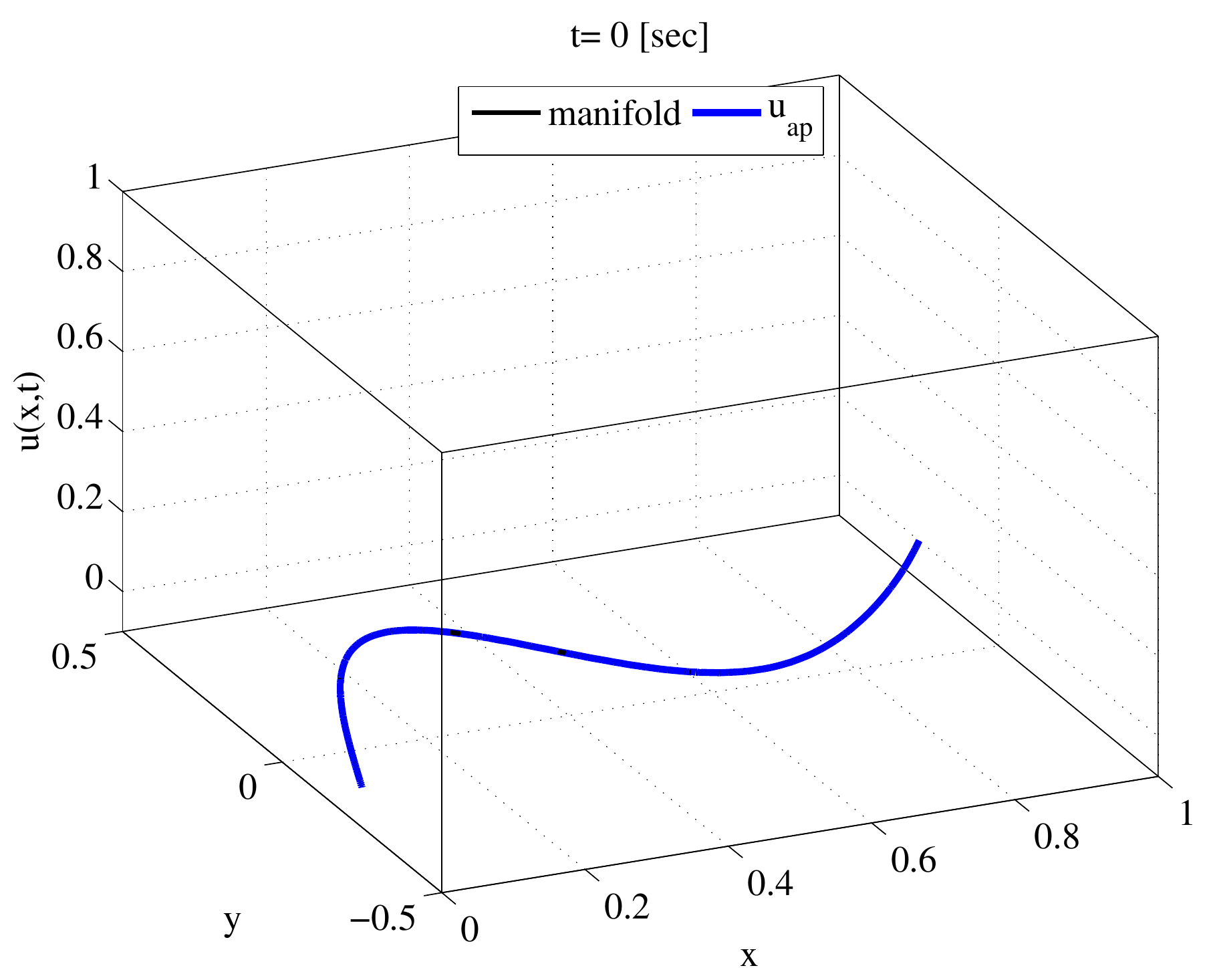}}\hfill\subfigure[State at $0.19$]{\includegraphics[width=0.32\textwidth]{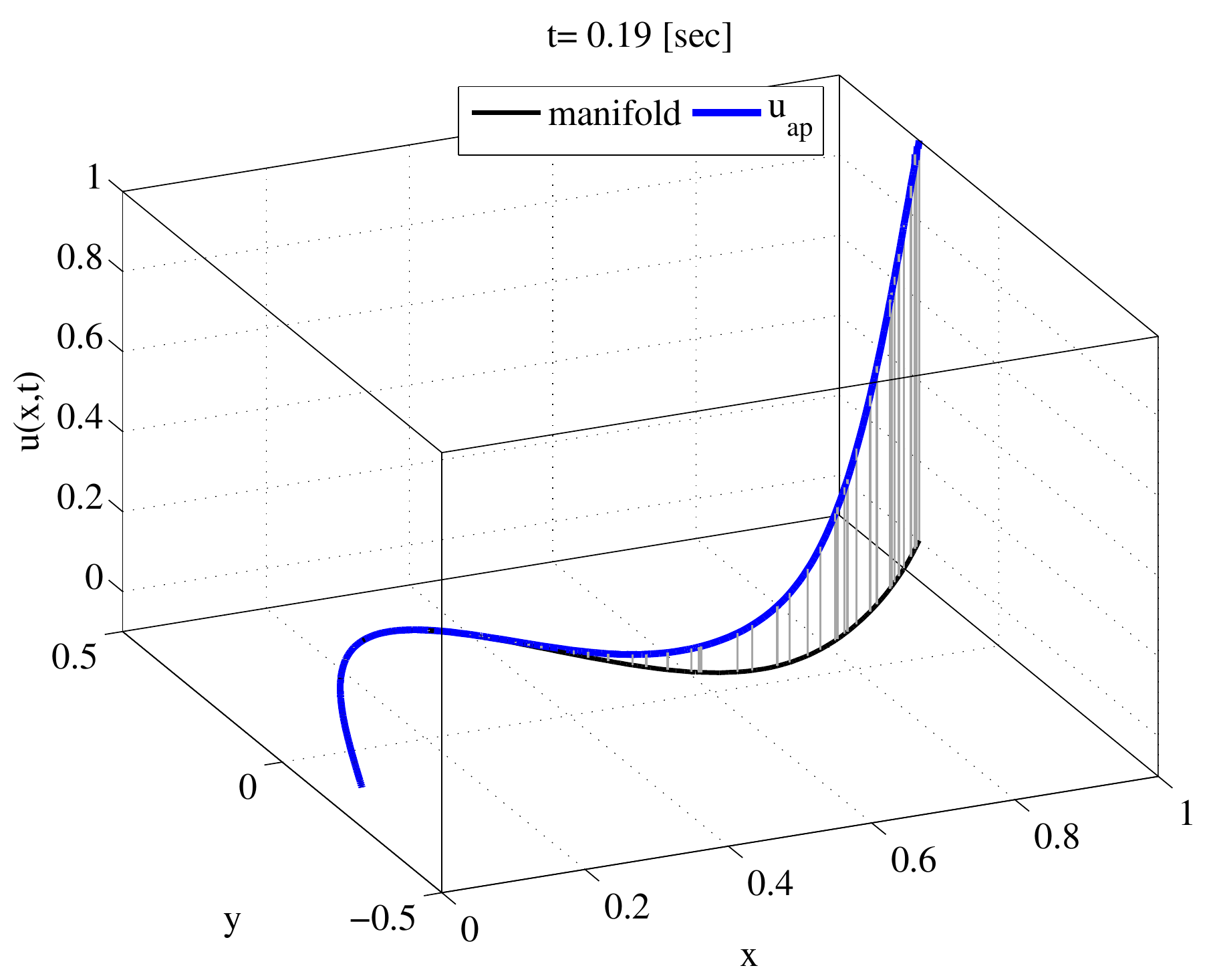}}\hfill\subfigure[State at $1.00$]{\includegraphics[width=0.32\textwidth]{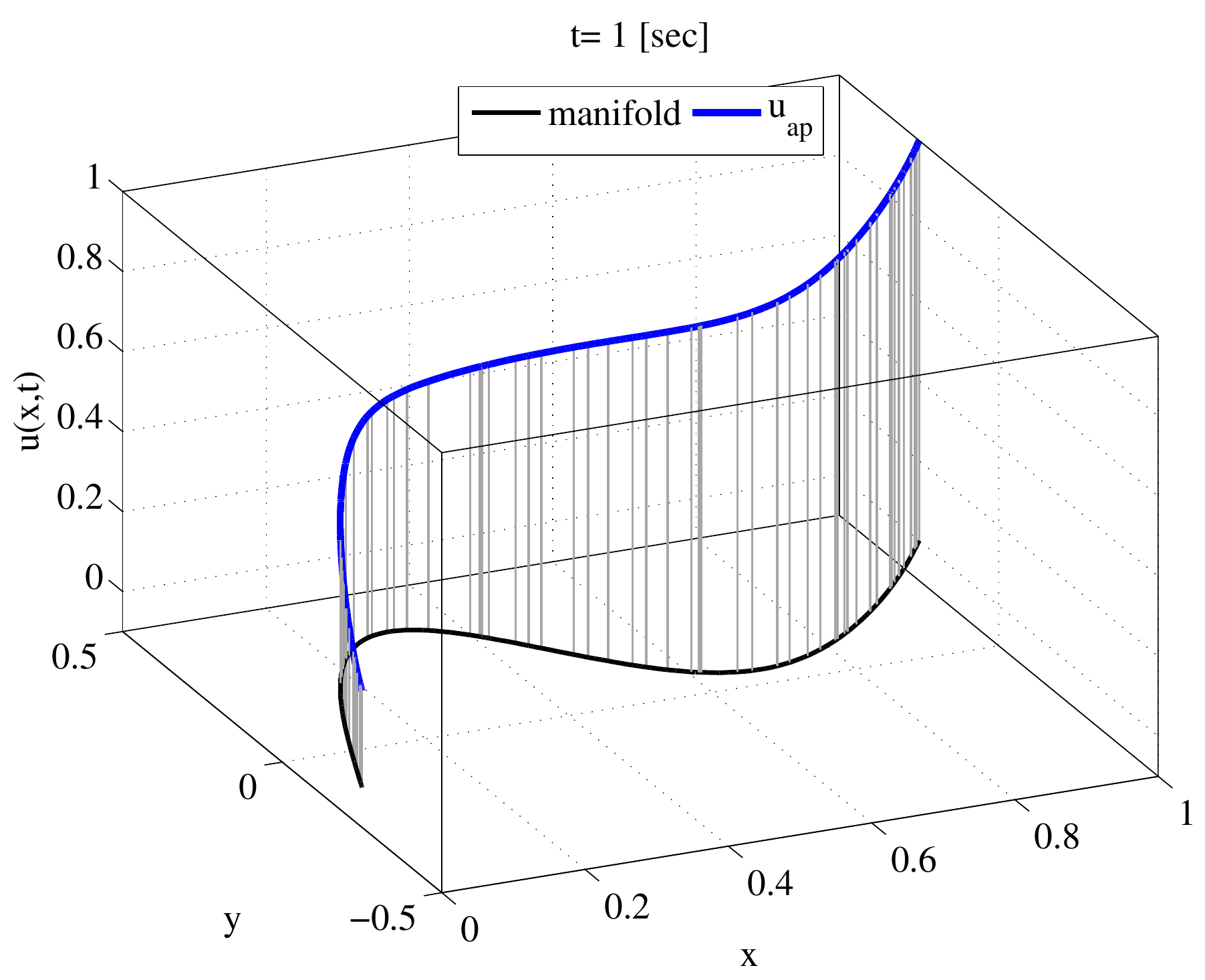}}

\caption{\label{fig:TranspSShapeRes} Solution of an instationary transport
problem on an S-shaped manifold.}
\end{figure}

Next, the instationary transport problem is solved on the S-shaped
manifold of Section \ref{sec:lap1dEx3} in the time interval $t\in(0,1)$.
The advection velocity is $\left\Vert \vek c_{\Gamma}\right\Vert =1$
and the diffusion coefficient $\lambda=0.15$. A Dirichlet boundary
condition of $u=1$ is prescribed at the inflow. The initial condition
is $u=0$ everywhere on $\Gamma$. As there is no analytical solution
available, no convergence study is performed and only a representative
approximation is shown in Fig.~\ref{fig:TranspSShapeRes}.

\subsubsection{Advection-diffusion on a hyperbolic paraboloid with bumps\label{sec:ex2dIn}}

\begin{figure}
\centering

\subfigure[Initial state]{\includegraphics[width=0.32\textwidth]{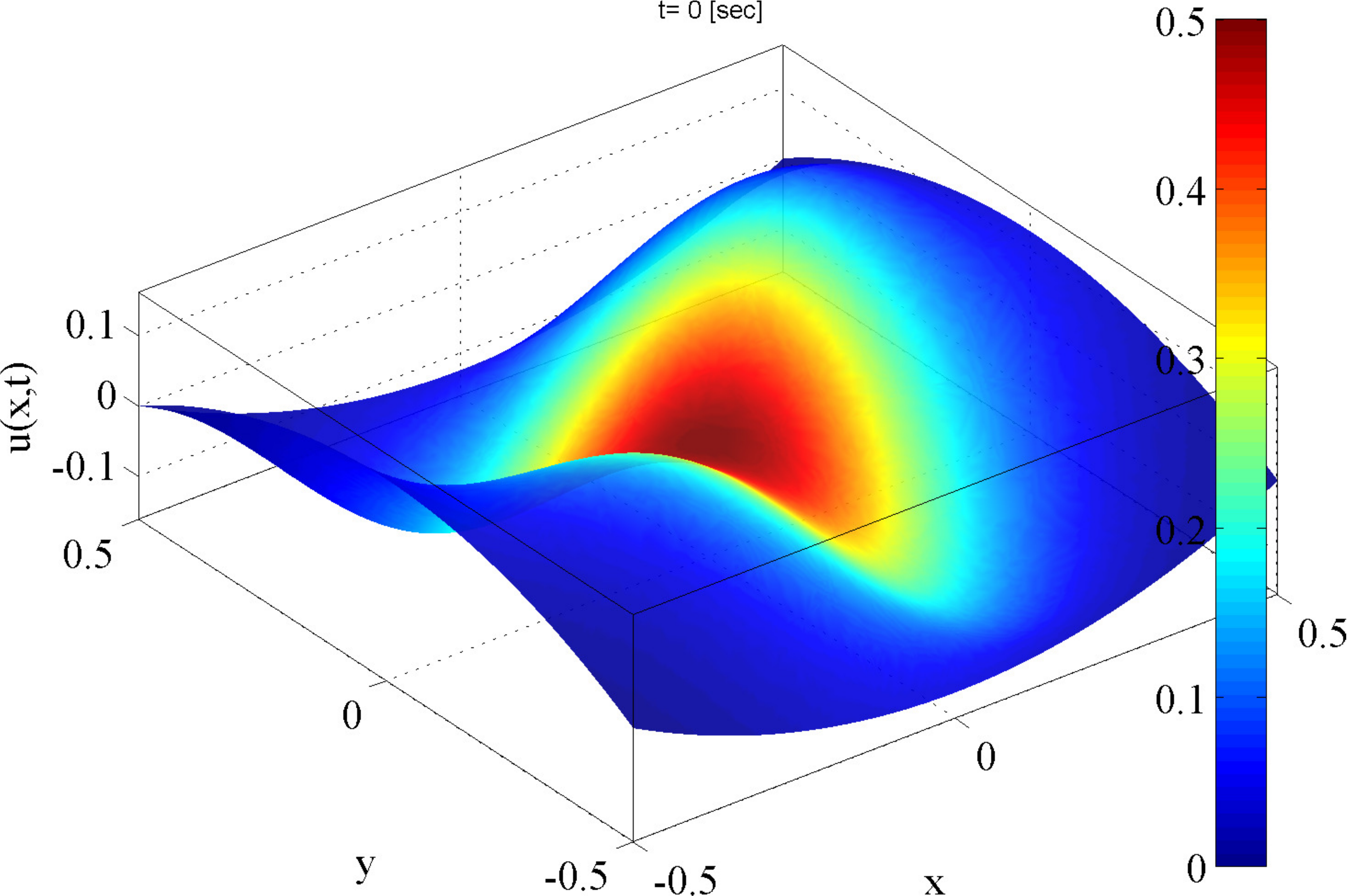}}\hfill\subfigure[State at $0.49$]{\includegraphics[width=0.32\textwidth]{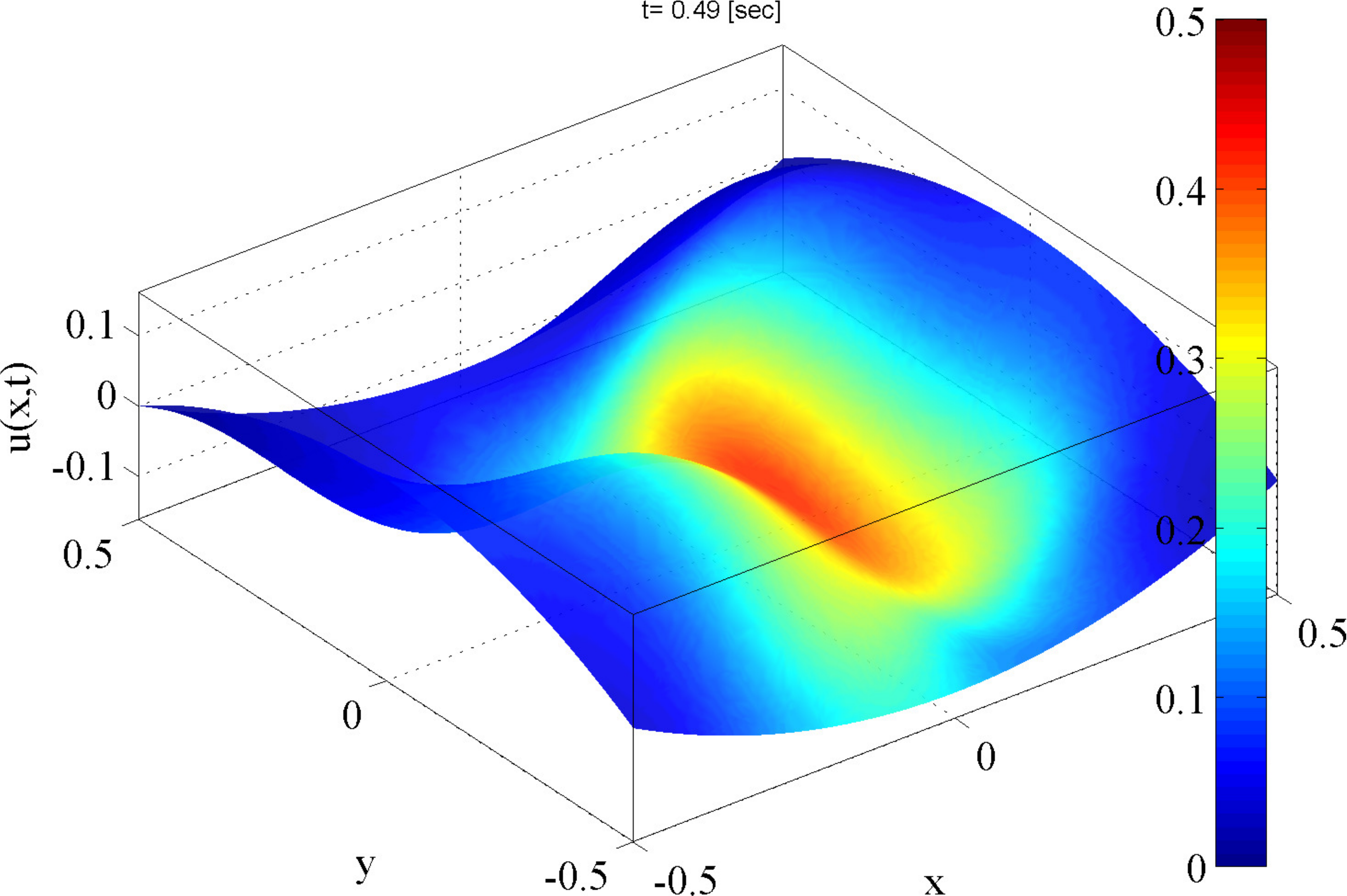}}\hfill\subfigure[State at $1.00$]{\includegraphics[width=0.32\textwidth]{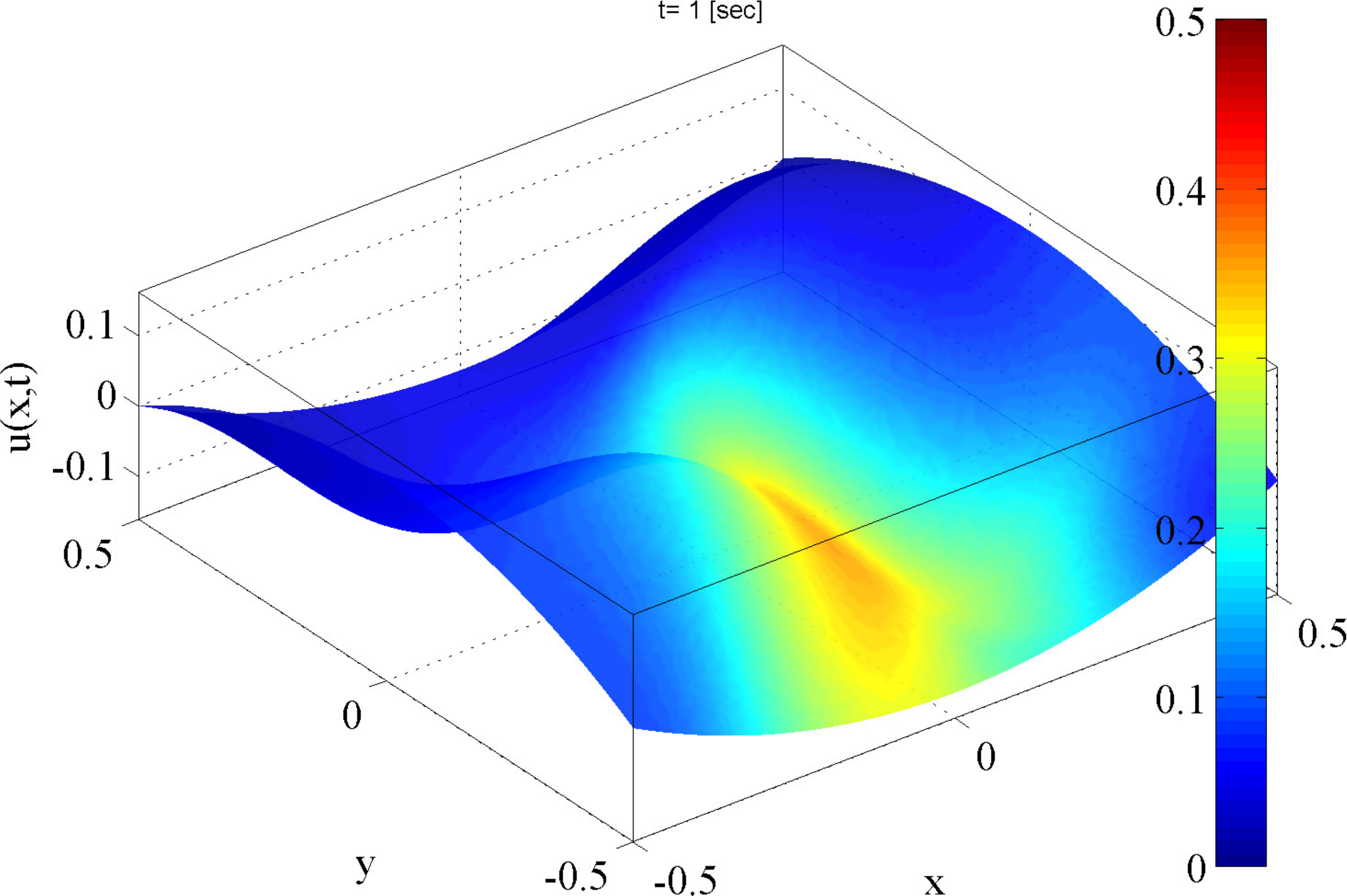}}

\caption{\label{fig:TranspHyperBolRes} Solution of an instationary transport
problem on a hyperbolic parabolid with bumps.}
\end{figure}

Finally, an advection-diffusion problem is considered on the hyperbolic
paraboloid with bumps introduced in Section \ref{sec:lap2dEx3}. The
time interval is again $t\in(0,1)$, the advection velocity $\left\Vert \vek c_{\Gamma}\right\Vert =1.25$
tangential to the manifold in the direction of $y$, and the diffusion
coefficient $\lambda=0.01$. The initial condition is given by $u(\vek x,0)=\dfrac{1}{2}\exp\left[-10\left(x^{2}+y^{2}\right)\right]$.
A representative result is seen in Fig.~\ref{fig:TranspHyperBolRes}.

\section{Conclusions\label{X_Conclusions}}

A higher-order accurate method for PDEs on surfaces is proposed. It
enables a completely automatic workflow from the geometric description
based on level-sets to the numerical analysis without any user-intervention.
This is an important advantage over other methods which are based
on \emph{handcrafted} surface meshes. Furthermore, these meshes are
often composed by flat triangles leading to a low-order representation
of the geometry and the resulting approximation of the BVP. Compared
to methods which solve BVPs on all zero-isosurfaces at once by using
volumetric elements, the proposed approach is more efficient as the
effort scales with standard planar, two-dimensional BVPs. Compared
to methods which employ shape functions of the background mesh as
in TraceFEM and CutFEM, it is an important advantage that boundary
conditions are enforced in the standard way without additional technqiues
for general constraints (Lagrange multipliers, Nitsche's method etc.).

The proposed method is characterized by the following key ingredients:
\begin{enumerate}
\item A geometry description of the bounded manifold based on several level-set
functions.
\item The automatic generation of surface elements which enable a $C_{0}$-continuous,
higher-order accurate representation of the level-set geometry. A
conforming surface mesh used for the finite element approximation
is automatically extracted from this set of surface elements. 
\item The manipulation of the background mesh by node movements in order
to ensure the shape regularity of the resulting surface elements and
a bounded condition number of the resulting system of equations. This
step may later be replaced by a stabilization similar to what is done
in the TraceFEM and CutFEM.
\end{enumerate}
The numerical results confirm that higher-order accurate solutions
of BVPs on curved surfaces in three dimensions are achieved. The next
steps will be to investigate the proposed technique for the solution
of more advanced BVPs on manifolds such as flow problems on curved
surfaces and the numerical analysis of membranes and shells.

\bibliographystyle{schanz}
\addcontentsline{toc}{section}{\refname}\bibliography{FriesRefs}
 
\end{document}